\documentclass[a4paper, 10pt, openany]{book}
\usepackage{graphicx,graphics}
\usepackage{calc}
\usepackage[cmtip,arrow,all,2cell]{xy}
\usepackage{pb-diagram,pb-xy}
\usepackage{amsmath,amsthm,amssymb,tabularx}
\usepackage[mathscr]{eucal}
\usepackage{hyperref}

\numberwithin{equation}{section}

\theoremstyle{plain}

\newtheorem{tm}{Theorem}[section]

\newtheorem{prop}[tm]{Proposition}
\newtheorem{lm}[tm]{Lemma}
\newtheorem{cor}[tm]{Corollary}

\theoremstyle{definition}

\newtheorem{ex}{Example}[section]

\newtheorem{rem}[ex]{Remark}

\newcommand{\btm}{\begin{tm}}
\newcommand{\etm}{\end{tm}}
\newcommand{\blm}{\begin{lm}}
\newcommand{\elm}{\end{lm}}
\newcommand{\bprop}{\begin{prop}}
\newcommand{\eprop}{\end{prop}}
\newcommand{\bcor}{\begin{cor}}
\newcommand{\ecor}{\end{cor}}
\newcommand{\bex}{\begin{ex}}
\newcommand{\eex}{\end{ex}}
\newcommand{\brem}{\begin{rem}}
\newcommand{\erem}{\end{rem}}

\newcommand{\bpr}{\begin{proof}}
\newcommand{\epr}{\end{proof}}
\newcommand{\beq}{\begin{equation}}
\newcommand{\eeq}{\end{equation}}
\newcommand{\bit}{\begin{itemize}}
\newcommand{\eit}{\end{itemize}}

\newcommand{\norm}[1]{\left\Vert#1\right\Vert}
\newcommand{\abs}[1]{\left\vert#1\right\vert}

\allowdisplaybreaks

\def \le {\leqslant}
\def \ge {\geqslant}

\def\e{\varepsilon}
\def\R{{\mathbb{R}}}
\def\C{{\Bbb C}}

\def \N {\mathbb{N}}

\def \Z {\mathbb{Z}}

\def\Bind{\operatorname{\sf Bind}}

\def \d {\operatorname{\sf{d}}}
\def \alt{\operatorname{\sf{alt}}}

\def\Spec{\operatorname{\sf{Spec}}}

\def\Ob{\operatorname{\sf Ob}}

\def\id{\operatorname{\sf id}}
\def\Mor{\operatorname{\sf Mor}}
\def\card{\operatorname{\sf card}}

\def\Sp{\operatorname{\sf Span}}

\def\d{\operatorname{d}}

\def\Mono{\operatorname{\sf Mono}}
\def\Epi{\operatorname{\sf Epi}}
\def\Bim{\operatorname{\sf Bim}}

\def\SMono{\operatorname{\sf SMono}}
\def\SEpi{\operatorname{\sf SEpi}}

\def\DEpi{\operatorname{\sf DEpi}}
\def\DBim{\operatorname{\sf DBim}}

\def\Iso{\operatorname{\sf Iso}}

\def\env{\operatorname{\sf env}}
\def\rf{\operatorname{\sf ref}}
\def\Env{\operatorname{\sf Env}}
\def\Rf{\operatorname{\sf Ref}}

\def\H{\operatorname{\sf H}}
\def\h{\operatorname{\sf h}}

\def\Coim{\operatorname{\sf Coim}}
\def\coim{\operatorname{\sf coim}}
\def\Im{\operatorname{\sf Im}}
\def\im{\operatorname{\sf im}}
\def\red{\operatorname{\sf red}}
\def\ker{\operatorname{\sf ker}}
\def\Ker{\operatorname{\sf Ker}}
\def\coker{\operatorname{\sf coker}}
\def\Coker{\operatorname{\sf Coker}}
\def\Dom{\operatorname{\sf Dom}}
\def\Ran{\operatorname{\sf Ran}}
\def\Sub{\operatorname{\sf Sub}}
\def\Quot{\operatorname{\sf Quot}}
\def\alg{\operatorname{\sf alg}}

\def\leftlim{\mathop{\varprojlim}\limits}
\def\rightlim{\mathop{\varinjlim}\limits}

\def \N {\mathbb{N}}
\def \Z {\mathbb{Z}}
\def \R {\mathbb{R}}
\def \T {\mathbb{T}}
\def \C {\mathbb{C}}

\def\id{\operatorname{\sf id}}

\def\d{\operatorname{\sf d}}

\def\exp{\operatorname{\sf exp}}
\def\card{\operatorname{\sf card}}
\def\alg{\operatorname{\sf alg}}

\def\d{\operatorname{\sf{d}}}

\def\i{\operatorname{\sf{i}}}
\def\h{\operatorname{\sf{h}}}

\def \e {\varepsilon}
\def \ph {\varphi}

\def\boxedast{\raise.8pt\hbox{\rlap{$\mkern2.5mu *$}}\square}

\def\subarr{\subset\kern-7pt\raise-2.4pt\hbox{\resizebox{5pt}{3.5pt}{$\to$}}\kern2pt}
\def\suparr{\reflectbox{$\subarr$}}

\def\osubarr{\subset\kern-6pt\raise-2.4pt\hbox{\resizebox{5pt}{3.5pt}{$\to$}}\kern-3pt\raise3.4pt\hbox{\resizebox{3pt}{3pt}{$\circ$}}\kern2pt}

\def\quarr{\gets\kern-3.7pt\raise0pt\hbox{\resizebox{2pt}{2.5pt}{$\backslash$}}\kern2pt}

\def\qparr{\kern-1pt\reflectbox{$\quarr$}\kern1pt}

\def\oquarr{\gets\kern-3.7pt\raise0pt\hbox{\resizebox{2pt}{2.5pt}{$\backslash$}}\kern-1.6pt\raise-3.5pt\hbox{\resizebox{3pt}{3pt}{$\circ$}}\kern5pt}

\def\ole{\le\kern-1.6pt\raise-3.5pt\hbox{\resizebox{3pt}{3pt}{$\circ$}}\kern5pt}

\def\ph{\varphi}
\def\l{\left(}
\def\r{\right)}

\def\ole{\le\kern-4.5pt\raise5.pt\hbox{\resizebox{3pt}{3pt}{$\circ$}}\kern2pt}

\numberwithin{equation}{section}

\makeindex

\addtocounter{section}{-1}

\begin{document}

\title{ENVELOPES AND REFINEMENTS IN CATEGORIES, WITH APPLICATIONS TO FUNCTIONAL ANALYSIS}

\author{S.S.Akbarov}

\maketitle

\tableofcontents

\def\theequation{\Alph{equation}}

\addcontentsline{toc}{section}{Introduction}
\section*{Introduction}\label{SEC-vvedenie}

In 1972 J.~L.~Taylor introduced in his paper \cite{Taylor-1} an operation, which associates to an arbitrary topological algebra $A$ a new topological algebra $\Env A$ called by A.~Ya.~Helemskii later ``the Arens-Michael envelope of $A$'' \cite{Helemskii}. Immediately after that in his next paper \cite{Taylor-2} Taylor gave an amusing formula\footnote{Taylor mentions this fact in passing on pages 207 and 251 in \cite{Taylor-2}.} which suggests an unexpectedly simple way to formalize the heuristically evident connection between algebraic geometry and complex analysis:
\beq\label{P^heartsuit(C^n)=O(C^n)}
\Env {\mathcal P}(\C^n)={\mathcal O}(\C^n)
\eeq
(here ${\mathcal P}(\C^n)$ and ${\mathcal O}(\C^n)$ are the algebras of polynomials and, respectively, of holomorphic functions on the complex space $\C^n$). Despite this promising application, up to the end of the century Taylor's construction did not manifest itself in mathematical literature, and only recently the interest to the operation $A\mapsto \Env A$ appeared again in A.Yu.Pirkovskii's papers on ``holomorphic non-commutative geometry'' \cite{Pirkovskii-1}, \cite{Pirkovskii-2}. In particular, in \cite{Pirkovskii-2} formula \eqref{P^heartsuit(C^n)=O(C^n)} was generalized to the case of arbitrary affine algebraic variety $M$:
\beq\label{P^heartsuit(M)=O(M)}
\Env {\mathcal P}(M)={\mathcal O}(M).
\eeq
This identity very soon was applied by the author in \cite{Akbarov-stein-groups} to the construction of a generalization of Pontryagin's duality from the category of commutative compactly generated Stein groups to the category of arbitrary (not necessarily, commutative) compactly generated Stein groups with the algebraic connected component of identity. The idea of the duality suggested in \cite{Akbarov-stein-groups} is illustrated by the diagram
 \beq\label{int:chetyrehugolnik-O-O*}
 \xymatrix @R=2.pc @C=4.pc @M=14pt
 {
 {\mathcal O}^\star(G)
  \ar@{|->}[r]^{\Env} &
 {\mathcal O}_{\exp}^\star(G)
 \ar@{|->}[d]^{\star}
 \\
 {\mathcal O}(G) \ar@{|->}[u]^{\star}
 &
 {\mathcal O}_{\exp}(G) \ar@{|->}[l]_{\Env}
 }
 \eeq
where $G$ is a group of the described class, ${\mathcal O}(G)$ the algebra of  holomorphic functions on $G$, ${\mathcal O}_{\exp}(G)$ its subalgebra, consisting of functions of exponential type, $A\mapsto \Env A$ the operation of taking Arens-Michael envelope, and $X\mapsto X^\star$ the operation of passage to the dual stereotype space in the sense of \cite{Akbarov}, i.e. to the space of linear continuous functionals with the topology of uniform convergence on totally bounded sets (in this case this is equivalent to the uniform convergence on compact sets).

One can call duality, presented in diagram \eqref{int:chetyrehugolnik-O-O*}, the {\it complex geometry duality}, having in mind the class of objects under consideration. The obtained theory for the described class of groups contrasts with the other existing theories in the following two points. First, its enveloping category (to which the group algebras belong) consists of Hopf algebras. And, second, the diagram \eqref{int:chetyrehugolnik-O-O*} suggests a natural way for constructing the analogous dualities for the ``other geometries'', in particular, for differential geometry and for  topology: one should just replace the Arens-Michael envelope in diagrams analogous to \eqref{int:chetyrehugolnik-O-O*} with some other envelopes (and this automatically leads to the replacing of the constructions in the corners of the diagram with some proper analogs from analysis).

This alleged connection between different dualities in geometry and different envelopes of topological algebras was recently vouched by other examples:
 \bit{
\item[1)] In the work by J.~N.~Kuznetsova \cite{Kuznetsova} the Arens-Michael envelope was replaced by the envelope generated by the $C^*$-quotient maps\footnote{Below on page \pageref{DEF:C*-obolochka} we define this construction as the Kuznetsova envelope.}, and this immediately led to a variant of {\it topological duality}, where the Stein groups are replaced by the Moore groups, and the algebras ${\mathcal O}(G)$ and ${\mathcal O}_{\exp}(G)$, respectively, by the algebra ${\mathcal C}(G)$ of continuous functions on $G$ and the algebra ${\mathcal K}(G)$ of coefficients of norm-continuous representations of $G$.

\item[2)] In author's work \cite{Akbarov-C-infty} a notion of smooth envelope $\Env_{\infty} A$ of a topological algebra $A$ was introduced. This construction replaces the Arens-Michael envelope in passage from complex analysis to the differential geometry, and an analogue of the Pirkovski theorem \eqref{P^heartsuit(M)=O(M)} was proved in the differential-geometric context: if a subalgebra $A$ in the algebra ${\mathcal C}^\infty(M)$ of smooth functions on a smooth variety $M$ has the same spectrum and the same tangent space in each point, then
$$
\Env_{\infty} A={\mathcal C}^\infty(M).
$$
This result gives a hope that a similar {\it duality theory in differential geometry} will be constructed in near future with a proper class of real Lie groups.
}\eit
It is interesting (and predictable), that in these theories the classical Fourier and Gelfand transforms are interpreted as envelopes with respect to the prescribed class of algebras (see e.g. below Theorems \ref{TH:Fourier-Stein=obolochka}, \ref{TH:Fourier-LCA=obolochka} and \ref{EX:C^*-obolochka}).

It is clear to the author that the obtained results are just first observations in the indicated field, but they already show the validity of the common philosophical idea which justifies and guides the investigations in this area: in each standard mathematical disciplines, where certain classes of symmetries play role (classes of groups, including those understood in generalized way, like quantum groups), a certain duality theory works (and apparently, not unique). This idea was suggested in the author's work \cite{Akbarov-stein-groups}, and among such disciplines the following four were mentioned:
 \bit{
\item[--] general topology,

\item[--] differential geometry,

\item[--] complex analysis,

\item[--] algebraic geometry.
}\eit

This paper is planned as a part of the program, described in \cite{Akbarov-stein-groups}. We discuss here the question (which remained open up to the last time), how one should define envelopes in general category theory, and under which conditions they exist and are functors? We suggest a natural definition (from our point of view) and establish some wide necessary and sufficient conditions for existence of envelopes and their dual constructions, which we call refinements. As applications, we show that in the categories ${\tt Ste}$ of stereotype spaces, and ${\tt Ste}^\circledast$ of stereotype algebras the envelopes and the refinements exist in a very wide class of situations.

\addcontentsline{toc}{section}{Acknowledgements}
\section*{Acknowledgements}

The author thanks A.~Blass, S.~Buschi, Ya.~A.~Kopylov, B.~V.~Novikov and A.~Yu.~Pirkovskii for useful consultations.

\addcontentsline{toc}{section}{Agreements and Notations}
\section*{Agreements and Notations}

Everywhere in category theory we use the terminology of textbooks \cite{Bucur-Deleanu}, \cite{Tsalenko-Shulgeifer} and of handbook  \cite{General-algebra}, and as a set-theoretic fundament for the notion of category we choose the Morse-Kelley theory \cite{38}.

The notations $\Mono({\tt K})$, $\Epi({\tt K})$, $\SMono({\tt K})$ and $\SEpi({\tt K})$ mean the classes of monomorphisms, epimorphisms, strong monomorphisms and strong epimorphisms (the last two are defined at p.\pageref{razbienie-kvadrata}) respectively in the category ${\tt K}$.
We say that a category ${\tt K}$ is
 \bit{
\item[---] {\it injectively (projectively) complete}, if each functor $K:{\tt M}\to{\tt K}$ from a small category $\tt M$  (i.e. a category where the class of morphisms is a set) has an injective (projective) limit,

\item[---] {\it complete}, if it is injectively and projectively complete,

\item[---] {\it finitely injectively (projectively) complete}, if each functor $K:{\tt M}\to{\tt K}$ from a finite category $\tt M$ (i.e. a category where the class of morphisms is a finite set) has an injective (projective) limit,

\item[---] {\it finitely complete}, if it is finitely injectively complete and finitely projectively complete,

 \item[---] {\it linearly complete}\label{DEF:lineino-polnaya-kategoriya}, if any functor from a linearly ordered set to ${\tt K}$ has injective and projective limits.
 }\eit

For any morphism $\ph:X\to Y$ in an arbitrary category the symbols
$\Dom\ph$ and $\Ran\ph$ mean respectively the domain and the range of $\ph$,
i.e. $\Dom\ph=X$ and $\Ran\ph=Y$. If $\tt L$ and $\tt M$ are two classes of objects in $\tt K$, then $\Mor({\tt L},{\tt M})$ means the class of morphisms with domains in $\tt L$ and ranges in $\tt M$:
$$
\Mor({\tt L},{\tt M})=\{\ph\in\Mor({\tt K}):\quad \Dom\ph\in{\tt L}\ \&\ \Ran\ph\in{\tt M}\}.
$$

Let $\varPhi$ be a class of morphisms and ${\tt L}$ a class of objects in a category ${\tt K}$. We say that
\bit{
\item[---]\label{DEF:goes-from} {\it $\varPhi$ goes from ${\tt L}$}, if for any object $X\in{\tt L}$ there is a morphism $\ph\in\varPhi$, going from $X$:
    $$
    \forall X\in{\tt L}\qquad \exists \ph\in\varPhi\qquad \Dom\ph=X;
    $$
in the special case, if $\tt L$ consists of only one object $X$, we say that $\varPhi$ {\it goes from $X$},

\item[---]\label{DEF:goes-to} {\it $\varPhi$ goes to ${\tt L}$}, if for any object $X\in{\tt L}$ there is a morphism $\ph\in\varPhi$, going to $X$:
    $$
    \forall X\in{\tt L}\qquad \exists \ph\in\varPhi\qquad \Ran\ph=X.
    $$
in the special case when $\tt L$ consists of only one object $X$, we say that $\varPhi$ {\it goes to $X$}.
    }\eit

In the theory of topological vector spaces we follow the textbook by H.~Schaefer \cite{Schaefer}, and in the theory of stereotype spaces and algebras author's papers \cite{Akbarov} and \cite{Akbarov-stein-groups}.
In particular, following \cite{Schaefer} we assume that {\it all locally convex spaces (LCS, shortly) are Hausdorff}.
By {\it topological algebra} we mean locally convex topological algebra in the spirit of textbook \cite{Mallios}, i.e. locally convex space $A$ over the field $\C$, endowed with associative multiplication which is {\it separately continuous} and has unit.

We use also the following notations.
First, for any locally convex space $X$ the
symbol  ${\mathcal U}(X)$ denotes the system of all neighborhoods of zero in $X$.
Second, for each neighborhood of zero $U$ in $X$ the set
$$\label{DEF:Ker(U)}
\Ker U=\bigcap_{\varepsilon>0}\varepsilon\cdot U
$$
will be called the {\it kernel} of this neighborhood of zero. If $U$ is an absolutely convex neighborhood of zero, then its kernel $\Ker U$ is a closed subspace in $X$. And, thurd, if a topological space $Y$ is imbedded into a topological space $X$ (injectively, but not necessarily in such a way that the topology of $Y$ is inherited from $X$),
and $A$ is a subset in $Y$, then to distinguish the closure of $A$ in $Y$ from its closure in $X$,
we denote the first one by $\overline{A}^Y$, and the second by $\overline{A}^X$\label{DEF:overline^X}.

Besides this we say that a subset $M$ in a locally convex space $X$ is {\it total (in $X$)},
if its linear span $\Sp M$ is dense in $X$:
$$
\overline{\Sp M}^X=X.
$$

\def\theequation{\arabic{section}.\arabic{equation}}

\section{Nodal decomposition and factorizations}

\subsection{Skeletally small graphs}

\paragraph{Graphs.} Recall that an {\it oriented graph} is a set $V$ with a given subset $\varGamma$ in its cartesian square $V\times V$. The elements of $V$ are called vertices, and the elements of $\varGamma$ edges of this graph. An oriented graph is said to be {\it reflexive}, if for each vertex  $x\in V$ the edge $(x,x)$ belongs to $\varGamma$, and {\it transitive}, if for any two edges $(x,y)$ and $(y,z)$ from $\varGamma$ the pair $(x,z)$ also belongs to $\varGamma$. Obviously, every reflexive transitive oriented graph is a (small) category, where objects and morphisms are respectively the vertices and the edges (the multiplication of edges $(x,y)$ and $(y,z)$ is the edge $(x,z)$, and local identities $1_x$ are  $(x,x)$). The characteristic property of such categories (apart from the requirement of being small), is that the sets of morphisms $\Mor(A,B)$ always contain at most one element. This justifies the following definition.

 \bit{
\item[$\bullet$] A {\it graph}\label{DEF:graph} is a category ${\tt K}$ (not necessarily small), where each set of morphisms $\Mor(A,B)$ contains at most one element:
\beq\label{DEF:harakt-sv-grafov}
\forall A,B\in\Ob({\tt K})\quad \card\Mor(A,B)\le 1.
\eeq
Clearly, this condition is equivalent to establishing the structure of (reflexive and transitive) oriented graph at the class $\Ob{\tt K}$ of objects of the category ${\tt K}$ (with the difference that $\Ob{\tt K}$ is not necessarily a set, but just a class).
 }\eit

\medskip
\centerline{\bf Properties of graphs:}

\bit{\it

\item[$1^\circ$.]\label{GRAPH:mu=mono-<=>-exists-<-} In any graph a morphism $\ph:A\to B$ is an isomorphism, iff there exists an arbitrary morphism to the reverse direction $\psi:A\gets B$,
\beq\label{ph-in-Iso<=>exists-psi-B->A}
\forall \ph\in\Mor(A,B)\qquad \Big( \ph\in\Iso\quad\Longleftrightarrow\quad \exists \psi\in\Mor(B,A)\Big).
\eeq

\item[$2^\circ$.]\label{GRAPH:psi-0-ph=1-<=>-ph-0-psi=1} In any graph a composition of morphisms is an identity iff the same remains true after replacing the factors:
\beq\label{psi-circ-ph=1<=>ph-circ-psi=1}
\psi\circ\ph=1\quad\Longleftrightarrow\quad \ph\circ\psi=1.
\eeq

\item[$3^\circ$.]\label{GRAPH:psi-circ-ph-in-Iso<=>psi-in-Iso-i-ph-in-Iso} In any graph a composition of morphisms $\psi\circ\ph$ is an isomorphism iff both $\psi$ and $\ph$ are isomorphisms:
\beq\label{psi-circ-ph-in-Iso<=>psi-in-Iso-i-ph-in-Iso}
\psi\circ\ph\in\Iso\quad\Longleftrightarrow\quad \psi\in\Iso\ \& \ \ph\in\Iso.
\eeq

}\eit

\bpr
1. If $\ph:A\to B$ and $\psi:A\gets B$, then $\psi\circ\ph$ acts from $A$ into $A$, so it must coincide with $1_A$. Similarly, $\ph\circ\psi$ acts from $B$ into $B$, so it must coincide with $1_B$.

2. From $\psi\circ\ph=1$ it follows that $\Ran\ph=\Dom\psi$ and $\Ran\psi=\Dom\ph$, and after that we apply the same reasoning as in step 1.

3. If $\omega=\psi\circ\ph\in\Iso$, then $\psi\circ\ph\circ\omega^{-1}=1$, so by \eqref{psi-circ-ph=1<=>ph-circ-psi=1},
$\ph\circ\omega^{-1}\circ\psi=1$, hence, $\psi\in\Iso$, and finally $\ph=\psi^{-1}\circ\omega\in\Iso$.
\epr

\paragraph{Partially ordered classes.}

Every partially ordered set $I$ can be considered as a category, where objects are elements of this set, and morphisms are pairs $(i,j)$, for which $i\le j$. Such categories ${\tt K}$, of course, are special cases of graphs, since every set of morphisms $\Mor(A,B)$ here contains at most one element (i.e. \eqref{DEF:harakt-sv-grafov} holds). But in addition (and this property distinguishes the partially ordered sets among all graphs), for $A\ne B$ the existence of a morphism $\ph:A\to B$ automatically makes impossible the existence of any morphisms $\psi:A\gets B$. This justifies the following definition.

 \bit{
\item[$\bullet$] A {\it partially ordered class}\label{DEF:chast-upor-klass} is a graph,
where the existence of opposite morphisms $\ph:A\to B$ and $\psi:A\gets B$ is possible only if $A=B$
(and then $\ph=\psi=1_A$). In other words,
\beq\label{DEF:harakt-sv-chast-upor-klassov}
\forall A\ne B\in\Ob({\tt K})\quad \Mor(A,B)\ne\varnothing\quad\Longrightarrow\quad \Mor(B,A)=\varnothing.
\eeq
Obviously, these requirements are equivalent to the establishing the structure of partial order at the class $\Ob{\tt K}$ of objects of the category ${\tt K}$ (again. like in the previous definition, with the difference that $\Ob{\tt K}$ is not necessarily a set, but just a class).
 }\eit

\bex {\it Category of ordinal numbers ${\tt Ord}$.} The class ${\tt Ord}$ of all ordinal numbers with its natural order (see e.g. \cite{38}) is an example of a partially ordered class which is not a set.
\eex

\bprop\label{PROP:iso-v-chast-upor-klasse} In a partially ordered class only local identities are isomorphisms:
$$
\forall \ph\in\Mor(A,B)\qquad\Big( \ph\in\Iso \quad\Longleftrightarrow\quad A=B\quad\&\quad \ph=1_A\Big).
$$
\eprop
 \bpr The identity $A=B$ follows from the fact that $\Mor(A,B)\ne\varnothing$ and $\Mor(B,A)\ne\varnothing$, and the identity $\ph=1_A$ from the fact that $\ph$ and $1_A$ are colinear arrows in a graph.
 \epr

\paragraph{Skeleton.}

A class $S$ of objects of a category ${\tt K}$ is called a {\it skeleton} of ${\tt K}$, if every object in ${\tt K}$ is isomorphic to an exactly one object of $S$. In other words, $S$ satisfies the following two requirements:
 \bit{
 \item[1)] elements of $S$ are isomorphic only if they coincide:
 $$
 \forall X,Y\in S\qquad (X\cong Y\quad\Leftrightarrow\quad X=Y);
 $$

 \item[2)] there exists a map $G:\Ob(\tt K)\to S$ from the class of objects of $\tt K$ to the class $S$ such that
$$
\forall X\in \Ob({\tt K})\quad X\cong G(X).
$$
}\eit

The skeleton $S$ is usually endowed with the structure of a full subcategory in ${\tt K}$.

\medskip
\centerline{\bf Properties of skeleton:}

\bit{\it

\item[$1^\circ$.]\label{1^0:skelet-sushestvuet} Each category ${\tt K}$ has a skeleton.

\item[$2^\circ$.] Each two skeletons in ${\tt K}$ are isomorphic (as categories).

\item[$3^\circ$.] Each category ${\tt K}$ is equivalent to its skeleton $S$.

\item[$4^\circ$.] Two categories ${\tt K}$ and ${\tt L}$ are equivalent if and only if their skeletons are isomorphic  (as categories).
}\eit
\bpr
Only the first proposition is not obvious here. It follows from the fact that the class $\tt Set$ of all sets can be well-ordered (see. \cite[V, 4.1]{Levy}): the class $\Ob(\tt K)$ of all objects of $\tt K$ is a subclass in the class $\tt Set$ of all sets, so $\Ob(\tt K)$ can also be well-ordered, and after that we can assign to each object $X\in\Ob(\tt K)$ the minimal object among all isomorphic to $X$ in $\tt K$ with respect to this order.
\epr

\bit{
\item[$\bullet$]
A category ${\tt K}$ is said to be
\bit{
\item[---] {\it skeletal}, if any two isomorphic objects coincide there (this is equivalent to the requirement that ${\tt K}$ is a skeleton for itself),

\item[---] {\it skeletally small}, if it has a skeleton, which is a set.
}\eit
}\eit

\bex Each partially ordered class is a skeletal category (since as we already noticed only local identities are isomorphisms there), but not vice versa. For instance, the category of all finite sets of the form $\{0,...,n\}$, $n\in\Z_+$, (with arbitrary maps as morphisms) is skeletal, but it is not a partially ordered class, since a set $\{0,...,n\}$ can have many bijections onto itself.
\eex

\paragraph{Transfinite chain condition.}

\bit{
\item[$\bullet$]
Let us say that a (covariant or contravariant) functor $F:{\tt Ord}\to {\tt K}$ {\it is stabilized}, if it satisfies the following two equivalent conditions:
    \bit{
\item[(i)] there exists an ordinal number $k\in{\tt Ord}$ such that
$$
\forall l\ge k\quad F(k,l)\in\Iso
$$

\item[(ii)] there exists an ordinal number $k\in{\tt Ord}$ such that
$$
\forall l,m\qquad\Big( k\le l\le m\quad\Longrightarrow\quad F(l,m)\in\Iso\Big)
$$
   }\eit
   }\eit
\bpr[Proof of equivalence] The implication  (i) $\Leftarrow$ (ii) is obvious, so we need to prove only (i) $\Rightarrow$ (ii). Let $F$ be a covariant functor (the case of a contravariant functor is considered similarly). If (i) holds, then for $k\le l\le m$ we have:
$$
\underbrace{F(k,m)}_{\scriptsize\begin{matrix}\text{\rotatebox{-90}{$\in$}}\\ \Iso\end{matrix}}=F(l,m)\circ \underbrace{F(k,l)}_{\scriptsize\begin{matrix}\text{\rotatebox{-90}{$\in$}}\\ \Iso\end{matrix}}\quad\Longrightarrow\quad
\underbrace{F(k,m)}_{\scriptsize\begin{matrix}\text{\rotatebox{-90}{$\in$}}\\ \Iso\end{matrix}}\circ \underbrace{F(k,l)^{-1}}_{\scriptsize\begin{matrix}\text{\rotatebox{-90}{$\in$}}\\ \Iso\end{matrix}}=F(l,m)\quad\Longrightarrow\quad F(l,m)\in\Iso
$$
\epr

\brem\label{REV:stabiliziruemost-v-chast-upor-klasse} If a category ${\tt K}$ is a partially ordered class, then by Proposition \ref{PROP:iso-v-chast-upor-klasse}, for a functor $F:{\tt Ord}\to {\tt K}$ the isomorphisms in (i) and (ii) become local identities:
    \bit{
\item[(i)$'$] there exists an ordinal number $k\in{\tt Ord}$ such that
$$
\forall l\ge k\quad F(k,l)=1_{F(k)}
$$

\item[(ii)$'$] there exists an ordinal number $k\in{\tt Ord}$ such that
$$
\forall l,m\qquad\Big( k\le l\le m\quad\Longrightarrow\quad F(l,m)=1_{F(l)}\Big)
$$
   }\eit
\erem

\btm[transfinite chain condition]\label{TH:obryv-transf-tsepei}
Every functor $F:{\tt Ord}\to {\tt K}$ into an arbitrary skeletally small graph ${\tt K}$ is stabilized.
\etm

We will need the following

\blm In the class ${\tt Ord}$ of ordinal numbers there is no a cofinal subclass, which is a set.
\elm
\bpr If $K$ is a cofinal subclass in ${\tt Ord}$, then ${\tt Ord}$ becomes a union of a family of sets, indexed by elements of $K$:
$$
{\tt Ord}=\bigcup_{k\in K}\{i\in {\tt Ord}:i\le k\}.
$$
Hence if $K$ is a set, then ${\tt Ord}$ must also be a set, but this is not true.
\epr

\bcor\label{COR:nepr-F:I->Ord} For any directed set $I$ each monotone map $F:I\to{\tt Ord}$ has a least upper bound in ${\tt Ord}$.
\ecor
\bpr It is sufficient to note here that the image $F(I)$ is bounded in ${\tt Ord}$. And this in its turn follows from the fact that $F(I)$ is a set, and thus cannot be a cofinal subclass in ${\tt Ord}$.
\epr

\bpr[Proof of Theorem \ref{TH:obryv-transf-tsepei}] Let $F:{\tt Ord}\to {\tt K}$ be a (covariant or contravariant) functor into a skeletally small graph ${\tt K}$. Suppose that it is not stabilized, i.e. for any ordinal number $i\in{\tt Ord}$ there is an ordinal number $j\in{\tt Ord}$ such that $F(i,j)\notin\Iso$. Let us construct a transfinite sequence of ordinal numbers $\{k_i;\ i\in{\tt Ord}\}\subseteq{\tt Ord}$ according to the following rules:
 \bit{

 \item[0)] We set $k_0=0$.

 \item[1)] If for some ordinal number $j\in{\tt Ord}$ all the ordinal numbers $k_i$ with the smaller indices $\{k_i;\ i<j\}$ are already chosen, then we consider two cases:

  \bit{

 \item[---] if $j$ is an isolated ordinal, i.e. $j=i+1$ for some $i<j$, then we take $k_j$ with the properties
     $$
     k_i<k_{i+1}=k_j,\qquad F(k_i,k_{i+1})=F(k_i,k_j)\notin\Iso
     $$
  ($k_j$ exists due to our assumption that $F$ is not stabilized),

 \item[---] if $j$ is a limit ordinal, i.e. $j\ne i+1$ for any $i<j$, then we take $k_j$ as the least upper bound of $k_i$:
     $$
    k_j=\lim_{i\to j}k_i=\sup_{i<j}k_i
     $$
   (it exists due to Corollary \ref{COR:nepr-F:I->Ord}).

 }\eit
 }\eit
We obtain a transfinite sequence $i\in{\tt Ord}\mapsto k_i\in{\tt Ord}$ with the following properties:

\bit{

 \item[(i)] It is cofinal in ${\tt Ord}$, since $i\le k_i$ for any $i\in{\tt Ord}$.

 \item[(ii)] For $i<j$ we have $F(k_i,k_j)\notin\Iso$, since
$$
i<j\quad\Longrightarrow\quad i+1\le j\quad\Longrightarrow\quad F(k_i,k_j)=F(k_{i+1},k_j)\circ \underbrace{F(k_i,k_{i+1})}_{\scriptsize\begin{matrix}\text{\rotatebox{-90}{$\notin$}}\\ \Iso\end{matrix}}\quad\underset{\tiny\eqref{psi-circ-ph-in-Iso<=>psi-in-Iso-i-ph-in-Iso}}{\Longrightarrow}\quad F(k_i,k_j)\notin\Iso
$$
(we assume here that $F$ is a covariant functor, but for a contravariant one the reasoning is the same).
 }\eit

Now let $S\subseteq{\tt K}$ be a skeleton of ${\tt K}$. For any $i\in{\tt Ord}$ we consider the object $G(i)\in S$ such that
$$
G(i)\cong F(k_i).
$$
Suppose now that $G(i)=G(j)$ for some $i\le j$. Then the morphism $F(k_i,k_j):G(i)\to G(j)$ must coincide with the local identity  $1_{G(i)}=1_{G(j)}$, since the category ${\tt S}$ is a graph, and therefore it cannot have two different colinear morphisms. Thus, $F(k_i,k_j)$ must be an isomorphism, and, by (ii), this is possible only if $i=j$. So we obtain that the map $G:{\tt Ord}\to S$ is injective. On the other hand, it turns the proper class ${\tt Ord}$ into the set $S$, and this is impossible.
\epr

\subsection{Some classes of monomorphisms and epimorphisms}\label{SUBSEC:mono-i-epi}

The widely used in the category theory notions of monomorphism and epimorphism have several variations,
and two of them, the so-called immediate and strong mono- and epimorphisms, will be important for us further.
As the reader will see, we will accentuate the analogy between mono/epimorphisms from the one hand and
strong mono/epimorphisms from the other. In the cases, where due to this analogy the proofs becomes identical
(up to the substitution of the epithet ``strong'' into the proper places, like in the results about
categories $\SMono_X$ and $\SEpi^X$\label{zamech-1}), as well as in the elementary propositions we omit the proofs.

\paragraph{Monomorphisms and epimorphisms.} Recall that a morphism $\ph:X\to Y$ is called
 \bit{
\item[---] a {\it monomorphism}, if any equality
$\ph\circ\alpha=\ph\circ\beta$ implies $\alpha=\beta$;

\item[---] an {\it epimorphism}, if any equality
$\alpha\circ\ph=\beta\circ\ph$ implies $\alpha=\beta$;

\item[---] a {\it bimorphism}, if it is a monomorphism and an epimorphism.
 }\eit

\bex\label{EX:v-grafe-morfizm=bimorfizm} {\it In any graph ${\tt K}$ every morphism is a bimorphism.} Indeed, if $\ph\circ\alpha=\ph\circ\beta$, then, since $\alpha$ and $\beta$ are colinear, they coincide, $\alpha=\beta$. So $\ph$ is a monomorphism. Similarly, it is an epimorphism.
\eex

\bprop\label{PROP:kompoz-monomorf} A composition of two monomorphisms (respectively, two epimorphisms) is a monomorphism (respectively, an epimorphism).
\eprop

\medskip
\centerline{\bf Properties of mono- and epimorphisms:}

\bit{\it

\item[$1^\circ$.]\label{PROP:ph-mu=mono-=>-mu=mono} If $\ph\circ\mu$ is a monomorphism, then $\mu$ is a monomorphism as well.

\item[$2^\circ$.]\label{mu-0-ph-iso-=>-ph-iso} If $\mu\circ\ph$ is an isomorphism, and $\mu$ a monomorphism, then $\mu$ and $\ph$ are isomorphisms.

\item[$3^\circ$.]\label{PROP:e-ph=epi-=>-e=epi} If $\e\circ\ph$ is an epimorphism, then $\e$ is an epimorphism as well.

\item[$4^\circ$.]\label{ph-0-e-iso-=>-ph-iso} If $\ph\circ\e$ is an isomorphism, and $\e$ an epimorphism, then $\ph$ and $\e$ are isomorphisms.

}\eit

By {\it covariant system} (respectively, by {\it contravariant system}) in a category ${\tt K}$ over a partially ordered set $(I,\le)$ we mean arbitrary covariant (respectively, contravariant) functor from $I$ into ${\tt K}$.

\bprop\label{PROP:proj-predel-monomorfizmov} If a covariant system $\{X^j;\iota_i^j\}$ over a directed set $(I,\le)$
 has projective limit $\{X;\pi^j\}$ and all the morphisms $\iota_i^j$ are monomorphisms, then all the morphisms $\pi^j$ are monomorphisms as well.
\eprop
\bpr Let us assume that $I$ is decreasingly directed. Take an index $k\in I$, and let  $Y\overset{\alpha}{\longrightarrow}X$ and $Y\overset{\beta}{\longrightarrow}X$ be two colinear morphisms such that
$$
\pi^k\circ\alpha=\pi^k\circ\beta.
$$
Then for any $j\le  k$ we have:
$$
\underbrace{\iota_j^k\circ\pi^j}_{\pi^k}\circ\ \alpha=\underbrace{\iota_j^k\circ\pi^j}_{\pi^k}\circ\ \beta.
$$
Here $\iota_j^k$ is a monomorphism, so we can cancel it:
$$
\pi^j\circ\alpha=\pi^j\circ\beta,\qquad j\le  k.
$$
Set $\sigma^j=\pi^j\circ\alpha=\pi^j\circ\beta$, then morphisms $Y\overset{\alpha}{\longrightarrow}X$ and $Y\overset{\beta}{\longrightarrow}X$ generate the same cone of the covariant system $\{X^j;\iota_i^j;\ i\le j\le  k\}$:
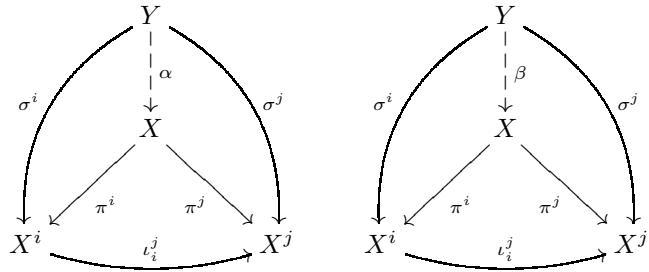
$$
 \xymatrix @R=2.5pc @C=2.5pc
 {
 & Y\ar@/_4ex/[ldd]_{\sigma^i}\ar@{-->}[d]^{\alpha} \ar@/^4ex/[rdd]^{\sigma^j} & \\
 & X\ar[ld]^{\pi^i}\ar[rd]_{\pi^j} & \\
 X^i \ar@/_2ex/[rr]^{\iota_i^j}  & & X^j
 }
 \qquad
  \xymatrix @R=2.5pc @C=2.5pc
 {
 & Y\ar@/_4ex/[ldd]_{\sigma^i}\ar@{-->}[d]^{\beta} \ar@/^4ex/[rdd]^{\sigma^j} & \\
 & X\ar[ld]^{\pi^i}\ar[rd]_{\pi^j} & \\
 X^i \ar@/_2ex/[rr]^{\iota_i^j}  & & X^j
 }
$$
(the projective limit of a covariant system over a cofinal interval $\{j\in I: \ j\le  k\}$ is $X$, the same as over $I$). This implies that $\alpha$ and $\beta$ coincide by the uniqueness of the corresponding arrow in the definition of projective limit:
$$
\alpha=\beta
$$
\epr

The dual proposition is the following:

\bprop\label{PROP:inj-predel-epimorfizmov} If a covariant system $\{X^j;\iota_i^j\}$ over a directed set $(I,\le)$ has injective limit $\{X;\rho_i\}$ and all the morphisms $\iota_i^j$ are epimorphisms, then all the morphisms $\rho_i$ are epimorphisms as well.
\eprop

\brem If the set of indices $I$ is not directed, then the projective (injective) limit of a covariant system of monomorphisms (epimorphisms) over it is not necessarily a cone of monomorphisms (epimorphisms). For example if the order in $I$ is discrete, i.e. $i\le j\ \Leftrightarrow \ i=j$, then the projective limit of any covariant system $\{X^i;\ \iota_i^j\}$ over $I$ is the direct product $\prod_{i\in I}X^i$, where the projections
$$
\prod_{i\in I}X^i\overset{\pi^k}{\longrightarrow}X^k
$$
as a rule are not monomorphisms (although the initial morphisms $\iota_i^i=1_{X^i}$ are monomorphisms). Similarly, injective limit of $\{X^i;\ \iota_i^j\}$ is a coproduct $\coprod_{i\in I}X_i$, and the corresponding injections
$$
X_k\overset{\rho_k}{\longrightarrow}\coprod_{i\in I} X_i
$$
as a rule are not epimorphisms here (although $\iota_i^i=1_{X_i}$ are epimorphisms).
\erem

\paragraph{Immediate monomorphisms and immediate epimorphisms.}\label{SUBSEC:neposr-mono-i-epi}

\bit{

\item[$\bullet$]
We call a {\it factorization} of a morphism $X\overset{\ph}{\longrightarrow}Y$ any its representation as a composition of epimorphism and a monomorphism, i.e. any commutative diagram
\beq\label{DEF:faktorizatsiya}
\begin{diagram}
\node{X}\arrow{se,b}{\e}\arrow[2]{e,t}{\ph}\node[2]{Y}\\
\node[2]{M}\arrow{ne,b}{\mu}
\end{diagram}
\eeq
where $\e$ is an epimorphism, and $\mu$ a monomorphism.

\item[$\bullet$]
A monomorphism $\mu:X\to Y$ is said to be {\it immediate}\label{DEF:immediate-mono}, if in any its factorization $\mu=\mu'\circ\e$ the epimorphism $\e$ is automatically an isomorphism. Note that for a monomorphism $\mu$ in any its factorization $\mu=\mu'\circ\e$ the epimorphism $\e$ is automatically a bimorphism. As a corollary, the condition of being immediate monomorphism for $\mu$ is equivalent to the requirement that in any decomposition $\mu=\mu'\circ\e$, where $\e$ is a bimorphism, and $\mu'$ a monomorphism, the morphism $\e$ must be an isomorphism. It is natural to call a monomorphism $\mu'$ in the factorization $\mu=\mu'\circ\e$ a {\it mediator} of the monomorphism $\mu$, then the epithet ``immediate'' for $\mu$ will mean that there are no non-trivial  mediators for $\mu$ (i.e. mediators, which are not isomorphic to $\mu$ in $\Mono_Y$ -- see below definition \eqref{DEF:varGamma(X)}, here $\varGamma=\Mono$).

\item[$\bullet$]
An epimorphism $\e:X\to Y$ is said to be {\it immediate}\label{DEF:immediate-epi}, if
if $\e$ is an immediate monomorphism in the dual category. In other words, in any factorization $\e=\mu\circ\e'$ the monomorphism $\mu$ must be automatically an isomorphism. Note that for an epimorphism $\e$ in any its factorization $\e=\mu\circ\e'$ the monomorphism $\mu$ is automatically a bimorphism. As a corollary, the condition of being immediate epimorphism for $\e$ is equivalent to the requirement that in any decomposition $\e=\mu\circ\e'$, where  $\mu$ is a bimorphism, and $\e'$ an epimorphism, the morphism $\mu$ must be an isomorphism. It is natural to call an epimorphism $\e'$ in the factorization $\e=\mu\circ\e'$ a {\it mediator} of the epimorphism $\e$, then the epithet ``immediate'' for  $\e$ will mean that there are no non-trivial mediators for $\e$ (i.e. mediators, which are not isomorphic to $\e$ in $\Epi^X$ -- see below definition \eqref{morphism-v-Epi(X)}, here $\varOmega=\Epi$).

}\eit

\brem If in the definition of the immediate monomorphism we omit the requirement that the morphism $\mu'$ in the representation $\mu=\mu'\circ\e$ is a monomorphism (i.e. if we claim only that each epimorphism $\e$ in such a representation must be an  isomorphism), then we obtain exactly the definition of the so-called {\it extremal monomorphism}. Similarly, if in the definition of the immediate epimorphism we omit the requirement that the morphism $\e'$ in the representation $\e=\mu\circ\e'$ is an epimorphism (i.e. if we claim only that each monomorphism $\mu$ in such a representation must be an isomorphism), then we obtain the definition of the {\it extremal epimorphism} \cite[Definition 4.3.2]{Borceux}. Certainly, each extremal monomorphism (respectively, extremal epimorphism) is an immediate monomorphism (respectively, immediate epimorphism). But the reverse implication is not true, and the following example shows this\footnote{This example was suggested to the author by B.~V.~Novikov.}. Consider a monoid $\langle a,b,c \ |\ ac=bc \rangle$ (generated by three elements $a,b,c$ with the equality $ac=bc$) as a category with the one object. In this category
 \bit{
 \item[1)] the morphisms $a,b,c$ are monomorphisms (since they can be canceled in the equalities like $a\cdot P=a\cdot Q$),

 \item[2)] the morphisms $a,b$ are epimorphisms (since they can be canceled in the equalities like $P\cdot a=Q\cdot a$),

 \item[3)] the morphism $c$ is not an epimorphism (since it cannot be canceled in the equality $a\cdot c=b\cdot c$),

 \item[4)] the morphism $ac=bc$ is
 \bit{

 \item[---] a monomorphism (since it can be canceled in the equalities like $ac\cdot P=ac\cdot Q$),

 \item[---] an epimorphism (since it can be canceled in the equalities like $P\cdot ac=Q\cdot ac$),

 \item[---] an immediate epimorphism (since there is only one possibility to write it in the form $(\text{mono})\circ(\text{epi})$, namely, $ac=1\cdot (ac)$, and then the first morphism in this decomposition, i.e. 1, is an isomorphism),

 \item[---] but not an extremal epimorphism (since it can be written in the form $(\text{mono})\circ(\text{...})$, namely, $ac=a\cdot c$, where the first morphism, i.e. $a$, is not an isomorphism).
 }\eit
 }\eit
In addition, the morphism $acac$ is not an immediate epimorphism, since it can be represented as
$$
acac=\underbrace{(ac)}_{\scriptsize\begin{matrix}\text{\rotatebox{90}{$\owns$}} \\ \Mono \end{matrix}}\cdot\underbrace{(ac)}_{\scriptsize\begin{matrix}\text{\rotatebox{90}{$\owns$}} \\ \Epi \end{matrix}}
$$
where the first morphism is not an isomorphism. This shows that {\it a composition of two immediate monomorphisms (respectively, of two immediate epimorphisms) is not necessarily an immediate monomorphism (respectively, an immediate epimorphism)}.
\erem

\bigskip

\begin{samepage}

\centerline{\bf Properties of immediate mono- and epimorphisms:}

\bit{\it

\item[$1^\circ$.]\label{ph-mu=Nmono-=>-mu=Nmono} If $\ph\circ\mu$ is an immediate monomorphism, then $\mu$ is an immediate monomorphism as well.

}\eit

\end{samepage}

\bit{\it

\item[$2^\circ$.]\label{mu-Nmono+epi-=>-mu-iso} If $\mu$ is an immediate monomorphism, and at the same time an epimorphism, then $\mu$ is an isomorphism.

\item[$3^\circ$.]\label{e-ph=Nepi-=>-e=Nepi} If $\e\circ\ph$ is an immediate epimorphism, then $\e$ is an immediate epimorphism as well.

\item[$4^\circ$.]\label{e-Nepi+mono-=>-e-iso} If $\e$ is an immediate epimorphism, and at the same time a monomorphism, then $\e$ is an isomorphism.

}\eit

\paragraph{Strong monomorphisms and strong epimorphisms.}\label{SUBSEC:strogie-mono-i-epi}

The following two definitions are due to M.~Sh.~Tsalenko and E.~G.~Shulgeifer \cite[Chapter 1 \S 7]{Tsalenko-Shulgeifer} and F.~Borceux \cite[4.3]{Borceux}.

\bit{
\item[$\bullet$] A monomorphism $C\overset{\mu}{\longrightarrow}D$ is said to be {\it strong}, if for any epimorphism $A\overset{\e}{\longrightarrow}B$ and for any morphisms $A\overset{\alpha}{\longrightarrow}C$ and $B\overset{\beta}{\longrightarrow}D$ such that
    $\beta\circ\e=\mu\circ\alpha$ there exists (the only possible) morphism $B\overset{\delta}{\longrightarrow}C$, such that the following diagram will be commutative:
\beq\label{razbienie-kvadrata}
\begin{diagram}
\node{A}\arrow{s,l}{\alpha}\arrow{e,t}{\e}\node{B}\arrow{s,r}{\beta}\arrow{sw,t,--}{\delta}
\\
\node{C}\arrow{e,b}{\mu}\node{D}
\end{diagram}
\eeq
\item[$\bullet$] Dually, an epimorphism $A\overset{\e}{\longrightarrow}B$ is said to be {\it strong}, if for any monomorphism  $C\overset{\mu}{\longrightarrow}D$ and for any morphisms $A\overset{\alpha}{\longrightarrow}C$ and $B\overset{\beta}{\longrightarrow}D$ such that
    $\beta\circ\e=\mu\circ\alpha$ there exists (the only possible) morphism $B\overset{\delta}{\longrightarrow}C$, such that diagram \eqref{razbienie-kvadrata} is commutative.
}\eit

\brem\label{REM:o-razbienii-kvadrata} The uniqueness of $\delta$ follows from monomorphity of $\mu$ (or from epimorphity of $\e$): if $\delta'$ is another morphism with the same property, then
$$
\mu\circ\delta=\beta=\mu\circ\delta'\quad\Longrightarrow\quad \delta=\delta'.
$$
Besides this, the commutativity of the upper triangle in \eqref{razbienie-kvadrata} imply the commutativity of the lower one, and vice versa. For example,
 \beq\label{VERH-treug<=>NIZHN-treug}
\alpha=\delta\circ\e\quad\Longrightarrow\quad \beta\circ\underset{\scriptsize\begin{matrix}\text{\rotatebox{90}{$\owns$}}\\ \Epi\end{matrix}}{\e}=\mu\circ\alpha=\mu\circ\delta\circ\underset{\scriptsize\begin{matrix}\text{\rotatebox{90}{$\owns$}}\\ \Epi\end{matrix}}{\e}\quad\Longrightarrow\quad \beta=\mu\circ\delta
 \eeq
\erem

The following propositions are proved in \cite[Proposition 4.3.6]{Borceux}:

\bprop\label{PROP:kompoz-S-monomorf} A composition of two strong monomorphisms (respectively, of two strong epimorphisms) is a strong monomorphism (respectively, a strong epimorphism).
\eprop

\bigskip

\centerline{\bf Properties of strong mono- and epimorphisms:}

\bit{\it

\item[$1^\circ$.]\label{ph-mu=Smono-=>-mu=Smono} If $\ph\circ\mu$ is a strong monomorphism, then $\mu$ is a strong monomorphism as well.

\item[$2^\circ$.]\label{mu-Smono-=>-mu-Nmono} Every strong monomorphism $\mu$ is an immediate monomorphism.

\item[$3^\circ$.]\label{e-ph=Sepi-=>-e=Sepi} If $\e\circ\ph$ is a strong epimorphism, then $\e$ is a strong epimorphism as well.

\item[$4^\circ$.]\label{e-Sepi-=>-e-Nepi} Every strong epimorphism $\e$ is an immediate epimorphism.

}\eit

\bprop\label{PROP:proj-predel-strogih-monomorfizmov}
If in a covariant system $\{X^j;\iota_i^j\}$ over a decreasingly directed set $(I,\le)$ the morphisms  $\iota_i^j$ are strong monomorphisms, then in its projective limit $\{X;\pi^j\}$ the morphisms $\pi^j$ are strong monomorphisms as well.
\eprop
\bpr Take an index $k\in I$. By Proposition \ref{PROP:proj-predel-monomorfizmov}, $\pi^k$ is a monomorphism, so we need only to show that it is strong. Consider a diagram
$$
\begin{diagram}
\node{A}\arrow{e,t}{\e}\arrow{s,l}{\alpha} \node{B}\arrow{s,r}{\beta} \\
\node{X}\arrow{e,b}{\pi^k} \node{X^k}
\end{diagram}
$$
where $\e$ is an epimorphism. For any index $j\le  k$ we can construct a diagram
$$
\begin{diagram}
\node{A}\arrow[2]{e,t}{\e}\arrow{se,t}{\pi^j\circ\alpha}\arrow[2]{s,l}{\alpha} \node[2]{B}\arrow[2]{s,r}{\beta} \\
\node[2]{X^j}\arrow{se,b}{\iota_j^k} \\
\node{X}\arrow[2]{e,b}{\pi^k}\arrow{ne,b}{\pi^j} \node[2]{X^k}
\end{diagram}
$$
and consider the following fragment:
$$
\begin{diagram}
\node{A}\arrow[2]{e,t}{\e}\arrow{se,b}{\pi^j\circ\alpha} \node[2]{B}\arrow[2]{s,r}{\beta} \\
\node[2]{X^j}\arrow{se,b}{\iota_j^k} \\
 \node[3]{X^k}
\end{diagram}
$$
Since $\e$ is an epimorphism, and $\iota_j^k$ is a strong monomorphism, there exists (a unique) morphism $\delta^j$ such that the following diagram is commutative:
$$
\begin{diagram}
\node{A}\arrow[2]{e,t}{\e}\arrow{se,b}{\pi^j\circ\alpha} \node[2]{B}\arrow[2]{s,r}{\beta}\arrow{sw,b,--}{\delta^j} \\
\node[2]{X^j}\arrow{se,b}{\iota_j^k} \\
 \node[3]{X^k}
\end{diagram}
$$
In particular,
$$
\iota_j^k\circ\delta^j=\beta, \qquad j\le  k
$$
As a corollary, if we take a new index $i\le  j$, then for the arising morphisms $\delta^j$ and $\delta^i$ we get
$$
\iota_j^k\circ\delta^j=\beta=\iota^i_k\circ\delta^i= \iota_j^k\circ\iota_i^j\circ\delta^i.
$$
Here $\iota_j^k$ is a monomorphism, so we can cancel it:
$$
\delta^j=\iota_i^j\circ\delta^i.
$$
Thus for any $i\le j\le k$ the following diagram is commutative:
$$
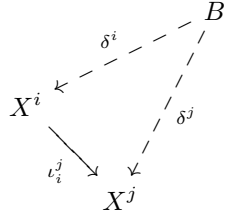
\begin{diagram}
\node[3]{B}\arrow{sww,t,--}{\delta^i}\arrow{ssw,b,--}{\delta^j} \\
\node{X^i}\arrow{se,b}{\iota_i^j} \\
 \node[2]{X^j}
\end{diagram}
$$
(for $j=k$ we have $\delta^k=\beta$).

This means that the system of morphisms $\{\delta^j:B\to X^j;\ j\le  k\}$ is a projective cone of a covariant system $\{\iota_i^j:X^i\to X^j;\ i\le  j\le  k\}$. Hence, there exists a unique morphism $\delta:B\to X$ such that all the following diagrams are commutative:
$$
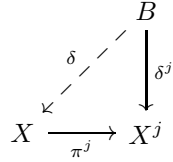
\begin{diagram}
\node[2]{B}\arrow{s,r}{\delta^j}\arrow{sw,t,--}{\delta} \\
\node{X}\arrow{e,b}{\pi^j} \node{X^j}
\end{diagram}
$$
(the limit along a cofinal interval $\{j\in I:\ j\le  k\}$ coincides with the limit along $I$).

In particular, for $j=k$ we get a commutative diagram
$$
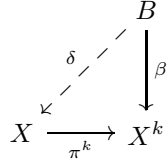
\begin{diagram}
\node[2]{B}\arrow{s,r}{\beta}\arrow{sw,t,--}{\delta} \\
\node{X}\arrow{e,b}{\pi^k} \node{X^k}
\end{diagram}
$$
It implies the following chain:
$$
\beta=\pi^k\circ\delta\quad\Longrightarrow\quad \underbrace{\pi^k}_{\scriptsize\begin{matrix}\text{\rotatebox{90}{$\owns$}}\\ \Mono\end{matrix}}\circ\ \alpha=\beta\circ\e=\underbrace{\pi^k}_{\scriptsize\begin{matrix}\text{\rotatebox{90}{$\owns$}}\\ \Mono\end{matrix}}\circ\ \delta\circ\e\quad\Longrightarrow\quad
\alpha=\delta\circ\e
$$
Thus, the following square is commutative:
$$
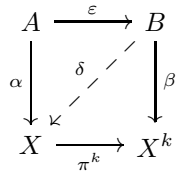
\begin{diagram}
\node{A}\arrow{e,t}{\e}\arrow{s,l}{\alpha} \node{B}\arrow{s,r}{\beta}\arrow{sw,t,--}{\delta} \\
\node{X}\arrow{e,b}{\pi^k} \node{X^k}
\end{diagram}
$$
\epr

The dual proposition is the following:

\bprop\label{PROP:inj-predel-strogih-epimorfizmov}
If in a covariant system $\{X^j;\iota_i^j\}$ over an increasingly directed set $(I,\le)$ the morphisms  $\iota_i^j$ are strong epimorphisms, then in its injective limit $\{X;\rho_i\}$ the morphisms $\rho_i$ are strong epimorphisms as well.
\eprop

\subsection{Categories of monomorphisms and epimorphisms.}

\paragraph{Categories of monomorphisms $\varGamma_X$ and systems of subobjects.}

Let $\varGamma$ be a class of monomorphisms in a category ${\tt K}$, and all local identities belong to it:
$$
\{1_X;\ X\in\Ob({\tt K})\}\subseteq\varGamma\subseteq\Mono({\tt K})
$$
(the key examples are the classes $\varGamma=\Mono$ and $\varGamma=\SMono$). For each object $X$ in ${\tt K}$ let us denote by $\varGamma_X$ the class of all morphisms in $\varGamma$ with $X$ as range:
 \beq\label{DEF:varGamma(X)}
 \varGamma_X=\{\sigma\in\varGamma:\quad \Ran\sigma=X\}.
 \eeq
It is a category, where a morphism $\rho\overset{\varkappa}{\longrightarrow}\sigma$ from an object $\rho\in\varGamma_X$ into an object $\sigma\in\varGamma_X$, i.e. a monomorphism $\rho:A\to X$ into a monomorphism $\sigma:B\to X$, is an arbitrary morphism $\varkappa:A\to B$ in ${\tt K}$ such that the following diagram is commutative:
\beq\label{morphism-v-Mono(X)}
\xymatrix @R=1pc @C=2pc
{
A\ar[rd]^{\rho}\ar[dd]_{\varkappa} &   \\
  & X \\
B\ar[ru]_{\sigma} &
}
\eeq
Actually, this diagram in the initial category ${\tt K}$ can be considered as a morphism $\rho\overset{\varkappa}{\longrightarrow}\sigma$ in the category  $\varGamma_X$. A composition of such morphisms $\rho\overset{\varkappa}{\longrightarrow}\sigma$ and $\sigma\overset{\lambda}{\longrightarrow}\tau$, i.e. of diagrams
$$
\xymatrix @R=1pc @C=2pc
{
A\ar[rd]^{\rho}\ar[dd]_{\varkappa} &   \\
  & X \\
B\ar[ru]_{\sigma} &
}
\qquad
\xymatrix @R=1pc @C=2pc
{
B\ar[rd]^{\sigma}\ar[dd]_{\lambda} &   \\
  & X \\
C\ar[ru]_{\tau} &
}
$$
is a morphism $\rho\overset{\lambda\circ\varkappa}{\longrightarrow}\tau$, i.e. a diagram
$$
\xymatrix @R=1pc @C=2pc
{
A\ar[rd]^{\rho}\ar[dd]_{\lambda\circ\varkappa} &   \\
  & X \\
C\ar[ru]_{\tau} &
}
$$
One can conceive it as a result of splicing of  the initial diagrams along the common edge $\sigma$, adding the arrow of composition $\varkappa\circ\lambda$, and then throwing away the vertex $B$ together with all its incidental edges:
$$
\xymatrix % @R=1pc @C=2pc
{
A\ar@/^2ex/[rrd]^{\rho}\ar[dd]_{\lambda\circ\varkappa}\ar@{-->}[rd]_{\varkappa} & &  \\
  & B\ar@{-->}[r]_{\sigma}\ar@{-->}[ld]_{\lambda} &  X \\
  C\ar@/_2ex/[rru]_{\tau} & &
}
$$
Of course, local identities in $\varGamma_X$ are diagrams of the form
$$
\xymatrix @R=1pc @C=2pc
{
A\ar[rd]^{\rho}\ar[dd]_{1_A} &   \\
  & X \\
A\ar[ru]_{\rho} &
}
$$

\brem
The composition of morphisms in $\varGamma_X$ can be defined in two ways. In our definition this operation is connected with the composition in ${\tt K}$ through the following identity:
$$
\lambda\kern-3pt\underset{\varGamma_X}{\circ}\kern-3pt\varkappa=\lambda\underset{\tt K}{\circ}\varkappa.
$$
\erem

\btm\label{card-Mono(a,b)-le-1}
For any object $X$ the category $\varGamma_X$ is a graph.
\etm
\bpr We should verify that for any two objects $\rho:A\to X$ and $\sigma:B\to X$ there exist at most one morphism  $\rho\overset{\varkappa}{\longrightarrow}\sigma$. Indeed, a morphism $\varkappa$ in diagram \eqref{morphism-v-Mono(X)} is unique, since the monomorphity of $\sigma$ gives the following implication: $\sigma\circ\varkappa=\rho=\sigma\circ\varkappa'$ $\Longrightarrow$ $\varkappa=\varkappa'$.
\epr

\brem\label{REM:svoistva-Mono(X)}
By Example \ref{EX:v-grafe-morfizm=bimorfizm} this means that {\it in the category $\varGamma_X$ all morphisms are bimorphisms}. The connection between the properties of a morphism $\rho\overset{\varkappa}{\longrightarrow}\sigma$ in $\varGamma_X$ and the properties of the same morphism $\varkappa:A\to B$ in the initial category ${\tt K}$, is expressed in the following observations:

\bit{\it

\item[---] \label{PROP:mono-v-Mono(X)}
every morphism $\rho\overset{\varkappa}{\longrightarrow}\sigma$ in $\varGamma_X$ is a monomorphism in ${\tt K}$,

\item[---] \label{PROP:iso-v-Mono(X)}
a morphism $\rho\overset{\varkappa}{\longrightarrow}\sigma$ in $\varGamma_X$ is an isomorphism in $\varGamma_X$ $\Longleftrightarrow$ $\varkappa$ is an isomorphism in ${\tt K}$.
}\eit
\erem
\bpr
1. A morphism $\varkappa$ in \eqref{morphism-v-Mono(X)} must be a monomorphism due to Property $1^\circ$ at the p.\pageref{PROP:ph-mu=mono-=>-mu=mono}, since $\sigma\circ\varkappa$ is a monomorphism.

2. If a morphism $\varkappa:A\to B$ in \eqref{morphism-v-Mono(X)} is an isomorphism in ${\tt K}$, then we can put $\lambda=\varkappa^{-1}:A\gets B$, and the diagrams
\beq\label{iso-v-Mono(X)}
\xymatrix % @R=1pc @C=2pc
{
A\ar[dd]_{1_A}\ar@/^2ex/[rrd]^{\rho}\ar@{-->}[rd]_{\varkappa} & &  \\
  & B\ar@{-->}[r]_{\sigma}\ar@{-->}[ld]_{\lambda} &  X \\
  A\ar@/_2ex/[rru]_{\rho} & &
}
\qquad
\xymatrix % @R=1pc @C=2pc
{
B\ar[dd]_{1_B}\ar@/^2ex/[rrd]^{\sigma}\ar@{-->}[rd]_{\lambda} & &  \\
  & A\ar@{-->}[r]_{\rho}\ar@{-->}[ld]_{\varkappa} &  X \\
B\ar@/_2ex/[rru]_{\sigma} & &
}
\eeq
will be commutative, since $\rho$ and $\sigma$ are monomorphisms. They mean that the morphisms $\rho\overset{\varkappa}{\longrightarrow}\sigma$ and $\sigma\overset{\lambda}{\longrightarrow}\rho$ in $\varGamma_X$ are inverse to each other. Conversely, if  morphisms $\rho\overset{\varkappa}{\longrightarrow}\sigma$ and $\sigma\overset{\lambda}{\longrightarrow}\rho$ are inverse to each other in $\varGamma_X$, then this means that diagrams \eqref{iso-v-Mono(X)} are commutative. Hence, morphisms $\varkappa$ and $\lambda$ are inverse to each other in ${\tt K}$, and thus $\varkappa$ must be an isomorphism in $\tt K$.
\epr

It is convenient to introduce a special notation, $\to$, for the pre-order in $\varGamma_X$:
 \beq\label{DEF:to-in-Mono(X)}
\rho\to\sigma\Longleftrightarrow\quad \exists \varkappa\in\Mor({\tt K})\quad
\rho=\sigma\circ\varkappa.
 \eeq
Here the morphism $\varkappa$, if it exists, must be unique, and besides this it is a monomorphism (this follows from the fact that $\sigma$ is a monomorphism). As a corollary, there is an operation, which to any pair of morphisms $\rho,\sigma\in\varGamma_X$ with the property $\rho\to\sigma$ assigns the morphism
$\varkappa=\varkappa^\sigma_\rho$ in \eqref{DEF:to-in-Mono(X)}:
 \beq\label{DEF:to-in-Mono(X)-*}
\rho=\sigma\circ\varkappa^\sigma_\rho.
 \eeq
If $\rho\to\sigma\to\tau$, then the chain
$$
\tau\circ\varkappa^\tau_\rho=\rho=\sigma\circ\varkappa^\sigma_\rho=
\tau\circ\varkappa^\tau_\sigma\circ\varkappa^\sigma_\rho,
$$
implies, due to monomorphy of $\tau$, the equality
\beq\label{varkappa^gamma_alpha=varkappa^gamma_beta-circ-varkappa^beta_alpha}
\varkappa^\tau_\rho=\varkappa^\tau_\sigma\circ\varkappa^\sigma_\rho.
 \eeq

\bit{

\item[$\bullet$] A {\it system of subobjects of the class $\varGamma$} in an object $X$ of a category ${\tt K}$ is an arbitrary skeleton $S$ of the category $\varGamma_X$, such that the morphism $1_X$ belongs to $S$. In other words, a subclass $S$ in $\varGamma_X$ is a system of subobjects in $X$, if
\bit{

\item[(a)]\label{DEF:toch-sist-podobj-a} the local identity of $X$ belongs to $S$:
$$
1_X\in S,
$$

\item[(b)]\label{DEF:toch-sist-podobj-b} every monomorphism $\mu\in\varGamma_X$ has an isomorphic monomorphism in the class $S$:
$$
\forall\mu\in\varGamma_X\qquad\exists\sigma\in S\qquad \mu\cong\sigma.
$$

\item[(c)]\label{DEF:toch-sist-podobj-c} in $S$ an isomorphism (in the sense of category $\varGamma_X$) is equivalent to the identity:
$$
\forall\sigma,\tau\in S\qquad \Big( \sigma\cong\tau\quad\Longleftrightarrow\quad \sigma=\tau \Big)
$$
}\eit
Due to property $1^\circ$ on page \pageref{1^0:skelet-sushestvuet}, such a class $S$ always exists.
The elements of $S$ are called {\it subobjects} of $X$ (of the class $\varGamma$). The class $S$ is endowed with the structure of  a full subcategory in $\varGamma_X$.
}\eit

\btm\label{S-chast-upor-klass} Any system of subobjects $S$ of an object $X$ is a partially ordered class.
\etm
\bpr
Let subobjects $\rho\in S$ and $\sigma\in S$ have two mutually inverse morphisms $\varkappa:A\gets B$ and $\lambda:A\to B$, i.e.
$$
\rho=\sigma\circ\varkappa,\qquad \sigma=\rho\circ\lambda.
$$
Then
$$
\rho\circ\lambda\circ\varkappa=\rho=\rho\circ 1_A,\qquad \sigma\circ\varkappa\circ\lambda=\sigma=\sigma\circ 1_B,
$$
and, since $\rho$ and $\sigma$ are monomorphisms in ${\tt K}$, one can cancel them:
$$
\lambda\circ\varkappa=1_A,\qquad \varkappa\circ\lambda=1_B,
$$
Thus, $\varkappa$ and $\lambda$ are isomorphisms. We obtain that $\rho\cong\sigma$, and by property (c), $\rho=\sigma$.
\epr

\btm\label{PROP:dost-mnozh-v-klasse-podobjektov} If $S$ is a system of subobjects in $X$, then for any subobject $\sigma\in S$, $\sigma:Y\to X$, the class of monomorphisms
$$
A=\{\alpha\in\varGamma_Y:\ \sigma\circ\alpha\in S\}
$$
is a system of subobjects in $Y$. If in addition $S$ is a set, then $A$ is a set as well.
\etm
\bpr
1. Property (a) is obvious: since $\sigma\circ 1_Y=\sigma\in S$, we have that $1_Y\in A$.

2. Property (b). Let $\beta:B\to Y$ be a monomorphism. The composition $\sigma\circ\beta:B\to X$ is a monomorphism from $\varGamma_X$, and since $S$ is a system of subobjects in $X$, there must exist $\tau\in S$ such that
$$
\tau\cong\sigma\circ\beta.
$$
This means that
$$
\tau=\sigma\circ\beta\circ\iota
$$
for some isomorphism $\iota$. Now we get that the monomorphism $\alpha=\beta\circ\iota$ is isomorphic to $\beta$
$$
\alpha\cong\beta
$$
and lies in $A$, since $\sigma\circ\alpha=\tau\in S$.

3. Property (c). Let $\alpha,\beta\in A$ be two isomorphic monomorphisms, i.e.
$$
\alpha=\beta\circ\iota
$$
for some isomorphism $\iota$. Then, first, the morphisms $\sigma\circ\alpha$ and $\sigma\circ\beta$ are isomorphic as well, since
$$
\sigma\circ\alpha=\sigma\circ\beta\circ\iota.
$$
And, second, they lay in $S$, since $\alpha$ and $\beta$ lay in $A$. But $S$ satisfies (c), hence the morphisms $\sigma\circ\alpha$ and $\sigma\circ\beta$ coincide:
$$
\sigma\circ\alpha=\sigma\circ\beta.
$$
In addition $\sigma$ is a monomorphism, so we have $\alpha=\beta$.

4. It remains to check that if $S$ is a set, then $A$ is a set as well. This follows from the fact that the map $\alpha\in A\mapsto \sigma\circ\alpha\in S$ is injective. Indeed, if for some $\alpha,\alpha'\in A$ we have
$$
\sigma\circ\alpha=\sigma\circ\alpha',
$$
then, since $\sigma$ is a monomorphism, we have $\alpha=\alpha'$.
\epr

 \bit{
\item[$\bullet$] We say that a category ${\tt K}$ is {\it well-powered in the class $\varGamma$}, if each object $X$ has a system of subobjects $S$ of the class $\varGamma$, which is a set (i.e. not a proper class); in other words, each category $\varGamma_X$ must be a skeletally small graph.
 }\eit

\bex The standard categories frequently used as examples, like the category of sets, groups, vector spaces, algebras (over a given field), topological spaces, topological vector spaces, topological algebras, etc., are, obviously, well-powered in the class $\Mono$.
\eex

\btm\label{TH:o-lok-malosti-v-podobjektah}
If a category ${\tt K}$ is well-powered in a class $\varGamma$, then there is a map $X\mapsto S_X$ which assigns to each object $X$ in ${\tt K}$ its system of subobjects $S_X$ of the class $\varGamma$ (and $S_X$ is a set).
\etm
\bpr
We must use here the theorem that the class of all sets can be well-ordered \cite[V, 4.1]{Levy}: this allows to assign to each $X$ the system of subobjects $S$, which is minimal with respect to this well-ordering.
\epr

\paragraph{Categories of epimorphisms $\varOmega^X$ and systems of quotient objects.}

Let $\varOmega$ be a class of epimorphisms in a category ${\tt K}$, and all local identities belong to it:
$$
\{1_X;\ X\in\Ob({\tt K})\}\subseteq\varOmega\subseteq\Epi({\tt K})
$$
(the key examples are the classes $\varOmega=\Epi$ and $\varOmega=\SEpi$). For each object $X$ in ${\tt K}$ we denote by $\varOmega^X$ the class of all morphisms in $\varOmega$ with the domain $X$:
 \beq\label{DEF:varOmega(X)}
 \varOmega^X=\{\sigma\in\varOmega:\quad \Dom\sigma=X\}.
 \eeq
This class forms a category where a morphism $\rho\overset{\varkappa}{\longrightarrow}\sigma$ from an object $\rho\in\varOmega^X$ into an object $\sigma\in\varOmega^X$, i.e. from an epimorphism $\rho:X\to A$ into an epimorphism $\sigma:X\to B$, is an arbitrary morphism $\varkappa:A\to B$ in ${\tt K}$ such that the following diagram is commutative
\beq\label{morphism-v-Epi(X)}
\xymatrix @R=1pc @C=2pc
{
 & A\ar[dd]^{\varkappa}    \\
X\ar[ru]^{\rho}\ar[rd]_{\sigma}  &  \\
 & B
}
\eeq
Actually, this diagram in the initial category ${\tt K}$ can be considered as a morphism $\rho\overset{\varkappa}{\longrightarrow}\sigma$ in $\varOmega^X$. A composition of two such morphisms $\rho\overset{\varkappa}{\longrightarrow}\sigma$ and $\sigma\overset{\lambda}{\longrightarrow}\tau$, i.e. diagrams $$
\xymatrix @R=1pc @C=2pc
{
 & A\ar[dd]^{\varkappa}    \\
X\ar[ru]^{\rho}\ar[rd]_{\sigma}  &  \\
 & B
}
\qquad
\xymatrix @R=1pc @C=2pc
{
 & B\ar[dd]^{\lambda}    \\
X\ar[ru]^{\sigma}\ar[rd]_{\tau}  &  \\
 & C
}
$$
is a morphism $\rho\overset{\lambda\circ\varkappa}{\longrightarrow}\tau$, i.e. a diagram
$$
\xymatrix @R=1pc @C=2pc
{
 & A\ar[dd]^{\lambda\circ\varkappa}    \\
X\ar[ru]^{\rho}\ar[rd]_{\tau}  &  \\
 & C
}
$$
One can conceive it as a result of splicing of  the initial diagrams along the common edge $\sigma$, adding the arrow of composition $\lambda\circ\varkappa$, and then throwing away the vertex $B$ together with all its incidental edges:
$$
\xymatrix % @R=1pc @C=2pc
{
 & & A\ar[dd]^{\lambda\circ\varkappa}\ar@{-->}[ld]^{\varkappa}    \\
X\ar@/^2ex/[rru]^{\rho}\ar@{-->}[r]_{\sigma}\ar@/_2ex/[rrd]_{\tau}  & B\ar@{-->}[rd]^{\lambda} &  \\
 & & C
}
$$
Of course, local identities in $\varOmega^X$ are diagrams of the form
$$
\xymatrix @R=1pc @C=2pc
{
 & A\ar[dd]^{1_A}    \\
X\ar[ru]^{\rho}\ar[rd]_{\sigma}  &  \\
 & A
}
$$

\brem
The composition of morphisms in $\varOmega^X$ can be defined in two ways. In our definition this operation is connected with the composition in ${\tt K}$ through the following identity:
$$
\lambda\kern-3pt\underset{\varOmega^X}{\circ}\kern-3pt\varkappa=\lambda\underset{\tt K}{\circ}\varkappa.
$$
\erem

By analogy with $\varGamma_X$ the following properties of $\varOmega^X$ are proved.

\btm\label{card-Epi(a,b)-le-1}
For any object $X$ the category $\varOmega^X$ is a graph.
\etm

\brem\label{REM:svoistva-Epi(X)}
By Example \ref{EX:v-grafe-morfizm=bimorfizm} this means that {\it in the category $\varOmega^X$ all the morphisms are bimorphisms}. The connection between the properties of a morphism $\rho\overset{\varkappa}{\longrightarrow}\sigma$ in $\varOmega^X$ and the properties of the same morphism $\varkappa:A\to B$ in the initial category ${\tt K}$, is expressed in the following observations:

\bit{\it

\item[---] \label{PROP:epi-v-Epi(X)}
every morphism $\rho\overset{\varkappa}{\longrightarrow}\sigma$ in $\varOmega^X$ is an epimorphism in ${\tt K}$,

\item[---] \label{PROP:iso-v-Epi(X)}
a morphism $\rho\overset{\varkappa}{\longrightarrow}\sigma$ in $\varOmega^X$ is an isomorphism in $\varOmega^X$ $\Longleftrightarrow$ $\varkappa$ is an isomorphism in ${\tt K}$.
}\eit
\erem

It is convenient to introduce a special notation, $\to$, for the pre-order in $\varOmega^X$:
 \beq\label{DEF:le-in-F_X}
\rho\to\sigma\Longleftrightarrow\quad \exists \iota\in\Mor({\tt K})\quad
\sigma=\iota\circ\rho.
 \eeq
Here the morphism $\iota$, if it exists, must be unique, and besides this it is an epimorphism (since $\rho$ and $\sigma$ are epimorphisms). As a corollary, there is an operation, which to each pair of morphisms $\rho,\sigma\in\varOmega^X$ with the property $\rho\to\sigma$ assigns the morphism $\iota=\iota^\sigma_\rho$ in
\eqref{DEF:le-in-F_X}:
 \beq\label{DEF:le-in-F_X-*}
\sigma=\iota^\sigma_\rho\circ\rho.
 \eeq
If $\pi\to\rho\to\sigma$, then the chain
$$
\iota^\sigma_\pi\circ\pi=\sigma=\iota^\sigma_\rho\circ\rho=\iota^\sigma_\rho\circ\iota^\rho_\pi\circ\pi,
$$
implies by epimorphy of $\pi$ the equality
\beq\label{iota_rho^tau=iota_rho^sigma-circ-iota_sigma^tau}
\iota^\sigma_\pi=\iota^\sigma_\rho\circ\iota^\rho_\pi.
 \eeq

\bit{

\item[$\bullet$] A {\it system of quotient objects of the class $\varOmega$} on an object $X$ in a category ${\tt K}$ is an arbitrary skeleton $Q$ of the category $\varOmega^X$, such that $1_X$ belongs to $Q$. In other words, a subclass $Q$ in $\varOmega^X$ is called a system of quotient objects on $X$, if
\bit{

\item[(a)] the local identity of $X$ belongs to $Q$:
$$
1_X\in Q,
$$

\item[(b)] every epimorphism $\e\in\varOmega^X$ has an isomorphic epimorphism in $Q$:
$$
\forall\e\in\varOmega^X\qquad\exists\pi\in Q\qquad \e\cong\pi,
$$

\item[(c)] in $Q$ an isomorphism (in the sense of category $\varOmega^X$) is equivalent to the identity:
$$
\forall\pi,\rho\in Q\qquad \Big( \pi\cong\rho\quad\Longleftrightarrow\quad \pi=\rho \Big)
$$
}\eit
By Property $1^\circ$ on page \pageref{1^0:skelet-sushestvuet} this class $Q$ always exists. The elements of $Q$ are called {\it quotient objects} on $X$. The class $Q$ is endowed with the structure of a full subcategory in $\varOmega^X$.
}\eit

By analogy with Theorems \ref{S-chast-upor-klass} and \ref{PROP:dost-mnozh-v-klasse-podobjektov} we have

\btm\label{Q-chast-upor-klass} Any system $Q$ of quotient objects of an object $X$ is a partially ordered class.
\etm

\btm\label{PROP:dost-mnozh-v-klasse-faktor-objektov} If $Q$ is a system of quotient objects of an object $X$, then for any quotient object $\pi\in Q$, $\pi:X\to Y$, the class of epimorphisms
$$
A=\{\alpha\in\varOmega^Y:\ \alpha\circ\pi\in Q\}
$$
is a system of quotient objects on $Y$. If in addition $Q$ is a set, then $A$ is a set as well.
\etm

 \bit{
\item[$\bullet$] We say that a category ${\tt K}$ is {\it co-well-powered in the class $\varOmega$}, if each object  $X$ has a system of quotient objects $Q$ of the class $\varOmega$, which is a set (i.e. not a proper class); in other words, each category $\varOmega^X$ must be a skeletally small graph.
 }\eit

\bex Among the standard categories -- the category of sets, groups, vector spaces, algebras over a given field, topological spaces, topological vector spaces, topological algebras -- some are co-well-powered in the class $\Epi$, but sometimes this is not easy to prove (see \cite{Adamek-Rosicky}). In contrast to this the co-well-poweredness in the class $\SEpi$ is verified much easier. \eex

By analogy with Theorem \ref{TH:o-lok-malosti-v-podobjektah} the following fact is proved:

 \btm\label{TH:o-lok-malosti-v-faktor-objektah}
If a category ${\tt K}$ is co-well-powered in the class $\varOmega$, then there exists a map $X\mapsto Q_X$ which assigns to any object $X$ in ${\tt K}$ a system of its quotient-objects $Q_X$ of the class $\varOmega$ (and $Q_X$ is a set).
 \etm

\subsection{Nodal decomposition}

\paragraph{Strong decompositions.}\label{paragraph:strogie-razlozhenija}
 \bit{
 \item[$\bullet$]
A representation of a morphism $\ph$ into a composition of three morphisms
$$
\ph=\iota\circ\rho\circ\gamma,
$$
where $\iota$ is a strong monomorphism, and $\gamma$ a strong epimorphism, will be called a {\it strong decomposition} of $\ph$.
 }\eit

\btm\label{TH:glav-sv-strog-razlozh} If $\ph=\iota\circ\rho\circ\gamma$ is a strong decomposition of $\ph$, then for any other decomposition
$$
\ph=\mu\circ\e
$$
 \bit{
 \item[---] the epimorphity of $\e$ implies the existence of a unique morphism $\mu'$ such that the following diagram is commutative:
\beq\label{strog-razlozh-mu}
\begin{diagram}
\node{X}\arrow[2]{s,l}{\gamma}\arrow{se,b}{\e}\arrow[2]{e,t}{\ph}\node[2]{Y}\\
\node[2]{M}\arrow{ne,b}{\mu}\arrow{se,t,--}{\mu'}\\
\node{X'}\arrow[2]{e,t}{\rho}\node[2]{Y'}\arrow[2]{n,r}{\iota}
\end{diagram}
\eeq
(in this case if $\mu$ is a monomorphism, then $\mu'$ is a monomorphism as well),

 \item[---] the monomorphity of $\mu$ implies the existence of a unique morphism $\e'$ such that the following diagram is commutative:
\beq\label{strog-razlozh-e}
\begin{diagram}
\node{X}\arrow[2]{s,l}{\gamma}\arrow{se,b}{\e}\arrow[2]{e,t}{\ph}\node[2]{Y}\\
\node[2]{M}\arrow{ne,b}{\mu}\\
\node{X'}\arrow[2]{e,t}{\rho}\arrow{ne,t,--}{\e'}\node[2]{Y'}\arrow[2]{n,r}{\iota}
\end{diagram}
\eeq
(in this case if $\e$ is an epimorphism, then $\e'$ is an epimorphism as well).
 }\eit
\etm
\bpr
Let $\e$ be an epimorphism. Consider the diagram
$$
\begin{diagram}
\node{X}\arrow[2]{s,l}{\gamma}\arrow{se,b}{\e}\node[2]{Y}\\
\node[2]{M}\arrow{ne,b}{\mu}\\
\node{X'}\arrow[2]{e,t}{\rho}\node[2]{Y'}\arrow[2]{n,r}{\iota}
\end{diagram}
$$
and transform it into the following one:
$$
\begin{diagram}
\node{X}\arrow{s,l}{\rho\circ\gamma}\arrow{e,t}{\e}\node{M}\arrow{s,r}{\mu}\\
\node{Y'}\arrow{e,b}{\iota}\node{Y}
\end{diagram}
$$
Here $\e$ is an epimorphism, and $\iota$ a strong monomorphism, hence there exists a (unique) morphism $\mu'$ such that
$$
\begin{diagram}
\node{X}\arrow{s,l}{\rho\circ\gamma}\arrow{e,t}{\e}\node{M}\arrow{s,r}{\mu}\arrow{sw,r,--}{\mu'}\\
\node{Y'}\arrow{e,b}{\iota}\node{Y}
\end{diagram}
$$
This is the morphism for \eqref{strog-razlozh-mu}. By Property $1^0$ on page \pageref{PROP:ph-mu=mono-=>-mu=mono}, if in addition $\mu=\iota\circ\mu'$ is a monomorphism, then $\mu'$ is also a monomorphism. The second case is dual.
\epr

Suppose we have two strong decompositions $\ph=\iota\circ\rho\circ\gamma$ and $\ph=\iota'\circ\rho'\circ\gamma'$ of one morphism $\ph$
$$
\xymatrix %@R=2.5pc @C=4.0pc
{
X\ar[r]^{\ph}\ar[d]_{\gamma} & Y \\
P\ar[r]_{\rho} & Q \ar[u]_{\iota}
}
\qquad
\xymatrix %@R=2.5pc @C=4.0pc
{
X\ar[r]^{\ph}\ar[d]_{\gamma'} & Y \\
P'\ar[r]_{\rho'} & Q' \ar[u]_{\iota'}
}
$$
If there exist (necessarily unique due to Theorem \ref{TH:glav-sv-strog-razlozh}) morphisms $\sigma:P\to P'$ and $\tau:Q'\to Q$ such that the following diagram is commutative
\beq\label{otnoshenie-mezhdu-stand-razlozh}
\xymatrix %@R=2.5pc @C=4.0pc
{
X\ar[rrr]^{\ph}\ar[rd]^{\gamma'}\ar[dd]_{\gamma} & & & Y \\
& P'\ar[r]^{\rho'} & Q'\ar[ru]^{\iota'}\ar@{-->}[rd]_{\tau} & \\
P\ar[rrr]_{\rho}\ar@{-->}[ru]_{\sigma} & & & Q \ar[uu]_{\iota}
}
\eeq
then we say that the strong decomposition $\ph=\iota\circ\rho\circ\gamma$ {\it is subordinated} to the strong decomposition $\ph=\iota'\circ\rho'\circ\gamma'$, and we write in this case
$$
 (\iota,\rho,\gamma)\le (\iota',\rho',\gamma').
$$
If in addition $\sigma$ and $\tau$ are isomorphisms here, then we say that the decompositions $\ph=\iota\circ\rho\circ\gamma$ and $\ph=\iota'\circ\rho'\circ\gamma'$ are {\it isomorphic}, and we write
$$
 (\iota,\rho,\gamma)\cong (\iota',\rho',\gamma').
$$

\bprop The two-sided inequality
$$
 (\iota,\rho,\gamma)\le (\iota',\rho',\gamma')\le (\iota,\rho,\gamma)
$$
is equivalent to the isomorphism of strong decompositions
$$
 (\iota,\rho,\gamma)\cong (\iota',\rho',\gamma').
$$
\eprop
\bpr
The first inequality here implies the existence of (unique) dotted arrows in \eqref{otnoshenie-mezhdu-stand-razlozh}, and the second one means that the reverse arrows exist as well (and again are unique). In addition the epimorphy of $\gamma$ and $\gamma'$ imply that $\sigma$ with its reverse arrow are mutually reverse isomorphisms, while the monomorphy of $\iota$ and $\iota'$ imply that the same is true for $\tau$ with its reverse arrow.
\epr

\paragraph{Nodal decomposition.}\label{paragraph:uzlovoe-razlozhenie}
If in a strong decomposition $\ph=\iota'\circ\rho'\circ\gamma'$ the middle morphism $\rho'$ is a bimorphism, then we call this a {\it nodal decomposition}. We say also that ${\tt K}$ is a {\it category with a nodal decomposition}, if every morphism $\ph$ in ${\tt K}$ has a nodal decomposition.

\bprop\label{PROP:uzlov-razlozh-podchin-ostalnye}
Each nodal decomposition $\ph=\iota'\circ\rho'\circ\gamma'$ subordinates each strong decomposition $\ph=\iota\circ\rho\circ\gamma$:
$$
(\iota,\rho,\gamma)\le (\iota',\rho',\gamma').
$$
As a corollary, a nodal decomposition is unique up to isomorphism.
\eprop
\bpr Let $\ph=\iota\circ\rho\circ\gamma$ be a strong decomposition. If we transform the diagram
\beq\label{proof:predelnoe-razlozhenie}
\xymatrix %@R=2.5pc @C=4.0pc
{
X\ar[rd]^{\gamma'}\ar[dd]_{\gamma} & & & Y \\
& P'\ar[r]^{\rho'} & Q'\ar[ru]^{\iota'} & \\
P\ar[rrr]_{\rho} & & & Q \ar[uu]_{\iota}
}
\eeq
into the diagram
$$
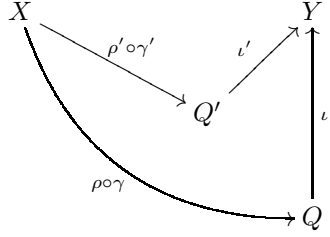
\xymatrix %@R=2.5pc @C=4.0pc
{
X\ar[rrd]^{\rho'\circ\gamma'}\ar@/_6ex/[rrrdd]_{\rho\circ\gamma} & & & Y \\
&  & Q'\ar[ru]^{\iota'} & \\
 & & & Q \ar[uu]_{\iota}
}
$$
then one can recognize here a quadrangle of the form \eqref{razbienie-kvadrata}, since $\iota$ is a strong monomorphism here, and $\rho'\circ\gamma'$ an epimorphism (as a composition of an epimorphism $\gamma'$ and a bimorphism $\rho'$). Hence, there is a unique morphism $\tau$ such that
$$
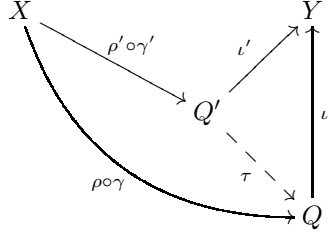
\xymatrix %@R=2.5pc @C=4.0pc
{
X\ar[rrd]^{\rho'\circ\gamma'}\ar@/_6ex/[rrrdd]_{\rho\circ\gamma} & & & Y \\
&  & Q'\ar[ru]^{\iota'}\ar@{-->}[rd]_{\tau} & \\
 & & & Q \ar[uu]_{\iota}
}
$$
Similarly, one can transform diagram \eqref{proof:predelnoe-razlozhenie} into
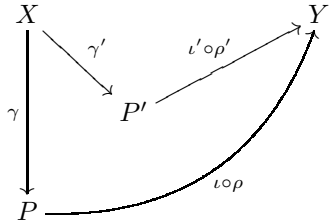
$$
\xymatrix %@R=2.5pc @C=4.0pc
{
X\ar[rd]^{\gamma'}\ar[dd]_{\gamma} & & & Y \\
& P'\ar[rru]^{\iota'\circ\rho'} &  & \\
P\ar@/_6ex/[rrruu]_{\iota\circ\rho} & & &
}
$$
and this again is a quadrangle of the form \eqref{razbienie-kvadrata}, since $\gamma$ is a strong epimorphism here, and $\iota'\circ\rho'$ a monomorphism (as a composition of a bimorphism $\rho'$ and a monomorphism $\iota'$). Hence, there exists a unique morphism $\sigma$ such that
$$
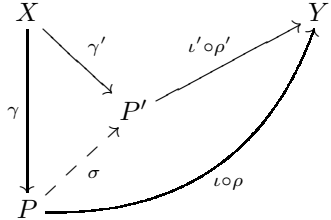
\xymatrix %@R=2.5pc @C=4.0pc
{
X\ar[rd]^{\gamma'}\ar[dd]_{\gamma} & & & Y \\
& P'\ar[rru]^{\iota'\circ\rho'} &  & \\
P\ar@{-->}[ru]_{\sigma}\ar@/_6ex/[rrruu]_{\iota\circ\rho} & & &
}
$$
These two morphisms together give diagram \eqref{otnoshenie-mezhdu-stand-razlozh}.
\epr

\bit{
\item[$\bullet$]
From the uniqueness (up to isomorphism) of the nodal decomposition $\ph=\iota'\circ\rho'\circ\gamma'$ it follows that one can assign notations to its components. We will further depict a nodal decomposition of a morphism $\ph:X\to Y$ as a diagram
\beq\label{DEF:oboznacheniya-dlya-uzlov-razlozh}
\begin{diagram}
\node{X}\arrow{s,l}{\coim_\infty\ph}\arrow{e,t}{\ph}\node{Y} \\
\node{\Coim_\infty\ph}\arrow{e,t}{\red_\infty\ph}\node{\Im_\infty\ph}\arrow{n,r}{\im_\infty\ph}
\end{diagram}
\eeq
(where elements are defined up to isomorphisms). The proof of Theorem \ref{TH:sush-uzlov-razlozh-v-polnoi-lok-maloi-kateg} below and Remark \ref{REM:struktura-uzlov-razlozh-v-kateg-s-nulem} justify these notations, since they show that $\coim_\infty$, $\red_\infty$ and $\im_\infty$ can be conceived as a sort of ``transfinite induction'' of the usual operation $\coim$, $\red$ and $\im$ in preabelian categories:
\begin{align*}
& \coim_\infty=\lim_{n\to\infty}\underbrace{\coim\circ\coim\circ...\circ\coim}_{\text{$n$ multipliers}} \\
& \red_\infty=\lim_{n\to\infty}\underbrace{\red\circ\red\circ...\circ\red}_{\text{$n$ multipliers}} \\
& \im_\infty=\lim_{n\to\infty}\underbrace{\im\circ\im\circ...\circ\im}_{\text{$n$ multipliers}}
\end{align*}
We will call
\bit{
\item[---] $\im_\infty\ph$ a {\it nodal image},

\item[---] $\red_\infty\ph$ a {\it nodal reduced part},

\item[---] $\coim_\infty\ph$  a {\it nodal coimage}
 }\eit
 of the morphism $\ph$.
}\eit

\brem
By Theorem \ref{TH:glav-sv-strog-razlozh},
\bit{
\item[---] for any decomposition $\ph=\mu\circ\e$, where $\e$ is an epimorphism, there is a unique morphism $\mu'$ such that
\beq\label{chastich-faktorizatsiya-v-uzlov-razlozhenii-1}
\xymatrix   @R=2.0pc @C=2.0pc
{
X\ar[rr]^{\ph}\ar[dd]_{\coim_\infty\ph}\ar[rd]^{\e} & & Y \\
 & M\ar[ru]^{\mu}\ar@{-->}[rd]^{\mu'} &  \\
 \Coim_\infty\ph\ar[rr]^{\red_\infty\ph} & & \Im_\infty\ph\ar[uu]_{\im_\infty\ph}
}
\eeq
(and if $\mu$ is a monomorphism, then $\mu'$ is a monomorphism),

\item[---] for any decomposition $\ph=\mu\circ\e$, where $\mu$ is a monomorphism, there is a unique morphism $\e'$ such that
\beq\label{chastich-faktorizatsiya-v-uzlov-razlozhenii-2}
\xymatrix  @R=2.0pc @C=2.0pc
{
X\ar[rr]^{\ph}\ar[dd]_{\coim_\infty\ph}\ar[rd]^{\e} & & Y \\
 & M\ar[ru]^{\mu} &  \\
 \Coim_\infty\ph\ar[rr]^{\red_\infty\ph}\ar@{-->}[ru]^{\e'} & & \Im_\infty\ph\ar[uu]_{\im_\infty\ph}
}
\eeq
(and if $\e$ is an epimorphism, then $\e'$ is an epimorphism).
}\eit
\erem

\paragraph{On existence of a nodal decomposition.}
Let us note that if $\mu$ is a monomorphism in a category ${\tt K}$, then for any its decomposition $\mu=\mu'\circ\e$, if $\e$ is a strong epimorphism, then $\e$ must be an isomorphism:
$$
\mu\in\Mono\quad\Longrightarrow\quad \Big(\forall \e\in\SEpi \quad \forall \mu'\quad \mu=\mu'\circ\e \quad\Longrightarrow\quad \e\in\Iso\Big).
$$
Indeed, by $1^\circ$ on p.\pageref{PROP:ph-mu=mono-=>-mu=mono}, the equality $\mu=\mu'\circ\e$ means that $\e$ must be a monomorphism, and, since in addition $\e$ is a strong epimorphism, so (by $4^\circ$ on p.\pageref{e-Sepi-=>-e-Nepi}), an immediate epimorphism, then by $4^\circ$ on p.\pageref{e-Nepi+mono-=>-e-iso} we obtain that $\e$ is an isomorphism.

\bit{
\item[$\bullet$]\label{DEF:strog-epi-razlich-mono} Let us say that in a category ${\tt K}$
 {\it strong epimorphisms discern monomorphisms}, if the reverse is true: from the fact that a morphism $\mu$ is not a monomorphism it follows that $\mu$ can be represented as a composition $\mu=\mu'\circ\e$, where $\e$ is a strong epimorphism, which is not an isomorphism.
}\eit

Dually, if $\e$ is an epimorphism in a category ${\tt K}$, then for any its decomposition  $\e=\mu\circ\e'$, if $\mu$ is a strong monomorphism, then  $\mu$ must be an isomorphism:
$$
\e\in\Epi\quad\Longrightarrow\quad \Big(\forall \mu\in\SMono \quad \forall \e'\quad \e=\mu\circ\e' \quad\Longrightarrow\quad \mu\in\Iso\Big).
$$

\bit{
\item[$\bullet$] Let us say that in a category ${\tt K}$ {\it strong monomorphisms discern epimorphisms}, if the reverse is true: from the fact that a morphism $\e$ is not an epimorphism it follows that $\e$ can be represented as a composition $\e=\mu\circ\e'$, where $\mu$ is a strong monomorphism, which is not an isomorphism.
 }\eit

Recall that the notion of linearly complete category was introduced on page \pageref{DEF:lineino-polnaya-kategoriya}.

\btm\label{TH:sush-uzlov-razlozh-v-polnoi-lok-maloi-kateg} Let ${\tt K}$ ba a linearly complete, well-powered in strong monomorphisms and co-well-powered in strong epimorphisms category, where strong epimorphisms discern monomorphisms, and, dually, strong monomorphisms discern epimorphisms. Then ${\tt K}$ is a category with nodal decomposition.
 \etm

Before proving this theorem let us introduce the following auxiliary construction. Take a morphism $\ph:X\to Y$ in a category ${\tt K}$. Since ${\tt K}$ is co-well-powered in strong epimorphisms, in the category $\SEpi^X$ of strong epimorphisms going from $X$ there exists a {\it set} of strong quotient objects  $Q\subseteq \SEpi^X$, and in the category $\SMono_Y$ of strong monomorphisms coming to $Y$ there is a {\it set} of strong subobjects $S\subseteq \SMono_Y$. We freeze these sets $Q$ and $S$.

\bit{

\item[$\bullet$]
A decomposition $\ph=\iota\circ\rho\circ\gamma$ of a morphism $\ph$ is said to be {\it admissible}\label{DEF:dop-razlozhenie}, if $\gamma\in Q$ and $\iota\in S$. Certainly, any strong decomposition $\ph=\iota'\circ\rho'\circ\gamma'$ of a morphism $\ph$ is isomorphic to some admissible decomposition $\ph=\iota\circ\rho\circ\gamma$.

\item[$\bullet$]
Let us call a {\it local basic decomposition}\label{DEF:loc-bas-razlozh} of a morphism $\ph$ in a category ${\tt K}$ an arbitrary map $\rho\mapsto(\coim\rho,\red\rho,\im\rho)$ that to each admissible decomposition $(\iota,\rho,\gamma)$ of the morphism $\ph$ assigns some strong decomposition $(\im\rho,\red\rho,\coim\rho)$ of $\rho$
\beq\label{DEF:lok-bazis-razlozh}
\begin{diagram}
\node{X}\arrow{s,l}{\gamma}\arrow{e,t}{\ph}\node{Y} \\
\node{\Dom\rho}\arrow{s,l}{\coim\rho}\arrow{e,t}{\rho}\node{\Ran\rho}\arrow{n,r}{\iota} \\
\node{\Coim\rho}\arrow{e,t}{\red\rho}\node{\Im\rho}\arrow{n,r}{\im\rho}
\end{diagram}
\eeq
in such a way that the following conditions are fulfilled:
 \bit{
 \item[(a)] the decomposition $(\iota\circ\im\rho,\red\rho,\coim\rho\circ\gamma)$ of $\ph$ is admissible (i.e. $\coim\rho\circ\gamma\in Q$ and $\iota\circ\im\rho\in S$),

 \item[(b)]\label{rho-mono<=>coim-rho=1} $\rho$ is a monomorphism $\Longleftrightarrow$ $\coim\rho$ is an isomorphism $\Longleftrightarrow$ $\coim\rho=1$,

 \item[(c)]\label{rho-epi<=>im-rho=1} $\rho$ is an epimorphism $\Longleftrightarrow$ $\im\rho$ is an epimorphism $\Longleftrightarrow$ $\im\rho=1$.
 }\eit

}\eit

\blm\label{LM:lok-baz-razlozh} Let ${\tt K}$ be a well-powered in strong monomorphisms and co-well-powered in strong epimorphisms category, where strong epimorphisms discern monomorphisms, and strong monomorphisms discern epimorphisms. Then each morphism $\ph$ in ${\tt K}$ has local basic decomposition.
\elm
\bpr
First of all, it is clear that adimssible decompositions always exist, for example one can take $\ph=1\circ\ph\circ 1$.
Let us then show that for any admissible decomposition $(\iota,\rho,\gamma)$ of $\ph$ a diagram \eqref{DEF:lok-bazis-razlozh} satisfying (a), (b), (c) exists. Let us freeze this decomposition $(\iota,\rho,\gamma)$ and consider several cases.
 \bit{
\item[1.] If $\rho$ is not a monomorphism, then there exists a decomposition
$\rho=\rho'\circ\e$, where $\e$ is a strong epimorphism, but not an isomorphism. Set $\coim\rho=\e$ and consider the morphism $\rho'$.
 \bit{
\item[1.1.] If $\rho'$ is not an epimorphism, then there exists a decomposition $\rho'=\mu\circ\rho''$, where $\mu$ is a strong monomorphism, but not an isomorphism. Then we set $\im\rho=\mu$ and $\red\rho=\rho''$.

\item[1.2.] If $\rho'$ is an epimorphism, then we set $\im\rho=1_{\Ran\rho}$ and $\red\rho=\rho'$.
 }\eit

\item[2.] If $\rho$ is a monomorphism, then we set $\coim\rho=1_{\Dom\rho}$ and again consider $\rho$.
 \bit{
\item[2.1.] If $\rho$ is not an epimorphism, then there exists a decomposition $\rho=\mu\circ\rho'$, where $\mu$ is a strong monomorphism, but not an isomorphism. We set $\im\rho=\mu$ and $\red\rho=\rho'$.

\item[2.2.] If $\rho$ is an epimorphism, then we set $\im\rho=1_Y$ and $\red\rho=\rho$.
 }\eit
 }\eit\noindent
In any case we obtain a decomposition $\rho=\im\rho\circ\red\rho\circ\coim\rho$, where $\im\rho$ is a strong monomorphism, $\coim\rho$ is a strong epimorphism, and (b) and (c) are fulfilled. Now to provide (a) we have to replace (if necessary) the epimorphism $\coim\rho$ with an isomorphic epimorphism $\pi\circ\coim\rho$ in such a way that $\pi\circ\coim\rho\circ\gamma\in Q$, and this can be done due to Theorem  \ref{PROP:dost-mnozh-v-klasse-faktor-objektov}. Similarly, the monomorphism $\im\rho$ should be replaced with an isomorphic monomorphism $\im\rho\circ\sigma$ in such a way that $\iota\circ\im\rho\circ\sigma\in S$, and this can be done due to Theorem \ref{PROP:dost-mnozh-v-klasse-podobjektov}.

Thus, for an arbitrary admissible decomposition $(\iota,\rho,\gamma)$ of $\ph$ diagram \eqref{DEF:lok-bazis-razlozh} satisfying (a), (b), (c), exists. Note now that from Theorems \ref{PROP:dost-mnozh-v-klasse-faktor-objektov} and \ref{PROP:dost-mnozh-v-klasse-podobjektov} it follows that for a given admissible decomposition $(\iota,\rho,\gamma)$ of morphism $\ph$ the class of decompositions $(\im\rho,\red\rho,\coim\rho)$ of $\rho$, which satisfy (a), (b), (c), is a set. Indeed, every such a decomposition $(\im\rho,\red\rho,\coim\rho)$ is uniquely defined by the morphisms $\im\rho$ and $\coim\rho$ (since from monomorphity of $\im\rho$ and epimorphity of $\coim\rho$ it follows that $\red\rho$, if exists, is unique). So the class of decompositions $(\im\rho,\red\rho,\coim\rho)$ can be conceived as a subclass in the cartesian product of sets $A\times B$, where
$A=\{\alpha\in\SMono_{\Ran\rho}:\ \iota\circ\alpha\in S\}$ is a class of monomorphisms where $\im\rho$ runs, and which is a set by Proposition \ref{PROP:dost-mnozh-v-klasse-podobjektov}, and $B=\{\beta\in\SEpi^{\Dom\rho}:\ \beta\circ\e\in Q\}$ is a class of epimorphisms, where $\coim\rho$ runs, and which is a set by Theorem \ref{PROP:dost-mnozh-v-klasse-faktor-objektov}).

We obtain that for any admissible decomposition $(\iota,\rho,\gamma)$ of $\ph$ the class of decompositions $(\coim\rho,\red\rho,\im\rho)$ satisfying \eqref{DEF:lok-bazis-razlozh} and (a), (b), (c), is a (non-empty) set. From this it follows that we can apply the axiom of choice and construct a map which to each admissible decomposition $(\iota,\rho,\gamma)$ of $\ph$ assigns a decomposition $(\coim\rho,\red\rho,\im\rho)$, satisfying \eqref{DEF:lok-bazis-razlozh} and (a), (b), (c). This is the required map $\rho\mapsto(\coim\rho,\red\rho,\im\rho)$.
\epr

\bpr[Proof of Theorem \ref{TH:sush-uzlov-razlozh-v-polnoi-lok-maloi-kateg}.] Take a morphism $\ph:X\to Y$, find a {\it set} of strong quotient objects $Q\subseteq \SEpi^X$ and a {\it set} of strong subobjects $S\subseteq \SMono_Y$, and construct a local basic decomposition like in Lemma \ref{LM:lok-baz-razlozh}. The proof consists in constructing a transfinite system of objects and morphisms, indexed by ordinal numbers $i\in{\tt Ord}$,
$$
X^i\overset{\ph^i}{\longrightarrow} Y^i,\quad X^i\overset{\e^i_j}{\longrightarrow} X^j, \quad Y^i\overset{\mu^i_j}{\longleftarrow} Y^j\quad (i\le j)
$$
the idea of which is illustrated by the following diagram (going infinitely below):
\beq\label{DIAGR:postroenie-uzlov-razlozhenija}
\xymatrix @R=2.5pc @C=5.0pc
{
X\ar[r]^{\ph}\ar@{=}[d]_{1_X} & Y \\
X^0\ar[r]^{\ph^0=\ph}\ar[d]_{\e^0_1=\coim\ph^0} & Y^0\ar@{=}[u]_{1_Y} \\
X^1\ar[r]^{\ph^1=\red\ph^0}\ar[d]_{\e^1_2=\coim\ph^1} & Y^1\ar[u]_{\mu^0_1=\im\ph^0} \\
X^2\ar[r]^{\ph^2=\red\ph^1}\ar[d]_{\e^2_3=\coim\ph^2} & Y^1\ar[u]_{\mu^1_2=\im\ph^1} \\
... & ... \ar[u]_{\mu^2_3=\im\ph^2}
}
\eeq
Here is how we do this.

0) Initially, we put
$$
 X^0=X,\qquad Y^0=Y,\qquad \ph^0=\ph, \qquad
 \e^0_1=\coim\ph^0,\qquad \mu^0_1=\im\ph^0,\qquad \ph^1=\red\ph^0.
$$

1) Then for an arbitrary ordinal number $k$ we put
$$
 \e^k_k=1_{X^k},\qquad \mu^k_k=1_{Y^k}
$$
 and
 \bit{
 \item[---] if $k$ is an isolated ordinal, i.e. $k=j+1$ for some $j$, then we set
 \begin{multline*}
 X^k=X^{j+1}=\Coim\ph^j,\qquad Y^k=Y^{j+1}=\Im\ph^j,\\
 \e^j_k=\e^j_{j+1}=\coim\ph^j,\qquad \mu^j_k=\mu^j_{j+1}=\im\ph^j,\qquad \ph^k=\ph^{j+1}=\red\ph^j
 \end{multline*}
 and after that for any other ordinal number $i<j$
 $$
 \e^i_k=\e^i_{j+1}=\e^j_{j+1}\circ\e^i_j,\qquad \mu^i_k=\mu^i_{j+1}=\mu^i_j\circ\mu^j_{j+1},
 $$

 \item[---] if $k$ is a limit ordinal, i.e. for any $j<k$ we have $j+1<k$, then $X^k$ is defined as the injective limit of the covariant system $\{X^j,\e^i_j;i\le j<k\}$, $Y^k$ as the projective limit of the contravariant system $\{Y^j,\mu^i_j;i\le j<k\}$,
$$
 X^k=\lim_{j\to k}X^j,\qquad Y^k=\lim_{k\gets j}Y^j,
$$
the system of morphisms $\{\e^i_k;i<k\}$ is the corresponding injective cone of morphism going to $X^k$, and the system of morphisms
$\{\mu^i_k;i<k\}$ is the corresponding projective cone of morphisms going from $Y^k$,
$$
 \e^i_k=\lim_{j\to k}\e^i_j,\qquad \mu^i_k=\lim_{k\gets j}\mu^i_j,\qquad i\le k.
$$
This automatically implies equalities
$$
 \e^i_k=\e^j_k\circ\e^i_j,\qquad  \mu^i_k=\mu^i_j\circ\mu^j_k,\qquad i\le j\le k
$$
and by Proposition \ref{PROP:inj-predel-strogih-epimorfizmov} all the morphisms $\e^i_k$ are strong epimorphisms, while by Proposition \ref{PROP:proj-predel-strogih-monomorfizmov} all the morphisms $\mu^i_j$ are strong monomorphisms. As a corollary, the object $X^k$ can be chosen in such a way that the epimorphism $\e^0_k$ lies in $Q$ (for this we just need to multiply from the left the system $\{\e^i_k;i<k\}$ of epimorphisms by a morphism, so that the property of being injective cone is preserved); similarly, the object $Y^k$ can be chosen in such a way that the monomorphism $\mu^0_k$ lies in the set $S$ (for this we just need to multiply from the right the system $\{\mu^i_k;i<k\}$ of monomorphisms, so that the property of being projective cone is preserved). That is what we will do, and after that the morphism $\ph^k$ can be defined by two equivalent formulas:
 $$
 \ph^k=\lim_{k\gets i}\lim_{j\to k}\mu^i_j\circ\ph^j=\lim_{i\to k}\lim_{k\gets j}\ph^j\circ\e^i_j
 $$
Here the first double limit should be understood as follows: for a given $i<k$ the family $\{\mu^i_j\circ\ph^j;\ i\le j<k\}$ is an injective cone of the covariant system $\{\e^l_j;\ i\le l,j<k\}$, so there exists a limit
 $$
 \lim_{j\to k}\mu^i_j\circ\ph^j;
 $$
after that the system $\{ \lim_{j\to k}\mu^i_j\circ\ph^j;\ i<k \}$ turns out to be a projective cone of the contravariant system $\{\mu^l_j;\ i\le l,j<k\}$, so there exists a limit
 $$
 \lim_{k\gets i}\lim_{j\to k}\mu^i_j\circ\ph^j.
 $$
Similarly, in the second double limit for a given $i<k$ the family $\{\ph^j\circ\e^i_j;\ i\le j<k\}$ is a projective cone of the contravariant system $\{\mu^l_j;\ i\le l,j<k\}$, so there exists a limit
 $$
 \lim_{k\gets j}\ph^j\circ\e^i_j;
 $$
after that the system $\{ \lim_{k\gets j}\ph^j\circ\e^i_j;\ i<k \}$ turns out to be an injective cone of the covariant system $\{\e^l_j;\ i\le l,j<k\}$, so there exists a limit
 $$
 \lim_{i\to k}\lim_{k\gets j}\ph^j\circ\e^i_j.
 $$
Each of these double limits gives an arrow from $X^k$ into $Y^k$ which makes all the necessary diagrams commutative, and since this arrow is unique (this follows from the fact that $\mu^i_k$ are monomorphisms and $\e^i_k$ are epimorphisms), those double limits (arrows) coincide.
 }\eit

Eventually we obtain a system of morphisms such that for any two ordinal numbers $i\le j$ the following diagram is commutative
$$
\begin{diagram}
\node{X^i}\arrow{s,l}{\e^i_j}\arrow{e,t}{\ph^i}\node{Y^i}\\
\node{X^j}\arrow{e,t}{\ph^j}\node{Y^j}\arrow{n,r}{\mu^i_j}
\end{diagram}
$$
and for any three ordinal numbers $i\le j\le k$ the following diagrams are commutative
$$
\begin{diagram}
\node[2]{X^i}\arrow{sw,l}{\e^i_j}\arrow[2]{s,r}{\e^i_k}  \\
\node{X^j}\arrow{se,b}{\e^j_k} \\
\node[2]{X^k}
\end{diagram}
\qquad
\begin{diagram}
\node[2]{Y^i} \\
\node{Y^j}\arrow{ne,l}{\mu^i_j}   \\
\node[2]{Y^k}\arrow{nw,b}{\mu^j_k}\arrow[2]{n,r}{\mu^i_k}
\end{diagram}
$$
and moreover, $\e^i_j$ are strong epimorphisms, and $\mu^i_j$ are strong monomorphisms. From the last two diagrams it follows that the formulas
$$
\begin{cases}
F(i)=\e_i^0,& i\in{\tt Ord}
\\
F(i,j)=\e^i_j,& i\le j\in{\tt Ord}
\end{cases}
\qquad
\begin{cases}
G(i)=\mu_i^0,& i\in{\tt Ord}
\\
G(i,j)=\mu^i_j,& i\le j\in{\tt Ord}
\end{cases}
$$
define a covariant functor $F:{\tt Ord}\to Q$ and a contravariant functor $G:{\tt Ord}\to S$. Since $Q$ and $S$ are sets, by Theorem  \ref{TH:obryv-transf-tsepei} these functors must stabilize, i.e. starting from some ordinal number $k$ (which can be chosen common for $F$ and $G$) the morphisms $F(i,j)$ and $G(i,j)$ become isomorphisms. Since in addition the categories $Q$ and $S$ are partially ordered classes (and as a corollary, only local identities are isomorphisms there, by Proposition \ref{PROP:iso-v-chast-upor-klasse}), we obtain (following Remark \ref{REV:stabiliziruemost-v-chast-upor-klasse}) that diagram \eqref{DIAGR:postroenie-uzlov-razlozhenija} is stabilized in the sense that, starting from some $k$,
 \bit{
\item[---] the objects $X^l$ become the same, and the morphisms $\e^l_m$ become local identities of $X^k$:
$$
\forall m>l\ge k\qquad X^m=X^l=X^k,\qquad \e^l_m=1_{X^k}
$$
\item[---] and the objects $Y^l$ become the same and the morphisms $\mu^l_m$ become local identities of $Y^k$:
$$
\forall m>l\ge k\qquad Y^m=Y^l=Y^k,\qquad \mu^l_m=1_{Y^k}
$$
 }\eit

Now let us consider the diagram
\beq\label{PROOF:sushestvovanie-uzlovogo-razlozhenija}
\begin{diagram}
\node{X}\arrow{s,l}{\e^0_k}\arrow{e,t}{\ph}\node{Y}\\
\node{X^k}\arrow{e,t}{\ph^k}\node{Y^k}\arrow{n,r}{\mu^0_k}
\end{diagram}
\eeq
Here $\e^0_k$ is a strong epimorphism, and $\mu^0_k$ a strong monomorphism. From the equality $\e^k_{k+1}=\coim\ph^k=1_{X^k}$ (which holds since the sequence $\e^0_j$ is stabilized for $j\ge k$) it follows by condition (b) on page \pageref{rho-mono<=>coim-rho=1}, that $\ph^k$ is a monomorphism. On the other hand, from the equality $\mu^k_{k+1}=\im\ph^k=1_{Y^k}$ (which holds since the sequence $\mu^0_j$ is stabilized for $j\ge k$) it follows by condition (c) on page \pageref{rho-epi<=>im-rho=1}, that $\ph^k$ is an epimorphism. Thus, $\ph^k$ is a bimorphism, hence \eqref{PROOF:sushestvovanie-uzlovogo-razlozhenija} is a nodal decomposition for $\ph$.
\epr

\paragraph{Connection with the basic decomposition in pre-Abelian categories}

Let us discuss the obvious analogy between nodal decomposition and the decomposition of a morphism $\ph$ in a pre-Abelian category ${\tt K}$ into a coimage $\coim\ph$, image $\im\ph$ and a morphism between them which we denote by $\red\ph$.

Recall (see definition in \cite{Bucur-Deleanu} or in \cite{General-algebra}) that {\it pre-Abelian category} is an enriched category ${\tt K}$ over the category ${\tt Ab}$ of Abelian groups, which is finitely complete and has zero object. In such a category every morphism $\ph:X\to Y$ has a kernel and a cokernel. From this it follows that $\ph$ can be represented as a composition
\beq\label{EX:bazis-razlozh}
\begin{diagram}
\node{X}\arrow{s,l}{\coim\ph}\arrow{e,t}{\ph}\node{Y} \\
\node{\Coim\ph}\arrow{e,t,--}{\red\ph}\node{\Im\ph}\arrow{n,r}{\im\ph}
\end{diagram}
\eeq
where the morphism $\coim\ph=\coker(\ker\ph)$ is called the {\it coimage} of $\ph$, the morphism $\im\ph=\ker(\coker\ph)$ the {\it image} of $\ph$, and the existence and uniqueness of the morphism $\red\ph$ is proved separately, and we will call it the {\it reduced part} of $\ph$.

\bit{
\item[$\bullet$] The representation of a morphism $\ph$ as a composition \eqref{EX:bazis-razlozh} we call the {\it basic decomposition} of $\ph$. }\eit

It is known (see \cite[Proposition 4.3.6(4)]{Borceux}) that in a pre-Abelian category (in fact, in a category with zero) every kernel $\ker\ph$ (and thus, every image $\im\ph$) is always a strong monomorphism, and every cokernel $\coker\ph$  (and thus, every coimage $\coim\ph$) is a strong epimorphism. As a corollary, we have

\btm\label{PROP:im-i-coim-strogie}
In a pre-Abelian category every basic decomposition is strong.
\etm

This implies that {\it if a category ${\tt K}$ is Abelian, then every basic decomposition in  ${\tt K}$ is nodal}. But if ${\tt K}$ is not Abelian, then these decompositions do not necessarily coincide, see below Example \ref{EX:reg(ph)-ne-Bim}.

The following two propositions are obvious:

\bprop\label{PROP:ph-mono=>ker-ph=0} In a pre-Abelian category for a morphism $\ph:X\to Y$ the following conditions are equivalent:
 \bit{
\item[(i)] $\ph$ is a monomorphism,

\item[(ii)] the zero morphism $0_{0,X}$ is a kernel for $\ph$: $0_{0,X}=\ker\ph$,

\item[(iii)] the identity morphism $1_X$ is a coimage for $\ph$: $1_X=\coim\ph$,

\item[(iv)] $\coim\ph$ is an isomorphism.
 }\eit
\eprop

\bprop\label{PROP:ph-epi=>coker-ph=0} In a pre-Abelian category for a morphism $\ph:X\to Y$ the following conditions are equivalent:
 \bit{
\item[(i)] $\ph$ is an epimorphism,

\item[(ii)] the zero morphism $0_{Y,0}$ is a cokernel for $\ph$: $0_{Y,0}=\coker\ph$,

\item[(iii)] the identity morphism $1_Y$ is an image for $\ph$: $1_Y=\im\ph$,

\item[(iv)] $\im\ph$ is an isomorphism.
 }\eit
\eprop

They imply

\bprop\label{PROP:bazis-razl=>strog-epi-razlich-mono} In a pre-Abelian category $\tt K$ the strong epimorphisms discern monomorphisms and the strong monomorphisms discern epimorphisms.
\eprop
\bpr Consider the basic decomposition of $\ph:X\to Y$:
$$
\ph=\im\ph\circ\red\ph\circ\coim\ph
$$
If $\ph:X\to Y$ is not a monomorphism, then by Proposition \ref{PROP:ph-mono=>ker-ph=0}, $\coim\ph$ is not an isomorphism. On the other hand, by Proposition \ref{PROP:im-i-coim-strogie}, $\coim\ph$ is a strong epimorphism. So, if we set $\ph'=\im\ph\circ\red\ph$, then in the decomposition $\ph=\ph'\circ\coim\ph$ the morphism $\coim\ph$ is a strong epimorphism, but not an isomorphism. This means that strong epimorphisms discern monomorphisms in $\tt K$. The statement about strong monomorphisms is proved similarly.
\epr

Proposition \ref{PROP:bazis-razl=>strog-epi-razlich-mono} implies that if a pre-Abelian category ${\tt K}$ is  well-powered in strong monomorphisms and co-well-powered in strong epimorphisms, then ${\tt K}$ has local basic decomposition (defined on page \pageref{DEF:loc-bas-razlozh}): the map $(\iota,\rho,\gamma)\mapsto (\coim\rho,\red\rho\im\rho)$ that to each admissible decomposition $(\iota,\rho,\gamma)$ (admissible decompositions were defined on page \pageref{DEF:dop-razlozhenie}) of a given morphism $\ph$ assigns the basic decomposition of $\rho$, is a local basic decomposition of $\ph$. Hence, the sufficient condition for existence of nodal decomposition (Theorem \ref{TH:sush-uzlov-razlozh-v-polnoi-lok-maloi-kateg}) becomes more simple:

\btm\label{TH:faktorizatsija-v-lok-maloi-kategorii} If a pre-Abelian category ${\tt K}$ is linearly complete, well-powered in strong monomorphisms and co-well-powered in strong epimorphisms, then every morphism $\ph:X\to Y$ in ${\tt K}$ has nodal decomposition \eqref{DEF:oboznacheniya-dlya-uzlov-razlozh}.
\etm

\brem\label{REM:struktura-uzlov-razlozh-v-kateg-s-nulem} From Proposition \ref{PROP:bazis-razl=>strog-epi-razlich-mono} and Diagram  \eqref{DIAGR:postroenie-uzlov-razlozhenija} it follows that
 \bit{
 \item[---] the nodal reduced part $\red_\infty\ph$ in diagram \eqref{DEF:oboznacheniya-dlya-uzlov-razlozh} can be conceived as a ``limit'' of transfinite sequence of ``usual'' reduced morphisms $\ph^{i+1}=\red\ph^i$,
 \item[---] the nodal coimage $\coim_\infty\ph$ is an injective limit of transfinite sequence of ``usual'' coimages $\coim\ph^i$ of this system of morphisms, and
 \item[---] the nodal image $\im_\infty\ph$ is a projective limit of transfinite sequence of ``usual'' images $\im\ph^i$ of this system of morphisms.
 }\eit
\erem

\brem
Since as we already noticed the basic decomposition $\ph=\im\ph\circ\red\ph\circ\coim\ph$ is strong, and thus, by Proposition \ref{PROP:uzlov-razlozh-podchin-ostalnye}, is subordinated to the nodal decomposition, there must exist unique morphisms $\sigma$ and $\tau$ such that the following diagram is commutative:
\beq\label{svayz-bazisnogo-i-uzlovogo-razlozheniya}
\xymatrix  @R=2.5pc @C=4.0pc
{
X\ar[rrr]^{\ph}\ar[rd]^{\coim_\infty\ph}\ar[dd]_{\coim\ph} & & & Y \\
& \Coim_\infty\ph\ar[r]_{\red_\infty\ph} & \Im_\infty\ph\ar[ru]^{\im_\infty\ph}\ar@{-->}[rd]_{\tau} & \\
\Coim\ph\ar[rrr]_{\red\ph}\ar@{-->}[ru]_{\sigma} & & & \Im\ph \ar[uu]_{\im\ph}
}
\eeq
At the same time, by Theorem \ref{TH:glav-sv-strog-razlozh},
\bit{
\item[---] for any decomposition $\ph=\mu\circ\e$, where $\e$ is an epimorphism, there exists a unique morphism $\mu'$ such that the following diagram is commutative:
\beq\label{1-sledstvie--v-bazisn-razl+uzlov-razlozh}
\xymatrix % @R=2.5pc @C=4.0pc
{
X\ar[rrrr]^{\ph}\ar[rdd]^{\coim_\infty\ph}\ar[ddd]_{\coim\ph}\ar[rrd]^{\e} & & & & Y \\
& & M\ar[rru]^{\mu}\ar@{-->}[rd]^{\mu'} & & \\
& \Coim_\infty\ph\ar[rr]_{\red_\infty\ph} & & \Im_\infty\ph\ar[ruu]^{\im_\infty\ph}\ar[rd]_{\tau} & \\
\Coim\ph\ar[rrrr]_{\red\ph}\ar[ru]_{\sigma} & & & & \Im\ph \ar[uuu]_{\im\ph}
}
\eeq
(in addition, if $\mu$ is a monomorphism, then $\mu'$ is a monomorphism as well);

\item[---] for any decomposition $\ph=\mu\circ\e$, where $\mu$ is a monomorphism, there exists a unique morphism $\e'$ such that the following diagram is commutative:
\beq\label{2-sledstvie-v-bazisn-razl+uzlov-razlozh}
\xymatrix  % @R=2.5pc @C=4.0pc
{
X\ar[rrrr]^{\ph}\ar[rdd]^{\coim_\infty\ph}\ar[ddd]_{\coim\ph}\ar[rrd]^{\e} & & & & Y \\
& & M\ar[rru]^{\mu} & & \\
& \Coim_\infty\ph\ar[rr]_{\red_\infty\ph}\ar@{-->}[ru]^{\e'} & & \Im_\infty\ph\ar[ruu]^{\im_\infty\ph}\ar[rd]_{\tau} & \\
\Coim\ph\ar[rrrr]_{\red\ph}\ar[ru]_{\sigma} & & & & \Im\ph \ar[uuu]_{\im\ph}
}
\eeq
(in addition, if $\e$ is an epimorphism, then $\e'$ is an epimorphism as well);

\item[---] in particular, for any factorization $\ph=\mu\circ\e$ of $\ph$ there exist unique morphisms $\Coim\ph\overset{\e'}{\longrightarrow}M$ and $M\overset{\mu'}{\longrightarrow}\Im\ph$ such that the following diagram is commutative:
\beq\label{DIAGR:faktorizatsija-v-bazisn-razl+uzlov-razlozh}
\xymatrix % @R=2.5pc @C=4.0pc
{
X\ar[rrrr]^{\ph}\ar[rdd]^{\coim_\infty\ph}\ar[ddd]_{\coim\ph}\ar[rrd]^{\e} & & & & Y \\
& & M\ar[rru]^{\mu}\ar@{-->}[rd]^{\mu'} & & \\
& \Coim_\infty\ph\ar[rr]_{\red_\infty\ph}\ar@{-->}[ru]^{\e'} & & \Im_\infty\ph\ar[ruu]^{\im_\infty\ph}\ar[rd]_{\tau} & \\
\Coim\ph\ar[rrrr]_{\red\ph}\ar[ru]_{\sigma} & & & & \Im\ph \ar[uuu]_{\im\ph}
}
\eeq
and in addition, $\e'$ is an epimorphism, and $\mu'$ a monomorphism.
}\eit
\erem

\subsection{Factorizations of a category}

\paragraph{Factorizations in a category with nodal decomposition.}\label{paragraph:faktorizatsija}
Recall that the notion of a factorization of a morphism was defined on page \pageref{DEF:faktorizatsiya}.
From \eqref{chastich-faktorizatsiya-v-uzlov-razlozhenii-1} and \eqref{chastich-faktorizatsiya-v-uzlov-razlozhenii-2} we immediately have

\bprop\label{PROP:vlozhenie-faktorizatsii-v-kanonicheskoe-razlozh} If $X\overset{\e}{\longrightarrow}M\overset{\mu}{\longrightarrow}Y$ is a factorization of a morphism $X\overset{\ph}{\longrightarrow}Y$ in a category ${\tt K}$ with a nodal decomposition, then there are unique morphisms  $\Coim_\infty\ph\overset{\e'}{\longrightarrow}M$ and $M\overset{\mu'}{\longrightarrow}\Im_\infty\ph$ such that the following diagram is commutative:
\beq\label{faktorizatsiya-v-uzlov-razlozhenii}
\xymatrix  @R=2.0pc @C=2.0pc
{
X\ar[rr]^{\ph}\ar[dd]_{\coim_\infty\ph}\ar[rd]^{\e} & & Y \\
 & M\ar[ru]^{\mu}\ar@{-->}[rd]^{\mu'} &  \\
 \Coim_\infty\ph\ar[rr]^{\red_\infty\ph}\ar@{-->}[ru]^{\e'} & & \Im_\infty\ph\ar[uu]_{\im_\infty\ph}
}
\eeq
Moreover, $\e'$ is an epimorphism, and $\mu'$ a monomorphism.
\eprop

Let $(\e,\mu)$ and $(\e',\mu')$ be two factorizations of $\ph$. We say that the factorization $(\e,\mu)$ {\it is subordinated} to the factorization  $(\e',\mu')$ (or $(\e',\mu')$ {\it subordinates} $(\e,\mu)$), and write
$$
(\e,\mu)\le (\e',\mu'),
$$
if there exists a morphism $\beta$ such that
$$
\e'=\beta\circ\e,\qquad \mu=\mu'\circ\beta
$$
i.e.
$$
\xymatrix @R=3.0pc @C=4.0pc
{
X\ar[d]_{\e}\ar[r]^{\ph}\ar[dr]^>>{\e'}|!{[d];[r]}\hole & Y
\\
M\ar@{-->}[r]_{\beta}\ar[ur]_>>>>{\mu} & M'\ar[u]_{\mu'}
}
$$
From Properties $1^\circ$ and $3^\circ$ on page \pageref{PROP:ph-mu=mono-=>-mu=mono} it follows that $\beta$, if exists, must be a bimorphism, and from the fact that $\mu'$ is a monomorphism (or from the fact that $\e$ is an epimorphism) that $\beta$ is unique.

\btm\label{TH:faktorizatsija-v-kategorii-s-uzlov-razlozh} In a category ${\tt K}$ with nodal decomposition
 \bit{
\item[(i)] every morphism $\ph$ has a factorization,

\item[(ii)] among all factorizations of $\ph$ there is a minimal $(\e_{\min},\mu_{\min})$ and a maximal $(\e_{\max},\mu_{\max})$, i.e. the factorizations that bound any other factorization $(\e,\mu)$:
    $$
    (\e_{\min},\mu_{\min})\le (\e,\mu)\le (\e_{\max},\mu_{\max})
    $$
 }\eit
\etm

\bpr Here (i) follows from (ii), so we prove (ii). Put
$$
\e_{\min}=\coim_\infty\ph,\qquad \mu_{\min}=\im_\infty\ph\circ\red_\infty\ph,\qquad \e_{\max}=\red_\infty\ph\circ\coim_\infty\ph,\qquad \mu_{\max}=\im_\infty\ph
$$
then these will be factorizations of $\ph$, and from \eqref{faktorizatsiya-v-uzlov-razlozhenii} it follows that the first one of them is minimal, and the second one is maximal.
\epr

\paragraph{Strong morphisms in a category with nodal decomposition.}

\btm\label{TH:harakter-strogih-morf-v-kat-s-uzlov-razlozh} In a category with nodal decomposition
 \bit{
\item[(a)]  $\mu$ is an immediate monomorphism $\Longleftrightarrow$ $\mu$ is a strong monomorphism  $\Longleftrightarrow$ $\mu\cong\im_\infty\mu$
 $\Longleftrightarrow$ $\coim_\infty\mu$ and $\red_\infty\mu$ are isomorphisms,

\item[(b)]  $\e$ is an immediate epimorphism $\Longleftrightarrow$ $\e$ is a strong epimorphism $\Longleftrightarrow$ $\e\cong\coim_\infty\e$
 $\Longleftrightarrow$ $\im_\infty\mu$ and $\red_\infty\mu$ are isomorphisms.
 }\eit
\etm
\bpr By the duality principle it is sufficient to prove (a).

1. If $\mu:X\to Y$ is an immediate monomorphism, then in its maximal factorization
$$
\mu=\mu_{\max}\circ\e_{\max}
$$
the morphism $\e_{\max}=\red_\infty\mu\circ\coim_\infty\mu$ must be an isomorphism. This implies formula
$$
1_X=(\e_{\max})^{-1}\circ\red_\infty\mu\circ\coim_\infty\mu
$$
from which one can conclude that $\coim_\infty\mu$ is a coretraction. On the other hand, $\coim_\infty\mu$ is an epimorphism, hence $\coim_\infty\mu$ is an isomorphism. This implies that $\red_\infty\mu=\e_{\max}\circ(\coim_\infty\mu)^{-1}$ is an isomorphism.

2. If $\coim_\infty\mu$ and $\red_\infty\mu$ are isomorphisms, then its composition $\chi=\red_\infty\mu\circ\coim_\infty\mu$ is an isomorphism as well, and at the same time $\mu=\im_\infty\mu\circ\chi$. This means that $\mu\cong\im_\infty\mu$.

3. If $\mu\cong\im_\infty\mu$, then, since $\im_\infty\mu$ is a strong monomorphism, $\mu$ is also a strong monomorphism.

4. If $\mu$ is a strong monomorphism, then by property $2^\circ$ on page \pageref{mu-Smono-=>-mu-Nmono}, $\mu$ is an immediate monomorphosm.
\epr

\paragraph{Factorization of a category.}

\bit{

\item
A pair of morphisms $(\mu,\e)$ is said to be {\it diagonizable}
\cite{General-algebra,Tsalenko-Shulgeifer}, if for all morphisms
$\alpha:\Dom\e\to\Dom\mu$ and $\beta:\Ran\e\to\Ran\mu$ such that $\mu\circ\alpha=\beta\circ\e$
there exists a morphism $\delta:B\to C$, such that the diagram \eqref{razbienie-kvadrata} is commutative:
$$
\begin{diagram}
\node{\Dom\e}\arrow{s,l}{\alpha}\arrow{e,t}{\e}\node{\Ran\e}\arrow{s,r}{\beta}\arrow{sw,t,--}{\delta}
\\
\node{\Dom\mu}\arrow{e,b}{\mu}\node{\Ran\mu}
\end{diagram}
$$
This is denoted by writing $\mu\downarrow\e$.
}\eit

\bex
The following example shows that in contrast to the situation considered above (in particular at the page \pageref{razbienie-kvadrata}), the relation $\mu\downarrow\e$ does not necessarily mean that $\mu\in\Mono$ and $\e\in\Epi$: in the category of vector spaces over the field $\C$ the pair of morphisms $\mu=0:\C\to 0$ and $\e=0:0\to\C$ is diagonizable:
$$
\begin{diagram}
\node{0}\arrow{s,l}{\alpha}\arrow{e,t}{\e=0}\node{\C}\arrow{s,r}{\beta}\arrow{sw,t,--}{\delta}
\\
\node{\C}\arrow{e,b}{\mu=0}\node{0}
\end{diagram}
$$
\eex

\bit{
\item For each class of morphisms $\varLambda$ in $\tt K$

 \bit{
\item[---] its {\it epimorphic conjugate class} is the class
$$
\varLambda^\downarrow=\{\e\in\Epi({\tt K}):\forall\lambda\in\varLambda\quad \lambda\downarrow\e\}.
$$
\item[---] its {\it monomorphic conjugate class} is the class
$$
{^\downarrow\kern-1pt\varLambda}=\{\mu\in\Mono({\tt K}):\forall\lambda\in\varLambda\quad \mu\downarrow\lambda\}.
$$
}\eit
}\eit
Clearly, for each class of morphisms $\varLambda$
\begin{align}
& \Iso\subseteq\varLambda^\downarrow\subseteq\Epi, && \Iso\circ\ \varLambda^\downarrow\subseteq\varLambda^\downarrow
\label{Iso-subseteq-varTheta^downarrow-subseteq-Epi}
\\
& \Iso\subseteq{^\downarrow\kern-1pt\varLambda}\subseteq\Mono, &&
{^\downarrow\kern-1pt\varLambda}\circ\Iso\subseteq{^\downarrow\kern-1pt\varLambda}
\label{Iso-subseteq-^downarrow-varTheta-subseteq-Mono}
\end{align}

\bit{

\item
Let us say that classes of morphisms $\varGamma$ and $\varOmega$ define a {\it factorization of the category\footnote{This construction is also called a {\it bicategory} \cite{General-algebra,Tsalenko-Shulgeifer}.}
$\tt K$}\label{DEF:faktorizatsija-v-kategorii}, if
 \bit{
 \item[F.1] $\varOmega$ is the epimorphic conjugate class for $\varGamma$:
 $$
 \varGamma^\downarrow=\varOmega
 $$

 \item[F.2] $\varGamma$ is the monomorphic conjugate class for $\varOmega$:
 $$
 \varGamma={^\downarrow\varOmega},
 $$

 \item[F.3] the composition of the class $\varGamma$ and $\varOmega$ covers the class of all morphisms:
 $$
 \varGamma\circ\varOmega=\Mor({\tt K})
 $$
 (this means that each morphism $\ph\in\Mor({\tt K})$ can be represented as a composition $\mu\circ\e$, where  $\mu\in\varGamma$, $\e\in\varOmega$).
 }\eit
If these conditions are fulfilled, we write
 \beq\label{K=varGamma-circledcirc-varOmega}
{\tt K}=\varGamma\circledcirc\varOmega.
 \eeq
}\eit

\bex\label{TH:faktorizatsija-v-kategorii-s-uzlov-razlozh-1}
In a category ${\tt K}$ with the nodal decomposition the following classes of morphisms define factorizations:
$$
{\tt K}=\Mono\circledcirc\SEpi=\SMono\circledcirc\Epi.
$$
\eex

The following fact is proved in \cite[Theorem 8.2]{Tsalenko-Shulgeifer}:

\btm\label{TH:o-faktorizatsii}
Classes $\varGamma$ and $\varOmega$ define a factorization of $\tt K$
$$
{\tt K}=\varGamma\circledcirc\varOmega
$$
if and only if the following conditions hold:
\bit{
 \item[(i)] $\varGamma\subseteq\Mono({\tt K})$ and $\varOmega\subseteq\Epi({\tt K})$,

 \item[(ii)] $\Iso({\tt K})\subseteq\varOmega\cap\varGamma$,

 \item[(iii)] for each morphism $\ph\in\Mor({\tt K})$ there is a decomposition
 \beq\label{faktorizatsiya-v-kat-s-faktoriz}
 \ph=\mu_{\ph}\circ\e_{\ph},\qquad \mu_{\ph}\in\varGamma,\quad \e_{\ph}\in\varOmega
 \eeq
 \item[(iv)] for any other decomposition with the same properties
$$
 \ph=\mu'\circ\e',\qquad \mu'\in\varGamma,\quad \e'\in\varOmega
$$
there is a morphism $\theta\in\Iso({\tt K})$ such that
$$
\mu'=\mu_{\ph}\circ\theta,\qquad \e'=\theta^{-1}\circ\e_{\ph}.
$$
 }\eit
\etm

 \bit{
\item Let us say that a class of morphisms $\varOmega$ in $\tt K$ is
{\it monomorphically complementable}\label{DEF:klass-monomorfno-dopolnyaem}, if
\beq\label{klass-monomorfno-dopolnyaem}
{\tt K}={^\downarrow\varOmega}\circledcirc\varOmega.
\eeq
In other words, $\varOmega$ must be epimorphic conjugate to its monomorphic conjugate class
$$
\varOmega=(^\downarrow\varOmega)^\downarrow,
$$
and the composition of ${^\downarrow\varOmega}$ and $\varOmega$ must cover the class of all morphisms:
$$
{^\downarrow\varOmega}\circ\varOmega=\Mor(\tt K).
$$
In this case the class ${^\downarrow\varOmega}$ will be called the {\it monomorphuc complement} to $\varOmega$.
}\eit

\brem
From \eqref{Iso-subseteq-varTheta^downarrow-subseteq-Epi} it follows that if a class of morphisms $\varOmega$ is monomorphically complementable, then
\beq\label{Iso-circ-varOmega-subseteq-varOmega}
\Iso\subseteq\varOmega\subseteq\Epi,\qquad \Iso\circ\varOmega\subseteq\varOmega
\eeq
\erem

 \bit{
\item Similarly, we say that the class of morphisms $\varGamma$ in $\tt K$ is
{\it epimorphically complementable}\label{DEF:klass-epimorfno-dopolnyaem}, if
\beq\label{klass-epimorfno-dopolnyaem}
{\tt K}=\varGamma\circledcirc\varGamma^\downarrow.
\eeq
In other words, $\varGamma$ must be the monomorphic conjugate to its epimorphic conjugate class
$$
\varGamma={^\downarrow(\varGamma^\downarrow)},
$$
and the composition of the classes $\varGamma$ and $\varGamma^\downarrow$ must cover the class of all morphisms:
$$
\varGamma\circ\varGamma^\downarrow=\Mor(\tt K).
$$
In this case the class $\varGamma^\downarrow$ will be called the {\it epimorphic complement} to $\varGamma$.
}\eit

\brem
From \eqref{Iso-subseteq-^downarrow-varTheta-subseteq-Mono} it follows that if a class $\varGamma$ is epimorphically complementable, then
\beq\label{varGamma-circ-Iso-subseteq-varGamma}
\Iso\subseteq\varGamma\subseteq\Mono,\qquad \varGamma\circ\Iso\subseteq\varGamma.
\eeq
\erem

\section{Envelope and refinement}\label{obolochka,otpechatok,uzl-razlozh}

\subsection{Envelope}

\paragraph{Envelope in a class of morphisms with respect to a class of morphisms.}
Suppose we have:
 \bit{

\item[---] a category ${\tt K}$ called an {\it enveloping category},

\item[---] a category ${\tt T}$ called an {\it attracting category},

\item[---] a covariant functor $F:{\tt T}\to{\tt K}$,

\item[---] two classes $\varOmega$ and $\varPhi$ of morphisms in ${\tt K}$, taking values in objects of the class $F({\tt T})$, and $\Omega$
is called the {\it class of realizing morphisms}, and $\varPhi$ the {\it class of test morphisms}. }\eit\noindent
Then

 \bit{
\item[$\bullet$] For given objects $X\in\Ob(\tt K)$ and $X'\in\Ob(\tt T)$ a morphism $\sigma:X\to F(X')$ is called an {\it extension of the object $X\in{\tt K}$ over the category ${\tt T}$ in the class of morphisms $\varOmega$ with respect to the class of morphisms $\varPhi$}, if $\sigma\in\varOmega$, and for any object $B$ in $\tt T$ and for each morphism $\ph:X\to F(B)$ from the class $\varPhi$ there exists a unique morphism $\ph':X'\to B$ in the category  ${\tt T}$ such that the following diagram is commutative:
 \beq\label{DEF:diagr-rasshirenie-T}
\begin{diagram}
\node[2]{X} \arrow{sw,t}{\varOmega\owns\sigma} \arrow{se,t}{\ph\in\varPhi}\\
\node{F(X')}  \arrow[2]{e,b,--}{F(\ph')} \node[2]{F(B)}
\end{diagram}
\eeq

\item[$\bullet$]\label{DEF:obolochka} An extension $\rho:X\to F(E)$ of an object $X\in{\tt K}$ over a category ${\tt T}$ in the class of morphisms  $\varOmega$ with respect to the class of morphisms $\varPhi$ is called an {\it envelope of the object $X\in{\tt K}$ over the category ${\tt T}$ in the class of morphisms $\varOmega$ with respect to the class of morphisms $\varPhi$},
if for each extension $\sigma:X\to F(X')$ (of the object $X\in{\tt K}$ over the category ${\tt T}$ in the class of morphisms $\varOmega$ with respect to the class of morphisms $\varPhi$) there exists a unique morphism $\upsilon:X'\to E$ in $\tt
T$ such that the following diagram is commutative
 \beq\label{DEF:diagr-obolochka-T}
\begin{diagram}
\node[2]{X} \arrow{sw,t}{\sigma} \arrow{se,t}{\rho}\\
\node{F(X')}  \arrow[2]{e,b,--}{F(\upsilon)} \node[2]{F(E)}
\end{diagram}
 \eeq
}\eit

In what follows we are almost exclusively interested in the case when ${\tt T}={\tt K}$, and $F:{\tt K}\to{\tt K}$ is the identity functor. It is useful to give the definitions for this case separately.

 \bit{
\item[$\bullet$] A morphism $\sigma:X\to X'$ in a category $\tt K$ is called an
{\it extension of the object $X\in\Ob({\tt K})$ in the class of morphisms $\varOmega$ with respect to the class of morphisms $\varPhi$}, if $\sigma\in\varOmega$, and for any morphism
$\ph:X\to B$ from the class $\varPhi$ there exists a unique morphism $\ph':X'\to B$ in
${\tt K}$ such that the following diagram is commutative:
 \beq\label{DEF:diagr-rasshirenie}
\begin{diagram}
\node[2]{X} \arrow{sw,t}{\varOmega\owns\sigma} \arrow{se,t}{\forall\ph\in\varPhi}\\
\node{X'}  \arrow[2]{e,b,--}{\exists!\ph'} \node[2]{B}
\end{diagram}
\eeq

\item[$\bullet$] An extension $\rho:X\to E$ of an object $X\in\Ob({\tt K})$ in the class of morphisms $\varOmega$
with respect to the class of morphisms $\varPhi$ is called an {\it envelope of $X$ in $\varOmega$ with respect to $\varPhi$}, if for any other extension $\sigma:X\to X'$ (of $X$ in $\varOmega$ with respect to $\varPhi$) there is a unique morphism $\upsilon:X'\to E$ in $\tt K$ such that the following diagram is commutative:
 \beq\label{DEF:diagr-obolochka}
\begin{diagram}
\node[2]{X} \arrow{sw,t}{\forall\sigma} \arrow{se,t}{\rho}\\
\node{X'}  \arrow[2]{e,b,--}{\exists!\upsilon} \node[2]{E}
\end{diagram}
 \eeq
For the morphism of envelope $\rho:X\to E$ we use the notation
    \beq\label{DEF:env_F^L(X)}
    \rho=\env_{\varPhi}^\varOmega X.
    \eeq
The very object $E$ is also called an {\it envelope} of $X$ (in $\varOmega$ with respect to $\varPhi$), and we use the following notation for it:
    \beq\label{DEF:Env_F^L(X)}
E=\Env_{\varPhi}^\varOmega X.
 \eeq
}\eit

\brem Clearly, the object $\Env_{\varPhi}^\varOmega X$ (if it exists) is defined up to an isomorphism. The question when the correspondence $X\mapsto\Env_{\varPhi}^\varOmega X$ can be defined as a functor is discussed below starting from page \pageref{DIAGR:funktorialnost-env-e-E}. \erem

\brem If $\varOmega=\varnothing$, then, of course neither extensions, nor envelopes in $\varOmega$ exist. So this construction can be interesting only when $\varOmega$ is a non-empty class. The following two situations will be of special interest
 \bit{
\item[---] $\varOmega=\Epi({\tt K})$ (i.e. $\varOmega$ coincides with the class of all epimorphisms in the category ${\tt K}$), then we will use the following notations
    \begin{align}\label{env_(varPhi)^Epi}
    &     \env_{\varPhi}^{\Epi}X:=\env_{\varPhi}^{\Epi({\tt K})}X, &&     \Env_{\varPhi}^{\Epi}X:=\Env_{\varPhi}^{\Epi({\tt K})}X.
    \end{align}

\item[---] $\varOmega=\Mor({\tt K})$ (i.e. $\varOmega$ coincides with the class of all morphisms in the category ${\tt K}$), in this case it is convenient to omit any mentioning about $\varOmega$ in the formulations and notations, so we will be speaking about the {\it envelope of object $X\in{\tt K}$ in the category ${\tt K}$ with respect to the class of morphisms $\varPhi$}, and the notations will be simplified as follows:
    \begin{align}\label{env_(varPhi)=env_(varPhi)^K}
    &     \env_{\varPhi}X:=\env_{\varPhi}^{\Mor({\tt K})}X, &&     \Env_{\varPhi}X:=\Env_{\varPhi}^{\Mor({\tt K})}X.
    \end{align}
 }\eit
\erem

\brem Another degenerate, but this time an informative case is when $\varPhi=\varnothing$. What is essential for a given object $X$, $\varPhi$ does not contain morphisms going from $X$:
$$
\varPhi^X=\{\ph\in\varPhi:\ \Dom\ph=X\}=\varnothing.
$$
Then, obviously, any morphism $\sigma\in\varOmega $, going from $X$, $\sigma:X\to X'$, is an extension for $X$ (in the class $\varOmega$ with respect to the class $\varnothing$). If in addition $\varOmega=\Epi$, then the envelope for $X$ is the terminal object in the category $\Epi^X$ (if it exists). This can be depicted by the formula
$$
\Env_\varnothing^\varOmega X=\max\Epi^X.
$$
In particular, if ${\tt K}$ is a category with zero $0$, and $\varOmega$ contains all morphisms going to $0$, then the envelope of any object with respect to the empty class of morphisms is $0$:
$$
\Env_\varnothing^\varOmega X=0.
$$
\erem

\brem Another extreme situation is when $\varPhi=\Mor({\tt K})$. For a given object $X$ the essential thing here is that the class $\varPhi$ contains the local identity of $X$:
$$
1_X\in\varPhi.
$$
For any extension $\sigma$ the diagram
$$
\xymatrix @R=2.pc @C=2.0pc % @M=14pt
{
X \ar[rr]^{\sigma} \ar[dr]_{1_X} & & X'\ar@{-->}[dl]\\
 & X &
}
$$
implies that $\sigma$ must be a co-retraction (moreover, the dotted arrow here must be unique). When $\varOmega\subseteq\Epi$ this is possible only if $\sigma$ is an isomorphism. As a corollary, in this case the envelope of $X$ coincides with $X$ (up to an isomorphism):
$$
\varOmega\subseteq\Epi\quad\Longrightarrow\quad\Env_{\Mor({\tt K})}^{\varOmega} X=X.
$$
\erem

\medskip
\centerline{\bf Properties of envelopes:}

\bit{\it

\item[$1^\circ$.]\label{LM:suzhenie-verh-klassa-morfizmov} Suppose that $\varSigma\subseteq\varOmega$, then for any object $X$ and for any class of morphisms $\varPhi$

\bit{

\item[(a)] each extension $\sigma:X\to X'$ in $\varSigma$ with respect to $\varPhi$ is an extension in $\varOmega$ with respect to $\varPhi$,

\item[(b)] if there are envelopes $\env_\varPhi^\varSigma X$ and $\env_\varPhi^\varOmega X$, then there is a unique morphism  $\rho:\Env_\varPhi^\varSigma X\to\Env_\varPhi^\varOmega X$ such that the following diagram is commutative:
\beq\label{suzhenie-verh-klassa-morfizmov}
\begin{diagram}
\node[2]{X} \arrow{sw,t}{\env_\varPhi^\varSigma X} \arrow{se,t}{\env^\varOmega _\varPhi X}\\
\node{\Env_\varPhi^\varSigma X}\arrow[2]{e,b,--}{\rho}   \node[2]{\Env^\varOmega _\varPhi X}
\end{diagram}
\eeq

\item[(c)] if there is an envelope $\env_\varPhi^\varOmega X$ (in a wider class), and it lies in $\varSigma$ (i.e. in a narrower class),
$$
\env_\varPhi^\varOmega X\in\varSigma,
$$
then it is an envelope $\env_\varPhi^\varSigma X$ (in a narrower class):
$$
\env_\varPhi^\varOmega X=\env_\varPhi^\varSigma X.
$$
}\eit

\item[$2^\circ$.]\label{LM:suzhenie-verh-klassa-morfizmov-2}
Let $\varSigma$, $\varOmega$, $\varPhi$ be classes of morphisms, and for an object $X$
 \bit{

\item[(a)] every extension $\sigma:X\to X'$ in $\varOmega$ with respect to $\varPhi$ belongs to $\varSigma$.
 }\eit
Then
 \bit{

\item[(b)] an envelope of $X$ with respect to $\varPhi$ in the class $\varOmega$ exists if and only if there exists an envelope of $X$ with respect to $\varPhi$ in the class $\varOmega\cap\varSigma$, and these envelopes coincide:
$$
\env_{\varPhi}^{\varOmega}=\env_{\varPhi}^{\varOmega\cap\varSigma},
$$

\item[(c)] if $\varSigma\subseteq\varOmega$, then an envelope of $X$ with respect to $\varPhi$ in the (narrower) class $\varSigma$ exists if and only if there exists an envelope of $X$ with respect to $\varPhi$ in the (wider) class $\varOmega$, and these envelopes coincide:
$$
\env_{\varPhi}^{\varOmega}X=\env_{\varPhi}^{\varSigma}X.
$$
 }\eit

\item[$3^\circ$.]\label{LM:suzhenie-klassa-morfizmov} Suppose $\varPsi\subseteq\varPhi$, then for any object $X$ and for any class of morphisms $\varOmega$

\bit{

\item[(a)] each extension $\sigma:X\to X'$ in $\varOmega$ with respect to $\varPhi$ is an extension in $\varOmega$ with respect to $\varPsi$,

\item[(b)] if there are envelopes $\env_\varPsi^\varOmega X$ and $\env_\varPhi^\varOmega X$, then there is a unique morphism $\alpha:\Env_\varPsi^\varOmega X\gets\Env_\varPhi^\varOmega X$ such that the following diagram is commutative:  \beq\label{suzhenie-klassa-morfizmov}
\begin{diagram}
\node[2]{X} \arrow{sw,t}{\env_\varPsi^\varOmega X} \arrow{se,t}{\env^\varOmega _\varPhi X}\\
\node{\Env_\varPsi^\varOmega X}   \node[2]{\Env^\varOmega _\varPhi X} \arrow[2]{w,b,--}{\alpha}
\end{diagram}
\eeq
}\eit

\item[$4^\circ$.]\label{TH:env_Psi=env_Phi} Suppose that $\varPhi\subseteq\Mor({\tt K})\circ\ \varPsi$ (i.e. each morphism $\ph\in\varPhi$ can be represented as a composition $\ph=\chi\circ\psi$, where $\psi\in\varPsi$), then for any object $X$ and for any class of morphisms $\varOmega$
\bit{\it

\item[(a)] if an extension $\sigma:X\to X'$ in $\varOmega$ with respect to $\varPsi$ is at the same time an epimorphism in ${\tt K}$, then it is an extension in $\varOmega$ with respect to $\varPhi$,

\item[(b)] if there are envelopes $\env_\varPsi^\varOmega X$ and $\env_\varPhi^\varOmega X$, and $\env_\varPsi^\varOmega X$ is at the same time an epimorphism in ${\tt K}$, then there exists a unique morphism $\beta:\Env_\varPsi^\varOmega X\to\Env_\varPhi^\varOmega X$ such that the following diagram is commutative:
 \beq\label{env_Psi=env_Phi-0}
\begin{diagram}
\node[2]{X} \arrow{sw,t}{\env_\varPsi^\varOmega X} \arrow{se,t}{\env^\varOmega _\varPhi X}\\
\node{\Env_\varPsi^\varOmega X}\arrow[2]{e,b,--}{\beta}   \node[2]{\Env^\varOmega _\varPhi X}
\end{diagram}
\eeq
}\eit

\item[$5^\circ$.]\label{PROP:deistvie-epimorfizma-na-Env} Suppose that $\varOmega$ and $\varPhi$ are some classes of morphisms, and $\e:X\to Y$ an epimorphism in ${\tt K}$ such that the following three conditions are fulfilled:
 \bit{
\item[(a)] there exists an envelope $\env_{\varPhi\circ\e}^\varOmega X$ with respect to the class of morphisms $\varPhi\circ\e=\{\ph\circ\e;\ \ph\in\varPhi\}$,

\item[(b)] there exists an envelope $\env_\varPhi^\varOmega Y$,

\item[(c)] the composition $\env_{\varPhi}^\varOmega Y\circ\e$ belongs to $\varOmega$.

    }\eit
Then there exists a unique morphism $\upsilon:\Env_{\varPhi\circ\e}^\varOmega X\gets \Env_{\varPhi}^\varOmega Y$ such that the following diagram is commutative:
\beq\label{deistvie-epimorfizma-na-Env}
\xymatrix @R=2.5pc @C=4.0pc
{
 X\ar@/_1ex/[rd]^{\ \env_{\varPhi}^\varOmega Y\circ\e}\ar[d]_{\env_{\varPhi\circ\e}^\varOmega X}\ar[r]^{\e} & Y\ar[d]^{\env_{\varPhi}^\varOmega Y} \\
 \Env_{\varPhi\circ\e}^\varOmega X &   \Env_{\varPhi}^\varOmega Y\ar@{-->}[l]^{\upsilon}
}
\eeq

}\eit

\bpr
1. If the morphism $\sigma$ satisfies \eqref{DEF:diagr-rasshirenie} with $\varSigma$ instead of $\varOmega$, then $\sigma$ satisfies the initial condition \eqref{DEF:diagr-rasshirenie}, since $\varSigma\subseteq\varOmega$. This proves (a). From this we have also that $\env_\varPsi^\varSigma X$ is an extension in $\varOmega$ with respect to $\varPhi$, so there must exist a unique dotted arrow in \eqref{suzhenie-verh-klassa-morfizmov}. This means that (b) is also true. Finally, if there exists an envelope $\env_\varPhi^\varOmega X$ (in a wider class), and it lies in $\varSigma$ (in a narrower class), $\env_\varPhi^\varOmega X\in\varSigma$, then $\env_\varPhi^\varOmega X$ is an extension in $\varSigma$. On the other hand any other extension $\sigma:X\to X'$ in $\varSigma$ is an extension in $\varOmega$ due to the property (a) which we have just proved, hence there is a unique morphism $\upsilon$ into the envelope in $\varOmega$:
$$
\begin{diagram}
\node[2]{X} \arrow{sw,t}{\sigma} \arrow{se,t}{\env_\varPhi^\varOmega X}\\
\node{X'}  \arrow[2]{e,b,--}{\upsilon} \node[2]{\Env_\varPhi^\varOmega X}
\end{diagram}
$$
This proves that $\env_\varPhi^\varOmega X$ is an envelope in $\varSigma$, and we have proved (c).

2. Suppose $2^\circ$(a) holds. If an object $X$ has an envelope
$\env_{\varPhi}^{\varOmega}X$ in the class $\varOmega$ with respect to $\varPhi$, then by (a),
this will be an extension in a narrower class $\varOmega\cap\varSigma$ with respect to $\varPhi$. Applying $1^\circ$(c) on page \pageref{LM:suzhenie-verh-klassa-morfizmov}, we obtain that $\env_{\varPhi}^{\varOmega}X$ is an envelope in the narrower class $\varOmega\cap\varSigma$, i.e.
$\env_{\varPhi}^{\varOmega}X=\env_{\varPhi}^{\varOmega\cap\varSigma}X$.

Conversely, suppose there exists an envelope $\env_{\varPhi}^{\varOmega\cap\varSigma}X$ with respect to $\varPhi$ in the class $\varOmega\cap\varSigma$. Then by $1^\circ$ (a) on p.\pageref{LM:suzhenie-verh-klassa-morfizmov}, it will be an envelope with respect to $\varPhi$ in the wider class $\varOmega$. Take another extension $\sigma:X\to X'$ with respect to $\varPhi$ in $\varOmega$. By (a), $\sigma$ is an extension with respect to $\varPhi$ in the class $\varOmega\cap\varSigma$. Hence, there exists a unique morphism $\upsilon:X'\to \Env_{\varPhi}^{\varOmega\cap\varSigma} X$ into the envelope in $\varOmega\cap\varSigma$, such that the following diagram is commutative:
$$
\begin{diagram}
\node[2]{X} \arrow{sw,t}{\sigma} \arrow{se,t}{\env_{\varPhi}^{\varOmega\cap\varSigma} X}\\
\node{X'}  \arrow[2]{e,b,--}{\upsilon} \node[2]{\Env_{\varPhi}^{\varOmega\cap\varSigma} X}
\end{diagram}
$$
This proves that $\env_{\varPhi}^{\varOmega\cap\varSigma}X$ is (not just an extension, but also) an envelope with respect to $\varPhi$ in $\varOmega$. We see that $2^\circ$ (b) is true, and $2^\circ$(c) is its corollary.

3. Suppose that $\varPsi\subseteq\varPhi$. Then (a) is obvious: each extension $\sigma:X\to X'$ with respect to $\varPhi$ is an extension with respect to the narrower class $\varPsi$. For (b) we have: since $\env_\varPhi^\varOmega X$ is an extension with respect to $\varPhi$, it must be an extension with respect to the narrower class $\varPsi$, so there exists a unique morphism from $\Env_\varPhi^\varOmega X$ into the envelope  $\Env_\varPsi^\varOmega X$ with respect to $\varPsi$ such that \eqref{suzhenie-klassa-morfizmov} is commutative.

4. Suppose that $\varPhi\subseteq\Mor({\tt K})\circ\ \varPsi$. For (a) our reasoning will be illustrated by the following diagram:
$$
\xymatrix @R=2.pc @C=5.0pc % @M=14pt
{
X \ar[rr]^{\sigma} \ar[dr]^{\psi}\ar@/_4ex/[ddr]_{\ph} & & X' \ar@{-->}[dl]_-{\psi'} \ar@{-->}@/^4ex/[ddl]^{\ph'}
\\
 & Y\ar[d]^{\chi} & \\
 & B &
}
$$
If $\sigma:X\to X'$ is an extension of $X$ in $\varOmega$ with respect to $\varPsi$, then for any morphism $\ph\in\varPhi$, $\ph:X\to B$, we take a decomposition $\ph=\chi\circ\psi$, where $\psi\in\varPsi$. Since $\sigma$ is an extension of $X$ in $\varOmega$ with respect to $\varPsi$, there is a morphism $\psi'$ such that $\psi=\psi'\circ\sigma$. After that we put $\ph'=\chi\circ\psi'$, and this will be a morphism such that
$$
\ph=\chi\circ\psi=\chi\circ\psi'\circ \sigma=\ph'\circ \sigma.
$$
The uniqueness of $\ph'$ follows from the epimorphy of $\sigma\in\varOmega$, and thus $\sigma$ is an extension of $X$ in $\varOmega$ with respect to  $\varPhi$. Once (a) is proved, (b) becomes its corollary: the morphism $\env^\varOmega_\varPsi X:X\to \Env^\varOmega _\varPsi X$ is an extension of  $X$ in $\varOmega$ with respect to $\varPsi$, hence, by (a), with respect to $\varPhi$ as well. So there must exist a morphism $\beta$ from  $\Env^\varOmega _\varPsi X$ into the envelope $\Env^\varOmega _\varPhi X$ with respect to $\varPhi$ such that \eqref{env_Psi=env_Phi-0} is commutative.

5. For any morphism $\ph:Y\to B$ lying in $\varPhi$ we have the following diagram:
$$
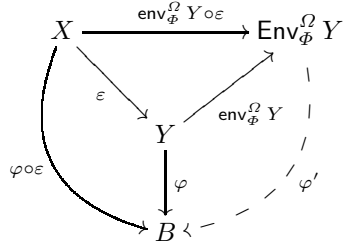
\xymatrix %@R=2.5pc @C=4.0pc
{
X\ar[rr]^{\env_{\varPhi}^\varOmega Y\circ\e}\ar[rd]_{\e}\ar@/_6ex/[rdd]_{\ph\circ\e} &  & \Env_{\varPhi}^\varOmega Y\ar@{-->}@/^6ex/[ldd]^{\ph'} \\
 & Y\ar[d]^{\ph}\ar[ru]_(.4){\env_{\varPhi}^\varOmega Y} & \\
 &  B  &
}
$$
It must be understood as follows. On the one hand, since $\env_{\varPhi}^\varOmega Y$ is an extension with respect to $\varPhi$, there exists a morphism $\ph'$, such that the lower right triangle is commutative, and as a corollary, the perimeter is commutative as well. On the other hand, if $\ph'$ is a morphism such that the perimeter is commutative, i.e.
$$
\ph'\circ\env_{\varPhi}^\varOmega Y\circ\e=\ph\circ\e,
$$
then, since $\e$ is an epimorphism, we can cancel it:
$$
\ph'\circ\env_{\varPhi}^\varOmega Y=\ph,
$$
So the lower right triangle must be commutative as well. This means that $\ph'$ must be unique (since by the definition of envelope, the dotted arrow in the lower right triangle must be unique).

We see that the perimeter has a unique dotted arrow $\ph'$. This is true for any $\ph\in\varPhi$, and in addition
$\env_{\varPhi}^\varOmega Y\circ\e\in\varOmega$. So we come to a conclusion that the morphism $\env_{\varPhi}^\varOmega Y\circ\e$ is an extension of  $X$ in $\varOmega$ with respect to the class of morphisms $\varPhi\circ\e$. As a corollary, there exists a unique morphism $\upsilon$ from
$\Env_{\varPhi}^\varOmega Y$ into the envelope $\Env_{\varPhi\circ\e}^\varOmega X$ with respect to $\varPhi\circ\e$ such that
\eqref{deistvie-epimorfizma-na-Env} is commutative. \epr

 \bit{
\item[$\bullet$] Let us say that in a category ${\tt K}$ {\it a class of morphisms $\varPhi$ is generated on the inside by a class of morphisms  $\varPsi$}, if
    \beq\label{DEF:morfizmy-porozhdayutsya-iznutri}
    \varPsi\subseteq\varPhi\subseteq\Mor({\tt K})\circ\ \varPsi.
    \eeq
  }\eit
\btm\label{TH:morfizmy-porozhdayutsya-iznutri} Suppose that in a category ${\tt K}$ a class of morphisms $\varPhi$ is generated on the inside by a class of morphisms  $\varPsi$. Then for any class of epimorphisms $\varOmega$ (it  is not necessary that $\varOmega$ contains all epimorphisms of ${\tt K}$) and for any object $X$ the existence of envelope $\env^\varOmega_\varPsi X$ is equivalent to the existence of envelope $\env^\varOmega_\varPhi X$, and these envelopes coincide:
 \beq\label{env_Psi=env_Phi}
\env^\varOmega _\varPsi X=\env^\varOmega _\varPhi X.
 \eeq
  \etm
\bpr 1. Suppose first that $\env^\varOmega_\varPsi X$ exists. Since it is an extension with respect to $\varPsi$, and at the same time an epimorphism, by $4^\circ$ (a) we have that it is an extension with respect to $\varPhi$ as well. If $\sigma:X\to X'$ is another extension with respect to
$\varPhi$, then by $3^\circ$ (a) it is an extension with respect to $\varPsi$ as well, so there exists a unique morphism
$\upsilon:\Env_\varPsi^\varOmega X\gets X'$ such that the following diagram is commutative:
$$
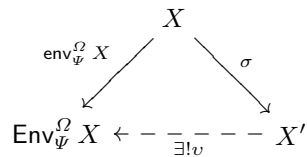
\begin{diagram}
\node[2]{X} \arrow{sw,t}{\env_\varPsi^\varOmega X} \arrow{se,t}{\sigma}\\
\node{\Env_\varPsi^\varOmega X}   \node[2]{X'}
\arrow[2]{w,b,--}{\exists!\upsilon}
\end{diagram}
$$
This means that $\env^\varOmega _\varPsi X$ is an envelope with respect to $\varPhi$, and \eqref{env_Psi=env_Phi} holds.

2. On the contrary, suppose that $\env^\varOmega _\varPhi X$ exists. It is an extension with respect to $\varPhi$, so by $3^\circ$ (a) it must be an extension with respect to $\varPsi$ as well. If $\sigma:X\to X'$ is another extension in $\varOmega$ with respect to $\varPsi$, then, since  $\sigma\in\Epi$, by $4^\circ$ (a), it must be an extension with respect to $\varPhi$, so there exists a unique morphism
$\upsilon:X'\to \Env_\varPhi^\varOmega X$ such that
$$
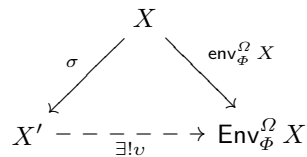
\begin{diagram}
\node[2]{X} \arrow{sw,t}{\sigma} \arrow{se,t}{\env_\varPhi^\varOmega X}\\
\node{X'} \arrow[2]{e,b,--}{\exists!\upsilon} \node[2]{\Env_\varPhi^\varOmega X}
\end{diagram}
$$
This means that $\env^\varOmega _\varPhi X$ is an envelope with respect to $\varPsi$, and again we have \eqref{env_Psi=env_Phi}.
\epr

\bit{
\item[$\bullet$]\label{DEF:varPhi-razlich-morfizmy-snaruzhi} Let us say that a class of morphisms $\varPhi$ in a category ${\tt K}$  {\it differs morphisms on the outside}, if for any two different parallel morphisms $\alpha\ne\beta:X\to Y$ there is a morphism $\ph:Y\to M$ from the class $\varPhi$ such that $\ph\circ\alpha\ne\ph\circ\beta$.
}\eit

\btm\label{TH:Phi-razdel-moprfizmy} If a class of morphisms $\varPhi$ differs morphisms on the outside, then for any class of morphisms $\varOmega$
 \bit{
\item[(i)] each extension in $\varOmega$ with respect to $\varPhi$ is a monomorphism,

\item[(ii)] an envelope with respect to $\varPhi$ in $\varOmega$ exists if and only if there exists an envelope with respect to $\varPhi$ in the class $\varOmega\cap\Mono$; in this case these envelopes coincide:
$$
\env_{\varPhi}^{\varOmega}=\env_{\varPhi}^{\varOmega\cap\Mono},
$$

\item[(iii)] if the class $\varOmega$ contains all monomorphisms,
 $$
 \varOmega\supseteq\Mono,
 $$
then the existence of the envelope with respect to $\varPhi$ in $\Mono$ automatically implies the existence of envelope with respect to $\varPhi$ in $\varOmega$, and the coincidence of these envelopes:
$$
\env_{\varPhi}^{\varOmega}=\env_{\varPhi}^{\Mono}.
$$

 }\eit
\etm
\bpr
1. Suppose that some extension $\sigma:X\to X'$ is not a monomorphism, i.e. there are two different parallel morphisms $\alpha\ne\beta:T\to X$ such that
\beq\label{sigma-circ-alpha=sigma-circ-beta}
\sigma\circ\alpha=\sigma\circ\beta.
\eeq
Since the class $\varPhi$ differs morphisms on the outside, there must exist a morphism $\ph:X\to M$, $\ph\in\varPhi$, such that
\beq\label{ph-circ-alpha-ne-ph-circ-beta}
\ph\circ\alpha\ne\ph\circ\beta.
\eeq
Since $\sigma:X\to X'$ is an extension with respect to $\varPhi$, there is a continuation $\ph':X'\to M$ of the morphism $\ph:X\to M$: $\ph=\ph'\circ\sigma$. Now we obtain
$$
\ph\circ\alpha=\ph'\circ\sigma\circ\alpha=\eqref{sigma-circ-alpha=sigma-circ-beta}=\ph'\circ\sigma\circ\beta=\ph\circ\beta,
$$
and this contradicts to \eqref{ph-circ-alpha-ne-ph-circ-beta}.

2. Suppose for an object $X$ there exists an envelope $\env_{\varPhi}^{\varOmega}X$ in $\varOmega$ with respect to $\varPhi$. Then, as we have already proved, it must be an extension in the narrower class $\varOmega\cap\Mono$ with respect to $\varPhi$. Applying property $1^\circ$ (c) on p.\pageref{LM:suzhenie-verh-klassa-morfizmov}, we obtain that $\env_{\varPhi}^{\varOmega}X$ is an envelope in the narrower class  $\varOmega\cap\Mono$, i.e. $\env_{\varPhi}^{\varOmega}X=\env_{\varPhi}^{\varOmega\cap\Mono}X$.

Conversely, suppose there is an envelope $\env_{\varPhi}^{\varOmega\cap\Mono} X$ with respect to $\varPhi$ in the class $\varOmega\cap\Mono$. By  $1^\circ$ (a) on p.\pageref{LM:suzhenie-verh-klassa-morfizmov}, it must be an extension with respect to $\varPhi$ in the wider class $\varOmega$. Consider another extension $\sigma:X\to X'$ with respect to $\varPhi$ in $\varOmega$. By the proposition (i) which we have already proved, $\sigma$ is an extension with respect to $\varPhi$ in $\varOmega\cap\Mono$. Hence, there is a unique morphism $\upsilon:X'\to \Env_{\varPhi}^{\varOmega\cap\Mono} X$ into the envelope in the class $\varOmega\cap\Mono$, such that
$$
\begin{diagram}
\node[2]{X} \arrow{sw,t}{\sigma} \arrow{se,t}{\env_{\varPhi}^{\varOmega\cap\Mono} X}\\
\node{X'}  \arrow[2]{e,b,--}{\upsilon} \node[2]{\Env_{\varPhi}^{\varOmega\cap\Mono} X}
\end{diagram}
$$
This proves that $\env_{\varPhi}^{\varOmega\cap\Mono} X$ is (not only an extension, but also) an envelope with respect to $\varPhi$ in the class $\varOmega$.

3. The proposition (iii) immediately follows from (ii).
\epr

\bit{
\item[$\bullet$] Let us remind that a class of morphisms $\varPhi$ in a category ${\tt K}$ is called a {\it right ideal}, if
$$
\varPhi\circ\Mor({\tt K})\subseteq\varPhi
$$
(i.e. for any $\ph\in\varPhi$ and for any morphism $\mu$ in ${\tt K}$ the composition $\ph\circ\mu$ belongs to $\varPhi$).
}\eit

\btm\label{TH:Phi-razdel-moprfizmy-*} If a class of morphisms $\varPhi$ differs morphisms on the outside and is a right ideal in the category ${\tt K}$, then for any class of morphisms $\varOmega$
 \bit{
\item[(i)] each extension in $\varOmega$ with respect to $\varPhi$ is a bimorphism,

\item[(ii)] an envelope with respect to $\varPhi$ in $\varOmega$ exists if an only if there exists an envelope with respect to $\varPhi$ in the class $\varOmega\cap\Bim$ of bimorphisms belonging to $\varOmega$; in this case these envelopes coincide:
$$
\env_{\varPhi}^{\varOmega}=\env_{\varPhi}^{\varOmega\cap\Bim}.
$$

\item[(iii)] if the class $\varOmega$ contains all bimorphisms,
 $$
 \varOmega\supseteq\Bim,
 $$
 then an envelope with respect to $\varPhi$ in $\varOmega$ exists if an only if there exists an envelope with respect to $\varPhi$ in $\Bim$, and these envelopes coincide:
$$
\env_{\varPhi}^{\varOmega}=\env_{\varPhi}^{\Bim}.
$$

 }\eit
\etm
\bpr
By property $2^\circ$ on p.\pageref{LM:suzhenie-verh-klassa-morfizmov-2}, (ii) and (iii) follow from (i), so we need to prove only (i). Let $\sigma:X\to X'$ be an extension in $\varOmega$ with respect to $\varPhi$. By Theorem \ref{TH:Phi-razdel-moprfizmy}(i), $\sigma$ must be a monomorphism. Suppose that it is not an epimorphism. This means that there are two different parallel morphisms  $\alpha\ne\beta:X'\to T$ such that
\beq\label{alpha-circ-sigma=beta-circ-sigma}
\alpha\circ\sigma=\beta\circ\sigma.
\eeq
Since $\varPhi$ differs morphisms on the outside, there must exist a morphism $\ph:T\to M$, $\ph\in\varPhi$, such that
$$
\ph\circ\alpha\ne\ph\circ\beta.
$$
In addition, by \eqref{alpha-circ-sigma=beta-circ-sigma},
$$
\ph\circ\alpha\circ\sigma=\ph\circ\beta\circ\sigma.
$$
If we now suppose that $\varPhi$ is a right ideal in ${\tt K}$, then the morphism $\ph\circ\alpha\circ\sigma=\ph\circ\beta\circ\sigma$ lies in $\varPhi$. So we can interpret this picture as follows: the test (i.e. lying in $\varPhi$) morphism  $\ph\circ\alpha\circ\sigma=\ph\circ\beta\circ\sigma:X\to M$ has two different continuations $\ph\circ\alpha\ne\ph\circ\beta:X'\to M$ along  $\sigma:X\to X'$. This means that $\sigma$ cannot be an extension with respect to $\varPhi$.
\epr

\paragraph{Envelope in a class of objects with respect to a class of objects.}

A  special case of the construction is a situation, where $\varOmega$ and/or $\varPhi$ are classes of all morphisms into the objects from some given subclasses in $\Ob({\tt K})$. The accurate formulation for the case, when both classes $\varOmega$ and $\varPhi$ are defined in such a way is the following. Suppose we have a category ${\tt K}$ and two subclasses ${\tt L}$ and ${\tt M}$ in the class $\Ob({\tt K})$ of objects in ${\tt K}$.

 \bit{
\item[$\bullet$] A morphism $\sigma:X\to X'$ is called an {\it extension of the object
$X\in{\tt K}$ in the class ${\tt L}$ with respect to the class ${\tt M}$}, if $X'\in{\tt L}$ and for any object
$B\in{\tt M}$ and any morphism $\ph:X\to B$ there exists a unique morphism $\ph':X'\to B$ such that the following diagram is commutative:
$$
\put(6,-35){$\begin{matrix}\text{\rotatebox{90}{$\owns$}}\\ {\tt L}
\end{matrix}$} \put(73,-35){$\begin{matrix}\text{\rotatebox{90}{$\owns$}}\\
{\tt M} \end{matrix}$}
\begin{diagram}
\node[2]{X} \arrow{sw,t}{\sigma} \arrow{se,t}{\forall\ph}\\
\node{X'}  \arrow[2]{e,b,--}{\exists!\ph'} \node[2]{B}
\end{diagram}
$$

\item[$\bullet$]\label{DEF:obolochka-otn-klassa} An extension $\rho:X\to E$ of an object $X\in{\tt K}$ in the class ${\tt L}$ with respect to the class ${\tt M}$ is called an {\it envelope of the object $X\in{\tt K}$ in the class ${\tt L}$ with respect to the class ${\tt M}$}, and we denote this by formula
    \beq\label{DEF:env_M^L(A)}
    \rho=\env_{\tt M}^{\tt L} X,
    \eeq
    if for any  other extension $\sigma:X\to X'$ (of the object $X$ in the class ${\tt L}$ with respect to the class ${\tt M}$) there exists a unique morphism $\upsilon:X'\to E$ such that the following diagram is commutative:
 \beq\label{diagr:obolochka}
\put(6,-35){$\begin{matrix}\text{\rotatebox{90}{$\owns$}}\\ {\tt L}
\end{matrix}$} \put(73,-35){$\begin{matrix}\text{\rotatebox{90}{$\owns$}}\\
{\tt L} \end{matrix}$}
\begin{diagram}
\node[2]{X} \arrow{sw,t}{\forall\sigma} \arrow{se,t}{\rho}\\
\node{X'}  \arrow[2]{e,b,--}{\exists!\upsilon} \node[2]{E}
\end{diagram}
 \eeq
The object $E$ is also called an {\it envelope} of the object $X$ (in the class of objects ${\tt L}$ with respect to the class of objects ${\tt M}$), and we will use the following notation for it:
 \beq\label{DEF:Env_M^L(A)}
E=\Env_{\tt M}^{\tt L} X.
 \eeq

}\eit

The following two extreme situations in the choice of the subclass ${\tt L}$ can occur:
 \bit{

\item[---] if ${\tt L} =\Ob({\tt K})$, then we will speak about an {\it envelope of the object $X\in{\tt K}$ in the category ${\tt K}$ with respect to the class of objects ${\tt M}$}, and the notations will be the following:
    \begin{align}\label{env_M=env_M^K}
    &     \env_{\tt M}X:=\env_{\tt M}^{\tt K}X, &&     \Env_{\tt M}X:=\Env_{\tt M}^{\tt K}X.
    \end{align}

\item[---]\label{L=M} if ${\tt L} ={\tt M}$, then the notions of the extension and the envelope coincide: {\it each extension of the object $X$ in the class ${\tt L}$ with respect to the same class ${\tt L}$ is an envelope of $X$ in the class ${\tt L}$ with respect to ${\tt L}$} (indeed, if  $\rho:X\to E$ and $\sigma:X\to X'$ are two extensions in ${\tt L}$ with respect to ${\tt L}$, then in diagram \eqref{diagr:obolochka} the morphism $\upsilon$ exists and is unique just because $\sigma$ is an extension); for simplicity, in case of ${\tt L} ={\tt M}$ we speak about the {\it envelope of $X$ in the class ${\tt L}$}, and our notations are simplified as follows:
    \begin{align}
    & \env^{\tt L} X=:\env_{\tt L}^{\tt L} X, &&  \Env^{\tt L} X=:\Env_{\tt L}^{\tt L} X.
    \end{align}

}\eit

\bit{
\item[$\bullet$] Let us say that a class of objects ${\tt M}$ in a category ${\tt K}$
 {\it differs morphisms on the outside}, if the class of morphisms with ranges in ${\tt M}$ possesses this property  (in the sense of definition on page \pageref{DEF:varPhi-razlich-morfizmy-snaruzhi}), i.e. for any two different parallel morphisms $\alpha\ne\beta:X\to Y$ there is a morphism $\ph:Y\to M\in{\tt M}$ such that $\ph\circ\alpha\ne\ph\circ\beta$.
}\eit

From Theorem \ref{TH:Phi-razdel-moprfizmy-*} we have

\btm\label{TH:M-razdel-moprfizmy} If a class of objects ${\tt M}$ differs morphisms on the outside, then for any class of objects ${\tt L}$
 \bit{
\item[(i)] each envelope in ${\tt L}$ with respect to ${\tt M}$ is a bimorphism,

\item[(ii)] an envelope in ${\tt L}$ with respect to ${\tt M}$ exists if and only if there exists an anvelope in the class of bimorphisms with the values in ${\tt L}$ with respect to ${\tt M}$; in this case these envelopes coincide:
$$
\env_{\tt M}^{\tt L}=\env_{\tt M}^{\Bim({\tt K},{\tt L})}.
$$
 }\eit
\etm

\paragraph{Examples of envelopes.}

\bex {\bf Universal enveloping algebra.} Let ${\tt K}={\tt LieAlg}$ be the category of Lie algebras (say, over the field $\C$), ${\tt T}={\tt Alg}$ the category of associative algebras (again, over $\C$) with the identity, and $F:{\tt Alg}\to{\tt LieAlg}$ the functor that every associative algebra $A$ represents as the Lie algebra with the Lie brackets
$$
[x,y]=x\cdot y-y\cdot x.
$$
Then the envelope of a Lie algebra $\frak{g}$ over the category ${\tt Alg}$ in the class $\Mor({\tt LieAlg},F({\tt Alg}))$ of all morphisms from  ${\tt LieAlg}$ into $F({\tt Alg})$ with respect to the same class $\Mor({\tt LieAlg},F({\tt Alg}))$ is exactly the universal enveloping algebra $U(\frak{g})$ of the Lie algebra $\frak{g}$
(cf.\cite{Bourbaki-LieAlg}):
$$
U(\frak{g})=\Env^{\Mor({\tt LieAlg},F({\tt Alg}))}
\frak{g}.
$$
 \eex

\bex {\bf Stone-\v{C}ech compactification.} In the category ${\tt Tikh}$ of Tikhonov spaces the Stone-\v{C}ech compactification $\beta:X\to\beta X$
is an envelope of the space $X$ in the class ${\tt Com}$ of compact spaces with respect to the same class ${\tt Com}$:
$$
\beta X=\Env^{\tt Com}X.
$$
\eex \bpr Here one uses the theorem \cite[Theorem 3.6.1]{Engelking}, which states that any continuous map $f:X\to K$ into an arbitrary compact space $K$ can be extended to a continuous map $F:\beta X\to K$. Since $\beta(X)$ is dense in $\beta X$, this extension $F$ is unique, and this means that  $\beta:X\to\beta X$ is an extension in the class ${\tt Com}$ with respect to ${\tt Com}$. By the Remark on p.\pageref{L=M}, in the case ${\tt
L}={\tt M}$ each extension is an envelope, so $\beta$ is an envelope. \epr

\bex {\bf Completion} $X^\blacktriangledown$ of a locally convex space $X$ is an envelope of $X$ in the category ${\tt LCS}$ of all locally convex spaces with respect to the class ${\tt Ban}$ of Banach spaces:
$$
X^\blacktriangledown=\Env_{\tt Ban}^{\tt LCS}X.
$$
\eex \bpr Let us denote the natural embedding of $X$ into its completion by $\blacktriangledown_X:X\to X^\blacktriangledown$ (we use the notations of \cite{Akbarov}).

First, each linear continuous map $f:X\to B$ into an arbitrary Banach space $B$ is uniquely extended to a linear continuous map  $F:X^\blacktriangledown\to B$ on the completion $X^\blacktriangledown$ of $X$ (one can refer here,  for instance, to the general theorem for all uniform spaces \cite[Theorem 8.3.10]{Engelking}). Hence, the completion $\blacktriangledown_X:X\to X^\blacktriangledown$ is an extension of the space $X$ in the category ${\tt LCS}$ of locally convex spaces with respect to the subclass ${\tt Ban}$ of Banach spaces.

Note then, that the class ${\tt Ban}$ of Banach spaces differs morphisms on the outside in the category ${\tt LCS}$. By Theorem
\ref{TH:M-razdel-moprfizmy} this means that any extension $\sigma:X\to X'$ with respect to ${\tt Ban}$ is a bimorphism in ${\tt LCS}$, i.e. is an injective map which image $\sigma(X)$ is dense in $X'$. Let us show that in addition it is an open map: for any neighborhood of zero  $U\subseteq X$ there is a neighborhood of zero $V\subseteq X'$ such that
 \beq\label{sigma(U)-supseteq-V-cap-sigma(X)}
\sigma(U)\supseteq V\cap\sigma(X)
 \eeq
We can think that $U$ is a closed convex neighborhood of zero in $X$. Then the set $\Ker U=\bigcap_{\varepsilon>0}\varepsilon\cdot U$
is a closed subspace in $X$. Consider the quotient space $X/\Ker U$ and endow it with the topology of normed space with the unit ball
$U+\Ker U$. Then $(X/\Ker U)^\blacktriangledown$ will be a Banach space, and we denote it by $A/U$. The natural map (the composition of the quotient map $X\to X/\Ker U$ and the completion $X/\Ker U\to (X/\Ker
U)^\blacktriangledown$) will be denoted by $\pi_U:X\to X/U$. Since $\sigma:X\to X'$ is an extension with respect to ${\tt Ban}$, the map $\pi_U:X\to X/U$ is extended to some linear continuous map $(\pi_U)':X'\to X/U$.
$$
\begin{diagram}
\node{X} \arrow[2]{e,t}{\sigma} \arrow{se,b}{\pi_U} \node[2]{X'}
\arrow{sw,b,--}{(\pi_U)'}
\\
\node[2]{X/U}
\end{diagram}
$$
If we denote by $W$ the unit ball in $X/U$, i.e the closure of the set $U+\Ker U$ in the space $(X/\Ker U)^\blacktriangledown=X/U$, then for the neighborhood of zero $V=((\pi_U)')^{-1}(W)$ we obtain the following chain, which proves \eqref{sigma(U)-supseteq-V-cap-sigma(X)}:
\begin{multline*}
y\in V\cap\sigma(X)\quad\Rightarrow\quad  \exists x\in X:\ y=\sigma(x)\ \& \ y\in V \quad \Rightarrow \\ \Rightarrow\quad \exists x\in X:\ y=\sigma(x)\ \&\ (\pi_U)'(y)=(\pi_U)'(\sigma(x))=\underbrace{\pi_U(x)\in W}_{\scriptsize\begin{matrix}\Updownarrow \\ x\in U\end{matrix}}
\quad\Rightarrow\quad \exists x\in U:\ y=\sigma(x)\quad\Rightarrow\quad y\in\sigma(U).
 \end{multline*}
Thus, $\sigma:X\to X'$ is an open and injective linear continuous map, and $\sigma(X)$ is dense in $X'$. This means that $X'$ can be perceived as a subspace in the completion $X^\blacktriangledown$ of the space $X$ with the topology induced from $X^\blacktriangledown$. I.e. there is a unique linear continuous map $\upsilon:X'\to X^\blacktriangledown$ such that
$$
\begin{diagram}
\node[2]{X} \arrow{sw,t}{\sigma} \arrow{se,t}{\blacktriangledown_X}\\
\node{X'}  \arrow[2]{e,b,--}{\upsilon} \node[2]{X^\blacktriangledown}
\end{diagram}
$$
We conclude that $\blacktriangledown_X:X\to X^\blacktriangledown$ is an envelope of $X$ in ${\tt LCS}$ with respect to ${\tt Ban}$. \epr

\subsection{Refinement}

\paragraph{Refinement in a class of morphisms by means of a class of morphisms.}

Suppose we have:
 \bit{

\item[---] a category ${\tt K}$, called {\it enveloping category},

\item[---] a category ${\tt T}$, called {\it repelling category},

\item[---] a covariant functor $F:{\tt T}\to{\tt K}$,

\item[---] two classes $\varGamma$ and $\varPhi$ of morphisms in ${\tt K}$, which have domains in objects of the class $F({\tt T})$, and $\varGamma$ is called a {\it class of realizing morphisms}, and $\varPhi$ a {\it class of test morphisms}.
     }\eit\noindent
Then

 \bit{
\item[$\bullet$] For given objects $X\in\Ob(\tt K)$ and $X'\in\Ob(\tt T)$ a morphism $\sigma:F(X')\to X$ is called an {\it enrichment of the object $X\in{\tt K}$ in the class of morphisms $\varGamma$ over the category ${\tt T}$ by means of the class of morphisms $\varPhi$}, if $\sigma\in\varGamma$, and for any object $B$ in $\tt T$ and for any morphism $\ph:F(B)\to X$, $\ph\in\varPhi$, there is a unique morphism $\ph':B\to X'$ in $\tt T$ such that the following diagram is commutative:
 \beq\label{DEF:suzhenie-T}
\begin{diagram}
\node[2]{X}  \\
\node{F(B)}\arrow{ne,t}{\varPhi\owns\ph}\arrow[2]{e,b,--}{F(\ph')}
\node[2]{F(X')}\arrow{nw,t}{\sigma\in\varGamma}
\end{diagram}
\eeq

\item[$\bullet$]\label{DEF:nachinka} An enrichment $\rho:F(E)\to X$ of the object $X\in{\tt K}$ in the class of morphisms $\varGamma$ over the category ${\tt T}$ by means of the class of morphisms $\varPhi$ is called a {\it refinement of the object $X\in{\tt K}$ in the class of morphisms $\varGamma$ over the category ${\tt T}$ by means of the class of morphisms $\varPhi$}, if for any other enrichment $\sigma:F(X')\to X$ (of the object $X\in{\tt K}$ in the class of morphisms $\varGamma$ over the category ${\tt T}$ by means of the class of morphisms $\varPhi$) there is a unique morphism $\upsilon:E\to X'$ in $\tt T$  such that the following diagram is commutative:
 \beq\label{DIAGR:otpechatok-T}
\begin{diagram}
\node[2]{X}  \\
\node{F(E)}\arrow{ne,t}{\rho}\arrow[2]{e,b,--}{F(\upsilon)}
\node[2]{F(X')}\arrow{nw,t}{\sigma}
\end{diagram}
 \eeq
}\eit

In what follows we are almost exclusively interested in the case when ${\tt T}={\tt K}$, and $F:{\tt K}\to{\tt K}$ is the identity functor. Like in the case of envelopes, we formulate the definitions for this situation separately.
 \bit{
\item[$\bullet$] A morphism $\sigma:X'\to X$ in the category  $\tt K$ is called an {\it enrichment of the object $X\in{\tt K}$ in the class of morphisms $\varGamma$ by means of the class of morphisms $\varPhi$}, if $\sigma\in\varGamma$, and for any morphism $\ph:B\to X$, $\ph\in\varPhi$, there exists a unique morphism $\ph':B\to X'$ in the category  $\tt K$, such that the following diagram is commutative:
 \beq\label{DEF:suzhenie}
\begin{diagram}
\node[2]{X}  \\
\node{B}\arrow{ne,t}{\forall\ph\in\varPhi}\arrow[2]{e,b,--}{\exists!\ph'}
\node[2]{X'}\arrow{nw,t}{\sigma\in\varGamma}
\end{diagram}
 \eeq

\item[$\bullet$] An enrichment $\rho:E\to X$ of the object $X\in{\sf Ob}(\tt K)$ in the class of morphisms $\varGamma$ by means of the class of morphisms $\varPhi$ is called a {\it refinement of $X$ in the class $\varGamma$ by means of $\varPhi$}, if for any other enrichment $\sigma:X'\to X$ (of $X$ in $\varGamma$ by means of  $\varPhi$) there exists a unique morphism $\upsilon:E\to X'$ in $\tt K$, such that the following diagram is commutative:
 \beq\label{DIAGR:otpechatok}
\begin{diagram}
\node[2]{X}  \\
\node{E}\arrow{ne,t}{\rho}\arrow[2]{e,b,--}{\exists!\upsilon}
\node[2]{X'}\arrow{nw,t}{\forall\sigma}
\end{diagram}
 \eeq
For the morphism of refinement $\rho:E\to X$ we use the notation
    \beq\label{DEF:tr_F^L(A)}
    \rho=\rf_{\varPhi}^\varGamma  X.
    \eeq
The very object $E$ is also called a {\it refinement} of $X$ in $\varGamma$ by means of $\varPhi$, and is denoted by
    \beq\label{DEF:Tr_F^L(A)}
E=\Rf_{\varPhi}^\varGamma  X.
 \eeq
}\eit

\brem Like in the case of envelope, the refinement $\Rf_{\varPhi}^\varGamma  X$ (if exists) is defined up to an isomorphism. The question when the correspondence $X\mapsto\Rf_{\varPhi}^\varGamma X$ can be defined as a functor is discussed below starting from page \pageref{DIAGR:funktorialnost-imp-i-I}. \erem

\brem If $\varGamma=\varnothing$, then, of course, neither enrichments, no refinements in the class $\varGamma$ exist in the category ${\tt K}$. So this construction is interesting only of $\varGamma$ is a non-empty class. The following two situations will be of special interest:
  \bit{

\item[---] $\varGamma=\Mono({\tt K})$ (i.e. $\varGamma$ coincides with the class of all monomorphisms of the category ${\tt K}$), then we will use the following notations:
    \begin{align}\label{imp_(varPhi)^Mono}
    &     \rf_{\varPhi}^{\Mono}X:=\rf_{\varPhi}^{\Mono({\tt K})}X, &&     \Rf_{\varPhi}^{\Mono}X:=\Rf_{\varPhi}^{\Mono({\tt K})}X.
    \end{align}

\item[---] $\varGamma=\Mor({\tt K})$ (i.e. $\varGamma$ coincides with the class of all morphisms of the category ${\tt K}$), in this case it is convenient to omit any mentioning about $\varGamma$ in the formulations and notations, so we will be speaking about {\it refinements of the object $X\in{\sf Ob}(\tt K)$ in the category ${\tt K}$ by means of the class of morphisms $\varPhi$}, and the notations will be simplified as follows:
    \begin{align}\label{imp_(varPhi)=imp_(varPhi)^K}
    &     \rf_{\varPhi}X:=\rf_{\varPhi}^{\Mor({\tt K})}X, &&     \Rf_{\varPhi}X:=\Rf_{\varPhi}^{\Mor({\tt K})}X.
    \end{align}
 }\eit
\erem

\brem Another degenerate, but this time an informative case is when $\varPhi=\varnothing$. What is essential for a given object $X$, $\varPhi$ does not contain morphisms coming to $X$:
$$
\varPhi_X=\{\ph\in\varPhi:\ \Ran\ph=X\}=\varnothing.
$$
Then, obviously, any morphism $\sigma\in\varGamma $ coming to $X$, $\sigma:X\gets X'$, is an enrichment of $X$ (in the class of morphisms $\varGamma$ by means of the class of morphisms $\varnothing$). If in addition $\varGamma=\Mono$, then the refinement will be the initial object of the category $\Mono_X$ (if it exists). This can be depicted by the formula
$$
\Rf_\varnothing^\varGamma X=\min\Mono_X.
$$
On the other hand, if ${\tt K}$ is a category with zero $0$, and $\varGamma$ contains all the morphisms going from $0$, then the refinement of each object in $\varGamma$ by means of the empty class of morphisms is $0$:
$$
\Rf_\varnothing^\varGamma X=0.
$$
\erem

\brem Another extreme situation is when $\varPhi=\Mor({\tt K})$. For a given object $X$ the essential thing here is that $\varPhi$ contains the local identity for $X$:
$$
1_X\in\varPhi.
$$
Then for any enrichment $\sigma$ the diagram
$$
\xymatrix @R=2.pc @C=2.0pc % @M=14pt
{
X   & & X'\ar[ll]_{\sigma}\\
 & B\ar[ul]^{1_X} \ar@{-->}[ur] &
}
$$
implies that $\sigma$ must be a coretraction (moreover, the dotted arrow here must be unique). In the special case if $\varGamma\subseteq\Mono$ this is possible only if $\sigma$ is an isomorphism. As a corollary, the refinement of $X$ in $\varGamma$ coincides here with $X$ (up to isomorphism):
$$
\varGamma\subseteq\Mono\quad\Longrightarrow\quad
\Rf_{\Mor({\tt K})}^{\varGamma} X=X.
$$
\erem

\medskip
\centerline{\bf Properties of refinements:}

\bit{\it

\item[$1^\circ$.]\label{LM:Imp,suzhenie-verh-klassa-morfizmov} Suppose $\varSigma\subseteq\varGamma$, then for any object $X$ and for any class of morphisms $\varPhi$

\bit{

\item[(a)] every enrichment $\sigma:X\gets X'$ in $\varSigma$ by means of $\varPhi$ is an enrichment in $\varGamma$ by means of  $\varPhi$,

\item[(b)] if there are refinements $\rf_\varPhi^\varSigma X$ and $\rf_\varPhi^\varGamma X$, then there is a unique morphism $\rho:\Rf_\varPhi^\varSigma X\gets\Rf_\varPhi^\varGamma X$ such that the following diagram is commutative: \beq\label{Imp:suzhenie-verh-klassa-morfizmov}
\begin{diagram}
\node[2]{X} \\
\node{\Rf_\varPhi^\varSigma X} \arrow{ne,t}{\rf_\varPhi^\varSigma X}
\node[2]{\Rf^\varGamma _\varPhi X}\arrow{nw,t}{\rf^\varGamma _\varPhi X}
\arrow[2]{w,b,--}{\rho}
\end{diagram}
\eeq

\item[(c)] if there is a refinement $\rf_\varPhi^\varGamma X$ (in the wider class), and it lies in the (narrower) class $\varSigma$,
$$
\rf_\varPhi^\varGamma X\in\varSigma,
$$
then it is a refinement $\rf_\varPhi^\varSigma X$ (in the narrower class):
$$
\rf_\varPhi^\varGamma X=\rf_\varPhi^\varSigma X.
$$

}\eit

\item[$2^\circ$.]\label{LM,Imp:suzhenie-verh-klassa-morfizmov-2}
Let $\varSigma$, $\varGamma$, $\varPhi$ be classes of morphisms, and for an object $X$
 \bit{

\item[(a)] every enrichment $\sigma:X\gets X'$ in $\varGamma$ by means of $\varPhi$
belongs to $\varSigma$.
 }\eit
Then
 \bit{

\item[(b)] the refinement of $X$ in $\varGamma$ by means of $\varPhi$ exists if and only if there exists the refinement of $X$ in $\varGamma\cap\varSigma$ by means of $\varPhi$; in this case the refinements coincide:
$$
\rf_{\varPhi}^{\varGamma}=\rf_{\varPhi}^{\varGamma\cap\varSigma},
$$

\item[(c)] if $\varSigma\subseteq\varGamma$, then the existence of the refinement of $X$ in the (narrower) class $\varSigma$ by means of $\varPhi$ automatically implies the existence of the refinement of $X$ in the (wider) class $\varGamma$ by means of $\varPhi$ and their coincidence:
$$
\rf_{\varPhi}^{\varGamma}X=\rf_{\varPhi}^{\varSigma}X.
$$
 }\eit

\item[$3^\circ$.]\label{LM,Imp:suzhenie-klassa-morfizmov} Suppose $\varPsi\subseteq\varPhi$, then for any object $X$ and for any class of morphisms  $\varGamma$

\bit{\it

\item[(a)] every enrichment $\sigma:X\gets X'$ of $X$ in $\varGamma$ by means of $\varPhi$ is an enrichment of $X$ in $\varGamma$ by means of $\varPsi$,

\item[(b)] if there are refinements $\Rf_\varPsi^\varGamma  X$ and $\Rf_\varPhi^\varGamma  X$, then there is a unique morphism $\alpha:\Rf_\varPsi^\varGamma  X\to\Rf_\varPhi^\varGamma  X$ such that the following diagram is commutative: \beq\label{Ipm:suzhenie-klassa-morfizmov}
\begin{diagram}
\node[2]{X}  \\
\node{\Rf_\varPsi^\varGamma  X}\arrow{ne,t}{\rf_\varPsi^\varGamma
X}\arrow[2]{e,b,--}{\alpha}   \node[2]{\Rf^\varGamma_\varPhi
X}\arrow{nw,t}{\rf^\varGamma_\varPhi X}
\end{diagram}
\eeq
}\eit

\item[$4^\circ$.]\label{TH:imp_Psi=imp_Phi}  Suppose $\varPhi\subseteq\varPsi\circ\Mor({\tt K})$ (i.e. each morphism $\ph\in\varPhi$ can be represented in the form $\ph=\psi\circ\chi$, where $\psi\in\varPsi$), then for any object $X$ and for any class of morphisms $\varGamma$

 \bit{\it

\item[(a)]  if an enrichment $\sigma:X\gets X'$ in $\varGamma$ by means of $\varPsi$ is at the same time a monomorphism in $\tt K$, then it is an enrichment in $\varGamma$ by means of $\varPhi$,

\item[(b)]  if there are refinements $\rf_\varPsi^\varGamma X$ and $\rf_\varPhi^\varGamma X$, and $\rf_\varPsi^\varGamma X$ is at the same time a monomorphism in $\tt K$, then there is a unique morphism $\beta:\Rf_\varPsi^\varGamma X\gets\Rf_\varPhi^\varGamma X$ such that the following diagram is commutative:
 \beq\label{imp_Psi=imp_Phi-0}
\begin{diagram}
\node[2]{X}  \\
\node{\Rf_\varPsi^\varGamma  X}\arrow{ne,t}{\rf_\varPsi^\varGamma  X}   \node[2]{\Rf^\varGamma  _\varPhi X}\arrow{nw,t}{\rf^\varGamma  _\varPhi X}\arrow[2]{w,b,--}{\beta}
\end{diagram}
\eeq
}\eit

\item[$5^\circ$.]\label{PROP:deistvie-monomorfizma-na-Imp} Let the classes of morphisms $\varGamma$, $\varPhi$ and a monomorphism $\mu:X\gets Y$ in ${\tt K}$ satisfy the following conditions:
 \bit{
\item[(a)]  there is a refinement $\Rf_{\mu\circ\varPhi}^\varGamma X$ by means of the class of morphisms $\mu\circ\varPhi=\{\mu\circ\ph;\ \ph\in\varPhi\}$,

\item[(b)]  there is a refinement $\Rf_\varPhi^\varGamma Y$,

\item[(c)]  the composition $\mu\circ\rf_{\varPhi}^\varGamma Y$ belongs to $\varGamma$.

    }\eit
Then there is a unique morphism $\upsilon:\Rf_{\mu\circ\varPhi}^\varGamma X\to \Rf_{\varPhi}^\varGamma Y$ such that the following diagram is commutative:
 \beq\label{deistvie-monomorfizma-na-Imp} \xymatrix @R=2.5pc @C=4.0pc {
 X & Y\ar[l]_{\mu} \\
 \Rf_{\mu\circ\varPhi}^\varGamma X\ar[u]^{\rf_{\mu\circ\varPhi}^\varGamma X}\ar@{-->}[r]_{\upsilon} &   \Rf_{\varPhi}^\varGamma Y\ar[u]_{\rf_{\varPhi}^\varGamma Y}\ar@/^1ex/[ul]_{\ \mu\circ\rf_{\varPhi}^\varGamma Y}
}
\eeq

}\eit

 \bit{
\item[$\bullet$]\label{DEF:morfizmy-porozhdayutsya-snaruzhi} Let us say that in a category ${\tt K}$ a {\it class of morphisms $\varPhi$ is generated on the outside by a class of morphisms $\varPsi$}, if
    $$
    \varPsi\subseteq\varPhi\subseteq\varPsi\circ\Mor({\tt K}).
    $$
    }\eit
The following fact is dual to Theorem \ref{TH:morfizmy-porozhdayutsya-iznutri} and is proved by analogy:

\btm\label{TH:morfizmy-porozhdayutsya-snaruzhi} Suppose in a category ${\tt K}$ a class of morphisms $\varPhi$ is generated on the outside by a class of morphisms $\varPsi$. Then for any class of monomorphisms $\varGamma$ (it is not necessary that $\varGamma$ contains all monomorphisms of the category ${\tt K}$) and for any object $X$ the existence of refinement $\rf^\varGamma_\varPsi X$ is equivalent to the existence of the refinement  $\rf^\varGamma_\varPhi X$, and these refinements coincide:
 \beq\label{imp_Psi=imp_Phi}
\rf^\varGamma_\varPsi X=\rf^\varGamma_\varPhi X.
 \eeq
\etm

\bit{
\item[$\bullet$]\label{DEF:varPhi-razlich-morfizmy-iznutri}  Let us say that a class of morphisms $\varPhi$ in a category ${\tt K}$  {\it differs morphisms on the inside}, if for any two different parallel morphisms $\alpha\ne\beta:X\to Y$ there is a morphism $\ph:M\to X$ from the class $\varPhi$ such that  $\alpha\circ\ph\ne\beta\circ\ph$.
}\eit

The following result is dual to Theorem \ref{TH:Phi-razdel-moprfizmy}:

\btm\label{TH:Phi-razdel-moprfizmy-iznutri} If the class of morphisms $\varPhi$ differs morphisms on the inside, then for any class of morphisms $\varGamma$
 \bit{
\item[(i)] every enrichment in $\varGamma$ by means of $\varPhi$ is an epimorphism,

\item[(ii)] the refinement in $\varGamma$ by means of $\varPhi$ exists if and only if there exists a refinement in $\varGamma\cap\Mono$ by means of $\varPhi$; in that case these refinements coincide:
$$
\rf_{\varPhi}^{\varGamma}=\rf_{\varPhi}^{\varGamma\cap\Epi},
$$

\item[(iii)] if the class $\varGamma$ contains all epimorphisms,
 $$
 \varGamma\supseteq\Epi,
 $$
 then the existence of a refinement in $\Epi$ by means of $\varPhi$ automatically implies the existence of a refinement in $\varGamma$ by means of  $\varPhi$, and the coincidence of these refinements:
$$
\rf_{\varPhi}^{\varGamma}=\rf_{\varPhi}^{\Epi}.
$$
 }\eit
\etm

\bit{
\item[$\bullet$] Let us remind that a class of morphisms $\varPhi$ in a category ${\tt K}$ is called a {\it left ideal}, if
$$
\Mor({\tt K})\circ\varPhi\subseteq\varPhi
$$
(i.e. for any $\ph\in\varPhi$ and for any morphism $\mu$ in ${\tt K}$ the composition $\mu\circ\ph$ belongs to $\varPhi$).
}\eit

The following is dual to Theorem \ref{TH:Phi-razdel-moprfizmy-*}

\btm\label{TH:Phi-razdel-moprfizmy-iznutri-*} If a class of morphisms $\varPhi$ differs morphisms on the inside and is a left ideal in the category  ${\tt K}$, then for any class of morphisms $\varGamma$
 \bit{
\item[(i)] every enrichment in $\varGamma$ by means of $\varPhi$ is a bimorphism,

\item[(ii)] a refinement in $\varGamma$ by means of $\varPhi$ exists if and only if there exists a refinement in  $\varGamma\cap\Bim$ by means of $\varPhi$; in that case these refinements coincide:
$$
\rf_{\varPhi}^{\varGamma}=\rf_{\varPhi}^{\varGamma\cap\Bim}.
$$

\item[(iii)] if $\varGamma$ contains all bimorphisms,
 $$
 \varGamma\supseteq\Bim,
 $$
 then a refinement in $\varGamma$ by means of $\varPhi$ exists if and only if there exists a refinement in $\Bim$ by means of $\varPhi$, and these refinements coincide:
$$
\rf_{\varPhi}^{\varGamma}=\rf_{\varPhi}^{\Bim}.
$$
 }\eit
\etm

\paragraph{Refinement in a class of objects by means of a class of objects.}

A special case is the situation when $\varGamma$ and/or $\varPhi$ are classes of all morphisms from some given subclass of objects in $\Ob({\tt K})$. An exact formulation for the case when both classes $\varGamma$ and $\varPhi$ are defined in this way is the following: suppose we have a category ${\tt K}$ and two subclasses ${\tt L}$ and ${\tt M}$ in the class $\Ob({\tt K})$ of objects of ${\tt K}$.

 \bit{
\item[$\bullet$] A morphism $\sigma:X'\to X$ is called an {\it enrichment of the object $X\in{\sf Ob}(\tt K)$ in the class of objects ${\tt L}$ by means of the class of objects ${\tt M}$}, if for any object $B\in{\tt M}$ and for any morphism $\ph:B\to X$ there is a unique morphism $\ph':B\to X'$ such that the following diagram is commutative:
$$
\put(5,-35){$\begin{matrix}\text{\rotatebox{90}{$\owns$}}\\ {\tt M}
\end{matrix}$} \put(71,-35){$\begin{matrix}\text{\rotatebox{90}{$\owns$}}\\
{\tt L} \end{matrix}$}
\begin{diagram}
\node[2]{X}  \\
\node{B}\arrow{ne,t}{\forall\ph}\arrow[2]{e,b,--}{\exists!\ph'}
\node[2]{X'}\arrow{nw,t}{\sigma}
\end{diagram}
$$

\item[$\bullet$]\label{DEF:nachinka-otn-klassa} An enrichment $\rho:E\to X$ of the object $X\in{\sf Ob}(\tt K)$ in the class of objects ${\tt L}$ by means of the class of objects ${\tt M}$ is called a {\it refinement of the object $X\in{\sf Ob}(\tt K)$ in the class of objects ${\tt L}$ by means of the class of objects ${\tt M}$}, and we write in this case
    \beq\label{DEF:tr_M^L(A)}
    \rho=\rf_{\tt M}^{\tt L} X,
    \eeq
    if for any other enrichment $\sigma:X'\to X$ (of the object $X\in{\sf Ob}(\tt K)$ in the class of objects ${\tt L}$ by means of the class of objects ${\tt M}$) there is a unique morphism $\upsilon:E\to X'$ such that the following diagram is commutative:
 \beq\label{diagr:sled}
\put(5,-35){$\begin{matrix}\text{\rotatebox{90}{$\owns$}}\\ {\tt L}
\end{matrix}$} \put(71,-35){$\begin{matrix}\text{\rotatebox{90}{$\owns$}}\\
{\tt L} \end{matrix}$}
\begin{diagram}
\node[2]{X}  \\
\node{E}\arrow{ne,t}{\rho}\arrow[2]{e,b,--}{\exists!\upsilon}
\node[2]{X'}\arrow{nw,t}{\forall\sigma}
\end{diagram}
 \eeq
The very object $E$ is also called a {\it refinement} of the object $X\in{\sf Ob}(\tt K)$ in the class of objects ${\tt L}$ by means of the class of objects ${\tt M}$, and we use the following notation for it:
    \beq\label{DEF:Tr_M^L(A)}
E=\Rf_{\tt M}^{\tt L} X.
 \eeq
}\eit

The following two extreme situations in the choice of ${\tt L}$ can occur:
 \bit{

\item[---] if ${\tt L} =\Ob({\tt K})$, then we speak about a {\it refinement of the object $X\in{\sf Ob}(\tt K)$ in the category ${\tt K}$ by means of the class of objects ${\tt M}$}, and the notations will be the following:
    \begin{align}\label{imp_M=imp_M^K}
    &     \rf_{\tt M}X:=\rf_{\tt M}^{\tt K}X, &&     \Rf_{\tt M}X:=\Rf_{\tt M}^{\tt K}X,
    \end{align}

\item[---] if ${\tt L} ={\tt M}$, then the notions of the enrichment and the refinement coincide: {\it every enrichment of an object $X\in{\tt K}$ in a class ${\tt L}$ by means of the very same class ${\tt L}$ is a refinement of $X$ in ${\tt L}$  by means of ${\tt L}$} (since if  $\rho:E\to X$ and $\sigma:X'\to X$ are two enrichments of $X$ in ${\tt L}$  by means of ${\tt L}$, then in diagram \eqref{diagr:sled} the morphism $\upsilon$ exists and is unique just because $\sigma$ is an enrichment); for simplicity, in this case we will be specking about {\it refinement of the object $X$ in the class ${\tt L}$}, and the notations will be simplified as follows:
    \begin{align}
    & \rf_{\tt L} ^{\tt L} X=:\rf^{\tt L} X, &&  \Rf_{\tt L} ^{\tt L} X=:\Rf^{\tt L} X
    \end{align}
}\eit

\bit{
\item[$\bullet$] Let us say that a class of objects ${\tt M}$ in the category ${\tt K}$  {\it differs morphisms on the inside}, if the class of all morphisms going from objects of ${\tt M}$ has this property (in the sense of definition on page \pageref{DEF:varPhi-razlich-morfizmy-iznutri}), i.e. for any two different parallel morphisms $\alpha\ne\beta:X\to Y$ there is a morphism $\ph:M\to X$ such that $\alpha\circ\ph\ne\beta\circ\ph$.
}\eit

Theorem \ref{TH:Phi-razdel-moprfizmy-iznutri-*} implies

\btm\label{TH:M-razdel-moprfizmy-iznutri} If a class of objects ${\tt M}$ differs morphisms on the inside, then for any class of objects ${\tt L}$
 \bit{
\item[(i)]  each domain of convergence in the class ${\tt L}$ by means of the class ${\tt M}$ is a bimorphism,

\item[(ii)] a refinement in the class ${\tt L}$ by means of the class ${\tt M}$ exists if and only if there exists a refinement in the class of bimorphisms going from ${\tt L}$ by means of the class ${\tt M}$; in that case these refinements coincide:
$$
\rf_{\tt M}^{\tt L}=\rf_{\tt M}^{\Bim({\tt L},{\tt K})}.
$$
 }\eit
\etm

\paragraph{Examples of refinements.}

\bex {\bf Simply connected covering} used in the theory of Lie groups is from the categorial point of view a refinement in the class of pointed simply connected coverings by means of empty class of morphisms in the category of connected locally connected and semilocally simply connected pointed topological spaces
(see definitions in \cite{Postnikov}). \eex

\bex {\bf Bornologification} (see definition in \cite{Kriegl-Michor}) $X_{\text{\rm born}}$ of a locally convex space $X$ is a refinement of $X$ in the category ${\tt LCS}$ of locally convex spaces by means of the subcategory ${\tt Norm}$ of normed spaces:
$$
X_{\text{\rm born}}=\Rf_{\tt Norm}^{\tt LCS}X
$$
\eex
\bpr This follows from the characterization of bornologification as the strongest locally convex topology on $X$, for which all the imbeddings  $X_B\to X$ are continuous, where $B$ runs over the system of bounded absolutely convex subsets in $X$, and $X_B$ is a normed space with the unit ball $B$ (see \cite[Chapter I, Lemma 4.2]{Kriegl-Michor}).
 \epr

\bex {\bf Saturation} $X^\blacktriangle$ of a pseudocomplete locally convex space $X$ is a refinement in the category ${\tt LCS}$ of locally convex spaces by means of the subcategory ${\tt Smi}$ of the Smith spaces (see definitions in \cite{Akbarov}):
$$
X^\blacktriangle=\Rf_{\tt Smi}^{\tt LCS}X
$$
 \eex

\subsection{Connection with factorizations and with nodal decomposition}\label{SUBSEC:obolochki<->uzl-razlozh}

\paragraph{Connection with projective and injective limits.}

The similarity between the notions of an envelope and a projective limit is formalized in the following

\blm\label{LM:obolochka-konusa} The projective limit $\rho=\leftlim \rho^i:X\to\leftlim X^i$ of any projective cone  $\{\rho^i:X\to X^i;\ i\in I\}$ from a given object $X$ into a covariant (or contravariant) system $\{X^i,\iota^j_i\}$ is an envelope of the object  $X$ in an arbitrary class $\varOmega$ containing $\rho$ with respect to the system of morphisms $\{\rho^i; i\in I\}$:
 \beq\label{obolochka-konusa-Omega}
\rho=\leftlim \rho^i\in\varOmega
\quad\Longrightarrow\quad
\Env_{\{\rho^i;\ i\in I\}}^\varOmega X=\leftlim X^i
 \eeq
In particular, this is always true for $\varOmega=\Mor({\tt K})$:
 \beq\label{obolochka-konusa}
\Env_{\{\rho^i;\ i\in I\}}^{\Mor({\tt K})}X=\leftlim X^i
 \eeq
\elm
 \bpr
1. First, the morphism $\rho$ is an extension of $X$ with respect to the system $\{\rho^i\}$, since the definition of projective limit guarantees that for any $\rho^j$ there exists a unique continuation $\pi^j$ on $\leftlim X^i$:
\beq\label{proof-obolochka=lim-0}
\begin{diagram}
\node{X}\arrow[2]{e,t}{\rho}\arrow{se,b}{\rho^j}\node[2]{\leftlim X^i}\arrow{sw,b,--}{\pi^j}
\\
\node[2]{X^j}
\end{diagram}
\eeq

2. Suppose then that $\sigma:X\to X'$ is another extension. Then for any morphism $\rho^j:X\to X^j$ there is a unique morphism $\upsilon^j:X'\to X^j$, such that
\beq\label{proof-obolochka=lim-1}
\begin{diagram}
\node{X}\arrow[2]{e,t}{\sigma}\arrow{se,b}{\rho^j}\node[2]{X'}\arrow{sw,b,--}{\upsilon^j}
\\
\node[2]{X^j}
\end{diagram}
\eeq
For each indices $i\le j$ in the following diagram
$$
  \xymatrix @R=2.5pc @C=2.5pc
 {
 & X\ar@/_4ex/[ldd]_{\rho^i}\ar[d]^{\sigma} \ar@/^4ex/[rdd]^{\rho^j} & \\
 & X'\ar@{-->}[ld]^{\upsilon^i}\ar@{-->}[rd]_{\upsilon^j} & \\
 X^i \ar@/_2ex/[rr]^{\iota_i^j}  & & X^j
 }
$$
the following elements will be commutative: two upper little triangles (each one has one dotted arrow) and the perimeter (without dotted arrows). This together with the uniqueness of $\upsilon^j$ in the upper right little triangle implies that the lower little triangle (with two dotted arrows) is commutative as well:
$$
\begin{cases}(\iota_i^j\circ\upsilon^i)\circ\sigma=\iota_i^j\circ(\upsilon^i\circ\sigma)=\iota_i^j\circ\rho^i=\rho^j
\\
\upsilon^j\circ\sigma=\rho^j
\end{cases}
\qquad\Longrightarrow\qquad \iota_i^j\circ\upsilon^i=\upsilon^j.
$$
The commutativity of the triangle with two dotted arrows means in its turn that $X'$ with the system of morphisms $\upsilon^i$ is a projective cone of the covariant system $\{X^i;\iota_i^j\}$. So there must exist a uniquely  defined morphism $\upsilon$ such that for any $j$ in the diagram
$$
  \xymatrix @R=2.5pc @C=2.5pc
 {
 & X\ar@/_4ex/[ldd]_{\rho}\ar[d]^{\sigma} \ar@/^4ex/[rdd]^{\rho^j} & \\
 & X'\ar@{-->}[ld]^{\upsilon}\ar[rd]_{\upsilon^j} & \\
 \leftlim X^i \ar@/_2ex/[rr]^{\pi^j}  & & X^j
 }
$$
the little lower triangle is commutative. On the other hand, the right little triangle here is also commutative, since this is diagram  \eqref{proof-obolochka=lim-1} turned around, and the perimeter is commutative, since this is diagram \eqref{proof-obolochka=lim-0} turned around. Together with the uniqueness of the morphism $\rho$ in the system of all those perimeters with different $j$ this implies that the left little triangle is also commutative:
$$
\left(
\forall j\qquad
\begin{cases}
\pi^j\circ\upsilon\circ\sigma=\upsilon^j\circ\sigma=\rho^j
\\
\pi^j\circ\rho=\rho^j
\end{cases}\right)\qquad\Longrightarrow\qquad \upsilon\circ\sigma=\rho
$$

We understood that there is a morphism $\upsilon$ such that diagram \eqref{DEF:diagr-obolochka} is commutative (with $E=\leftlim X^i$). It remains to verify that such a morphism is unique. Let $\upsilon'$ be another morphism with the same property: $\rho=\upsilon'\circ\sigma$. Consider the following diagram:
$$
  \xymatrix @R=2.5pc @C=2.5pc
 {
 & X\ar@/_4ex/[ldd]_{\rho}\ar[d]^{\sigma} \ar@/^4ex/[rdd]^{\rho^j} & \\
 & X'\ar@{-->}[ld]^{\upsilon'}\ar[rd]_{\upsilon^j} & \\
 \leftlim X^i \ar@/_2ex/[rr]^{\pi^j}  & & X^j
 }
$$
Here (apart from the upper left triangle) the upper right triangle will be commutative (since this is the turned around diagram  \eqref{proof-obolochka=lim-1}) and the perimeter as well (since this is diagram \eqref{proof-obolochka=lim-0} turned around). Together with the uniqueness of the arrow $\upsilon^j$ in the upper right triangle this implies that the lower little triangle is also commutative:
$$
\begin{cases}
\pi^j\circ\upsilon'\circ\sigma=\pi^j\circ\rho=\rho^j
\\
\upsilon^j\circ\sigma=\rho^j
\end{cases}
\qquad\Longrightarrow\qquad
\pi^j\circ\upsilon'=\upsilon^j.
$$
This is true for each index $j$, so the morphism $\upsilon'$ must coincide with the morphism $\upsilon$ which we constructed before: $\upsilon'=\upsilon$.
\epr

\blm\label{LM:obolochka-konusa-v-kat-s-uzl-razl} Let $\varOmega$ be a monomorphically complemented class in a category ${\tt K}$, $\{X^i,\iota^j_i\}$ a covariant (or contravariant) system, and $\{\rho^i:X\to X^i;\ i\in I\}$ a projective cone from a given object $X$ into $\{X^i,\iota^j_i\}$. If there exists the projective limit $\rho=\leftlim\rho^i:X\to\leftlim X^i$, then in its factorization
$$
\rho=\mu_{\rho}\circ\e_{\rho},\qquad \mu_{\rho}\in{^\downarrow\varOmega},\quad \e_{\rho}\in\varOmega
$$
the epimorphism $\e_{\rho}$ is an envelope of the object $X$ with respect to the system of morphisms $\{\rho^i; i\in I\}$ in the class $\varOmega$:
\begin{align}\label{obolochka-konusa-v-kat-s-uzl-razl}
&
\e_{\leftlim \rho^i}=
\e_{\rho}=
\env_{\{\rho^i;\ i\in I\}}^{\varOmega}X, &&
\Ran\e_{\leftlim \rho^i}=\Ran\e_{\rho}=
\Env_{\{\rho^i;\ i\in I\}}^{\varOmega}X
\end{align}
\elm
 \bpr
1. By definition of projective limit every morphism $\rho^j$ has an extension $\pi^j$ to $\leftlim X^i$.
The restriction of $\pi^j$ to $\Ran\e_{\rho}$, i.e. the composition $\tau^j=\pi^j\circ \mu_{\rho}$ is an extension of $\rho^j$ to $\Ran\e_{\rho}$ along $\e_{\rho}$:
\beq\label{proof-obolochka=lim-0-v-kat-s-uzl-razl}
\xymatrix @R=2.pc @C=5.0pc % @M=14pt
{
X\ar[r]_{\e_{\rho}}\ar@/_1ex/[dr]_{\rho^j}\ar@/^4ex/[rr]^{\rho}
 &  \Ran\e_{\rho}\ar[r]_{\mu_\rho}\ar@{-->}[d]^{\tau^j}   & \leftlim X^i \ar@/^1ex/@{-->}[dl]^{\pi^j}
\\
& X^j &
}
\eeq
Such an extension $\tau^j$ is unique since $\e_{\rho}\in\Epi$, and we can say that $\e_{\rho}$ is an extension of $X$ in $\varOmega$ with respect to the system $\{\rho^i\}$.

2. Further, let $\sigma:X\to X'$ be another extension of $X$ in $\varOmega$ with respect to $\{\rho^i\}$. Like in the proof of Lemma \ref{LM:obolochka-konusa}, we find a morphism $\upsilon$ such that $\upsilon\circ\sigma=\rho$. We have
$$
\upsilon\circ\sigma=\rho=\mu_\rho\circ\e_\rho
$$
and since $\sigma\in\varOmega$, $\mu_\rho\in{^\downarrow\varOmega}$, there exists a diagonal morphism
$\delta$ such that
$$
\delta\circ\sigma=\e_\rho.
$$
This morphism is unique since $\sigma\in\varOmega\subseteq\Epi$.
\epr

The dual results are as follows.

\blm\label{LM:nachinka-konusa}
The injective limit $\rho=\rightlim\rho^i:X\gets\rightlim X^i$ of any injective cone $\{\rho^i:X\gets X^i;\ i\in I\}$ into a given object $X$ from a covariant (or contravariant) system $\{X^i,\iota^j_i\}$ is a refinement of the object $X$ in an arbitrary class of objects $\varGamma$ containing $\rho$ by means of the system of morphisms $\{\rho^i; i\in I\}$:
 \beq\label{otpechatok-konusa-Omega}
\rho=\rightlim\rho^i\in\varGamma
\quad\Longrightarrow\quad
\Rf_{\{\rho^i;\ i\in I\}}^\varGamma X=\rightlim X^i
 \eeq
In particular, this is true for $\varGamma=\Mor({\tt K})$:
 \beq\label{otpechatok-konusa}
\Rf_{\{\rho^i;\ i\in I\}}^{\Mor({\tt K})}X=\rightlim X^i
 \eeq
 \elm

\blm\label{LM:otpechatok-konusa-v-kat-s-uzl-razl}
Let $\varGamma$ be an epimorphically complementable class in a category ${\tt K}$, $\{X^i,\iota^j_i\}$ a covariant  (or contravariant) system, and $\{\rho^i:X\gets X^i;\ i\in I\}$ an injective cone from $\{X^i,\iota^j_i\}$ into a given object $X$. If there esists the injective limit $\rho=\rightlim\rho^i:X\gets\rightlim X^i$, then in its factorization
$$
\rho=\mu_{\rho}\circ\e_{\rho},\qquad \mu_{\rho}\in\varGamma,\quad \e_{\rho}\in\varGamma^\downarrow
$$
the monomorphism $\mu_{\rho}$ is a refinement of the object $X$ in the class $\varGamma$ by means of the system of morphisms $\{\rho^i; i\in I\}$:
 \begin{align}\label{otpechatok-konusa-v-kat-s-uzl-razl}
& \rf_{\{\rho^i;\ i\in I\}}^{\varGamma}X=\mu_{\rho}=\mu_{\rightlim \rho^i},
&& \Rf_{\{\rho^i;\ i\in I\}}^{\varGamma}X=\Dom\mu_{\rho}=\Dom\mu_{\rightlim \rho^i}
 \end{align}
\elm

\paragraph{Existence of envelopes and refinements for complementable classes.}

\blm\label{LM:env_varPhi^varOmega_X=env_(e_ph;ph_in_varPhi)^varOmega_X}
Let $\varOmega$ be a monomorphically complementable class in a category ${\tt K}$. Then for each object $X$ and for any class of morphisms $\varPhi$
\beq\label{env_varPhi^varOmega_X=env_(e_ph;ph_in_varPhi)^varOmega_X}
\env_{\varPhi}^{\varOmega}X=\env_{\{\e_{\ph};\ \ph\in \varPhi\}}^{\varOmega}X
\eeq
(this means that if one of these envelopes exists then the other one exists as well and they coincide).
\elm
\bpr Let $\ph=\mu_{\ph}\circ\e_{\ph}$ be the factorization with $\mu_{\ph}\in{^\downarrow\varOmega}$ and
$\e_{\ph}\in\varOmega$. We need to verify that the extensions with respect to classes $\varPhi$ and $\{\e_{\ph};\ \ph\in \varPhi\}$ are the same. Let $\sigma:X\to X'$ be an extension of $X$ in $\varOmega$ with respect to morphisms  $\{\e_{\ph};\ \ph\in \varPhi\}$. Then in the diagram
$$
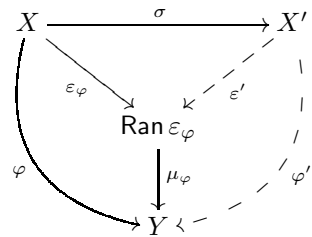
\xymatrix %@R=2.5pc @C=4.0pc
{
X\ar[rr]^{\sigma}\ar[rd]_{\e_{\ph}}\ar@/_6ex/[rdd]_{\ph} &  & X'\ar@{-->}[ld]^{\e'}\ar@{-->}@/^6ex/[ldd]^{\ph'} \\
 & \Ran\e_{\ph}\ar[d]^{\mu_{\ph}} & \\
 &  Y  &
}
$$
the existence of morphism $\e'$, for which the upper little triangle is commutative, implies the existence of morphism $\ph'$, for which the lower right little triangle is commutative, and since the last (left) little triangle is commutative, we conclude that the big triangle (the perimeter) is commutative as well. In addition, $\ph'$ is unique since $\sigma$ is an epimorphism. Hence, $\sigma:X\to X'$ is an extension of $X$ with respect to morphisms $\varPhi$.

Conversely, suppose that $\sigma:X\to X'$ is an extension of $X$ with respect to $\varPhi$. Then for any $\ph\in\varPhi$ there exists a morphism  $\ph'$ such that in the diagram
$$
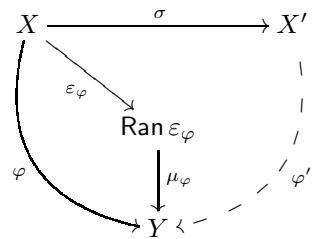
\xymatrix %@R=2.5pc @C=4.0pc
{
X\ar[rr]^{\sigma}\ar[rd]_{\e_{\ph}}\ar@/_6ex/[rdd]_{\ph} &  & X'\ar@{-->}@/^6ex/[ldd]^{\ph'} \\
 & \Ran\e_{\ph}\ar[d]^{\mu_{\ph}} & \\
 &  Y  &
}
$$
the big triangle (perimeter) is commutative. The lower left little triangle here is commutative as well due to \eqref{faktorizatsiya-v-kat-s-faktoriz}, hence the following quadrangle is also commutative:
$$
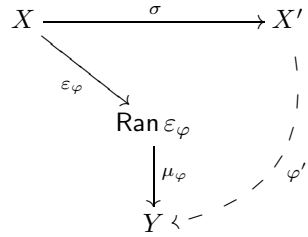
\xymatrix %@R=2.5pc @C=4.0pc
{
X\ar[rr]^{\sigma}\ar[rd]_{\e_{\ph}} &  & X'\ar@{-->}@/^6ex/[ldd]^{\ph'} \\
 & \Ran\e_{\ph}\ar[d]^{\mu_{\ph}} & \\
 &  Y  &
}
$$
Here $\sigma\in\varOmega$ and $\mu_{\ph}\in\varOmega^\downarrow$. Thus, there exists a diagonal $\e'$:
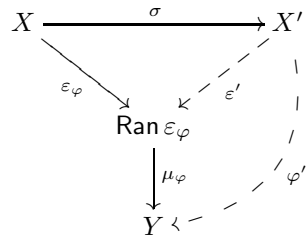
$$
\xymatrix %@R=2.5pc @C=4.0pc
{
X\ar[rr]^{\sigma}\ar[rd]_{\e_{\ph}} &  & X'\ar@{-->}@/^6ex/[ldd]^{\ph'}\ar@{-->}[ld]^{\e'} \\
 & \Ran\e_{\ph}\ar[d]^{\mu_{\ph}} & \\
 &  Y  &
}
$$
In particular, the upper triangle is commutative, and, since this is true for any $\ph\in\varPhi$, $\sigma:X\to X'$ is an extension of $X$ with respect to morphisms $\{\e_{\ph};\ \ph\in\varPhi\}$.
\epr

\medskip
\centerline{\bf Properties of envelopes in monomorphically complementable classes:}

\bit{\it

\item[] Let $\varOmega$ be a monomorphically complementable class in a category $\tt K$.

\item[$1^\circ$.]\label{1^0:obolochka-otn-1-morphizma}
For each morphism $\ph:X\to Y$ in ${\tt K}$ the epimorphism $\e_{\ph}$ in the factorization  $\ph=\mu_{\ph}\circ\e_{\ph}$ (defined by classes ${^\downarrow\varOmega}$ and $\varOmega$) is an envelope of $X$ in $\varOmega$ with respect to $\ph$:
 \begin{align}\label{obolochka-otn-morfizma}
 & \env_{\ph}^\varOmega X=\e_{\ph}, && \Env_{\ph}^\varOmega X=\Ran\e_{\ph}
 \end{align}

\item[$2^\circ$.]\label{2^0:obolochka-otn-konechnogo-mnozhestva-morphizmov}
If ${\tt K}$ is a category with finite products, then each object $X$ in ${\tt K}$ has an envelope in $\varOmega$ with respect to arbitrary finite set of morphisms $\varPhi$
going from $X$.

\item[$3^\circ$.]\label{3^0:obolochka-otn-mnozhestva-morphizmov} If ${\tt K}$ is a category with products\footnote{In propositions $3^\circ$-$5^\circ$ we assume that ${\tt K}$ has products over arbitrary index sets, not necessarily finite.}, then every its object $X$ has an envelope in $\varOmega$ with respect to an arbitrary set of morphisms $\varPhi$ going from $X$.

\item[$4^\circ$.]\label{4^0:obolochka-otn-klassa-morphizmov-porozhd-mnozhestvom}
If ${\tt K}$ is a category with products, then every its object $X$ has an envelope in $\varOmega$ with respect to an arbitrary class of morphisms $\varPhi$ going from $X$ and having a subset which generates $\varPhi$ on the inside (see page\pageref{DEF:morfizmy-porozhdayutsya-iznutri}).

\item[$5^\circ$.]\label{5^0:obolochka-otn-klassa-morphizmov}
If ${\tt K}$ has products, and is co-well-powered in the class $\varOmega$, then in ${\tt K}$ every object $X$ has envelope in the class $\varOmega$ with respect to arbitrary class of morphisms $\varPhi$, going from $X$.
}\eit

\bpr 1. The morphism $\e_{\ph}$ is an extension of $X$ in $\varOmega$ with respect to $\ph$, as it is seen from  diagram
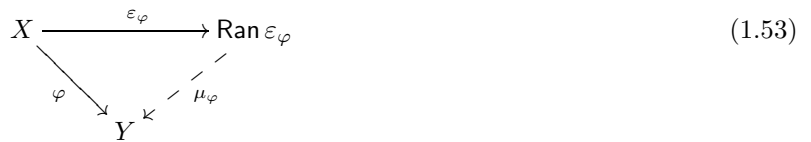
\beq\label{e_ph=obolochka}
\xymatrix %@R=2.5pc @C=4.0pc
{
X \ar[rr]^{\e_{\ph}}  \ar[dr]_{\ph} & & \Ran\e_{\ph} \ar@{-->}[dl]^{\mu_{\ph}}\\
& Y &
}
\eeq
Let $\sigma:X\to N$ be another extension of $X$ in $\varOmega$ with respect to $\ph$:
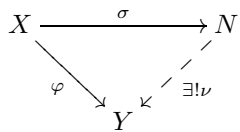
$$
\xymatrix %@R=2.5pc @C=4.0pc
{
X \ar[rr]^{\sigma} \ar[dr]_{\ph} & & N \ar@{-->}[dl]^{\exists!\nu}
\\
& Y &
}
$$
We have a commutative diagram
$$
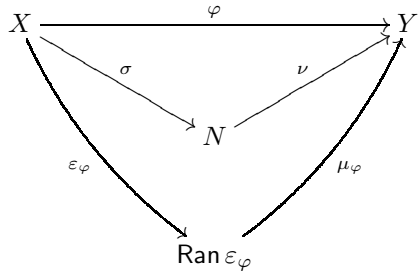
\xymatrix @R=2.5pc @C=4.0pc
{
X\ar[rr]^{\ph}\ar[rd]^{\sigma}\ar@/_2ex/[rdd]_{\e_{\ph}} & & Y \\
& N\ar[ru]^{\nu} & \\
& \Ran\e_{\ph}\ar@/_2ex/[ruu]_{\mu_{\ph}} &
}
$$
Here $\sigma\in\varOmega$ and $\mu_{\ph}\in{^\downarrow\varOmega}$, hence there exists a diagonal of the lower quadrangle:
$$
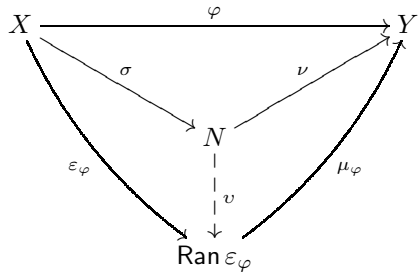
\xymatrix @R=2.5pc @C=4.0pc
{
X\ar[rr]^{\ph}\ar[rd]^{\sigma}\ar@/_2ex/[rdd]_{\e_{\ph}} & & Y \\
& N\ar[ru]^{\nu}\ar@{-->}[d]^{\upsilon} & \\
& \Ran\e_{\ph}\ar@/_2ex/[ruu]_{\mu_{\ph}} &
}
$$
The morphism $\upsilon$ is the very same morphism in Diagram \eqref{DEF:diagr-obolochka} which connects the extension $\sigma$ with the envelope $\e_{\ph}$. Its uniqueness follows from epimorphy of $\sigma$.

2. Let $X$ be an object and $\varPhi$ a finite set of morphisms. Certainly, it is sufficient to pick out in $\varPhi$ a subset $\varPhi^X=\{\ph:X\to Y_{\ph};\ \ph\in\varPhi^X\}$ of morphisms going from $X$,
$$
\ph\in\varPhi^X\quad\Longleftrightarrow\quad \ph\in\varPhi\quad\&\quad \Dom\ph=X.
$$
Then the envelope with respect to $\varPhi$ is the same as the envelope with respect to $\varPhi^X$. Consider the product of objects $\prod_{\ph\in\varPhi^X}Y_{\ph}$ and the product of morphisms $\prod_{\ph\in\varPhi^X}\ph:X\to \prod_{\ph\in\varPhi^X}Y_{\ph}$. The envelope of $X$ with respect to the $\varPhi^X$ is exactly the envelope of $X$ with respect to one morphism, $\prod_{\ph\in\varPhi^X}\ph$. After that we apply $1^\circ$.

3. Let ${\tt K}$ be a category with products over arbitrary (not necessarily finite) index set. Then the above reasoning work in the case when $\varPhi$ is a set (not necessarily finite) of morphisms.

4. Let $\varPsi\subseteq\varPhi$ be a subset (not a proper class), generating $\varPhi$ on the inside. By the just  proven property $3^\circ$, every object $X$ has an envelope with respect to $\varPsi$. And by \eqref{env_Psi=env_Phi} this envelope coincides with the envelope with respect to $\varPhi$.

5. Let ${\tt K}$ be a category with products (over arbitrary set of indices), $A$ an object in ${\tt K}$, and  $\varPhi$ a class of morphisms (not necessarily a set). The idea of the proof is to replace the {\it class} $\varPhi$ by a {\it set} of morphisms $M$, such that the envelope will be the same. Like in 2, we can think that $\varPhi$ consists of morphisms going from $X$:
$$
\forall\ph\in\varPhi\quad \Dom\ph=X.
$$
For any $\ph\in\varPhi$ we consider the morphism $\e_{\ph}$. By Lemma \ref{LM:env_varPhi^varOmega_X=env_(e_ph;ph_in_varPhi)^varOmega_X}, we can replace $\varPhi$ by the class   $\{\e_{\ph};\ \ph\in\varPhi\}$:
$$
\env_{\varPhi}^{\Epi}X=\env_{\{\e_{\ph};\ \ph\in\varPhi\}}^{\Epi}X.
$$
After that we need to recall that all $\e_{\ph}$ belong to $\varOmega$, and since our category is co-well-powered in the class $\varOmega$, we can choose among $\e_{\ph}$ a set $M$ such that every $\e_{\ph}$ will be isomorphisc to some $\e\in M$, i.e. $\e_{\ph}=\iota\circ\e$ for some isomorphism $\iota$. The set $M$ now replaces the class  $\{\e_{\ph};\ \ph\in\varPhi\}$ (and hence tha class $\varPhi$), and after that $3^\circ$ works.
\epr

The dual results for refinements look as follows.

\blm\label{LM:imp_varPhi^varGamma_X=imp_(mu_ph;ph_in_varPhi)^varGamma_X}
Let $\varGamma$ be an epimorphically complementable class in a category ${\tt K}$. then for every object $X$ and for every class of morphisms $\varPhi$
\beq\label{imp_varPhi^varGamma_X=imp_(mu_ph;ph_in_varPhi)^varGamma_X}
\rf_{\varPhi}^{\varGamma}X=\rf_{\{\mu_{\ph};\ \ph\in \varPhi\}}^{\varGamma}X
\eeq
(this means that if one of these refinements exists then the other one exists as well and they coincide).
\elm

\medskip
\centerline{\bf Properties of refinements in epimorphically complementable classes:}

\bit{\it

\item[] Let $\varGamma$ be an epimorphically complementable class of morphisms in a category $\tt K$.

\item[$1^\circ$.]\label{1^0:otpechatok-posr-1-morphizma}
For each morphism $\ph:X\gets Y$ in ${\tt K}$ the monomorphism $\mu_{\ph}$ in the factorization  $\ph=\mu_{\ph}\circ\e_{\ph}$ (defined by classes $\varGamma$ and $\varGamma^\downarrow$) is a refinement of $X$ in $\varGamma$ by means of the morphism $\ph$:
 \begin{align}\label{otpechatok-posr-morfizma}
 & \rf_{\ph}^\varGamma X=\mu_{\ph}, && \Rf_{\ph}^\varGamma X=\Dom\mu_{\ph}
 \end{align}

\item[$2^\circ$.]\label{2^0:otpechatok-posr-konechnogo-mnozhestva-morphizmov}
If ${\tt K}$ is a category with finite co-products, then every object $X$ in ${\tt K}$ has a refinement in $\varGamma$ by means of arbitrary finite set of morphisms $\varPhi$ going to  $X$.

\item[$3^\circ$.]\label{3^0:otpechatok-posr-mnozhestva-morphizmov}
If ${\tt K}$ is a category with co-products\footnote{\label{FOOT:lyuboe-I}In propositions $3^\circ$-$5^\circ$ we assume that ${\tt K}$ has co-products over arbitrary index set, not necessarily finite.}, then every object $X$ in ${\tt K}$ has a refinement in $\varGamma$ by means of an arbitrary set of morphisms $\varPhi$ going to $X$.

\item[$4^\circ$.]\label{4^0:otpechatok-posr-klassa-morphizmov-porozhd-mnozhestvom} If ${\tt K}$ is a category with co-products, then every object $X$ in ${\tt K}$ has a refinement in $\varGamma$ by means of an arbitrary set of morphisms $\varPhi$ going to $X$, and having a set that generates it on the inside.

\item[$5^\circ$.]\label{5^0:otpechatok-posr-klassa-morphizmov} If ${\tt K}$ has co-products and is well-powered in the class $\varGamma$, then every object $X$ in ${\tt K}$ has a refinement in $\varGamma$ by means of an arbitrary class of morphisms $\varPhi$ going to $X$.
}\eit

\paragraph{Existence of envelopes and refinements in categories with nodal decomposition.} \label{obolochki-i-nachinki-otn-1-morfizma}

The general properties on page \pageref{1^0:obolochka-otn-1-morphizma}, being applied to the cases $\varOmega=\Epi$ and $\varOmega=\SEpi$, give the following:

\medskip
\centerline{\bf Properties of envelopes in $\Epi$ and in $\SEpi$ in a category with nodal decomposition:}

\bit{\it

\item[] Let ${\tt K}$ be a category with nodal decomposition.

\item[$1^\circ$.]\label{1^0:obolochka-env_ph^Epi} For each morphism $\ph:X\to Y$ in ${\tt K}$
 \bit{
 \item[---] the epimorphism $\red_\infty\ph\circ\coim_\infty\ph$ in the nodal decomposition of $\ph$ is an envelope of the object $X$ in the class $\Epi$ of all epimorphisms with respect to $\ph$:
 \begin{align}\label{obolochka--env_ph^Epi}
 & \env_{\ph}^{\Epi} X=\red_\infty\ph\circ\coim_\infty\ph, && \Env_{\ph}^{\Epi}X=\Im_\infty\ph
 \end{align}

 \item[---] the epimorphism $\coim_\infty\ph$ in the nodal decomposition of $\ph$ is an envelope of the object $X$ in the class $\SEpi$ of strong epimorphisms with respect to $\ph$:
 \begin{align}\label{obolochka--env_ph^SEpi}
 & \env_{\ph}^{\SEpi} X=\coim_\infty\ph, && \Env_{\ph}^{\SEpi}X=\Coim_\infty\ph
 \end{align}

}\eit

\item[$2^\circ$.]\label{2^0:obolochka-env_Phi^Epi-otn-konechnogo-mnozhestva-morphizmov} If ${\tt K}$ has finite products, then every object $X$ in $\tt K$ has envelopes in the classes $\Epi$ and $\SEpi$ with respect to an arbitrary finite set of morphisms $\varPhi$ going from $X$.

\item[$3^\circ$.]\label{3^0:obolochka-env_Phi^Epi-otn-mnozhestva-morphizmov} If ${\tt K}$ is a category with products\footnote{Similarly to footnote \ref{FOOT:lyuboe-I}.}, then every object $X$ in ${\tt K}$ has envelopes in the classes $\Epi$ and $\SEpi$ with respect to arbitrary set of morphisms $\varPhi$ going from $X$.

\item[$4^\circ$.]\label{4^0:obolochka-env_Phi^Epi-otn-klassa-morphizmov-porozhd-mnozhestvom} If ${\tt K}$ is a category with products, then every object $X$ in $\tt K$ has envelopes in the classes $\Epi$ and $\SEpi$ with respect to an arbitrary class of morphisms $\varPhi$ going from $X$ and having a set that generates $\varPhi$ in the inside.

\item[$5^\circ$.]\label{5^0:obolochka-env_Phi^Epi-otn-klassa-morphizmov} If ${\tt K}$ is a category with products, co-well-powered in the class $\Epi$ (respectively, in $\SEpi$), then in ${\tt K}$ every object $X$ has envelope in the class $\Epi$ (respectively, in the class $\SEpi$) with respect ot an arbitrary class of morphisms $\varPhi$ going from $X$.
}\eit

\bprop\label{PROP:obolochka-v-Omega,Bim,Epi}
If ${\tt K}$ is a category with products, with nodal decomposition, and co-well-powered in the class $\Epi$, then in ${\tt K}$ every object $X$ has envelope in each class $\varOmega$ which contains bimorphisms,
    $$
    \varOmega\supseteq\Bim,
    $$
with respect to an arbitrary right ideal of morphisms $\varPhi$, which differs morphisms on the outside\footnote{See definition on p.\pageref{DEF:varPhi-razlich-morfizmy-snaruzhi}.} and goes from $X$, and the envelope of $X$ in $\varOmega$ with respect to $\varPhi$ coincides with the envelopes in $\Bim$ and in $\Epi$ with respect to $\varPhi$:
$$
\env_{\varPhi}^{\varOmega}X=\env_{\varPhi}^{\Bim}X=\env_{\varPhi}^{\Epi}X.
$$
\eprop
\bpr
By property $5^\circ$, there exists an envelope $\env_{\varPhi}^{\Epi}X$. By Theorem \ref{TH:Phi-razdel-moprfizmy}(i) this envelope is a monomorphism, and hemce, a bimorphism: $\env_{\varPhi}^{\Epi}X\in\Bim$. Then by Property  $1^\circ$(c) on p.\pageref{LM:suzhenie-verh-klassa-morfizmov} the envelope in $\Epi$ must be an envelope in $\Bim$:
$\env_{\varPhi}^{\Epi}X=\env_{\varPhi}^{\Bim}X$. Now by Theorem \ref{TH:Phi-razdel-moprfizmy-*} the envelope in  $\Bim$ must be an envelope in $\varOmega$: $\env_{\varPhi}^{\Bim}X=\env_{\varPhi}^{\varOmega}X$.
\epr

The dual results for refinements look as follows.

\medskip
\centerline{\bf Properties of refinements of $\Mono$ and $\SMono$ in a category with nodal decomposition:}

\bit{\it

\item[] Let ${\tt K}$ be a category with nodal decomposition.

\item[$1^\circ$.]\label{1^0:otpechatok-imp_ph^Mono} For each morphism $\ph:X\gets Y$ in ${\tt K}$
 \bit{
 \item[---] the monomorphism $\im_\infty\ph\circ\red_\infty\ph$ in the nodal decomposition of $\ph$ is a refinement in the class $\Mono$ of all monomorphisms in $X$ by means of $\ph$:
 \begin{align}\label{otpechatok-imp_ph^Mono}
 & \rf_{\ph}^{\Mono} X=\im_\infty\ph\circ\red_\infty\ph, && \Rf_{\ph}^{\Mono}X=\Coim_\infty\ph
 \end{align}

 \item[---] the monomorphism $\im_\infty\ph$ in the nodal decomposition of $\ph$ is a refinement in the class $\SMono$ of strong monomorphisms in $X$ by means of $\ph$:
 \begin{align}\label{otpechatok-imp_ph^SMono}
 & \rf_{\ph}^{\SMono} X=\im_\infty\ph, && \Rf_{\ph}^{\SMono}X=\Im_\infty\ph
 \end{align}

}\eit

\item[$2^\circ$.]\label{2^0:otpechatok-imp_Phi^Mono-posr-konechnogo-mnozhestva-morphizmov} If ${\tt K}$ is a category with finite co-products, then every object $X$ in $\tt K$ has refinements in the classes $\Mono$ and $\SMono$ by means of an arbitrary finite set of morphisms $\varPhi$ going to $X$.

\item[$3^\circ$.]\label{3^0:otpechatok-imp_Phi^Mono-posr-mnozhestva-morphizmov} If ${\tt K}$ is a category with co-products\footnote{In propositions $3^\circ$-$5^\circ$ we assume that ${\tt K}$ has co-products over arbitrary sets of indices, not necessarily finite.}, then every object $X$ in ${\tt K}$ has refinements in the classes $\Mono$ and $\SMono$ by means of an arbitrary set of morphisms $\varPhi$ going to $X$.

\item[$4^\circ$.]\label{4^0:otpechatok-imp_Phi^Mono-posr-klassa-morphizmov-porozhd-mnozhestvom} If ${\tt K}$ is a category with co-products, then every object $X$ in ${\tt K}$ has refinements in the classes $\Mono$ and $\SMono$ by means of an arbitrary class of morphisms $\varPhi$ coming to $X$, and having a set which generates $\varPhi$ on the outside.

\item[$5^\circ$.]\label{5^0:otpechatok-imp_Phi^Mono-posr-klassa-morphizmov} If ${\tt K}$ is a category with co-products, and well-powered in the class $\Mono$ (respectively, in $\SMono$), then in ${\tt K}$ each object $X$ has a refinement in the class $\Mono$ (respectively, of $\SMono$) by means of an arbitrary class of morphisms $\varPhi$ going to $X$.
}\eit

\bprop\label{PROP:otpechatok-v-Omega,Bim,Epi} If ${\tt K}$ is a category with co-products, with nodal decomposition and well-powered in the class $\Mono$, then in ${\tt K}$ every object $X$ has a refinement in an arbitrary class $\varGamma$, which contains bimorphisms,
    $$
    \varGamma\supseteq\Bim,
    $$
by means of arbitrary left ideal of morphisms $\varPhi$, which differs morphisms on the inside\footnote{See Definition on page \pageref{DEF:varPhi-razlich-morfizmy-iznutri}.} and goes to $X$, and the refinement of $X$ in $\varGamma$ by means of $\varPhi$ coincides with refinements in $\Bim$ and in $\Mono$ by means of $\varPhi$:
$$
\rf_{\varPhi}^{\varGamma}X=\rf_{\varPhi}^{\Bim}X=\rf_{\varPhi}^{\Mono}X.
$$
\eprop

\paragraph{Existence of nodal decomposition in categories with envelopes and refinements.}
By analogy with definitions on p.\pageref{DEF:strog-epi-razlich-mono} we will say that in a category ${\tt K}$
 \bit{
 \item[---]\label{DEF:epi-raspoznayut-mono}  {\it epimorphisms discern monomorphisms}, if from the fact that a morphism $\mu$ is not a monomorphism it follows that $\mu$ can be represented as a composition $\mu=\mu'\circ\e$, where $\e$ is an epimorphism, which is not an isomorphism,

 \item[---] {\it monomorphisms discern epimorphisms}, if from the fact that a morphism $\e$ is not an epimorphism it follows that $\e$ can be represented as a composition $\e=\mu\circ\e'$, where $\mu$ is a monomorphism, which is not an isomorphism.
 }\eit

\btm\label{TH:env+imp=>uzl=razl} Suppose that in a category ${\tt K}$
 \bit{
 \item[\rm (a)] epimorphisms discern monomorphisms, and, dually, monomorphisms discern epimorphisms,

 \item[\rm (b)] every immediate monomorphism is a strong monomorphism, and, dually, every immediate epimorphism is a strong epimorphism,

 \item[\rm (c)] every object $X$ has an envelope in the class $\Epi$ of all epimorphisms with respect to any morphism, starting from $X$, and, dually, in every object $X$ there is a refinement in the class $\Mono$ of all monomorphisms with respect to any morphism coming to $X$.
   }\eit
Then ${\tt K}$ is a category with nodal decomposition.
 \etm
\bpr Consider a morphism $\ph:X\to Y$.

1. Suppose $\e:X\to N$ is an envelope of $X$ in $\Epi$ with respect to $\ph$, and denote by $\beta$ the dashed arrow in \eqref{DEF:diagr-rasshirenie}:
$$
\ph=\beta\circ\e
$$
Note first that $\beta$ is a monomorphism. Indeed, if $\beta$ is not a monomorphism, then by (a), there exists a decomposition
$\beta=\beta'\circ\pi$, where $\pi$ is an epimorphism, but not an isomorphism. If we denote by $N'$ the range of $\pi$, then we get a diagram
\beq\label{PROOF:sushestv-uzlov-razlozhenija}
\xymatrix @R=4.0pc @C=4.0pc
{
X\ar[d]_{\e}\ar[r]^{\ph}\ar[dr]^(.7){\e'}|!{[d];[r]}\hole & Y
\\
N\ar[r]_{\pi}\ar[ur]^(.7){\beta} & N'\ar[u]_{\beta'}
}
\eeq
where by definition $\e'=\pi\circ\e$, and this will be is epimorphism, as a composition of two epimorphisms. Thus, $\e'$ is another extension of $X$ with respect to $\ph$. Hence, there exists a unique morphism $\upsilon$ such that the following diagram is commutative:
$$
\begin{diagram}\dgARROWLENGTH=3em
\node[2]{X}\arrow{sw,t}{\e} \arrow{se,t}{\e'} \\
\node{N}  \node[2]{N'}\arrow[2]{w,b,--}{\upsilon}
\end{diagram}
$$
Here we have:
$$
\pi\circ\e=\e'\quad\Longrightarrow\quad\upsilon\circ\pi\circ\e=\upsilon\circ\e'=\e=1_N\circ\e\quad\Longrightarrow\quad \upsilon\circ\pi=1_N
$$
and
$$
\upsilon\circ\e'=\e\quad\Longrightarrow\quad\pi\circ\upsilon\circ\e'=\pi\circ\e=\e'=1_{N'}\circ\e'\quad\Longrightarrow\quad \pi\circ\upsilon=1_{N'}.
$$
I.e. $\pi$ must be an isomorphism, and this contradicts our assumption that $\pi$ is not an isomorphism.

2. Similarly one can prove that $\beta$ is an immediate monomorphism. Indeed, any its factorization $\beta=\beta'\circ\pi$ leads again to diagram \eqref{PROOF:sushestv-uzlov-razlozhenija}, and the same reasoning gives that $\pi$ is an isomorphism.

3. The fact that $\beta$ is an immediate monomorphism together with condition (b) imply that $\beta$ is a strong monomorphism.

4. Denote by $\mu:M\to Y$ the refinement in the class $\Mono$ in $Y$ by means of the morphism $\ph$, and by $\alpha$ the dashed arrow in the corresponding diagram \eqref{DEF:suzhenie}, i.e.
$$
\ph=\mu\circ\alpha
$$
Using the dual reasoning to what we used when proving that $\beta$ is a strong monomorphism, one can show that $\alpha$ is a strong epimorphism.

5. Consider now a diagram
$$
\begin{diagram}\dgARROWLENGTH=3em
\node{X}\arrow[2]{r,t}{\ph}\arrow[2]{s,t}{\alpha}\arrow{se,t,-}{} \node[2]{Y} \\
 \node[2]{} \arrow{se,t}{\e} \\
\node{M}\arrow[2]{ne,t,3}{\mu}  \node[2]{N}\arrow[2]{n,r}{\beta}
\end{diagram}
$$
As we already understood, here $\alpha$ is an epimorphism, hence $\alpha$ is an extension of $X$ in $\Epi$ with respect to $\ph$. At the same time $\e$ is an envelope of $X$ in $\Epi$ with respect to $\ph$. Hence there exists a morphism $\upsilon$ such that the following diagram is commutative:
$$
\begin{diagram}\dgARROWLENGTH=3em
\node{X}\arrow{s,t}{\alpha} \arrow{se,t}{\e} \\
\node{M}\arrow{e,b}{\upsilon}  \node{N}
\end{diagram}
$$
As a corollary, the following diagram is commutative as well:
\beq\label{PROOF:sushestv-uzlov-razlozhenija-2}
\begin{diagram}\dgARROWLENGTH=3em
\node{X}\arrow{e,t}{\ph}\arrow{s,t}{\alpha} \arrow{se,t}{\e}\node{Y} \\
\node{M}\arrow{e,b}{\upsilon}  \node{N}\arrow{n,r}{\beta}
\end{diagram}
\eeq
Similarly, $\beta$ is a monomorphism, so it is an enrichment of $Y$ in $\Mono$ by means of $\ph$. At the same time, $\mu$ is a refinement of $Y$ in $\Mono$ by means of morphism $\ph$. Hence, there exists a morphism $\upsilon'$ such that the following diagram is commutative:
$$
\begin{diagram}\dgARROWLENGTH=3em
\node[2]{Y} \\
\node{M}\arrow{ne,t}{\mu}\arrow{e,b}{\upsilon'}  \node{N}\arrow{n,r}{\beta}
\end{diagram}
$$
As a corollary, the following diagram is commutative as well:
\beq\label{PROOF:sushestv-uzlov-razlozhenija-3}
\begin{diagram}\dgARROWLENGTH=3em
\node{X}\arrow{e,t}{\ph}\arrow{s,t}{\alpha} \node{Y} \\
\node{M}\arrow{e,b}{\upsilon'}\arrow{ne,t}{\mu}  \node{N}\arrow{n,r}{\beta}
\end{diagram}
\eeq
From \eqref{PROOF:sushestv-uzlov-razlozhenija-2} and \eqref{PROOF:sushestv-uzlov-razlozhenija-3} we have:
$$
\underset{\scriptsize\begin{matrix}\text{\rotatebox{90}{$\owns$}}\\ \Mono\end{matrix}}{\beta}\circ\upsilon\circ\underset{\scriptsize\begin{matrix}\text{\rotatebox{90}{$\owns$}}\\ \Epi\end{matrix}}{\alpha}=\ph=\underset{\scriptsize\begin{matrix}\text{\rotatebox{90}{$\owns$}}\\ \Mono\end{matrix}}{\beta}\circ\upsilon'\circ\underset{\scriptsize\begin{matrix}\text{\rotatebox{90}{$\owns$}}\\ \Epi\end{matrix}}{\alpha}
\quad\Longrightarrow\quad
\upsilon=\upsilon'
$$
I.e. the following diagram is commutative:
$$
\begin{diagram}\dgARROWLENGTH=3em
\node{X}\arrow[2]{r,t}{\ph}\arrow[2]{s,t}{\alpha}\arrow{se,t,-}{} \node[2]{Y} \\
 \node[2]{} \arrow{se,t}{\e} \\
\node{M}\arrow[2]{ne,t,3}{\mu}\arrow[2]{e,b}{\upsilon}  \node[2]{N}\arrow[2]{n,r}{\beta}
\end{diagram}
$$
Here $\e=\upsilon\circ\alpha$ is an epimorphism, hence $\upsilon$ is an epimorphism as well. On the other hand, $\mu=\beta\circ\upsilon$ is a monomorphism, so $\upsilon$ is a monomorphism as well. Thus, $\upsilon$ is a bimorphism, and $\ph=\beta\circ\upsilon\circ\alpha$ is a nodal decomposition of $\ph$.
\epr

\btm\label{TH:env^SEpi=>uzlo-razl} Suppose that in a category ${\tt K}$
 \bit{
 \item[\rm (a)] strong epimorphisms discern monomorphisms and strong monomorphisms discern epimorphisms\footnote{See  definitions on p.\pageref{DEF:strog-epi-razlich-mono}.},

 \item[\rm (b)] each object $X$ has an envelope in the class $\SEpi$ of all strong epimorphisms with respect to an arbitrary morphism that goes from $X$, and, dually, in each object $X$ there is a refinement in the class $\SMono$ of all strong monomorphisms by means of an arbitrary morphisms that comes to $X$.
   }\eit
Then ${\tt K}$ is a category with nodal decomposition.
 \etm
\bpr Take a morphism $\ph:X\to Y$.

1. By condition (b), there is an envelope $\env_{\ph}^{\SEpi}X:X\to \Env_{\ph}^{\SEpi}X$ of the object $X$ in the class $\SEpi$ of all strong epimorphisms with respect to the morphism $\ph$. Denote by $\alpha$ the morphism that continues $\ph$ at $\Env_{\ph}^{\SEpi}X$:
$$
\xymatrix %@R=4.0pc
@C=4.0pc
{
X\ar[d]_{\env_{\ph}^{\SEpi}X}\ar[rd]^{\ph} &
\\
\Env_{\ph}^{\SEpi}X\ar@{-->}[r]_{\alpha} & Y
}
$$

2. Similarly, by (b) there is a refinement $\rf_{\ph}^{\SMono}Y:\Rf_{\ph}^{\SMono}Y\to Y$ of the object $Y$ in the class $\SMono$ of all strong monomorphisms by means of $\ph$. Denote by $\beta$ the morphism that lifts $\ph$ to $\Rf_{\ph}^{\SMono}X$:
$$
\xymatrix %@R=4.0pc
 @C=4.0pc
{
X\ar[rd]_{\ph}\ar@{-->}[r]^{\beta} & \Rf_{\ph}^{\SEpi}Y\ar[d]^{\rf_{\ph}^{\SMono}Y}
\\
 & Y
}
$$

3. Pasting these triangles together by the common side $\ph$, and throwing away this side, we obtain a quadrangle:
$$
\xymatrix %@R=4.0pc
@C=4.0pc
{
X\ar[d]_{\env_{\ph}^{\SEpi}X}\ar[r]^{\beta} & \Rf_{\ph}^{\SEpi}Y\ar[d]^{\rf_{\ph}^{\SMono}Y}
\\
\Env_{\ph}^{\SEpi}X\ar[r]_{\alpha} & Y
}
$$
Here $\env_{\ph}^{\SEpi}X$ is a strong epimorphism, and $\rf_{\ph}^{\SMono}Y$ a monomorphism, so there is a diagonal $\delta$
\beq\label{PROOF:env^SEpi=>uzlo-razl}
\xymatrix %@R=4.0pc
@C=4.0pc
{
X\ar[d]_{\env_{\ph}^{\SEpi}X}\ar[r]^{\beta} & \Rf_{\ph}^{\SEpi}Y\ar[d]^{\rf_{\ph}^{\SMono}Y}
\\
\Env_{\ph}^{\SEpi}X\ar[r]_{\alpha}\ar@{-->}[ru]_{\delta} & Y
}
\eeq
Let us show that $\delta$ is a bimorphism.

4. Suppose first that $\delta$ is not a monomorphism. Then, since the strong epimorphisms discern monomorphisms (by (a)), there is a decomposition  $\delta=\delta'\circ\e$, where $\e$ is a strong epimorphism, which is not an isomorphism. As a corollary, the following diagram is commutative:
$$
\xymatrix %@R=4.0pc
@C=4.0pc
{
X\ar[d]_{\env_{\ph}^{\SEpi}X}\ar[rd]_{\beta}\ar[r]^{\ph} & Y
\\
\Env_{\ph}^{\SEpi}X\ar[r]_{\delta}\ar@{-->}[d]_{\e} & \Rf_{\ph}^{\SEpi}Y\ar[u]_{\rf_{\ph}^{\SMono}Y} \\
M\ar@/_2ex/@{-->}[ru]_{\delta'}
}
$$
We see here that the composition $\rf_{\ph}^{\SMono}Y\circ\delta'$ is a continuation of $\ph$ along $\e\circ\env_{\ph}^{\SEpi}X$, which in its turn is a strong epimorphism (as a comopsition of two strong epimorphisms). This means that $\e\circ\env_{\ph}^{\SEpi}X$ is an extension of $X$ in the class $\SEpi$ with respect to morphism $\ph$. Hence, there is a morphism $\upsilon$ from the extension $M$ to the envelope $\Env_{\ph}^{\SEpi}X$, such that diagram \eqref{DEF:diagr-obolochka} is commutative:
$$
\xymatrix %@R=4.0pc
@C=4.0pc
{
 & X\ar@/_2ex/[dl]_{\env_{\ph}^{\SEpi}X}\ar@/^2ex/[rd]^{\e\circ\env_{\ph}^{\SEpi}X} &
\\
\Env_{\ph}^{\SEpi}X & & M\ar@{-->}[ll]_{\upsilon} \\
}
$$
We have now $\upsilon\circ\e\circ\env_{\ph}^{\SEpi}X=\env_{\ph}^{\SEpi}X=1_M\circ\env_{\ph}^{\SEpi}X$, and, since $\env_{\ph}^{\SEpi}X$ is an epimorphism, this implies the equality $\upsilon\circ\e=1_M$, which means that $\e$ is a coretraction. On the other hand, this is an epimorphism, and together this means that $\e$ must be an isomorphism. This contradicts to the choice of $\e$.

5. Thus, $\delta$ must be a monomorphism. By analogy we prove that this is an epimorphism. Let us now add $\ph$ to Diagram  \eqref{PROOF:env^SEpi=>uzlo-razl} and twist it as follows:
$$
\xymatrix @R=4.0pc @C=4.0pc
{
X\ar[d]_{\env_{\ph}^{\SEpi}X}\ar[r]^{\ph}\ar[dr]^(.7){\beta}|!{[d];[r]}\hole & Y
\\
\Env_{\ph}^{\SEpi}X\ar[r]_{\delta}\ar[ur]^(.7){\alpha} & \Rf_{\ph}^{\SEpi}Y\ar[u]_{\rf_{\ph}^{\SMono}Y}
}
$$
We see now that $\ph=\rf_{\ph}^{\SMono}Y\circ\delta\circ\env_{\ph}^{\SEpi}X$ is a nodal decomposition of $\ph$.
\epr

\subsection{Nets and functoriality.}

Apparently, in the general case the operations of taking envelopes and refinements are not functors. But under some assumptions they are, and in the last part of this section we discuss this. Let us use the following definition. Suppose $\varOmega$ and $\varPhi$ are classes of morphisms in a category ${\tt K}$.

\bit{
    \item
Let us say that {\it the envelope $\Env^\varOmega_\varPhi$ can be defined as a functor}, if there exist
        \bit{
\item[E.1] a map $X\mapsto (E(X),e_X)$, that to each object $X$ in $\tt K$ assigns a morphism $e_X:X\to E(X)$ in ${\tt K}$, which is an envelope in $\varOmega$ with respect to $\varPhi$:
$$
E(X)=\Env_{\varPhi}^\varOmega X,\qquad e_X=\env_{\varPhi}^\varOmega X
$$
    \item[E.2] a map $\alpha\mapsto E(\alpha)$, that each morphism $\alpha:X\to Y$ in $\tt K$ turns into a morphism  $E(\alpha):E(X)\to E(Y)$ in $\tt K$ in such a way that the following diagram is commutative
    \beq\label{DIAGR:funktorialnost-env-e-E}
\xymatrix @R=2.pc @C=5.0pc % @M=14pt
{
X\ar[d]^{\alpha}\ar[r]^{e_X} & E(X)\ar@{-->}[d]^{E(\alpha)} \\
Y\ar[r]^{e_Y} & E(Y) \\
}
\eeq
}\eit
and the following identities hold
\beq\label{tozhdestva:funktorialnost-env-e-E}
E(1_X)=1_{E(X)},\qquad E(\beta\circ\alpha)=E(\beta)\circ E(\alpha)
\eeq
Clearly, in this case the map $(X,\alpha)\mapsto(E(X),E(\alpha))$ is a covariant functor from ${\tt K}$ into ${\tt K}$, and the map $X\mapsto e_X$ is a natural transformation of the identity functor $(X,\alpha)\mapsto(X,\alpha)$ into the functor $(X,\alpha)\mapsto(E(X),E(\alpha))$.

    \item Let us say that {\it the envelope $\Env^\varOmega_\varPhi$ can be defined as an idempotent functor}, if in addition to E.1 and E.2 one can ensure the condition

 \bit{
 \item[E.3] for each object $X\in\Ob({\tt K})$ the morphism $e_{E(X)}:E(X)\to E(E(X))$ is the local identity:
\beq\label{e_(E(X))=1_(E(X))}
E(E(X))=E(X),\qquad e_{E(X)}=1_{E(X)}\qquad X\in \Ob({\tt K}).
\eeq
 }\eit
 }\eit

\brem\label{REM:E(e_X)=1_(E(X))}
If $\varOmega\subseteq\Epi$, then \eqref{e_(E(X))=1_(E(X))} implies
\beq\label{E(e_X)=1_(E(X))}
E(e_X)=1_{E(X)}\qquad X\in \Ob({\tt K}).
\eeq
Indeed, if we put $\alpha=e_X$ into \eqref{DIAGR:funktorialnost-env-e-E}, we obtain
$$
\xymatrix @R=2.pc @C=8.0pc % @M=14pt
{
X\ar[d]^{e_X}\ar[r]^{e_X} & E(X)\ar[d]^{E(e_X)} \\
E(X)\ar[r]^{e_{E(X)}=1_{E(X)}} & E(E(X))=E(X) \\
}
$$
i.e. $E(e_X)\circ e_X=1_{E(X)}\circ e_X$, and, since $e_X\in\varOmega\subseteq\Epi$, we can cancel it: $E(e_X)=1_{E(X)}$.
\erem

   \bit{
    \item
Let us say that {\it the refinement $\Rf^\varGamma_\varPhi$ can be defined as a functor}, if there exist
        \bit{
\item[R.1] a map $X\mapsto (I(X),i_X)$, that to each object $X$ in $\tt K$ assigns a morphism $i_X:I(X)\to X$ in ${\tt K}$, which is a refinement in $\varGamma$ by means of $\varPhi$:
$$
I(X)=\Rf_{\varPhi}^\varGamma X,\qquad i_X=\rf_{\varPhi}^\varGamma X
$$
    \item[R.2] a map $\alpha\mapsto I(\alpha)$, that each morphism $\alpha:X\gets Y$ in $\tt K$ turns into a morphism $I(\alpha):I(X)\gets I(Y)$ in $\tt K$ in such a way that the following diagram is commutative
    \beq\label{DIAGR:funktorialnost-imp-i-I}
\xymatrix @R=2.pc @C=5.0pc % @M=14pt
{
X & I(X)\ar[l]_{i_X} \\
Y\ar[u]_{\alpha} & I(Y)\ar[l]_{i_Y}\ar@{-->}[u]_{I(\alpha)} \\
}
\eeq
}\eit
and the following identities hold
\beq\label{tozhdestva:funktorialnost-imp-i-I}
I(1_X)=1_{I(X)},\qquad I(\beta\circ\alpha)=I(\beta)\circ I(\alpha)
\eeq
In this case the map $(X,\alpha)\mapsto(I(X),I(\alpha))$ is a covariant functor from ${\tt K}$ into ${\tt K}$, and the map $X\mapsto i_X$ is a natural transformation of the identity functor $(X,\alpha)\mapsto(X,\alpha)$ into the functor $(X,\alpha)\mapsto(I(X),I(\alpha))$.

     \item  Let  us say that {\it the refinement $\Rf^\varGamma_\varPhi$ can be defined as an idempotent functor}, if in addition to R.1 and R.2 one can ensure the condition

 \bit{
 \item[R.3] for each object $X\in\Ob({\tt K})$ the morphism $i_{I(X)}:I(X)\gets I(I(X))$ is the local identity:
\beq\label{i_(I(X))=1_(I(X))}
I(I(X))=I(X),\qquad i_{I(X)}=1_{I(X)}\qquad X\in \Ob({\tt K}).
\eeq
 }\eit
 }\eit

\brem\label{REM:I(i_X)=1_(I(X))}
If $\varGamma\subseteq\Mono$, then \eqref{i_(I(X))=1_(I(X))} implies
\beq\label{I(i_X)=1_(I(X))}
I(i_X)=1_{I(X)}\qquad X\in \Ob({\tt K}).
\eeq
\erem

\brem For envelopes in the most important cases when the class $\varOmega$ consists of epimorphisms,  $\varOmega\subseteq\Epi$, the identities \eqref{tozhdestva:funktorialnost-env-e-E} automatically follow from E.1 and E.2. Dually, for refinements, when $\varGamma$ consists of monomorphisms, $\varGamma\subseteq\Mono$, the identities  \eqref{tozhdestva:funktorialnost-imp-i-I} automatically follow from R.1 and R.2.
\erem

\paragraph{Nets of epimorphisms.}

\bit{ \item[$\bullet$] Suppose that to each object $X\in\Ob({\tt K})$ in a category ${\tt K}$ it is assigned a subset ${\mathcal N}^X$ in the class
$\Epi^X$ of all epimorphisms of the category ${\tt K}$, going from $X$, and the following three requirements are fulfilled:
    \bit{
\item[(a)]\label{AX:set-Epi-a} for each object $X$ the set ${\mathcal N}^X$ is non-empty and is directed to the left with respect to the pre-order \eqref{DEF:le-in-F_X} inherited from $\Epi^X$:
$$
\forall \sigma,\sigma'\in {\mathcal N}^X\quad \exists\rho\in{\mathcal N}^X\quad
\rho\to\sigma\ \& \ \rho\to\sigma',
$$

\item[(b)]\label{AX:set-Epi-b} for each object $X$ the covariant system of morphisms generated by ${\mathcal N}^X$  \beq\label{DEF:kategoriya-svyazyv-morpfizmov}
    \Bind({\mathcal N}^X):=\{\iota_\rho^\sigma;\ \rho,\sigma\in{\mathcal N}^X,\ \rho\to\sigma\}
    \eeq
    (the morphisms $\iota_\rho^\sigma$ were defined in \eqref{DEF:le-in-F_X-*}; by \eqref{iota_rho^tau=iota_rho^sigma-circ-iota_sigma^tau} this system is a covariant functor from the set ${\mathcal N}^X$ considered as a full subcategory in $\Epi^X$ into $\tt K$) has a projective limit in $\tt K$;

\item[(c)]\label{AX:set-Epi-c}  for each morphism $\alpha:X\to Y$ and for each element $\tau\in{\mathcal N}^Y$ there are an element $\sigma\in{\mathcal N}^X$ and a morphism $\alpha_\sigma^\tau:\Ran\sigma\to\Ran\tau$ such that the following diagram is commutative
    \beq\label{DIAGR:set} \xymatrix @R=2.5pc @C=4.0pc {
 X\ar[r]^{\alpha}\ar@{-->}[d]_{\sigma} & Y\ar[d]^{\tau} \\
 \Ran\sigma\ar@{-->}[r]_{\alpha_\sigma^\tau} & \Ran\tau
    } \eeq
    (a remark: for given $\alpha$, $\sigma$ and $\tau$ the morphism $\alpha_\sigma^\tau$, if exists, must be unique, since $\sigma$ is an epimorphism).

 }\eit
Then
 \bit{
\item[---] we call the family of sets ${\mathcal N}=\{{\mathcal N}^X;\ X\in\Ob({\tt K})\}$ a {\it net of epimorphisms}\label{DEF:set-epimorf} in the category ${\tt K}$, and the elements of the sets ${\mathcal N}^X$ {\it elements of the net} ${\mathcal N}$,

\item[---] for each object $X$ the system of morphisms $\Bind({\mathcal N}^X)$ defined by equalities \eqref{DEF:kategoriya-svyazyv-morpfizmov} will be called the {\it system of binding morphisms of the net ${\mathcal N}$ over the vertex $X$}, its projective limit (which exists by condition (b)) is a projective cone, whose vertex will be denoted by $X_{\mathcal N}$, and the morphisms going from it by $\sigma_{\mathcal N}=\leftlim_{\rho\in{\mathcal N}^X} \iota^\sigma_\rho: X_{\mathcal N}\to \Ran\sigma$:
 \beq\label{X_F-proektiv-sistema} \xymatrix @R=2.5pc @C=2.0pc {
 & X_{\mathcal N}\ar[dr]^{\sigma_{\mathcal N}}\ar[dl]_{\rho_{\mathcal N}} &  \\
 \Ran\rho\ar[rr]^{\iota^\sigma_\rho} & &  \Ran\sigma
}\qquad (\rho\to\sigma);
 \eeq
in addition, by \eqref{DEF:le-in-F_X-*}, the system of epimorphisms ${\mathcal N}^X$ is also a projective cone of the system $\Bind({\mathcal N}^X)$:
 \beq\label{F_X-proektiv-sistema} \xymatrix @R=2.5pc @C=2.0pc {
 & X\ar[dr]^{\sigma}\ar[dl]_{\rho} &  \\
 \Ran\rho\ar[rr]^{\iota^\sigma_\rho} & &  \Ran\sigma
}\qquad (\rho\to\sigma),
 \eeq
 so there must exist a natural morphism from $X$ into the vertex $X_{\mathcal N}$ of the projective limit of the system $\Bind({\mathcal N}^X)$. We  denote this morphism by $\leftlim{\mathcal N}^X$ and call it the {\it local limit of the net ${\mathcal N}$ of epimorphisms at the object $X$}: \beq\label{DIAGR:sigma-sigma_F}
\xymatrix @R=2.5pc @C=2.0pc {
   X\ar[dr]_{\sigma}\ar@{-->}[rr]^{\leftlim{\mathcal N}^X} & & X_{\mathcal N}\ar[ld]^{\sigma_{\mathcal N}}  \\
   & \Ran\sigma  &
}\qquad (\sigma\in{\mathcal N}^X). \eeq

\item[---] the element $\sigma$ of the net in diagram \eqref{DIAGR:set} will be called a {\it counterfort} of the element $\tau$ of the net.
 }\eit
 }\eit

The examples of net of epimorphisms will be given in \ref{SEC:proj-ster-algebry}\ref{SUBSEC;Arens-Michael} and in \ref{SEC:proj-ster-algebry}\ref{SUBSEC:C^*-envelopes}.

\btm\label{TH:funktorialnost-obolochki_F} Let ${\mathcal N}$ be a net of epimorphisms in a category ${\tt K}$. Then \bit{
\item[(a)] for each object $X$ in ${\tt K}$ the local limit
$\leftlim{\mathcal N}^X:X\to X_{\mathcal N}$ is an envelope $\env_{\mathcal N} X$ in the category ${\tt
K}$ with respect to the class of morphisms ${\mathcal N}$:
 \beq\label{lim F_X=env_F-X} \leftlim {\mathcal N}^X=\env_{\mathcal N} X,
 \eeq

\item[(b)] for each morphism $\alpha:X\to Y$ in ${\tt K}$ and for any choice of local limits $\leftlim {\mathcal N}^X$ and $\leftlim {\mathcal N}^Y$ the formula
 \beq\label{DEF:alpha_F}
\alpha_{\mathcal N}=\leftlim_{\tau\in{\mathcal
N}_Y}\leftlim_{\sigma\in{\mathcal N}^X}\alpha_\sigma^\tau\circ\sigma_{\mathcal
N}
 \eeq
defines a morphism $\alpha_{\mathcal N}:X_{\mathcal N}\to Y_{\mathcal N}$ such that the following diagram is commutative:
 \beq\label{DIAGR:funktorialnost-lim-F}
\xymatrix @R=2.pc @C=10.0pc % @M=14pt
{
X\ar[d]^{\alpha}\ar[r]^{\leftlim {\mathcal N}^X=\env_{\mathcal N}X} &  X_{\mathcal N}
=\Env_{\mathcal N}X\ar@{-->}[d]^{\alpha_{\mathcal N}} \\
Y\ar[r]^{\leftlim {\mathcal N}^Y=\env_{\mathcal N}Y} &  Y_{\mathcal N}=\Env_{\mathcal N}Y
},
 \eeq

\item[(c)] the envelope $\Env_{\mathcal N}$ can be defined as a functor.
 }\eit
 \etm
\bpr
1. By Lemma \ref{LM:obolochka-konusa} the projective limit $\leftlim
{\mathcal N}^X$ is an envelope of $X$ in $\tt K$ with respect to the cone of morphisms
${\mathcal N}^X$:
$$
\leftlim {\mathcal N}^X=\env_{{\mathcal N}^X} X
$$
Here one can replace ${\mathcal N}^X$ by ${\mathcal N}$, since ${\mathcal N}^X$ is exactly the subclass in ${\mathcal N}$ consisting of morphisms with $X$ as domain:
$$
\leftlim {\mathcal N}^X=\env_{{\mathcal N}^X} X=\env_{\mathcal N} X.
$$

2. Let us explain first the sense of formula \eqref{DEF:alpha_F}. Take a morphism $\alpha:X\to Y$. For each element $\tau\in{\mathcal N}^Y$ of the net denote
 \beq\label{opredelenie-alpha^tau} \alpha^\tau=\tau\circ\alpha.
 \eeq
Clearly, the family of morphisms $\{\alpha^\tau:X\to\Ran\tau;\ \tau\in{\mathcal N}^Y\}$ is a projective cone of the system of binding morphisms  $\Bind({\mathcal N}^Y)$:
 \beq\label{X-alpha^tau=konus}
 \xymatrix @R=2.5pc @C=2.0pc {
 & X\ar[dr]^{\alpha^\upsilon}\ar[dl]_{\alpha^\tau} &  \\
 \Ran\tau\ar[rr]^{\iota_\tau^\upsilon} & &  \Ran\upsilon
}\qquad (\tau\to\upsilon).
 \eeq

By property (c) for each element $\tau\in{\mathcal N}^Y$ there are an element $\sigma\in{\mathcal N}^X$ and a morphism $\alpha_\sigma^\tau:\Ran\sigma\to\Ran\tau$ such that diagram \eqref{DIAGR:set} is commutative, and we have already denoted by $\alpha^\tau$ the diagonal there:
\beq\label{alpha^tau=tau-circ-alpha=alpha_sigma^tau-circ-sigma}
\alpha^\tau=\tau\circ\alpha=\alpha_\sigma^\tau\circ\sigma.
\eeq
Put
 \beq\label{opredelenie-beta^tau} \alpha_{\mathcal
N}^\tau=\alpha_\sigma^\tau\circ\sigma_{\mathcal N},
 \eeq
then we  obtain a diagram
 \beq\label{DIAFR:alpha^tau} \xymatrix @R=2.5pc @C=2.0pc {
   X\ar@/_4ex/[ddr]_{\alpha^\tau}\ar[dr]_{\sigma}\ar[rr]^{\leftlim{\mathcal N}^X} & & X_{\mathcal N}\ar[ld]^{\sigma_{\mathcal N}}\ar@{-->}@/^4ex/[ddl]^{\alpha_{\mathcal N}^\tau}  \\
   & \Ran\sigma\ar[d]^{\alpha_\sigma^\tau}  & \\
   & \Ran\tau &
 }\qquad (\sigma\in{\mathcal N}^X).
 \eeq

Note then that for any other element $\rho\in{\mathcal N}^X$ such that $\rho\to\sigma$ the following equality analogous to \eqref{opredelenie-beta^tau} is true:
 \beq\label{opredelenie-beta^tau-1}
\alpha_{\mathcal N}^\tau=\alpha_\rho^\tau\circ\rho_{\mathcal N},\qquad
\rho\to\sigma.
 \eeq
Indeed, for $\rho\to\sigma$ diagram \eqref{DIAGR:set} can be added to the diagram
 \beq\label{DIAGR:set-3} \xymatrix
@R=2.5pc @C=4.0pc {
 & X\ar[r]^{\alpha}\ar[d]_{\sigma}\ar@/_3ex/[ddl]_{\rho} & Y\ar[d]^{\tau}  \\
 & \Ran\sigma\ar[r]_{\alpha_\sigma^\tau} & \Ran\tau  \\
 \Ran\rho\ar@{-->}@/_3ex/[rru]_{\alpha_\rho^\tau}\ar[ru]_{\iota_\rho^\sigma} & &
}
 \eeq
(here the dotted arrow is initially defined as composition $\alpha_\sigma^\tau\circ\iota_\rho^\sigma$, and then, since such an arrow if it exists is unique, we deduce that this is the morphism $\alpha_\rho^\tau$). After that we have:
$$
\alpha_{\mathcal N}^\tau=\alpha_\sigma^\tau\circ\sigma_{\mathcal
N}=\eqref{X_F-proektiv-sistema}=
\alpha_\sigma^\tau\circ\iota_\rho^\sigma\circ\rho_{\mathcal
N}=\eqref{DIAGR:set-3}= \alpha_\rho^\tau\circ\rho_{\mathcal N}.
$$

From \eqref{opredelenie-beta^tau-1} it follows that the definition of $\alpha_{\mathcal N}^\tau$ by \eqref{opredelenie-beta^tau} does not depend on the choice of element $\sigma\in{\mathcal N}^X$, since if $\sigma'\in{\mathcal N}^X$ is another element such that there exists a morphism  $\alpha_{\sigma'}^\tau:\Ran\sigma'\to\Ran\tau$ for which diagram \eqref{DIAGR:set} is commutative (where $\sigma$ is replaced by
$\sigma'$), then we can take $\rho\in{\mathcal N}^X$ standing from the left of $\sigma$ and $\sigma'$,
$$
\rho\to\sigma,\qquad \rho\to\sigma',
$$
(at this moment we use Axiom (a) of the net of epimorphisms) and we have the chain
$$
\alpha_{\mathcal N}^\tau=\alpha_\sigma^\tau\circ\sigma_{\mathcal
N}=\eqref{opredelenie-beta^tau-1}= \alpha_\rho^\tau\circ\rho_{\mathcal
N}=\eqref{opredelenie-beta^tau-1}=\alpha_{\sigma'}^\tau\circ\sigma'_{\mathcal
N}.
$$

We can deduce now that formula \eqref{opredelenie-beta^tau} correctly defines a map $\tau\in{\mathcal N}^Y\mapsto\alpha_{\mathcal N}^\tau$. Let us show that the family of morphisms $\{\alpha_{\mathcal N}^\tau:X_{\mathcal N}\to \Ran\tau;\ \tau\in{\mathcal N}^Y\}$ is a projective cone of the system of binding morphisms $\Bind({\mathcal N}^Y)$:
 \beq\label{X_F-proektiv-sistema-dllya-F_Y} \xymatrix
@R=2.5pc @C=2.0pc {
 & X_{\mathcal N}\ar[dr]^{\alpha_{\mathcal N}^\upsilon}\ar[dl]_{\alpha_{\mathcal N}^\tau} &  \\
 \Ran\tau\ar[rr]^{\iota_\tau^\upsilon} & &  \Ran\upsilon
}\qquad (\tau\to\upsilon\in{\mathcal N}^Y).
 \eeq
For $\tau\to\upsilon$ diagram \eqref{DIAGR:set} can be added to the diagram
 \beq\label{DIAGR:set-4}
\xymatrix @R=2.5pc @C=4.0pc
{
  X\ar[r]^{\alpha}\ar[d]_{\sigma} & Y\ar[d]^{\tau}\ar@/^3ex/[ddr]^{\upsilon} & \\
  \Ran\sigma\ar[r]_{\alpha_\sigma^\tau}\ar@{-->}@/_3ex/[rrd]_{\alpha_\sigma^\upsilon} & \Ran\tau\ar[rd]_{\iota_\tau^\upsilon} & \\
  & & \Ran\upsilon
}
 \eeq
(where the dotted arrow is initially defined as the composition $\iota_\tau^\upsilon\circ\alpha_\sigma^\tau$, and then, since such an arrow, if it exists, is unique, we deduce that this is the morphism $\alpha_\sigma^\upsilon$). Using this diagram we have:
$$
\iota_\tau^\upsilon\circ\alpha_{\mathcal
N}^\tau=\eqref{opredelenie-beta^tau}=\iota_\tau^\upsilon\circ\alpha_\sigma^\tau\circ\sigma_{\mathcal
N}= \eqref{DIAGR:set-4}=\alpha_\sigma^\upsilon\circ\sigma_{\mathcal
N}=\eqref{opredelenie-beta^tau}=\alpha_{\mathcal N}^\upsilon.
$$
From diagram \eqref{X_F-proektiv-sistema-dllya-F_Y} it follows now that there must exist a natural morphism $\alpha_{\mathcal N}$ from $X_{\mathcal N}$ into the projective limit $Y_{\mathcal N}$ of the system $\Bind({\mathcal N}^Y)$:
 \beq\label{sushestvovanie-alpha_F} \xymatrix @R=2.5pc
@C=2.0pc {
   X_{\mathcal N}\ar[dr]_{\alpha_{\mathcal N}^\tau}\ar@{-->}[rr]^{\alpha_{\mathcal N}} & & Y_{\mathcal N}\ar[ld]^{\tau_{\mathcal N}}  \\
   & \Ran\tau  &
}\qquad (\tau\in{\mathcal N}^Y).
 \eeq
Recall now that by Axiom (b) of the net the passage from $X$ to the projective limit $\leftlim\Bind({\mathcal N}^X)$ can be organized as a map.
The further steps on building $\alpha_{\mathcal N}$ (the choice of the vertex $X_{\mathcal N}$ of the cone $\leftlim\Bind({\mathcal N}^X)$, and then the choice of the arrow $\alpha_{\mathcal N}$ such that all the diagrams \eqref{sushestvovanie-alpha_F} are commutative) are also unambiguous, so the correspondence $\alpha\mapsto\alpha_{\mathcal N}$ can also be treated as a map.

3. Note further that for the morphisms $\alpha_{\mathcal N}$ the diagrams of the form \eqref{DIAGR:funktorialnost-lim-F} are commutative. In the diagram $$
\xymatrix @R=2.pc @C=5.0pc % @M=14pt
{
X\ar[dd]^{\alpha}\ar[rr]^{\leftlim{\mathcal N}^X}\ar[dr]_{\alpha^\tau} &  & X_{\mathcal N}\ar@{-->}[dd]^{\alpha_{\mathcal N}}\ar[dl]^{\alpha_{\mathcal N}^\tau} \\
 & \Ran\tau & \\
Y\ar[rr]^{\leftlim {\mathcal N}^Y}\ar[ur]^{\tau} & &  Y_{\mathcal N}\ar[ul]_{\tau_{\mathcal N}} \\
}
$$
all the inner triangles are commutative: the upper inner triangle is commutative because this is the perimeter of \eqref{DIAFR:alpha^tau}, the left inner triangle because this is a variant of formula \eqref{opredelenie-alpha^tau}, the lower inner triangle because this is up to notations diagram  \eqref{DIAGR:sigma-sigma_F}, and the right inner triangle because this is a rotated diagram \eqref{sushestvovanie-alpha_F}. So the following equalities are true:
$$
\tau_{\mathcal N}\circ\leftlim{\mathcal
N}_Y\circ\alpha=\alpha^\tau=\tau_{\mathcal N}\circ\alpha_{\mathcal
N}\circ\leftlim{\mathcal N}^X \qquad (\tau\in{\mathcal N}^Y)
$$
One can interpret this as follows: each of the morphisms $\leftlim{\mathcal N}^Y\circ\alpha$ and $\alpha_{\mathcal N}\circ\leftlim{\mathcal N}^X$
is a lifting of the projective cone $\{\alpha^\tau:X\to\Ran\tau;\ \tau\in{\mathcal N}^Y\}$ for the system of binding morphisms $\Bind({\mathcal
N}_Y)$ which we were talking about in diagram \eqref{X-alpha^tau=konus} to the projective limit of this system. I.e. $\leftlim{\mathcal N}^Y\circ\alpha$ and $\alpha_{\mathcal N}\circ\leftlim{\mathcal N}^X$ are the very same dotted arrow in the definition of projective limit, for which all the diagrams of the form
$$
\xymatrix @R=2.5pc @C=2.0pc
{
   X\ar[dr]_{\alpha^\tau}\ar@{-->}[rr] & & Y_{\mathcal N}\ar[ld]^{\tau_{\mathcal N}}  \\
   & \Ran\tau  &
}\qquad (\tau\in{\mathcal N}^Y).
$$
are commutative. But such an arrow is unique, so these morphisms must coincide:
$$
\leftlim{\mathcal N}^Y\circ\alpha=\alpha_{\mathcal N}\circ\leftlim{\mathcal
N}_X.
$$
This is exactly diagram \eqref{DIAGR:funktorialnost-lim-F}.

4. The theorem on well-ordering of the class of all sets \cite[V, 4.1]{Levy} allows to define the operation of taking local limit as a map:
$$
    X\mapsto \leftlim\Bind({\mathcal N}^X)
$$
(i.e. there is a map that assigns to each object $X\in\Ob(\tt K)$ a concrete projective limit of the subcategory $\Bind({\mathcal N}^X)$ among all its projective limits in ${\tt K}$). Let us show that in this case the arising map $(X,\alpha)\mapsto(X_{\mathcal N},\alpha_{\mathcal N})$ is a functor, i.e. the following identities hold:
 \beq\label{funktorialnost-alpha_F}
 (1_X)_{\mathcal N}= 1_{X_{\mathcal N}},\qquad
 (\beta\circ\alpha)_{\mathcal N}=\beta_{\mathcal N}\circ\alpha_{\mathcal N}.
 \eeq
Suppose first that $\alpha=1_X:X\to X$. Then
\begin{multline*}
\alpha^\tau=\eqref{opredelenie-alpha^tau}=\tau\circ\alpha=\tau\circ 1_X=\tau\quad\Longrightarrow\quad
\alpha_\sigma^\tau\circ\sigma=\eqref{alpha^tau=tau-circ-alpha=alpha_sigma^tau-circ-sigma}=\alpha^\tau=\tau=\eqref{DEF:le-in-F_X-*}=\iota_\sigma^\tau\circ\sigma
\quad\Longrightarrow\\ \Longrightarrow\quad \alpha_\sigma^\tau=\iota_\sigma^\tau \quad\Longrightarrow\quad  \alpha_{\mathcal N}^\tau=\iota_\sigma^\tau\circ\sigma_{\mathcal N}=\tau_{\mathcal N}
\end{multline*}
So in diagrams \eqref{sushestvovanie-alpha_F} we can replace $\alpha_{\mathcal N}^\tau$ by $\tau_{\mathcal N}$:
 $$
  \xymatrix @R=2.5pc
@C=2.0pc {
   X_{\mathcal N}\ar[dr]_{\tau_{\mathcal N}}\ar@{-->}[rr]^{\alpha_{\mathcal N}} & & X_{\mathcal N}\ar[ld]^{\tau_{\mathcal N}}  \\
   & \Ran\tau  &
}\qquad (\tau\in{\mathcal N}^X).
 $$
These diagrams are commutative for all $\tau\in{\mathcal N}^X$, and the dotted arrow $\alpha_{\mathcal N}$ is defined here as the lifting of the projective cone $\{\alpha_{\mathcal N}^\tau=\tau_{\mathcal N}:X_{\mathcal N}\to\Ran\tau\}$ to the projective limit $\{\tau_{\mathcal N}:X_{\mathcal N}\to\Ran\tau\}$. Such an arrow is unique, so it must coincide with the morphism $1_{X_{\mathcal N}}$, for which all these diagrams are trivially commutative: $\alpha_{\mathcal N}=1_{X_{\mathcal N}}$.

Let us now prove the second identity in \eqref{funktorialnost-alpha_F}. Consider the sequence of morphisms $X\overset{\alpha}{\longrightarrow}Y\overset{\beta}{\longrightarrow}Z$. Take an element $\upsilon\in{\mathcal N}^Z$ and, using Axiom (c), let us choose first an element $\tau\in{\mathcal N}^Y$ and a morphism $\beta_\tau^\upsilon$ such that
$$
\upsilon\circ\beta=\beta_\tau^\upsilon\circ\tau,
$$
And then again using Axiom (c) choose an element $\sigma\in{\mathcal N}^X$ and a morphism $\alpha_\sigma^\tau$ such that
$$
\tau\circ\alpha=\alpha_\sigma^\tau\circ\sigma.
$$
We obtain the following diagram:
$$
\xymatrix @R=2.5pc @C=2.0pc
{
   X\ar[d]^{\sigma}\ar[r]^{\alpha} & Y\ar[d]^{\tau}\ar[r]^{\beta} & Z\ar[d]^{\upsilon}  \\
   \Ran\sigma\ar[r]^{\alpha_\sigma^\tau} & \Ran\tau\ar[r]^{\beta_\tau^\upsilon}  & \Ran\upsilon
}.
$$
If we remove here the middle arrow, then we obtain a diagram
$$
\xymatrix @R=2.5pc @C=4.0pc
{
   X\ar[d]^{\sigma}\ar[r]^{\beta\circ\alpha} & Z\ar[d]^{\upsilon}  \\
   \Ran\sigma\ar[r]^{\beta_\tau^\upsilon\circ\alpha_\sigma^\tau} & \Ran\upsilon
},
$$
which can be understood in such a way that the morphism $\beta_\tau^\upsilon\circ\alpha_\sigma^\tau$ is exactly the unique dotted arrow from diagram \eqref{DIAGR:set}, but the difference is that $Y$ is replaced here by $Z$, $\alpha$ by $\beta\circ\alpha$, and $\tau$ by $\upsilon$. Hence we can deduce that there exists a morphism $(\beta\circ\alpha)_\sigma^\upsilon$ which coincide with $\beta_\tau^\upsilon\circ\alpha_\sigma^\tau$:
\beq\label{beta_tau^upsilon-circ-alpha_sigma^tau=(beta-circ-alpha)_sigma^upsilon}
\beta_\tau^\upsilon\circ\alpha_\sigma^\tau=(\beta\circ\alpha)_\sigma^\upsilon
\eeq
This equality is used in the following chain:
$$
\underbrace{\upsilon_{\mathcal N}\circ\beta_{\mathcal N}}_{\scriptsize
\begin{matrix}\phantom{\tiny \eqref{sushestvovanie-alpha_F}} \ \text{\rotatebox{90}{$=$}} \ {\tiny \eqref{sushestvovanie-alpha_F}}\\ \beta_{\mathcal N}^\upsilon \\ \phantom{\tiny \eqref{opredelenie-beta^tau}}  \ \text{\rotatebox{90}{$=$}} \ {\tiny \eqref{opredelenie-beta^tau}} \\
\beta_\tau^\upsilon\circ\tau_{\mathcal N} \end{matrix}} \circ\ \alpha_{\mathcal
N}= \beta_\tau^\upsilon\circ\kern-3pt \underbrace{\tau_{\mathcal
N}\circ\alpha_{\mathcal N}}_{\scriptsize
\begin{matrix}\phantom{\tiny \eqref{sushestvovanie-alpha_F}} \ \text{\rotatebox{90}{$=$}} \ {\tiny \eqref{sushestvovanie-alpha_F}}\\
\alpha_{\mathcal N}^\tau \\ \phantom{\tiny \eqref{opredelenie-beta^tau}}  \ \text{\rotatebox{90}{$=$}} \ {\tiny \eqref{opredelenie-beta^tau}} \\
\alpha_\sigma^\tau\circ\sigma_{\mathcal N} \end{matrix}} =
\underbrace{\beta_\tau^\upsilon\circ\alpha_\sigma^\tau}_{\scriptsize
\begin{matrix}\phantom{\tiny \eqref{beta_tau^upsilon-circ-alpha_sigma^tau=(beta-circ-alpha)_sigma^upsilon}} \ \text{\rotatebox{90}{$=$}} \ {\tiny \eqref{beta_tau^upsilon-circ-alpha_sigma^tau=(beta-circ-alpha)_sigma^upsilon}}\\
(\beta\circ\alpha)_\sigma^\upsilon \end{matrix}} \kern-3pt\circ\
\sigma_{\mathcal N}=
\underbrace{(\beta\circ\alpha)_\sigma^\upsilon\circ\sigma_{\mathcal
N}}_{\scriptsize
\begin{matrix}\phantom{\tiny \eqref{opredelenie-beta^tau}} \ \text{\rotatebox{90}{$=$}} \ {\tiny \eqref{opredelenie-beta^tau}}\\
(\beta\circ\alpha)_{\mathcal N}^\upsilon  \end{matrix}}
=(\beta\circ\alpha)_{\mathcal
N}^\upsilon=\eqref{sushestvovanie-alpha_F}=\upsilon_{\mathcal N}\circ
(\beta\circ\alpha)_{\mathcal N}.
$$
If we omit the intermediate calculations, we arrive at the following double equality:
$$
\upsilon_{\mathcal N}\circ(\beta_{\mathcal N}\circ\alpha_{\mathcal
N})=(\beta\circ\alpha)_{\mathcal N}^\upsilon=\upsilon_{\mathcal N}\circ
(\beta\circ\alpha)_{\mathcal N}.
$$
This is true for each $\upsilon\in{\mathcal N}^Z$. So this can be treated as if both $\beta_{\mathcal N}\circ\alpha_{\mathcal
N}$ and $(\beta\circ\tau)_{\mathcal N}$ were liftings of the projective cone $\{(\beta\circ\alpha)_{\mathcal N}^\upsilon:X_{\mathcal
N}\to\Ran\upsilon;\ \upsilon\in{\mathcal N}^Z\}$ for the system of binding morphisms $\Bind({\mathcal N}^Z)$ (and this family is indeed a projective cone due to diagram \eqref{X_F-proektiv-sistema-dllya-F_Y} where one should replace $Y$ by $Z$, and $\alpha$ by $\beta\circ\alpha$) to the projective limit of this system. Thus, $\beta_{\mathcal N}\circ\alpha_{\mathcal N}$ and $(\beta\circ\tau)_{\mathcal N}$ are exactly the dotted arrow in the definition of projective limit, for which all the diagrams of the form
$$
\xymatrix @R=2.5pc @C=2.0pc
{
   X_{\mathcal N}\ar[dr]_{(\beta\circ\alpha)_{\mathcal N}^\upsilon}\ar@{-->}[rr] & & Z_{\mathcal N}\ar[ld]^{\upsilon_{\mathcal N}}  \\
   & \Ran\upsilon  &
}\qquad (\upsilon\in{\mathcal N}^Z).
$$
are commutative. But this dotted arrow is unique, so these morphisms must coincide:
$$
\beta_{\mathcal N}\circ\alpha_{\mathcal N}=(\beta\circ\tau)_{\mathcal N}.
$$
This is the identity \eqref{funktorialnost-alpha_F}.
\epr

\btm\label{TH:funktorialnost-obolochki} Let ${\mathcal N}$ be a net of epimorphisms in a category ${\tt K}$ that generates a class of morphisms $\varPhi$ on the inside:
    $$
    {\mathcal N}\subseteq\varPhi\subseteq\Mor({\tt K})\circ {\mathcal N}.
    $$
Then for any class of epimorphisms $\varOmega$ in $\tt K$, which contains all local limits
$\leftlim{\mathcal N}^X$,
\beq\label{lim_N-subseteq-varOmega}
\{\leftlim{\mathcal N}^X; \ X\in\Ob(\tt K)\}\subseteq\varOmega\subseteq\Epi(\tt
K),
\eeq
the following holds:
 \bit{
\item[(a)] for each object $X$ in ${\tt K}$ the local limit
$\leftlim{\mathcal N}^X$ is an envelope $\env_\varPhi^\varOmega X$ in
$\varOmega$ with respect to $\varPhi$: \beq\label{lim F_X=env_M^varOmega-X}
\leftlim {\mathcal N}^X=\env_\varPhi^\varOmega X, \eeq

\item[(b)] the envelope $\Env_\varPhi^\varOmega$ can be defined as a functor.
 }\eit
 \etm
 \bpr 1. By Theorem \ref{TH:funktorialnost-obolochki_F} the local limit of the net $\leftlim {\mathcal N}^X$ is an envelope of $X$ in the class
$\Mor({\tt K})$ of all morphisms of the category $\tt K$ with respect to the class of morphisms ${\mathcal N}$:
$$
\leftlim {\mathcal N}^X=\env_{\mathcal N} X:=\env_{\mathcal N}^{\Mor(\tt K)} X.
$$
On the other hand, by (i) $\leftlim {\mathcal N}^X$ belongs to a narrower class $\varOmega$, so by $1^\circ$ (c) on page \pageref{LM:suzhenie-verh-klassa-morfizmov}, $\leftlim {\mathcal N}^X$ must be an envelope in this narrower class $\varOmega$:
$$
\leftlim {\mathcal N}^X=\env_{\mathcal N} X=\env_{\mathcal N}^{\Mor(\tt K)}
X=\env_{\mathcal N}^\varOmega X.
$$
Further, since ${\mathcal N}$ generates $\varPhi$ on the inside, and $\varOmega$ consists of epimorphisms, by \eqref{env_Psi=env_Phi} the envelope with respect to ${\mathcal N}$ must coincide with the envelope with respect to $\varPhi$:
$$
\leftlim {\mathcal N}^X=\env_{\mathcal N} X=\env_{\mathcal N}^{\Mor(\tt K)}
X=\env_{\mathcal N}^\varOmega X=\env_\varPhi^\varOmega X.
$$
This proves \eqref{lim F_X=env_M^varOmega-X}. After that (b) follows from Theorem
\ref{TH:funktorialnost-obolochki_F}(c). \epr

One can get rid of the left side of \eqref{lim_N-subseteq-varOmega}, if the class $\varOmega$ is monomorphically complementable:

\btm\label{TH:funktorialnost-pri-seti-Epi-i-dolonyaemosti}
Let ${\mathcal N}$ be a net of epimorphisms in a category ${\tt K}$, that generates a class of morphisms $\varPhi$ on the inside:
    $$
    {\mathcal N}\subseteq\varPhi\subseteq\Mor({\tt K})\circ {\mathcal N}.
    $$
Then for each monomorphically complementable\footnote{See definition on p.\pageref{DEF:klass-monomorfno-dopolnyaem}.} class of epimorphisms $\varOmega$,
$$
{^\downarrow\varOmega}\circledcirc\varOmega={\tt K},
$$
the following holds:
 \bit{
\item[(a)] for each object $X$ in ${\tt K}$ the morphism $\e_{\leftlim{\mathcal N}^X}$ in the factrization  \eqref{faktorizatsiya-v-kat-s-faktoriz} defined by the classes ${^\downarrow\varOmega}$ and $\varOmega$, is an envelope  $\env_\varPhi^{\varOmega} X$ in $\varOmega$ with respect to $\varPhi$:
\beq\label{im_infty-lim-F_X=env_M^Epi-X-1}
\e_{\leftlim{\mathcal N}^X}=\env_\varPhi^{\varOmega} X,
\eeq

\item[(b)] for each morphism $\alpha:X\to Y$ in ${\tt K}$ and for any choice of envelopes $\env_\varPhi^{\varOmega}X$ and $\env_\varPhi^{\varOmega}Y$ there exists a unique morphism $\Env_\varPhi^{\varOmega} \alpha:\Env_\varPhi^{\varOmega}  X\to \Env_\varPhi^{\varOmega}Y$ in ${\tt K}$ such that the following diagram is commutative:
\beq\label{DIAGR:funktorialnost-env_varPhi^Epi-v-kat-s-uzl-razl-1}
\xymatrix @R=2.pc @C=5.0pc % @M=14pt
{
X\ar[d]^{\alpha}\ar[r]^{\env_\varPhi^{\varOmega}  X} & \Env_\varPhi^{\varOmega}  X\ar@{-->}[d]^{\Env_\varPhi^{\varOmega} \alpha} \\
Y\ar[r]^{\env_\varPhi^{\varOmega}  Y} & \Env_\varPhi^{\varOmega}  Y \\
}
\eeq

\item[(c)] if in addition $\tt K$ is co-well-powered in the class $\varOmega$, then the envelope  $\Env_\varPhi^{\varOmega}$ can be defined as a functor.
}\eit
 \etm \bpr
1. Since ${\mathcal N}$ generates $\varPhi$, and $\varOmega$ consists of epimorphisms, by \eqref{env_Psi=env_Phi} the envelope with respect to ${\mathcal N}$ coincides with the envelope with respect to $\varPhi$:
$$
\env_{\mathcal N}^\varOmega X=\env_\varPhi^\varOmega X.
$$
After that the equality \eqref{obolochka-konusa-v-kat-s-uzl-razl} from Lemma
\ref{LM:obolochka-konusa-v-kat-s-uzl-razl} implies \eqref{im_infty-lim-F_X=env_M^Epi-X-1}:
$$
\env_{\varPhi}^{\varOmega}X=\env_{\mathcal N}^{\varOmega}X=\e_{\leftlim{\mathcal
N}_X}.
$$

2. The property \eqref{DIAGR:funktorialnost-env_varPhi^Epi-v-kat-s-uzl-razl-1} is proved as follows. First we add Diagram \eqref{DIAGR:funktorialnost-lim-F} by decomposing limits $\leftlim{\mathcal
N}_X$ and $\leftlim{\mathcal N}^Y$ as follows:
$$
\xymatrix @R=3.pc @C=5.0pc % @M=14pt
{
X\ar@/^5ex/[rr]^{\leftlim{\mathcal N}^X}\ar[d]^{\alpha}\ar[r]_(.3){\env_\varPhi^{\varOmega}  X} & \Env_\varPhi^{\varOmega}X=\Dom\mu_{\leftlim{\mathcal N}^X}\ar[r]_(.7){\mu_{\leftlim{\mathcal N}^X}} & X_{\mathcal N}\ar[d]^{\alpha_{\mathcal N}} \\
Y\ar@/_5ex/[rr]_{\leftlim{\mathcal N}^Y}\ar[r]^(.3){\env_\varPhi^{\varOmega}  Y} & \Env_\varPhi^{\varOmega}Y=\Dom\mu_{\leftlim{\mathcal N}^Y}\ar[r]^(.7){\mu_{\leftlim{\mathcal N}^Y}} & Y_{\mathcal N} \\
}
$$
Then the inner quadrangle we represent as follows
$$
\xymatrix @R=3.pc @C=5.0pc % @M=14pt
{
X\ar[d]^{\alpha}\ar[r]^(.3){\env_\varPhi^{\varOmega}  X} & \Env_\varPhi^{\varOmega}X=\Dom\mu_{\leftlim{\mathcal N}^X} \ar@/^2ex/[dr]^(.7){\alpha_{\mathcal N}\circ\mu_{\leftlim{\mathcal N}^X}} & \\
Y\ar[r]^(.3){\env_\varPhi^{\varOmega}  Y} & \Env_\varPhi^{\varOmega}Y=\Dom\mu_{\leftlim{\mathcal N}^Y} \ar[r]_(.6){\mu_{\leftlim{\mathcal N}^Y}} & Y_{\mathcal N} \\
}
$$
Here the upper horizontal arrow, $\env_\varPhi^{\varOmega}X$, belongs to $\varOmega$, and the second lower horizontal arrow, $\mu_{\leftlim{\mathcal N}^Y}$, belongs to $\varGamma=\varOmega^\downarrow$. Hence there exists a morphism $\xi$ such that the following diagram is commutative:
$$
\xymatrix @R=3.pc @C=5.0pc % @M=14pt
{
X\ar[d]^{\alpha}\ar[r]^(.3){\env_\varPhi^{\varOmega}  X} & \Env_\varPhi^{\varOmega}X=\Dom\mu_{\leftlim{\mathcal N}^X} \ar@{-->}[d]^{\xi}\ar@/^2ex/[dr]^(.7){\alpha_{\mathcal N}\circ\mu_{\leftlim{\mathcal N}^X}} & \\
Y\ar[r]^(.3){\env_\varPhi^{\varOmega}  Y} & \Env_\varPhi^{\varOmega}Y=\Dom\mu_{\leftlim{\mathcal N}^Y} \ar[r]_(.6){\mu_{\leftlim{\mathcal N}^Y}} & Y_{\mathcal N} \\
}
$$
This $\xi$ will be the vertical arrow in \eqref{DIAGR:funktorialnost-env_varPhi^Epi-v-kat-s-uzl-razl-1} that we need.

3. Let $\tt K$ be co-well-powered in $\varOmega$, i.e. for each object $X$ the category  $\varOmega^X=\varOmega\cap\Epi^X$ is skeletally small. Let $S_X$ be its skeleton, which is a set. Using Theorem
\ref{TH:o-lok-malosti-v-faktor-objektah}, we can take a map $X\mapsto S_X$, which assigns to each object $X$ a skeleton $S_X$ in $\varOmega^X$. Let us fix this map. To define the envelope $\Env_\varPhi^{\varOmega}$ as a functor, we now define (by the axiom of choice) a map $X\in\Ob({\tt K})\mapsto \env_\varPhi^{\varOmega}X\in S_X$. Then the object $\Env_\varPhi^{\varOmega}X$ is defined as the domain of the morphism
$\env_\varPhi^{\varOmega}X$, and the morphism $\Env_\varPhi^{\varOmega}\alpha$
\eqref{DIAGR:funktorialnost-env_varPhi^Epi-v-kat-s-uzl-razl-1} arises automatically (as the unique possible morphism).
 \epr

\paragraph{Nets of monomorphisms.}

\bit{

\item[$\bullet$] Suppose that to each object $X\in\Ob({\tt K})$ in a category ${\tt K}$ it is assigned a subset ${\mathcal N}_X$ in the class $\Mono_X$ of all monomorphisms of ${\tt K}$ coming to $X$, and the following three requirements are fulfilled:
    \bit{

\item[(a)]\label{AX:set-Mono-a} for each object $X$ the set ${\mathcal N}_X$ is non-empty and is directed to the right with respect to pre-order \eqref{DEF:to-in-Mono(X)} inherited from $\Mono_X$:
$$
\forall \rho,\rho'\in {\mathcal N}_X\quad \exists\sigma\in{\mathcal N}_X\quad
\rho\to\sigma\ \& \ \rho'\to\sigma,
$$

\item[(b)]\label{AX:set-Mono-b} for each object $X$ the covariant system of morphisms generated by the set ${\mathcal N}_X$
 \beq\label{DEF:kategoriya-svyazyv-morpfizmov-*}
\Bind({\mathcal N}_X):=\{\varkappa^\sigma_\rho;\ \rho,\sigma\in{\mathcal N}_X,\
\rho\to\sigma\}
 \eeq
(the morphisms $\varkappa^\sigma_\rho$ were defined in \eqref{DEF:to-in-Mono(X)-*}; according to \eqref{varkappa^gamma_alpha=varkappa^gamma_beta-circ-varkappa^beta_alpha}, this system is a covariant functor from the set ${\mathcal N}_X$ considered as a full subcategory in $\Mono_X$ into $\tt K$) has an injective limit in $\tt K$;

\item[(c)]\label{AX:set-Mono-c} for each morphism $\alpha:X\to Y$ and for each element $\sigma\in{\mathcal N}_X$ there is an element  $\tau\in{\mathcal N}_Y$ and a morphism $\alpha_\sigma^\tau:\Dom\sigma\to\Dom\tau$ such that the following diagram is commutative:
 \beq\label{DIAGR:set-mono} \xymatrix @R=2.5pc @C=4.0pc {
 X\ar[r]^{\alpha} & Y \\
 \Dom\sigma\ar[u]^{\sigma}\ar@{-->}[r]_{\alpha_\sigma^\tau} & \Dom\tau\ar@{-->}[u]_{\tau}
 } \eeq
(a remark: for given $\alpha$, $\sigma$ and $\tau$ the morphism $\alpha_\sigma^\tau$, if exists, must be unique, since $\tau$ is a monomorphism).

 }\eit
Then
 \bit{
\item[---] we call the family of sets ${\mathcal N}=\{{\mathcal N}_X;\ X\in\Ob({\tt K})\}$ a {\it net of monomorphisms}\label{DEF:set-monomorf} in the category ${\tt K}$, and the elements of the sets ${\mathcal N}_X$ {\it elements of the net} ${\mathcal N}$,

\item[---] for each object $X$ the system of morphisms $\Bind({\mathcal N}_X)$ defined by equalities \eqref{DEF:kategoriya-svyazyv-morpfizmov-*} will be called a {\it system of binding morphisms of the net ${\mathcal N}$ over the vertex $X$}, its injective limit (which exists by condition (b)) is an injective cone whose vertex will be denoted by $X_{\mathcal N}$, and the morphisms coming to it by $\rho_{\mathcal N}=\rightlim_{\sigma\in{\mathcal N}_X} \varkappa^\sigma_\rho: X_{\mathcal N}\gets\Ran\sigma$:
 \beq\label{X_F-injektiv-sistema} \xymatrix @R=2.5pc @C=2.0pc {
 & X_{\mathcal N} &  \\
 \Dom\rho\ar[ur]^{\rho_{\mathcal N}}\ar[rr]^{\varkappa^\sigma_\rho} & &  \Dom\sigma\ar[ul]_{\sigma_{\mathcal N}}
}\qquad (\rho\to\sigma);
 \eeq
in addition, by \eqref{DEF:to-in-Mono(X)-*}, the system of monomorphisms ${\mathcal N}_X$ is also called an injective cone of the system  $\Bind({\mathcal N}_X)$:
 \beq\label{F_X-injektiv-sistema} \xymatrix @R=2.5pc @C=2.0pc {
 & X &  \\
 \Dom\rho\ar[rr]^{\varkappa^\sigma_\rho}\ar[ur]^{\rho} & &  \Dom\sigma\ar[ul]_{\sigma}
}\qquad (\rho\to\sigma),
 \eeq
so there must exist a natural morphism into $X$ from the vertex $X_{\mathcal N}$ of the injective limit of the system $\Bind({\mathcal N}_X)$. This morphism will be denoted by $\rightlim{\mathcal N}_X$ and will be called a {\it local limit of the net of monomorphisms ${\mathcal N}$ at the object  $X$}:
 \beq\label{DIAGR:sigma-sigma_F-mono}
 \xymatrix @R=2.5pc @C=2.0pc {
   X_{\mathcal N}\ar@{-->}[rr]^{\rightlim{\mathcal N}_X} & & X  \\
   & \Dom\sigma\ar[ur]_{\sigma}\ar[ul]^{\sigma_{\mathcal N}}  &
}\qquad (\sigma\in{\mathcal N}_X). \eeq

\item[---] the element $\tau$ of the net in diagram \eqref{DIAGR:set-mono} will be called a {\it shed} for the element $\sigma$ of the net.
 }\eit
 }\eit

The following propositions are dual to Theorems \ref{TH:funktorialnost-obolochki_F},
 \ref{TH:funktorialnost-obolochki} and \ref{TH:funktorialnost-pri-seti-Epi-i-dolonyaemosti}.

\btm\label{TH:funktorialnost-otpechatka_F} Let ${\mathcal N}$ be a net of monomorphisms in a category ${\tt K}$. Then \bit{
\item[(a)] for each object $X$ in ${\tt K}$ the local limit $\rightlim{\mathcal N}_X:X_{\mathcal N}\to X$ is a refinement $\rf_{\mathcal N} X$ of $X$ in the category ${\tt K}$ by means of the class of morphisms ${\mathcal N}$:
 \beq\label{lim F_X=imp_F-X}
\rightlim {\mathcal N}_X=\rf_{\mathcal N} X,
 \eeq

\item[(b)] for each morphism $\alpha:X\to Y$ in ${\tt K}$ and for any choice of local limits $\rightlim {\mathcal N}_X$ and $\rightlim {\mathcal N}_Y$ the formula
 \beq\label{DEF:alpha_F-mono}
\alpha_{\mathcal N}=\rightlim_{\sigma\in{\mathcal
N}_X}\rightlim_{\tau\in{\mathcal N}_Y}\tau_{\mathcal N}\circ\alpha_\sigma^\tau
 \eeq
defines a morphism $\alpha_{\mathcal N}:X_{\mathcal N}\to Y_{\mathcal N}$ such that the following diagram is commutative:
 \beq\label{DIAGR:funktorialnost-lim-F-mono}
\xymatrix @R=2.pc @C=10.0pc % @M=14pt
{
X\ar[d]^{\alpha} &  X_{\mathcal N}=\Rf_{\mathcal N} X \ar@{-->}[d]^{\alpha_{\mathcal N}}\ar[l]_{\rightlim {\mathcal N}_X=\rf_{\mathcal N} X} \\
Y &  Y_{\mathcal N}=\Rf_{\mathcal N} Y\ar[l]_{\rightlim {\mathcal N}_Y=\rf_{\mathcal N} Y} \\
}
 \eeq

\item[(c)] the refinement $\Rf_{\mathcal N}$ can be defined as a functor.
 }\eit
\etm

\btm\label{TH:funktorialnost-otpechatka}
Let ${\mathcal N}$ be a net of monomorphisms in a category ${\tt K}$, that generates a class of morphisms $\varPhi$ on the outside:
    $$
    {\mathcal N}\subseteq\varPhi\subseteq{\mathcal N}\circ\Mor({\tt K}).
    $$
Then for every class of monomorphisms $\varGamma$ in ${\tt K}$, that contains all local limits
$\rightlim{\mathcal N}_X$,
$$
\{\rightlim{\mathcal N}_X; \ X\in\Ob(\tt
K)\}\subseteq\varGamma\subseteq\Mono(\tt K),
$$
the following holds:
\bit{

\item[(a)] for each object $X$ in ${\tt K}$ the local limit $\rightlim{\mathcal N}_X$ is a refinement $\rf_\varPhi^\varGamma X$ in $\varGamma$ by means of $\varPhi$:
 \beq\label{lim F_X=imp_M^varOmega-X}
\rightlim {\mathcal N}_X=\rf_\varPhi^\varGamma X, \eeq

\item[(b)] the refinement $\Rf_\varPhi^\varGamma$ can be defined as a functor.
 }\eit
 \etm

\btm\label{TH:funktorialnost-pri-seti-Mono-i-dolonyaemosti}
Let ${\mathcal N}$ be a net of monomorphisms in a category ${\tt K}$, that generates the class $\varPhi$ on the outside:
    $$
    {\mathcal N}\subseteq\varPhi\subseteq{\mathcal N}\circ\Mor({\tt K}).
    $$
Then for every epimorphically complementable\footnote{See definition on p.\pageref{DEF:klass-epimorfno-dopolnyaem}.} class of monomorphisms $\varGamma$,
$$
\varGamma\circledcirc\varGamma^\downarrow={\tt K},
$$
the following holds:
 \bit{
\item[(a)] in ${\tt K}$ there exists a net of monomorphisms ${\mathcal N}$ such that for any object $X$ in ${\tt K}$ the morphism $\mu_{\rightlim{\mathcal N}_X}$ in the factorization \eqref{faktorizatsiya-v-kat-s-faktoriz} is a refinement $\rf_\varPhi^{\varGamma} X$ in $\varGamma$ by means of $\varPhi$:
\beq\label{mu_lim N_X=imp_Phi^Gamma-X}
\mu_{\rightlim{\mathcal N}_X}=\rf_\varPhi^{\varGamma} X,
\eeq

\item[(b)] for each morphism $\alpha:X\to Y$ in ${\tt K}$ and for any choice of refinements $\rf_\varPhi^{\varGamma}X$ and $\rf_\varPhi^{\varGamma}Y$ there is a unique morphism $\Rf_\varPhi^{\varGamma} \alpha:\Rf_\varPhi^{\varGamma}  X\to \Rf_\varPhi^{\varGamma}  Y$ in ${\tt K}$ such that the following diagram is commutative:
\beq\label{DIAGR:funktorialnost-imp_varPhi^Mono-v-kat-s-uzl-razl-1}
\xymatrix @R=2.pc @C=5.0pc % @M=14pt
{
X\ar[d]^{\alpha} & \Rf_\varPhi^{\varGamma}  X\ar[l]_{\rf_\varPhi^{\varGamma}X} \ar@{-->}[d]^{\Rf_\varPhi^{\varGamma} \alpha} \\
Y & \Rf_\varPhi^{\varGamma}  Y\ar[l]_{\rf_\varPhi^{\varGamma}Y} \\
}
\eeq

\item[(c)] if a category ${\tt K}$ is well-powered in $\varGamma$, then the refinement $\Rf_\varPhi^{\varGamma}$ can be defined as a functor.
}\eit
\etm

\paragraph{Existence of nets of epimorphisms and semi-regular envelopes.}

\btm\label{TH:sushestvovanie-seti-pri-faktorizatsii} Suppose a category $\tt K$ and classes of morphisms $\varOmega$ and $\varPhi$ in it satisfy the following conditions:
\bit{

\item[RE.1:] $\tt K$ is projectively complete,

\item[RE.2:] $\varOmega$ is monomorphically complementable: ${^\downarrow\varOmega}\circledcirc\varOmega={\tt K}$,

\item[RE.3:] ${\tt K}$ is co-well-powered in the class $\varOmega$,

\item[RE.4:] $\varPhi$ goes from\footnote{In the sense of definition on p.\pageref{DEF:goes-from}.} $\Ob({\tt K})$ and is a right ideal in ${\tt K}$:
$$
\Ob(\tt K)=\{\Dom\ph;\ \ph\in\varPhi\},\qquad \varPhi\circ\Mor({\tt K})\subseteq\varPhi.
$$
}\eit
Then
 \bit{
\item[(a)] there is a net of epimorphisms ${\mathcal N}$ in ${\tt K}$ such that for each object $X$ in ${\tt K}$ the morphism $\e_{\leftlim{\mathcal N}^X}$ in the factorization \eqref{faktorizatsiya-v-kat-s-faktoriz} is an envelope $\env_\varPhi^{\varOmega} X$ in $\varOmega$ with respect to $\varPhi$:
\beq\label{im_infty-lim F_X=env_M^Epi-X-1-*}
\e_{\leftlim{\mathcal N}^X}=\env_\varPhi^{\varOmega} X,
\eeq

\item[(b)] for each morphism $\alpha:X\to Y$ in ${\tt K}$ and for any choice of envelopes   $\env_\varPhi^{\varOmega}X$ and $\env_\varPhi^{\varOmega}Y$ there exists a unique morphism  $\Env_\varPhi^{\varOmega} \alpha:\Env_\varPhi^{\varOmega}  X\to \Env_\varPhi^{\varOmega}Y$ in ${\tt K}$, such that the following diagram is commutative:
\beq\label{DIAGR:funktorialnost-env_varPhi^Epi-v-kat-s-uzl-razl-1-*}
\xymatrix @R=2.pc @C=5.0pc % @M=14pt
{
X\ar[d]^{\alpha}\ar[r]^{\env_\varPhi^{\varOmega}  X} & \Env_\varPhi^{\varOmega}  X\ar@{-->}[d]^{\Env_\varPhi^{\varOmega} \alpha} \\
Y\ar[r]^{\env_\varPhi^{\varOmega}  Y} & \Env_\varPhi^{\varOmega}  Y \\
}
\eeq

\item[(c)] the envelope $\Env_\varPhi^{\varOmega}$ can be defined as a functor.
}\eit
\etm

\bit{

\item If the conditions RE.1-RE.4 are fulfilled, then we say that {\it the classes $\varOmega$ and $\varPhi$ define a semiregular envelope in $\tt K$}\label{DEF:polureg-obolochka}, or {\it the envelope  $\Env_\varPhi^\varOmega$ is semiregular}.
}\eit

\bpr
1. By RE.3, for each $X$ the category $\varOmega^X=\varOmega\cap\Epi^X$ is skeletally small. Let $S_X$ be its skeleton (which is a set). Using Theorem \ref{TH:o-lok-malosti-v-faktor-objektah}, we can choose a map $X\mapsto S_X$ that to each object $X$ assigns a skeleton $S_X$ in the category $\varOmega^X$. Let us fix this map.

For every object $X$ in ${\tt K}$ let us denote by $\varPhi^X$ the subclass of morphisms from $\varPhi$ which have  $X$ as a domain:
$$
\varPhi^X=\{\ph\in\varPhi:\ \Dom\ph=X\}
$$
(from RE.4 it follows that $\varPhi^X\ne\varnothing$). Denote by $2_{\varPhi^X}$ the class of finite subsets in $\varPhi^X$:
$$
\varPsi\in 2_{\varPhi^X}\qquad\Longleftrightarrow\qquad \varPsi\subseteq\varPhi^X\quad\&\quad \card\varPsi<\infty.
$$
To each object $X$ in ${\tt K}$ and to each morphism $\varPsi\in 2_{\varPhi^X}$ let us assign a morphism
$$
\overline{\varPsi}=\prod_{\psi\in\varPsi}\psi:X\to\prod_{\psi\in\varPsi}\Ran\psi,\qquad
$$
and morphisms $\mu_{\varPsi}\in\varGamma$ and $\e_{\varPsi}\in S_X$ such that
\beq\label{overline(varPsi)=mu_varPsi-circ-e_(varPsi)}
\overline{\varPsi}=\mu_{\varPsi}\circ\e_{\varPsi}.
\eeq
(since $S_X$ is a skeleton in $\varOmega^X$, such morphisms are unique). Put
$$
{\mathcal N}^X=\{\e_{\varPsi};\ \varPsi\in 2_{\varPhi^X}\}.
$$
Since ${\mathcal N}^X\subseteq S_X$, this is a set, and since the correspondence $X\mapsto S_X$ is a map, we obtain a map $X\mapsto {\mathcal N}^X$.

2. Let us check that the system $\mathcal N$ satisfies that axioms of a net of epimorphisms (p.\pageref{AX:set-Epi-a}). First, let us show that ${\mathcal N}^X$ is directed to the left with respect to the pre-order \eqref{DEF:le-in-F_X}, inherited from $\Epi^X$. For each two sets $\varPsi,\varPsi'\in 2_{\varPhi^X}$ let us consider the diagram
$$
\xymatrix @R=3.pc @C=5.0pc % @M=14pt
{
& X\ar[d]^{\overline{\varPsi\cup\varPsi'}}\ar@/_3ex/[ld]_{\overline{\varPsi}}
\ar@/^3ex/[rd]^{\overline{\varPsi'}} & \\
\Ran\overline{\varPsi}\ar@{=}[d] & \Ran\overline{\varPsi\cup\varPsi'}\ar[l]_{\pi}\ar[r]^{\pi'}\ar@{=}[d] & \Ran\overline{\varPsi'}\ar@{=}[d]\\
\prod\limits_{\psi\in\varPsi}\Ran\psi & \prod\limits_{\psi\in\varPsi\cup\varPsi'}\Ran\psi & \prod\limits_{\psi\in\varPsi'}\Ran\psi
}
$$
where $\pi$ and $\pi'$ are natural projections. Let us decompose the arrows going from $X$ by  \eqref{overline(varPsi)=mu_varPsi-circ-e_(varPsi)}:
\beq\label{DIAGR:sushestvovanie-seti-Epi-pri-uzlovom-razlozhenii}
\xymatrix @R=3.pc @C=5.0pc % @M=14pt
{
& X\ar[d]^{\e_{\varPsi\cup\varPsi'}}\ar@/_3ex/[ld]_{\e_{\varPsi}}
\ar@/^3ex/[rd]^{\e_{\varPsi'}} & \\
\Ran\e_{\varPsi}\ar[d]^{\mu_{\varPsi}} & \Ran\e_{\varPsi\cup\varPsi'}\ar[d]^{\mu_{\varPsi\cup\varPsi'}} & \Ran\e_{\varPsi'}\ar[d]^{\mu_{\varPsi'}}
\\
\Ran\overline{\varPsi} & \Ran\overline{\varPsi\cup\varPsi'}\ar[l]_{\pi}\ar[r]^{\pi'} & \Ran\overline{\varPsi'}
}
\eeq
The left side of the diagram we represent as a quadrangle:
$$
\xymatrix @R=3.pc @C=5.0pc % @M=14pt
{
& X\ar[d]^{\e_{\varPsi\cup\varPsi'}}\ar@/_3ex/[ld]_{\e_{\varPsi}}
\\
\Ran\e_{\varPsi}\ar[d]^{\mu_{\varPsi}} & \Ran\e_{\varPsi\cup\varPsi'}\ar@/^3ex/[dl]^{\pi\circ\mu_{\varPsi\cup\varPsi'}}
\\
\Ran\overline{\varPsi} &
}
$$
Here $\e_{\varPsi\cup\varPsi'}$ is an epimorphism, and $\mu_{\varPsi}$ a strong monomorphism, hence there exists a horizontal arrow to the left:
$$
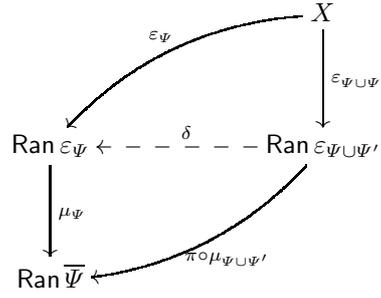
\xymatrix @R=3.pc @C=5.0pc % @M=14pt
{
& X\ar[d]^{\e_{\varPsi\cup\varPsi'}}\ar@/_3ex/[ld]_{\e_{\varPsi}}
\\
\Ran\e_{\varPsi}\ar[d]^{\mu_{\varPsi}} & \Ran\e_{\varPsi\cup\varPsi'}\ar@/^3ex/[dl]^{\pi\circ\mu_{\varPsi\cup\varPsi'}}
\ar@{-->}[l]_{\delta}
\\
\Ran\overline{\varPsi} &
}
$$
By the same reason in \eqref{DIAGR:sushestvovanie-seti-Epi-pri-uzlovom-razlozhenii} there is an arrow to the right, and we obtain a diagram
$$
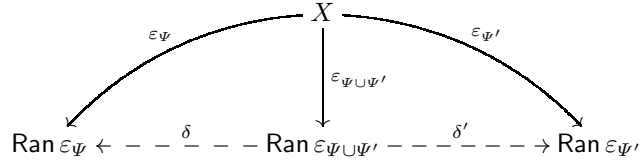
\xymatrix @R=3.pc @C=5.0pc % @M=14pt
{
& X\ar[d]^{\e_{\varPsi\cup\varPsi'}}\ar@/_3ex/[ld]_{\e_{\varPsi}}
\ar@/^3ex/[rd]^{\e_{\varPsi'}} & \\
\Ran\e_{\varPsi} & \Ran\e_{\varPsi\cup\varPsi'}
\ar@{-->}[l]_{\delta}\ar@{-->}[r]^{\delta'} & \Ran\e_{\varPsi'}
}
$$
It means that in the category $\varOmega^X$ the morphism $\e_{\varPsi\cup\varPsi'}$ majorates morphisms  $\e_{\varPsi}$ and $\e_{\varPsi'}$:
$$
\e_{\varPsi\cup\varPsi'}\to \e_{\varPsi}\quad\&\quad \e_{\varPsi\cup\varPsi'}\to\e_{\varPsi'}.
$$

The second condition in the definition of the net of epimorphisms is fulfilled automatically: since the category $\tt K$ is projectively complete, the system of binding morphisms $\Bind({\mathcal N}^X)$, defined in  \eqref{DEF:kategoriya-svyazyv-morpfizmov}, always has a projective limit.

Let us check the third condition. Let $\alpha:X\to Y$ be a morphism in $\varPsi\in 2_{\varPhi_Y}$. By RE.4, $\varPhi$ is a right ideal, hence for each $\psi\in\varPsi\subseteq\varPhi$ the composition $\psi\circ\alpha$ belongs to  $\varPhi$, and we can consider the set $\varPsi\circ\alpha\in 2_{\varPhi^X}$. We obtain the following diagram:
 $$
  
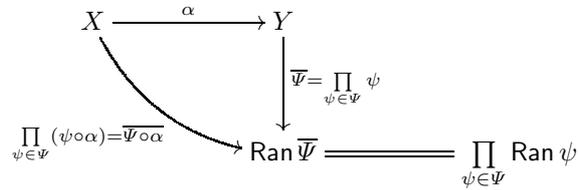
\xymatrix @R=2.5pc @C=4.0pc {
 X\ar[r]^{\alpha}\ar@/_3ex/[dr]_{\prod\limits_{\psi\in\varPsi}(\psi\circ\alpha)=\overline{\varPsi\circ\alpha}\quad} & Y\ar[d]^{\overline{\varPsi}=\prod\limits_{\psi\in\varPsi}\psi} & \\
  & \Ran\overline{\varPsi}\ar@{=}[r] & \prod\limits_{\psi\in\varPsi}\Ran\psi
 }
 $$
Let us represent morphisms coming to $\Ran\overline{\varPsi}$ as their factorizations \eqref{overline(varPsi)=mu_varPsi-circ-e_(varPsi)}:
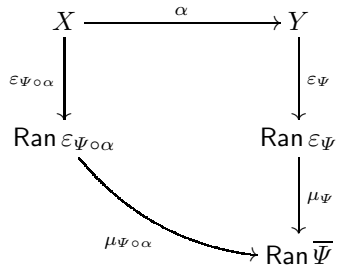
 $$
  \xymatrix @R=2.5pc @C=4.0pc {
 X\ar[r]^{\alpha}\ar[d]_{\e_{\varPsi\circ\alpha}} & Y\ar[d]^{\e_{\varPsi}}  \\
 \Ran\e_{\varPsi\circ\alpha}\ar@/_3ex/[dr]_{\mu_{\varPsi\circ\alpha}\quad} & \Ran\e_{\varPsi}\ar[d]^{\mu_{\varPsi}}  \\
  & \Ran\overline{\varPsi}
 }
 $$
This diagram can be represented as a quadrangle
 $$
  
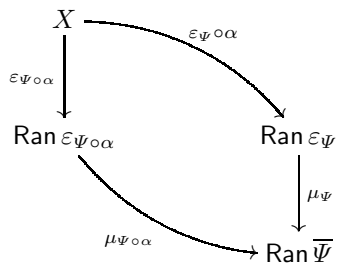
\xymatrix @R=2.5pc @C=4.0pc {
 X\ar@/^3ex/[dr]^{\e_{\varPsi}\circ\alpha}\ar[d]_{\e_{\varPsi\circ\alpha}} &   \\
 \Ran\e_{\varPsi\circ\alpha}\ar@/_3ex/[dr]_{\mu_{\varPsi\circ\alpha}\quad} & \Ran\e_{\varPsi}\ar[d]^{\mu_{\varPsi}}  \\
  & \Ran\overline{\varPsi}
 }
 $$
where $\e_{\varPsi\circ\alpha}\in\varOmega$ and $\mu_{\varPsi}\in\varGamma=\varOmega^\downarrow$. So there must exist a horizontal arrow to the right:
 $$
  \xymatrix @R=2.5pc @C=4.0pc {
 X\ar@/^3ex/[dr]^{\e_{\varPsi}\circ\alpha}\ar[d]_{\e_{\varPsi\circ\alpha}} &   \\
 \Ran\e_{\varPsi\circ\alpha}\ar@/_3ex/[dr]_{\mu_{\varPsi\circ\alpha}\quad}\ar@{-->}[r]^{\delta} & \Ran\e_{\varPsi}\ar[d]^{\mu_{\varPsi}}  \\
  & \Ran\overline{\varPsi}
 }
 $$
This will be the horizontal arrow that we need in \eqref{DIAGR:set}:
 $$
  \xymatrix @R=2.5pc @C=4.0pc {
 X\ar[r]^{\alpha}\ar[d]_{\e_{\varPsi\circ\alpha}} & Y\ar[d]^{\e_{\varPsi}}  \\
 \Ran\e_{\varPsi\circ\alpha}\ar@{-->}[r]^{\delta} & \Ran\e_{\varPsi}
 }
 $$

3. Note further that in the class $\varOmega$ the envelopes with respect to the classes $\varPhi$, $2_{\varPhi}$, $\{\e_{\varPsi};\ \varPsi\in 2_{\varPhi}\}$ and $\mathcal N$ coincide:
$$
\env_{\varPhi}^{\varOmega}X=\env_{2_{\varPhi}}^{\varOmega}X=\eqref{env_varPhi^varOmega_X=env_(e_ph;ph_in_varPhi)^varOmega_X}
=\env_{\{\e_{\varPsi};\ \varPsi\in 2_{\varPhi}\}}^{\varOmega}X=\env_{\mathcal N}^{\varOmega}X.
$$
After that the proof of Theorem \ref{TH:funktorialnost-pri-seti-Epi-i-dolonyaemosti} works.
\epr

\paragraph{Existence of nets of monomophisms and semi-regular refinements.}

The dual proposition for refinements look as follows:

\btm\label{TH:sushestvovanie-seti-Mono-pri-faktorizatsii}
Suppose a category $\tt K$ and classes of morphisms $\varGamma$ and $\varPhi$ satisfy the following conditions:
\bit{

\item[RR.1:] $\tt K$ is injectively complete,

\item[RR.2:] $\varGamma$ is epimorphically complementable in $\tt K$:
$\varGamma\circledcirc\varGamma^\downarrow={\tt K}$,

\item[RR.3:] ${\tt K}$ is well-powered in the class $\varGamma$,

\item[RR.4:] $\varPhi$ goes to\footnote{In the sence of definition on p.\pageref{DEF:goes-to}.} $\Ob(\tt K)$ and is a left ideal in ${\tt K}$:
$$
\Ob(\tt K)=\{\Ran\ph;\ \ph\in\varPhi\},\qquad \Mor({\tt K})\circ\varPhi\subseteq\varPhi.
$$
}\eit
Then
 \bit{
\item[(a)] there exists a net of monomorphisms ${\mathcal N}$ in ${\tt K}$ such that for each object $X$ in ${\tt K}$ the morphism $\mu_{\rightlim{\mathcal N}_X}$ in the factorization \eqref{faktorizatsiya-v-kat-s-faktoriz} is a refinement $\rf_\varPhi^{\varGamma} X$ in $\varGamma$ by means of $\varPhi$:
\beq\label{mu_lim N_X=imp_Phi^Gamma-X-*}
\mu_{\rightlim{\mathcal N}_X}=\rf_\varPhi^{\varGamma} X,
\eeq

\item[(b)] for each morphism $\alpha:X\to Y$ in ${\tt K}$ and for any choice of refinements $\rf_\varPhi^{\varGamma}X$ and $\rf_\varPhi^{\varGamma}Y$ there is a unique morphism $\Rf_\varPhi^{\varGamma} \alpha:\Rf_\varPhi^{\varGamma}  X\to \Rf_\varPhi^{\varGamma}  Y$ in ${\tt K}$ such that the following diagram is commutative:
\beq\label{DIAGR:funktorialnost-imp_varPhi^Mono-v-kat-s-uzl-razl-1-*}
\xymatrix @R=2.pc @C=5.0pc % @M=14pt
{
X\ar[d]^{\alpha} & \Rf_\varPhi^{\varGamma}  X\ar[l]_{\rf_\varPhi^{\varGamma}X} \ar@{-->}[d]^{\Rf_\varPhi^{\varGamma} \alpha} \\
Y & \Rf_\varPhi^{\varGamma}  Y\ar[l]_{\rf_\varPhi^{\varGamma}Y} \\
}
\eeq

\item[(c)] the refinement $\Rf_\varPhi^{\varGamma}$ can be defined as a functor.
}\eit \etm

\bit{

\item If RR.1-RR.4 are fulfilled, then we say that {\it the classes of morphisms $\varGamma$ and $\varPhi$ define a semi-regular refinement $\Rf_\varPhi^\varGamma$ in $\tt K$}\label{DEF:polureg-otpechatok}, or that {\it the refinement $\Rf_\varPhi^\varGamma$ is semiregular}.
}\eit

\paragraph{Pushing, regular envelope and complete objects.}

\bit{
\item Let us say that {\it a class of morphisms $\varOmega$ pushes a class of morphisms $\varPhi$}, if
\beq\label{Omega-podderzhivaet-Phi}
\forall\psi\in\Mor({\tt K})\qquad \forall \sigma\in\varOmega\qquad \big(\psi\circ\sigma\in\varPhi\quad\Longrightarrow\quad \psi\in\varPhi\big).
\eeq
}\eit

\brem\label{Omega-podtalkivaet-Mor(K,M)} Obviously, \eqref{Omega-podderzhivaet-Phi} holds, if $\varPhi$ is a class of morphisms with ranges in some class of objects $\tt M$ in $\tt K$:
$$
\varPhi=\{\ph\in\Mor(\tt K):\ \Ran\ph\in{\tt M}\},
$$
\erem

\blm\label{LM:kompozitsiya-rasshirenij} If a class $\varOmega$ pushes a class $\varPhi$, then the composition  $\sigma\circ\rho:X\to X''$ of any two extensions $\rho:X\to X'$ and $\sigma:X'\to X''$ (in $\varOmega$ with respect to $\varPhi$) is an extension (in $\varOmega$ with respect to $\varPhi$).
\elm
\bpr
This is seen from in following diagram:
$$
\xymatrix @R=2.pc @C=2.0pc % @M=14pt
{
X\ar@/_2ex/[dr]_{\ph}\ar[r]^{\rho} & X'\ar[r]^{\sigma}\ar@{-->}[d]^{\ph'} & X''\ar@/^2ex/@{-->}[dl]^{\ph''} \\
 & {\tt M} &
}
$$
Since $\rho$ is an extension, for any $\ph\in\varPhi$ there exists $\ph'$, and since $\varOmega$ pushes $\varPhi$, we have $\ph'\in\varPhi$. Then since $\sigma$ is an extension, there exists $\ph''$. Therewith every next arrow is uniquely defined by the previous one.
\epr

\bprop\label{PROP:harakterizatsiya-polnoty}
Suppose $\varOmega\subseteq\Epi$, then for each object $A\in\Ob({\tt K})$ the following conditions are equivalent:
 \bit{
\item[(i)] each extension $\sigma:A\to A'$ in $\varOmega$ with respect to $\varPhi$ is an isomorphism;

\item[(ii)] the local identity $1_A:A\to A$ is an envelope of $A$ in $\varOmega$ with respect to $\varPhi$;

\item[(iii)] there exists an envelope of $A$ in $\varOmega$ with respect to $\varPhi$, which is an isomorphism: $\env^\varOmega_\varPhi A\in\Iso$.
 }\eit
If in addition $\varOmega$ pushes $\varPhi$, then these conditions are equivalent to the following one:
 \bit{
\item[(iv)] $A$ is isomorphic to an envelope of some object $X\in\Ob({\tt K})$: $A\cong\Env^\varOmega_\varPhi X$.
 }\eit
\eprop

 \bit{
\item We will say that an object $A$ in $\tt K$ is {\it complete} in the class $\varOmega\subseteq\Epi$ with respect to the class $\varPhi$, if it satisfies the properties (i)-(iii) of this proposition.
 }\eit

\bpr
1. (i)$\Longrightarrow$(ii). Suppose that each extension $\sigma:A\to A'$ is an isomorphism. Then for the local identity $1_A:A\to A$ (which is also an extension) we have the diagram
$$
\xymatrix @R=2.pc @C=2.0pc % @M=14pt
{
& A \ar[dl]_{\sigma} \ar[dr]^{1_A}\\
A'\ar[rr]_{\sigma^{-1}} & & A
}
$$
which can be considered as the special case of \eqref{DEF:diagr-obolochka}, and this means that $1_A$ is an envelope.

2. The implication (ii)$\Longrightarrow$(iii) is obvious.

3. (iii)$\Longrightarrow$(i). Let $\rho:A\to E$ be an envelope, and at the same time an isomorphism. Then for any extension $\sigma:A\to A'$ we can take a morphism $\upsilon$ in \ref{DEF:diagr-obolochka} and we get  $\upsilon\circ\sigma=\rho\in\Iso$, hence $\sigma$ is a coretraction. On the other hand, $\sigma\in\varOmega\subseteq\Epi$, hence $\sigma\in\Iso$.

4. The implication (iii)$\Longrightarrow$(iv) is also obvious: if $\env^\varOmega_\varPhi A\in\Iso$, then $A\cong\Env^\varOmega_\varPhi A$.

5. Now it is sufficient to prove (iv)$\Longrightarrow$(i) in the case when $\varOmega$ pushes $\varPhi$. Suppose that  $A\cong\Env^\varOmega_\varPhi X$ for some $X\in\Ob({\tt K})$. Then $A$ can be considered as an envelope of $X$, i.e. there exists a morphism $\rho:X\to A$, which is an envelope. Take any extension $\sigma:A\to A'$ of $A$. By Lemma  \ref{LM:kompozitsiya-rasshirenij}, the composition $\sigma\circ\rho:X\to A'$ is an extension for $X$, so there is a morphism $\upsilon$ such that \eqref{DEF:diagr-obolochka} is commutative:
$$
\xymatrix @R=2.pc @C=2.0pc % @M=14pt
{
& X\ar[dl]_{\sigma\circ \rho}\ar[dr]^{\rho} & \\
A'\ar@{-->}[rr]_{\upsilon} & & A
}
$$
Now we have:
$$
\upsilon\circ\sigma\circ \rho=\rho=1_A\circ \underset{\scriptsize\begin{matrix}\text{\rotatebox{90}{$\owns$}}\\ \Epi\end{matrix}}{\rho}\qquad\Longrightarrow\qquad \upsilon\circ\sigma=1_A.
$$
In the last equality the morphism $\upsilon$ must be unique, since $\sigma$ is an epimorphism. We have that the extension $\sigma$ is subordinated to the extension $1_A$, and since this is true for each $\sigma$, the morphism  $\rho$ muts be an envelope for $A$.
\epr

Let us denote by $\tt L$ the class of complete object in $\tt K$ (in $\varOmega$ with respect to $\varPhi$). We consider $\tt L$ as a full subcategory in $\tt K$.

\bprop\label{PROP:obolochka-na-L}
In the conditions of Theorem \ref{TH:sushestvovanie-seti-pri-faktorizatsii} the functor of envelope  $(X,\alpha)\mapsto (E(X),E(\alpha))$ on the subcategory of complete objects ${\tt L}\subseteq {\tt K}$ is isomorphic to the identity functor:
\beq\label{E:K->L=funktor}
\forall A\in {\tt L}\qquad E(A)\cong A, \qquad \forall \alpha:\underset{\scriptsize\begin{matrix}\text{\rotatebox{90}{$\owns$}}\\ {\tt L}\end{matrix}}{A}\to \underset{\scriptsize\begin{matrix}\text{\rotatebox{90}{$\owns$}}\\ {\tt L}\end{matrix}}{A'}\qquad
E(\alpha)=e_{A'}\circ\alpha\circ e_A^{-1}.
\eeq
\eprop
\bpr
Take an arbitrary morphism $\alpha:A\to A'$ in $\tt L$, i.e. a morphism in $\tt K$, whose domain and range belong to $\tt L$. Then in diagram \eqref{DIAGR:funktorialnost-env-e-E} the horizontal arrows are isomorphisms, so
$$
\xymatrix @R=2.pc @C=5.0pc % @M=14pt
{
A\ar[d]^{\alpha} & A\ar[d]^{E(\alpha)}\ar[l]_{e_A^{-1}} \\
A'\ar[r]^{e_{A'}} & A' \\
}
$$
\epr

 \bit{
\item We say that {\it the classes $\varOmega$ and $\varPhi$ define a regular envelope in the category $\tt K$}, or  {\it the envelope $\Env_\varPhi^\varOmega$ is regular}, if in addition to the conditions RE.1-RE.4 of Theorem \ref{TH:sushestvovanie-seti-pri-faktorizatsii} the class $\varOmega$ pushes the class $\varPhi$.
}\eit

 \btm\label{TH:regulyarnaya-obolochka}
If the classes $\varOmega$ and $\varPhi$ define a regular envelope in $\tt K$, then $\Env_\varPhi^\varOmega$ can be defined as an idempotent functor.
 \etm
\bpr
Consider the functor of envelope $E$ built in Theorem \ref{TH:sushestvovanie-seti-pri-faktorizatsii}, and denote by ${\tt L}_0$ the class of all objects, which are values of the map $X\mapsto E(X)$:
\beq\label{DEF:L_0}
A\in{\tt L}_0\qquad\Longleftrightarrow\qquad \exists X\in\Ob({\tt K})\quad A=E(X).
\eeq
Define a system of isomorphisms
$$
\forall X\in\Ob({\tt K})\qquad \zeta_X=\begin{cases}1_X,& X\notin{\tt L}_0\\
e_X^{-1},& X\in{\tt L}_0 \end{cases}
$$
(this is a correct definition by Proposition \ref{PROP:harakterizatsiya-polnoty}). After that consider the maps $X\mapsto F(X)$, $X\mapsto f_X$, $\alpha\mapsto F(\alpha)$, defined by rules
 \begin{align*}
& \forall X\in\Ob({\tt K})\qquad F(X)=\begin{cases}E(X),& X\notin{\tt L}_0\\
X,& X\in{\tt L}_0 \end{cases}, \qquad f_X=\begin{cases}e_X,& X\notin{\tt L}_0\\
1_X,& X\in{\tt L}_0 \end{cases} \\
& \forall\alpha\in\Mor({\tt K})\qquad F(\alpha)=\zeta_{\Ran E(\alpha)}\circ E(\alpha)\circ\zeta_{\Dom E(\alpha)}^{-1}
 \end{align*}
The connection with the functor $E$ is reflected in the diagram
\beq\label{DIAGR:E(E(X))=E(X)}
\xymatrix @R=2.pc @C=5.0pc % @M=14pt
{
 &  & F(X)\ar[ddd]^{F(\alpha)} \\
X\ar[d]^{\alpha}\ar[r]^{e_X} \ar@/^3ex/[rru]^{f_X}  & E(X)\ar[d]^{E(\alpha)}\ar[ru]_{\zeta_X} & \\
Y\ar[r]^{e_Y}\ar@/_3ex/[rrd]_{f_Y} & E(Y)\ar[rd]^{\zeta_Y} & \\
 & & F(Y)
}
\eeq
For any $X$ the morphism $f_X:X\to F(X)$ is an envelope of $X$, since $f_X$ and $e_X$ are connected by the isomorphism $\zeta_X$. The map $(X,\alpha)\mapsto (F(X),F(\alpha))$ is a functor, since, first,
 \begin{multline*}
F(\beta\circ\alpha)=\zeta_{\Ran\beta}\circ E(\beta\circ\alpha)\circ\zeta_{\Dom\alpha}^{-1}=
\zeta_{\Ran\beta}\circ E(\beta)\circ E(\alpha)\circ\zeta_{\Dom\alpha}^{-1}=\\=
\zeta_{\Ran\beta}\circ E(\beta)\circ\zeta_{\Dom\beta}^{-1}\circ \zeta_{\Ran\alpha}\circ E(\alpha)\circ\zeta_{\Dom\alpha}^{-1}=F(\beta)\circ F(\alpha),
 \end{multline*}
and, second, for $X\notin{\tt L}_0$ diagram \eqref{DIAGR:E(E(X))=E(X)} has the form
$$
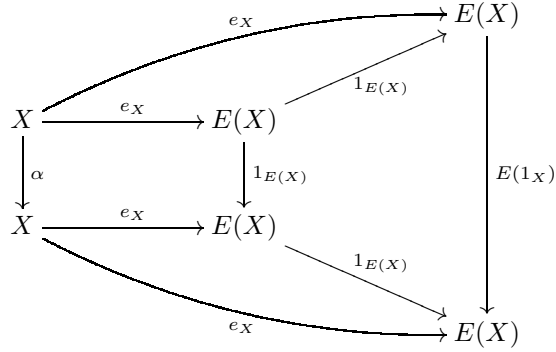
\xymatrix @R=2.pc @C=5.0pc % @M=14pt
{
 &  & E(X)\ar[ddd]^{E(1_X)} \\
X\ar[d]^{\alpha}\ar[r]^{e_X} \ar@/^3ex/[rru]^{e_X}  & E(X)\ar[d]^{1_{E(X)}}\ar[ru]_{1_{E(X)}} & \\
X\ar[r]^{e_X}\ar@/_3ex/[rrd]_{e_X} & E(X)\ar[rd]^{1_{E(X)}} & \\
 & & E(X)
}
$$
hence
$$
F(1_X)=\zeta_X\circ E(1_X)\circ\zeta_X^{-1}=1_{E(X)}^{-1}\circ 1_{E(X)}\circ 1_{E(X)}=1_{E(X)}=1_{F(X)},
$$
and for $X\in{\tt L}_0$ diagram \eqref{DIAGR:E(E(X))=E(X)} turns into
$$
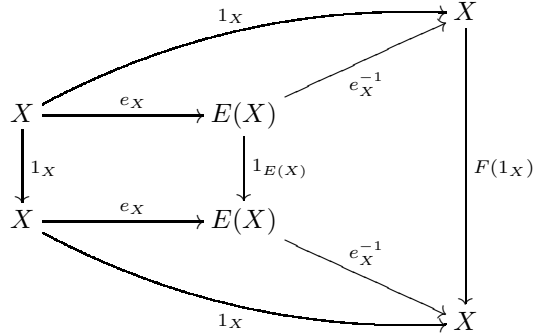
\xymatrix @R=2.pc @C=5.0pc % @M=14pt
{
 &  & X\ar[ddd]^{F(1_X)} \\
X\ar[d]^{1_X}\ar[r]^{e_X} \ar@/^3ex/[rru]^{1_X}  & E(X)\ar[d]^{1_{E(X)}}\ar[ru]_{e_X^{-1}} & \\
X\ar[r]^{e_X}\ar@/_3ex/[rrd]_{1_X} & E(X)\ar[rd]^{e_X^{-1}} & \\
 & & X
}
$$
If we replace $F(1_X)$ by $1_X$, then the perimeter will still be a commutative diagram. Since this arrow is unique we have
$$
F(1_X)=1_X=1_{F(X)}.
$$

The condition \eqref{e_(E(X))=1_(E(X))} holds for the functor $F$ by definition: since always
$F(X)\in{\tt L}_0$, we have $f_{F(X)}=1_{F(X)}$.
\epr

\btm[description of envelope in terms of complete objects]\label{TH:obolochki-i-polnye-obyekty}
Suppose that the classes $\varOmega$ and $\varPhi$ define a regular envelope in $\tt K$. Then a given morphism $\rho:X\to A$ is an envelope (in $\varOmega$ with respect to $\varPhi$) if and only if the following condition are fulfilled:
  \bit{
\item[(i)] $\rho:X\to A$ is an epimorphism,

\item[(ii)] $A$ is a complete object (in $\varOmega$ with respect to $\varPhi$),

\item[(iii)] for any complete object $B$ (in $\varOmega$ with respect to $\varPhi$) and for any morphism $\xi:X\to B$ there is a unique morphism $\xi':A\to B$ such that the following diagram is commutative:
 \beq\label{DIAGR:obolochki-i-polnye-obyekty}
\xymatrix @R=2.pc @C=2.0pc % @M=14pt
{
X\ar[rr]^{\rho}\ar[dr]_{\xi} & & A\ar@{-->}[dl]^{\xi'}\\
& B &
}
 \eeq
 }\eit
 \etm
\bpr Let $\rho:X\to A$ be an envelope. Then, first, this is an epimorphism, since $\varOmega\subseteq\Epi$. Second, by Proposition \ref{PROP:harakterizatsiya-polnoty}, $A\cong \Env^\varOmega_\varPhi X$ is a complete object. Third, if  $\xi:X\to B$ is a morphism into a complete object $B$, then we can consider diagram \eqref{DIAGR:funktorialnost-env_varPhi^Epi-v-kat-s-uzl-razl-1-*} which in this situation has the form
 $$
\xymatrix @R=2.pc @C=6.0pc % @M=14pt
{
X\ar[r]^{\rho=\env^\varOmega_\varPhi X}\ar[d]_{\xi} & A=\Env^\varOmega_\varPhi X\ar[d]^{\Env^\varOmega_\varPhi\xi}\\
B\ar[r]^{\env^\varOmega_\varPhi B} & \Env^\varOmega_\varPhi B
}
$$
Here $\env^\varOmega_\varPhi B$ must be an isomorphism, since $B$ is a complete object, and as a corollary, there exists a morphism
$$
\xi'=(\env^\varOmega_\varPhi B)^{-1}\circ \Env^\varOmega_\varPhi\xi.
$$
It is the dotted arrow in \eqref{DIAGR:obolochki-i-polnye-obyekty}.

On the contrary, suppose (i)--(iii) hold. In our circumstances Theorem
\ref{TH:sushestvovanie-seti-pri-faktorizatsii} works, so we can consider diagram
\eqref{DIAGR:funktorialnost-env_varPhi^Epi-v-kat-s-uzl-razl-1-*}:
 $$
\xymatrix @R=2.pc @C=5.0pc % @M=14pt
{
X\ar[d]_{\rho}\ar[r]^{\env^\varOmega_\varPhi X} & \Env^\varOmega_\varPhi X\ar[d]^{\Env^\varOmega_\varPhi\rho}\\
A\ar[r]^{\env^\varOmega_\varPhi {A}} & \Env^\varOmega_\varPhi A
}
$$
Here $\env^\varOmega_\varPhi {A}$ is an isomorphism (since $A$ is a complete object). Hence if we take  $\zeta=\env^\varOmega_\varPhi {A}^{-1}\circ \Env^\varOmega_\varPhi(\rho)$, we obtain a commutative diagram
 $$
\xymatrix @R=2.pc @C=5.0pc % @M=14pt
{
X\ar[d]_{\rho}\ar[r]^{\env^\varOmega_\varPhi X} & \Env^\varOmega_\varPhi X\ar@{-->}[dl]^{\zeta}\\
A &
}
$$
On the other hand, by Proposition \ref{PROP:harakterizatsiya-polnoty}, $\Env^\varOmega_\varPhi X$ is a complete object, so by (iii), there exists a morphism $\eta$ such that the following diagram is commutative:
 $$
\xymatrix @R=2.pc @C=5.0pc % @M=14pt
{
X\ar[d]_{\rho}\ar[r]^{\env^\varOmega_\varPhi X} & \Env^\varOmega_\varPhi X\\
A\ar@{-->}[ur]_{\eta} &
}
$$
In these diagrams both $\rho$ and $\env^\varOmega_\varPhi X$ are epimorphisms, so $\zeta$ and $\eta$ are mutially inverse morphisms. Thus, $\rho=\zeta\circ\env^\varOmega_\varPhi X$, where $\zeta\in\Iso$. By  \eqref{Iso-circ-varOmega-subseteq-varOmega}, we have that $\rho\in\varOmega$, and thus it is an envelope.
\epr

\paragraph{Pulling, regular refinement and saturated objects.}

\bit{
\item Let us say that {\it a class of morphisms $\varGamma$ pulls a class of morphisms $\varPhi$}, if
\beq\label{Omega-podtyagivaet-Phi}
\forall\psi\in\Mor({\tt K})\qquad \forall \sigma\in\varGamma\qquad \big(\sigma\circ\psi\in\varPhi\quad\Longrightarrow\quad \psi\in\varPhi\big).
\eeq
}\eit

\brem Obviously, \eqref{Omega-podtyagivaet-Phi} holds if $\varPhi$ is a class of morphisms with domains in a subclass $\tt M$ of objects in $\tt K$:
$$
\varPhi=\{\ph\in\Mor(\tt K):\ \Dom\ph\in{\tt M}\},
$$
\erem

\blm\label{LM:kompozitsiya-oblastej-vliyaniya} If $\varGamma$ pulls $\varPhi$, then the composition  $\sigma\circ\rho:X\gets X''$ of any two enrichments $\sigma:X\gets X'$ and $\rho:X'\gets X''$ (in $\varGamma$ by means of $\varPhi$) is an enrichment (in $\varGamma$ by means of $\varPhi$).
\elm
\bpr
This is seen from the diagram
$$
\xymatrix @R=2.pc @C=2.0pc % @M=14pt
{
X & X'\ar[l]_{\sigma}  & X''\ar[l]_{\rho}  \\
 & {\tt M} \ar@/^2ex/[ul]^{\ph}\ar@{-->}[u]_{\ph'} \ar@/_2ex/@{-->}[ur]_{\ph''}&
}
$$
\epr

\bprop\label{PROP:harakterizatsiya-nasyshennosti}
Suppose $\varGamma\subseteq\Mono$, then for an object $A\in\Ob({\tt K})$ the following conditions are equivalent:
 \bit{
\item[(i)] every enrichment $\sigma:A\gets A'$ in $\varGamma$ by means of $\varPhi$ is an isomorphism;

\item[(ii)] the local identity $1_A:A\to A$ is a refinement of $A$ in $\varGamma$ by means of $\varPhi$;

\item[(iii)] there exists a refinement of $A$ in $\varGamma$ by means of $\varPhi$, which is an isomorphism: $\rf^\varGamma_\varPhi A\in\Iso$.
 }\eit
If in addition $\varGamma$ pulls $\varPhi$, then these conditions are equivalent to the following one:
 \bit{
\item[(iv)] $A$ is isomorphic to a refinement of some object $X\in\Ob({\tt K})$: $A\cong\Rf^\varGamma_\varPhi X$.
 }\eit
\eprop

 \bit{
\item We say that an object $A$ in $\tt K$ is {\it saturated} in the class $\varGamma\subseteq\Mono$ by means of the class $\varPhi$, if it satisfies the conditions (i)-(iii) of this proposition.
 }\eit

Denote by $\tt L$ the class of all saturated objects in $\tt K$ (in $\varGamma\subseteq\Mono$ by means of $\varPhi$). We consider $\tt L$ as a full subcategory in $\tt K$.

\bprop\label{PROP:otpechatok-na-L}
In the conditions of Theorem \ref{TH:sushestvovanie-seti-Mono-pri-faktorizatsii} the functor of refinement $(X,\alpha)\mapsto (I(X),I(\alpha))$ on a subcategory of saturated objects ${\tt L}\subseteq {\tt K}$ is isomorphic to the identity functor:
\beq\label{I:K->L=funktor}
\forall A\in {\tt L}\qquad I(A)\cong A, \qquad \forall \alpha:\underset{\scriptsize\begin{matrix}\text{\rotatebox{90}{$\owns$}}\\ {\tt L}\end{matrix}}{A}\gets \underset{\scriptsize\begin{matrix}\text{\rotatebox{90}{$\owns$}}\\ {\tt L}\end{matrix}}{A'}\qquad
E(\alpha)=i_A^{-1}\circ\alpha\circ i_{A'}.
\eeq
\eprop

 \bit{
\item We say that {\it the classes $\varGamma$ and $\varPhi$ define a regular refinement in $\tt K$}, or
{\it the refinement $\Rf_\varPhi^\varGamma$ is regular}, if in addition to the conditions RR.1-RR.5 of Theorem
\ref{TH:sushestvovanie-seti-Mono-pri-faktorizatsii} the class $\varGamma$ pulls the class $\varPhi$.
}\eit

 \btm\label{TH:regulyarnayi-otpechatok}
If the classes $\varGamma$ and $\varPhi$ define a regular refinement in $\tt K$, then $\Rf_\varPhi^\varGamma$ can be defined as an idempotent functor.
 \etm

\btm[description of refinement in terms of saturated objects]\label{TH:otpechatki-i-nasyshennye-obyekty}
Suppose the classes $\varGamma$ and $\varPhi$ define a regular refinement in $\tt K$. Then a given morphism $\rho:X\gets A$ is a refinement (in $\varGamma$ by means of $\varPhi$) if and only if the following conditions hold:
 \bit{

 \item[(i)] $\rho:X\gets A$ is a monomorphism,

 \item[(ii)] $A$ is a saturated object (in $\varGamma$ by means of $\varPhi$),

 \item[(iii)] for any saturated object $B$ (in $\varGamma$ by means of $\varPhi$) and for any morphism $\xi:X\gets B$ there is a unique morphism $\xi':A\gets B$ such that the following diagram is commutative:
 \beq\label{DIAGR:otpechatki-i-nasyshennye-obyekty}
\xymatrix @R=2.pc @C=2.0pc % @M=14pt
{
X & & A\ar[ll]_{\rho}\\
& B\ar[ul]^{\xi} \ar@{-->}[ur]_{\xi'} &
}
 \eeq
 }\eit
\etm

\paragraph{Functoriality on epimorphisms and monomorphisms.}

Let us denote by ${\tt K}^{\Epi}$ the subcategory in $\tt K$ with the same class of objects as in $\tt K$, but with epimorphisms from $\tt K$ as morphisms:
$$
\Ob({\tt K}^{\Epi})=\Ob({\tt K}),\qquad \Mor({\tt K}^{\Epi})=\Epi({\tt K}).
$$

\btm\label{TH:funktorialnost-v-K^Epi}
Let ${\tt K}$ be a category with products (over arbitrary index sets), and classes of morphisms $\varOmega$ and $\varPhi$ in ${\tt K}$ satisfy the following conditions:
\bit{

\item[---] $\varOmega$ is monomorphically complementable in ${\tt K}$,

\item[---] ${\tt K}$ is co-well-powered in the class $\varOmega$,

\item[---] $\varPhi$ goes from\footnote{In the sense of definition on p.\pageref{DEF:goes-from}.} ${\tt K}$,

\item[---] $\varPhi\circ\varOmega\subseteq\varPhi$.

}\eit
Then
\bit{
\item[(a)] each object $X$ in ${\tt K}$ has an envelope $\Env_{\varPhi}^{\varOmega}X$ in $\varOmega$ with respect to $\varPhi$,

\item[(b)] for each epimorphism $\pi:X\to Y$ there is a unique epimorphism  $\Env_{\varPhi}^{\varOmega}\pi:\Env_{\varPhi}^{\varOmega}X\to \Env_{\varPhi}^{\varOmega}Y$ such that the following diagram is commutative:
\beq\label{DIAGR:funktorialnost-v-K^Epi}
\xymatrix @R=2.pc @C=5.0pc % @M=14pt
{
X\ar[d]^{\pi}\ar[r]^{\env_{\varPhi}^{\varOmega} X} & \Env_{\varPhi}^{\varOmega} X\ar@{-->}[d]^{\Env_{\varPhi}^{\varOmega}\pi} \\
Y\ar[r]^{\env_{\varPhi}^{\varOmega} Y} & \Env_{\varPhi}^{\varOmega} Y \\
}
\eeq

\item[(c)] the envelope $\Env_{\varPhi}^\varOmega$ can be defined as a functor from ${\tt K}^{\Epi}$ into ${\tt K}^{\Epi}$.
}\eit
\etm

We'll need the following

\blm\label{LM:Env_varPhi-circ-e^Omega-X=Env_varPhi^Omega-Y} If ${\tt K}$ is a category with products (over arbitrary index sets), co-well-powered in the class $\varOmega$, and $\varOmega$ is monomorphically complementable in ${\tt K}$, then for any class of morphisms $\varPhi$ and for any epimorphism $\pi:X\to Y$ the following formula holds:
    \beq\label{Env_varPhi-circ-e^Omega-X=Env_varPhi^Omega-Y}
    \Env_{\varPhi\circ\pi}^\varOmega X=\Env_{\varPhi}^\varOmega Y
    \eeq
(the envelope of $X$ in ${\varOmega}$ with respect to $\varPhi\circ\pi=\{\ph\circ\pi;\ \ph\in\varPhi\}$ coincides with the envelope of $Y$ in ${\varOmega}$ with respect to $\varPhi$).
\elm
\bpr Note first that the existence of envelopes in \eqref{Env_varPhi-circ-e^Omega-X=Env_varPhi^Omega-Y} is guaranteed by property $5^\circ$ on p.\pageref{5^0:obolochka-otn-klassa-morphizmov}. In addition, by
$5^\circ$ on p.\pageref{PROP:deistvie-epimorfizma-na-Env}, there exists a morphism $\upsilon$ such that
\eqref{deistvie-epimorfizma-na-Env} is commutative:
$$
\xymatrix %@R=2.5pc @C=4.0pc
{
 & X\ar[ld]_{\env_{\varPhi}^{\varOmega} Y\circ\pi}\ar[rd]^{\env_{\varPhi\circ\pi}^{\varOmega} X} & \\
 \Env_{\varPhi}^{\varOmega} Y\ar@{-->}[rr]^{\upsilon} &   & \Env_{\varPhi\circ\pi}^{\varOmega} X
}
$$
Let us show that there is a reverse morphism. Consider the envelope $\env_{\varPhi}^{\varOmega} Y: Y\to \Env_{\varPhi}^{\varOmega} Y$ and represet it as an envelope with respect to a set of morphisms $M$, like in the proof of property $5^\circ$ on p.\pageref{5^0:obolochka-otn-klassa-morphizmov}. Then, like in the proof of $3^\circ$ on p.\pageref{3^0:obolochka-otn-mnozhestva-morphizmov}, let us replace $M$ by a unique morphism $\psi=\prod_{\chi\in M}\chi$. By property $1^\circ$ on p.\pageref{1^0:obolochka-otn-1-morphizma}, the envelope with respect to $\psi$ will be described as epimorphism $\e_\psi$ in the factorization $\psi$:
$$
\env_{\varPhi}^{\varOmega} Y=\env_{M}^{\varOmega} Y=\env_{\psi}^{\varOmega} Y=\e_\psi.
$$
We obtain a diagram
$$
\xymatrix %@R=2.5pc @C=4.0pc
{
X\ar@/^2ex/[rd]^(.7){\env_{\varPhi}^{\varOmega} Y\circ\pi=\e_\psi\circ\pi}
\ar@/^2ex/[rrr]^{\env_{\varPhi\circ\pi}^{\varOmega} X}\ar[dd]_{\pi} & &  & \Env_{\varPhi\circ\pi}^{\varOmega} X\ar[dd]^{(\psi\circ\pi)'}\ar@{-->}@/_2ex/[dl]_(.7){\delta} \\
 & \Env_{\varPhi}^{\varOmega} Y\ar@{=}[r] & \Ran\e_\psi\ar@/_2ex/[dr]_(.3){\mu_\psi} & \\
Y\ar@/_2ex/[rrr]_{\psi}\ar@/_2ex/[ru]_(.7){\env_{\varPhi}^{\varOmega} Y=\e_\psi}  & & & B\\
}
$$
where $(\psi\circ\pi)'$ is an extenstion of $\psi\circ\pi\in\varPhi\circ\pi$ along the envelope $\env_{\varPhi\circ\pi}^{\varOmega} X$. Here the existence of morphism $\delta$ follows froom the fact that  $\env_{\varPhi\circ\pi}^{\varOmega} X\in\varOmega$, and $\mu_\psi\in{^\downarrow\varOmega}$. We have now the diagram
$$
\xymatrix %@R=2.5pc @C=4.0pc
{
 & X\ar[ld]_{\env_{\varPhi}^{\varOmega} Y\circ\pi}\ar[rd]^{\env_{\varPhi\circ\pi}^{\varOmega} X} & \\
 \Env_{\varPhi}^{\varOmega} Y &   & \Env_{\varPhi\circ\pi}^{\varOmega} X\ar@{-->}[ll]_{\delta}
}
$$

It remains to verify that $\upsilon$ and $\delta$ are mutually reverse. First,
$$
\delta\circ\upsilon\circ\underbrace{\env_{\varPhi}^{\varOmega} Y\circ\pi}_{\scriptsize\begin{matrix}\text{\rotatebox{90}{$\owns$}}\\ \Epi\end{matrix}}=\delta\circ\env_{\varPhi\circ\pi}^{\varOmega} X=\env_{\varPhi}^{\varOmega} Y\circ\pi=1_{\Env_{\varPhi}^{\varOmega} Y}\circ\underbrace{\env_{\varPhi}^{\varOmega} Y\circ\pi}_{\scriptsize\begin{matrix}\text{\rotatebox{90}{$\owns$}}\\ \Epi\end{matrix}}
\quad\Longrightarrow\quad
\delta\circ\upsilon=1_{\Env_{\varPhi}^{\varOmega} Y}.
$$
And, second,
$$
\upsilon\circ\delta\circ\underbrace{\env_{\varPhi\circ\pi}^{\varOmega} X}_{\scriptsize\begin{matrix}\text{\rotatebox{90}{$\owns$}}\\ \Epi\end{matrix}}=
\upsilon\circ\env_{\varPhi}^{\varOmega} Y\circ\pi=\env_{\varPhi\circ\pi}^{\varOmega} X=1_{\Env_{\varPhi\circ\pi}^{\varOmega} X}\circ
\underbrace{\env_{\varPhi\circ\pi}^{\varOmega} X}_{\scriptsize\begin{matrix}\text{\rotatebox{90}{$\owns$}}\\ \Epi\end{matrix}}
\quad\Longrightarrow\quad
\upsilon\circ\delta=1_{\Env_{\varPhi\circ\pi}^{\varOmega} X}.
$$
\epr

\bpr[Proof of Theorem \ref{TH:funktorialnost-v-K^Epi}]
Proposition (a) follows from property $5^\circ$ on p.\pageref{5^0:obolochka-otn-klassa-morphizmov}. Let us prove (b). By Lemma \ref{LM:Env_varPhi-circ-e^Omega-X=Env_varPhi^Omega-Y}, $\Env_{\varPhi}^{\varOmega} Y=\Env_{\varPhi\circ\pi}^{\varOmega} X$, and by property $3^\circ$ on p.\pageref{LM:suzhenie-klassa-morfizmov}, when we pass to a narrower class of morphisms $\varPhi\circ\pi\subseteq\varPhi$ a dotted arrow arises in the upper triangle of the diagram
$$
\xymatrix @R=2.pc @C=5.0pc % @M=14pt
{
X\ar[dr]_{\env_{\varPhi}^{\varOmega} Y\circ\pi}\ar[dd]_{\pi}\ar[r]^{\env_{\varPhi}^{\varOmega} X} &  \Env_{\varPhi}^{\varOmega} X\ar@{-->}[d]^{\Env_{\varPhi}^{\varOmega}\pi} \\
     &             \Env_{\varPhi\circ\pi}^{\varOmega} X \ar@{=}[d] \\
Y\ar[r]^{\env_{\varPhi}^{\varOmega} Y} & \Env_{\varPhi}^{\varOmega} Y \\
}
$$
It will be the dotted arrow in \eqref{DIAGR:funktorialnost-v-K^Epi}, but we need to verify that it is an epimorphism (so that this will be a morphism in ${\tt K}^{\Epi}$). This follows from property $3^\circ$ on p.\pageref{PROP:e-ph=epi-=>-e=epi}: since $\Env_{\varPhi}^{\varOmega}\pi\circ\env_{\varPhi}^{\varOmega} X=\env_{\varPhi}^{\varOmega} Y\circ\pi\in\Epi$, we have $\Env_{\varPhi}^{\varOmega}\pi\in\Epi$.

When (a) and (b) are proven, (c) becomes its corollary due to Theorem \ref{TH:o-lok-malosti-v-podobjektah}: $\tt K$ is co-well-powered in $\varOmega$, hence we can choose a map $X\mapsto S_X$, which assigns to each object a skeleton $S_X$ in the category $\varOmega\cap\Epi^X$. After that it becomes possible to choose a map $X\mapsto\env_\varPhi^\varOmega X$, and for any epimorphism $\pi:X\to Y$ the arrow $\Env_\varPhi^\varOmega\pi$ automeatically appears from diagram \eqref{DIAGR:funktorialnost-v-K^Epi}.
\epr

The dual results for refinements look as follows.

Denote by ${\tt K}^{\Mono}$ the subcategory in $\tt K$ with the same class of objects as in $\tt K$, but with monomorphisms from $\tt K$ as morphisms:
$$
\Ob({\tt K}^{\Mono})=\Ob({\tt K}),\qquad \Mor({\tt K}^{\Mono})=\Mono({\tt K}).
$$

\btm\label{TH:funktorialnost-v-K^Mono}
Let ${\tt K}$ be a category with coproducts (over arbitraty index sets), and classes of morphisms $\varGamma$ and $\varPhi$ in ${\tt K}$ satisfy the following conditions:
\bit{

\item[---] $\varGamma$ epimorphically complementable in ${\tt K}$,

\item[---] ${\tt K}$ well-powered in the class $\varGamma$,

\item[---] $\varPhi$ goes to\footnote{In the sense of definition on p.\pageref{DEF:goes-to}} ${\tt K}$,

\item[---] $\varGamma\circ\varPhi\subseteq\varPhi$.

}\eit
Then
\bit{
\item[(a)] each object $X$ in ${\tt K}$ has a refinement $\Rf_{\varPhi}^{\varGamma}X$ in $\varGamma$ by means of $\varPhi$,

\item[(b)] for each monomorphism $\pi:X\to Y$ there exists a unique monomorphism  $\Rf_{\varPhi}^{\varGamma}\pi:\Rf_{\varPhi}^{\varGamma}X\to \Rf_{\varPhi}^{\varGamma}Y$, such that the following diagram is commutative:
\beq\label{DIAGR:funktorialnost-v-K^Mono}
\xymatrix @R=2.pc @C=5.0pc % @M=14pt
{
X\ar[d]^{\pi} & \Rf_{\varPhi}^{\varGamma} X\ar[l]_{\Rf_{\varPhi}^{\varGamma} X}\ar@{-->}[d]^{\Rf_{\varPhi}^{\varGamma}\pi} \\
Y & \Rf_{\varPhi}^{\varGamma} Y\ar[l]_{\rf_{\varPhi}^{\varGamma} Y} \\
}
\eeq

\item[(c)] the refinement $\Rf_{\varPhi}^\varGamma$ can be defined as a functor from ${\tt K}^{\Mono}$ into ${\tt K}^{\Mono}$.
}\eit
\etm

The following lemma is used in the proof:

\blm\label{LM:Imp_varPhi-circ-i^Gama-X=Imp_varPhi^Gamma-Y} If ${\tt K}$ is a category with coproducts, well-powered in the class $\varGamma$, and $\varGamma$ is epimorphically complemented in ${\tt K}$, then for each class of morphisms $\varPhi$ and for each monomorphism $\pi:X\gets Y$ the following formula holds:
    \beq\label{Imp_varPhi-circ-e^Omega-X=Imp_varPhi^Omega-Y}
    \Rf_{\varPhi}^\varGamma X=\Rf_{\pi\circ\varPhi}^\varGamma Y
    \eeq
(the refinement of $X$ in ${\varGamma}$ by means of $\varPhi$ coincides with the refinement of $Y$ in ${\varGamma}$  by means of $\pi\circ\varPhi=\{\pi\circ\ph;\ \ph\in\varPhi\}$).
\elm

\paragraph{The case of $\Env_{\tt L}^{\tt L}$ and $\Rf_{\tt L}^{\tt L}$.}

Theorem \ref{TH:sushestvovanie-seti-pri-faktorizatsii} has important corollaries in the case when the classes of test morphisms and realizing morphisms coincide $\varPhi=\varOmega$, and are the class of all morphisms with ranges in a given class of objects $\tt L$ (this is the special case of the situation described on p.\pageref{DEF:obolochka-otn-klassa}, where ${\tt L}={\tt M}$).

\btm\label{TH:env^L-funktor} Suppose a category ${\tt K}$ and a class of objects ${\tt L}$ have the following properties:
 \bit{
\item[(i)] ${\tt K}$ is projectively complete,

\item[(ii)] ${\tt K}$ has nodal decomposition,

\item[(iii)] ${\tt K}$ is co-well-powered in the class $\Epi$,

\item[(iv)] $\Mor({\tt K},{\tt L})$ goes from $\tt K$:
$\forall X\in\Ob({\tt K})\quad \exists \ph\in\Mor({\tt K})\quad \Dom\ph=X\quad\&\quad\Ran\ph\in{\tt L}$,

\item[(v)] ${\tt L}$ differs morphisms on the outside,

\item[(vi)] ${\tt L}$ is closed with respect to passage to projective limits,

\item[(vii)] ${\tt L}$ is closed with respect to passage from the range of a morphism to its nodal image: if $\Ran\alpha\in{\tt L}$, then $\Im_\infty\alpha\in{\tt L}$.
 }\eit
Then
 \bit{

\item[(a)] each object $X$ has an envelope $\env_{\tt L}^{\tt L}X$ in the class of objects ${\tt L}$ (with respect to the same class ${\tt L}$)

\item[(b)] each envelope $\env_{\tt L}^{\tt L}X$ is a bimorphism,

\item[(c)] the envelope $\Env_{\tt L}^{\tt L}$ can be defined as a functor.
}\eit \etm
 \bpr
Conditions (i)-(v) mean that the classes $\Epi$ and $\varPhi=\Mor({\tt K},{\tt L})$ satisfy the promises of Theorem
\ref{TH:sushestvovanie-seti-pri-faktorizatsii}, i.e. define a semiregular envelope $\Env_\varPhi^{\Epi}=\Env_{\tt L}^{\Epi}$. In the proof of Theorem \ref{TH:sushestvovanie-seti-pri-faktorizatsii} this envelope is constructed by passing from the spaces $\Ran\ph\in{\tt L}$ ($\ph\in\varPhi$) to their projective limits, which belong to $\tt K$ by (vi), and then to the nodal image, which belongs to $\tt L$ by (vii). Therefore, $\Env_{\tt L}^{\Epi}\in{\tt L}$, hence by property $1^\circ$ on p.\pageref{LM:suzhenie-verh-klassa-morfizmov}
$\Env_{\tt L}^{\Epi}=\Env_{\tt L}^{\Epi({\tt K},{\tt L})}$. By construction, the class $\varPhi$ is a right ideal, and by (v), $\varPhi$ differs morphisms on the outside. So by Theorem \ref{TH:Phi-razdel-moprfizmy-*},
$\Env_{\tt L}^{\Epi({\tt K},{\tt L})}=\Env_{\tt L}^{\Bim({\tt K},{\tt L})}$. Further, by the same theorem
\ref{TH:Phi-razdel-moprfizmy-*} the envelope in the class $\Bim({\tt K},{\tt L})=\Mor({\tt K},{\tt L})\cap\Bim$
exists if and only if there exist the envelope in the class $\Mor({\tt K},{\tt L})$, and these envelopes coincide:
$\Env_{\tt L}^{\Bim({\tt K},{\tt L})}=\Env_{\tt L}^{\Mor({\tt K},{\tt L})}$. We obtain the following logical chain:
$$
\Env_{\tt L}^{\Epi}=\Env_{\tt L}^{\Epi({\tt K},{\tt L})}=
\Env_{\tt L}^{\Bim({\tt K},{\tt L})}=\Env_{\tt L}^{\Mor({\tt K},{\tt L})}=
\Env^{\tt L}_{\tt L}
$$
This proves (a) and (c), and incidentally (b).
 \epr

The dual result is as follows:

\btm\label{TH:imp^L-funktor} Suppose a category ${\tt K}$ and a class of objects ${\tt L}$ satisfy the following conditions:
 \bit{
\item[(i)] ${\tt K}$ is injectively complete,

\item[(ii)] ${\tt K}$ has nodal decomposition,

\item[(iii)] ${\tt K}$ is well-powered in the class $\Mono$,

\item[(iv)] $\Mor({\tt L},{\tt K})$ goes to $\tt K$:
$\forall X\in\Ob({\tt K})\quad \exists \ph\in\Mor({\tt K})\quad \Dom\ph\in{\tt L}\quad\&\quad\Ran\ph=X$,

\item[(v)] ${\tt L}$ differs morphisms on the inside,

\item[(vi)] ${\tt L}$ is closed with respect to the operation of taking injective limits,

\item[(vii)] ${\tt L}$ is closed with respect to passage from domain of a morphism to its nodal coimage: if $\Dom\alpha\in{\tt L}$, then $\Coim_\infty\alpha\in{\tt L}$.
 }\eit
Then
 \bit{

\item[(a)] each object $X$ has refinement $\rf_{\tt L}^{\tt L}X$ in the class ${\tt L}$ (by means of the same class ${\tt L}$)

\item[(b)] each refinement $\rf_{\tt L}^{\tt L}X$ is a bimorphism,

\item[(c)] the refinement $\Rf_{\tt L}^{\tt L}$ can be defined as a functor.
}\eit \etm

\subsection{Envelopes in monoidal categories}

\paragraph{Envelopes coherent with tensor product.}

Let $\tt K$ be a monoidal category \cite{MacLane} with the tensor product $\otimes$ and the unit object $I$.
 \bit{\label{DEF:obolochka-soglasovana-s-tenz-proizv}
\item Let us say that {\it the envelope $\Env^\varOmega_\varPhi$ is coherent with the tensor product $\otimes$} in $\tt K$, if the following conditions are fulfilled:

 \bit{

\item[T.1]\label{DEF:T.1} The tensor product $\rho\otimes\sigma:X\otimes Y\to X'\otimes Y'$ of any two extensions  $\rho:X\to X'$ and $\sigma:Y\to Y'$ (in $\varOmega$ with respect to $\varPhi$) is an extension (in $\varOmega$ with respect to $\varPhi$).

\item[T.2]\label{DEF:T.2} The local identity $1_I:I\to I$ of the unit object $I$ is the envelope (in $\varOmega$ with respect to $\varPhi$):
\beq\label{E(I)=I}
\env_\varPhi^\varOmega I=1_I
\eeq
 }\eit
 }\eit

Everywhere below in this section we consider the case when the classes $\varOmega$ and $\varPhi$ define a regular envelope in $\tt K$. By Theorem \ref{TH:regulyarnaya-obolochka} this means that $\Env_\varPhi^\varOmega$ can be defined as an idempotent functor. We denote it by $E:{\tt K}\to{\tt K}$, and the natural transformation of the identity functor into $E$ we denote by $e$:
$$
E(X):=\Env_\varPhi^\varOmega X, \qquad E(\ph):=\Env_\varPhi^\varOmega\ph,\qquad e_X:=\env_\varPhi^\varOmega X.
$$
The class of all complete objects in $\tt K$ (in $\varOmega$ with respect to $\varPhi$) we denote by $\tt L$.

\blm\label{PROP:E(e_X-otimes-e_Y)-in-Iso} Let $\Env_\varPhi^\varOmega$ be a regular envelope coherent with tensor product in ${\tt K}$. Then
 \bit{

\item[(i)] for any objects $A\in{\tt L}$ and $X\in\Ob({\tt K})$ the envelope $E(1_A\otimes e_X)$ of morphism $1_A\otimes e_X:A\otimes X\to A\otimes E(X)$ is an isomorphism (in $\tt K$ and in $\tt L$):
 \beq\label{E(1_A-otimes-e_X)-in-Iso}
 E(1_A\otimes e_X)\in\Iso.
 \eeq

\item[(ii)] for any objects $X,Y\in\Ob({\tt K})$ the envelope $E(e_X\otimes e_Y)$ of the morphism $e_X\otimes e_Y:X\otimes Y\to E(X)\otimes E(Y)$ is an isomorphism (in $\tt K$ and in $\tt L$):
 \beq\label{E(e_X-otimes-e_Y)-in-in-Iso}
 E(e_X\otimes e_Y)\in\Iso.
 \eeq
 }\eit
\elm
\bpr
1. Take $A\in{\tt L}$ and $X\in\Ob({\tt K})$. The product of morphisms $1_A:A\to A$ and $e_X:X\to E(X)$ is $1_A\otimes e_X:A\otimes X\to A\otimes E(X)$. If we put it instead of $\alpha$ into \eqref{DIAGR:funktorialnost-env-e-E}, we obtain:
\beq\label{DIAGR:svyaz-otimes-s-env-0}
\xymatrix @R=2.pc @C=5.0pc % @M=14pt
{
A\otimes X \ar[d]^{1_A\otimes e_X}\ar[r]^{e_{A\otimes X}} & E(A\otimes X)\ar[d]^{E(1_A\otimes e_X)} \\
A\otimes E(X)\ar[r]^{e_{A\otimes E(X)}} & E(A\otimes E(X)) \\
}
\eeq
From diagram
$$
\xymatrix @R=2.pc @C=5.0pc % @M=14pt
{
A\otimes X \ar@/_2ex/[dr]_{\ph}\ar[r]^{1_A\otimes e_X} & A\otimes E(X)\ar[r]^{e_{A\otimes E(X)}}\ar@{-->}[d]^{\ph'} & E(A\otimes E(X))\ar@/^2ex/@{-->}[dl]^{\ph''} \\
 & B &
}
$$
it is seen that the composition $e_{A\otimes E(X)}\circ 1_A\otimes e_X$ is an extension for $A\otimes X$ (here in the left triangle we use T.1). Hence $e_{A\otimes E(X)}\circ 1_A\otimes e_X$ is subordinated to the envelope of $A\otimes X$:
\beq\label{DIAGR:svyaz-otimes-s-env-00}
\xymatrix @R=2.pc @C=5.0pc % @M=14pt
{
A\otimes X \ar[d]^{1_A\otimes e_X}\ar[r]^{e_{A\otimes X}} & E(A\otimes X) \\
A\otimes E(X)\ar[r]^{e_{A\otimes E(X)}} & E(A\otimes E(X))\ar@{-->}[u]_{\upsilon} \\
}
\eeq
for some (unique) $\upsilon$. In addition, $\varOmega\subseteq\Epi$, hence morphisms
$e_{A\otimes E(X)}\circ 1_A\otimes e_X$ and $e_{A\otimes X}$, being extensions, are epimorphisms. As a corollary,
\eqref{DIAGR:svyaz-otimes-s-env-0} and \eqref{DIAGR:svyaz-otimes-s-env-00} together give
$$
\upsilon=E(1_A\otimes e_X)^{-1}.
$$

2. For any two objects $X$ and $Y$ the product of morphisms $e_X:X\to E(X)$ and $e_Y:Y\to E(Y)$ is $e_X\otimes e_Y:X\otimes Y\to E(X)\otimes E(Y)$. If we put it instead of $\alpha$ into \eqref{DIAGR:funktorialnost-env-e-E}, we get
\beq\label{DIAGR:svyaz-otimes-s-env}
\xymatrix @R=2.pc @C=5.0pc % @M=14pt
{
X\otimes Y \ar[d]^{e_X\otimes e_Y}\ar[r]^{e_{X\otimes Y}} & E(X\otimes Y)\ar[d]^{E(e_X\otimes e_Y)} \\
E(X)\otimes E(Y)\ar[r]^{e_{E(X)\otimes E(Y)}} & E(E(X)\otimes E(Y)) \\
}
\eeq
From the diagram
$$
\xymatrix @R=2.pc @C=5.0pc % @M=14pt
{
X\otimes Y \ar@/_2ex/[dr]_{\ph}\ar[r]^{e_X\otimes e_Y} & E(X)\otimes E(Y)\ar[r]^{e_{E(X)\otimes E(Y)}}\ar@{-->}[d]^{\ph'} & E(E(X)\otimes E(Y))\ar@/^2ex/@{-->}[dl]^{\ph''} \\
 & B &
}
$$
we see that the composition $e_{E(X)\otimes E(Y)}\circ e_X\otimes e_Y$ is an extension for $X\otimes Y$ (in the left triangle we use T.1). Hence $e_{E(X)\otimes E(Y)}\circ e_X\otimes e_Y$ is subordinated to $X\otimes Y$:
\beq\label{DIAGR:svyaz-otimes-s-env-*}
\xymatrix @R=2.pc @C=5.0pc % @M=14pt
{
X\otimes Y \ar[d]^{e_X\otimes e_Y}\ar[r]^{e_{X\otimes Y}} & E(X\otimes Y) \\
E(X)\otimes E(Y)\ar[r]^{e_{E(X)\otimes E(Y)}} & E(E(X)\otimes E(Y))\ar@{-->}[u]_{\upsilon} \\
}
\eeq
for some (unique) $\upsilon$. And like in the previous case, the morphisms $e_{E(X)\otimes E(Y)}\circ e_X\otimes e_Y$ and $e_{X\otimes Y}$, being extensions, are epimorpisms, so diagrams
\eqref{DIAGR:svyaz-otimes-s-env} and \eqref{DIAGR:svyaz-otimes-s-env-*} together give
$$
\upsilon=E(e_X\otimes e_Y)^{-1}.
$$
\epr

\paragraph{Monoidal structure on the class of complete objects.}
Let $\Env_\varPhi^\varOmega$ be a regular envelope, coherent with the tensor product in $\tt K$,
$E=\Env_\varPhi^\varOmega$ the idempotent functor built in Theorem \ref{TH:regulyarnaya-obolochka}, and ${\tt L}$ the (full) subcategory of complete objects in $\tt K$. For any objects $A,B\in{\tt L}$ and for any morphisms
$\ph,\psi\in{\tt L}$ we put
\beq\label{A-overset(E)(otimes)B:=E(A-otimes-B)}
A\overset{E}{\otimes}B:=E(A\otimes B),\qquad \ph\overset{E}{\otimes}\psi:=E(\ph\otimes\psi).
\eeq
Let us notice the following identity:
\beq\label{E(E(X)-otimes-E(Y))=E(X)-overset(E)(otimes)E(Y)}
E(X)\overset{E}{\otimes}E(Y)=E(E(X)\otimes E(Y)),\qquad X,Y\in \Ob(\tt K)
\eeq
(this is the equality of objects, since by Proposition \ref{PROP:harakterizatsiya-polnoty}, always
$E(X),E(Y)\in{\tt L}$).

\btm\label{TH:sushestvovanie-tenz-proizv-v-L}
Suppose $\Env_\varPhi^\varOmega$ is a regular envelope, coherent with the tensor product in $\tt K$. Then the formulas \eqref{A-overset(E)(otimes)B:=E(A-otimes-B)} define a structure of monoidal category on $\tt L$ (with  $\overset{E}{\otimes}$ as tensor product and $I$ as unit object).
\etm
\bpr
1. The tensor product of local identities must be a local identity. Let us put $1_{A\otimes B}$ instead of $\alpha$ into \eqref{DIAGR:funktorialnost-env-e-E} :
$$
\xymatrix @R=2.pc @C=5.0pc % @M=14pt
{
A\otimes B\ar[d]^{1_{A\otimes B}}\ar[r]^{e(A\otimes B)} & A\overset{E}{\otimes}B\ar[d]^{E(1_{A\otimes B})} \\
A\otimes B\ar[r]^{e(A\otimes B)} & A\overset{E}{\otimes}B \\
}
$$
If we replace here $E(1_{A\otimes B})$ by $1_{A\overset{E}{\otimes}B}$, then the diagram will be also commutative. But this arrow is unique (since $e(A\otimes B)$ is an epimorphism), so these arrows must coincide, and this is used in the last euqlity in the following chain:
$$
1_A\overset{E}{\otimes} 1_B=\eqref{A-overset(E)(otimes)B:=E(A-otimes-B)}=E(1_A\otimes 1_B)=E(1_{A\otimes B})=1_{A\overset{E}{\otimes}B}.
$$

2. Tensor product of commutative diagrams must be a commutative diagram. Suppose we have tow commutative diagrams in ${\tt L}$:
$$
\begin{diagram}
\node[2]{B}\arrow{se,l}{\chi} \\
\node{A}\arrow{ne,l}{\ph}\arrow[2]{e,r}{\psi}\node[2]{C}
\end{diagram}\qquad
\begin{diagram}
\node[2]{B'}\arrow{se,l}{\chi'} \\
\node{A'}\arrow{ne,l}{\ph'}\arrow[2]{e,r}{\psi'}\node[2]{C'}
\end{diagram}
$$
If we multiply them in ${\tt K}$, we obtain a commutative diagram
$$
\begin{diagram}
\node[2]{B\otimes B'}\arrow{se,l}{\chi\otimes \chi'} \\
\node{A\otimes A'}\arrow{ne,l}{\ph\otimes
\ph'}\arrow[2]{e,r}{\psi\otimes\psi'}\node[2]{C\otimes C'}
\end{diagram}
$$
Then we apply the functor $E$ and again obtain a commutative diagram:
$$
\begin{diagram}
\node[2]{E(B\otimes B')}\arrow{se,l}{E(\chi\otimes \chi')} \\
\node{E(A\otimes A')}\arrow{ne,l}{E(\ph\otimes
\ph')}\arrow[2]{e,r}{E(\psi\otimes\psi')}\node[2]{E(C\otimes C')}
\end{diagram}
$$
By \eqref{A-overset(E)(otimes)B:=E(A-otimes-B)} this is the diagram that we need:
$$
\begin{diagram}
\node[2]{B\overset{E}{\otimes} B'}\arrow{se,l}{\chi\overset{E}{\otimes}\chi'} \\
\node{A\overset{E}{\otimes}A'}\arrow{ne,l}{\ph\overset{E}{\otimes}\ph'}
\arrow[2]{e,r}{\psi\overset{E}{\otimes}\psi'}\node[2]{C\overset{E}{\otimes}C'}
\end{diagram}
$$

3. Notice that from what we already proved it follows that the tensor product of isomorphisms in $\tt L$ is also an isomorphism:
\beq\label{Iso-E/otimes-Iso-subseteq-Iso}
\ph,\psi\in\Iso \qquad\Longrightarrow\qquad \ph\overset{E}{\otimes}\psi:=E(\ph\otimes\psi)\in\Iso
\eeq
Indeed,
$$
(\ph\otimes\psi)\circ(\ph^{-1}\otimes\psi^{-1})=(\ph\circ\ph^{-1})\otimes(\psi\circ\psi^{-1})=1\otimes 1=1,
$$
so
$$
(\ph\overset{E}{\otimes}\psi)\circ (\ph^{-1}\overset{E}{\otimes}\psi^{-1})=E(\ph\otimes\psi)\circ E(\ph^{-1}\otimes\psi^{-1})=E((\ph\otimes\psi)\circ(\ph^{-1}\otimes\psi^{-1}))=E(1)=1.
$$
And similarly,
$$
(\ph^{-1}\overset{E}{\otimes}\psi^{-1})\circ (\ph\overset{E}{\otimes}\psi)=1.
$$

4. If $\alpha_{A,B,C}: (A\otimes B)\otimes C\to A\otimes (B\otimes C)$ is the associativity transform in $\tt K$, then the associativity transform $\alpha^E_{A,B,C}: (A\overset{E}{\otimes}B)\overset{E}{\otimes}C\to
A\overset{E}{\otimes}(B\overset{E}{\otimes}C)$ in $\tt L$ is defined by the diagram
\beq\label{DIAGR:opred-alpha-v-L}
\xymatrix @R=2.pc @C=5.0pc % @M=14pt
{
E((A\otimes B)\otimes C) \ar[r]^{E(\alpha_{A,B,C})} & E(A\otimes (B\otimes C))\ar[d]^{E(1_A\otimes e_{B\otimes C})} \\
E(E(A\otimes B)\otimes C)\ar[u]^{E(e_{A\otimes B}\otimes 1_C)^{-1}} & E(A\otimes E(B\otimes C)) \\
E(A\otimes B)\overset{E}{\otimes}C\ar@{=}[u] & A\overset{E}{\otimes}E(B\otimes C)\ar@{=}[u] \\
(A\overset{E}{\otimes}B)\overset{E}{\otimes}C\ar@{=}[u]\ar[r]^{\alpha^E_{A,B,C}} & A\overset{E}{\otimes}
(B\overset{E}{\otimes}C)\ar@{=}[u] \\
}
\eeq
(here we use \eqref{E(E(X)-otimes-E(Y))=E(X)-overset(E)(otimes)E(Y)} and Lemma \ref{PROP:E(e_X-otimes-e_Y)-in-Iso}, which implies that the morphism $E(e_{A\otimes B}\otimes 1_C)$ is invertible).

5. Let us show that the transform $\alpha^E$ is natural with respect to the tensor product:
$$
\alpha^E:\Big((A,B,C)\mapsto (A\overset{E}{\otimes}
B)\overset{E}{\otimes} C\Big)\rightarrowtail\Big((A,B,C)\mapsto A\overset{E}{\otimes} (B\overset{E}{\otimes} C)\Big)
$$
Take morphisms $\ph:A\to A'$, $\chi:B\to B'$, $\psi: C\to C'$ in ${\tt L}$, and consider the diagram of naturality for $\alpha$:
 \beq\label{DIAGR:estestvennost-associativnosti}
\xymatrix @R=2.pc @C=5.0pc % @M=14pt
{
(A\otimes B)\otimes C\ar[d]_{(\ph\otimes \chi)\otimes\psi}
\ar[r]^{\alpha_{A,B,C}} & A\otimes (B\otimes C)\ar[d]_{\ph\otimes
(\chi\otimes\psi)}
\\
(A'\otimes B')\otimes C'\ar[r]_{\alpha_{A',B',C'}}& A'\otimes(B'\otimes C')
}
 \eeq
After applying the functor $E$ we have
$$
\xymatrix @R=2.pc @C=5.0pc % @M=14pt
{
E((A\otimes B)\otimes C)\ar[d]_{E((\ph\otimes \chi)\otimes\psi)}
\ar[r]^{E(\alpha_{A,B,C})} & E(A\otimes (B\otimes C))\ar[d]_{E(\ph\otimes(\chi\otimes\psi))}
\\
E((A'\otimes B')\otimes C')\ar[r]_{E(\alpha_{A',B',C'})}& E(A'\otimes(B'\otimes C'))
}
$$
Let us add this diagram as follows:
$$
\xymatrix @R=3.pc @C=4.0pc % @M=14pt
{
(A\overset{E}{\otimes}B)\overset{E}{\otimes}C \ar@{-->}[r]_{E(e_{A\otimes B}\otimes 1_C)^{-1}}
\ar@{-->}[d]_{(\ph\overset{E}{\otimes}\chi)\overset{E}{\otimes}\psi}\ar@{-->}@/^4ex/[rrr]^{\alpha^E_{A,B,C}}
 & E((A\otimes B)\otimes C)\ar[d]_{E((\ph\otimes \chi)\otimes\psi)}
\ar[r]^{E(\alpha_{A,B,C})} & E(A\otimes (B\otimes C))\ar[d]_{E(\ph\otimes(\chi\otimes\psi))}
\ar@{-->}[r]_{E(1_A\otimes e_{B\otimes C})}
& A\overset{E}{\otimes}(B\overset{E}{\otimes}C)
\ar@{-->}[d]_{\ph\overset{E}{\otimes}(\chi\overset{E}{\otimes}\psi)}
\\
(A'\overset{E}{\otimes}B')\overset{E}{\otimes}C'\ar@{-->}[r]^{E(e_{A'\otimes B'}\otimes 1_{C'})^{-1}}
\ar@{-->}@/_4ex/[rrr]_{\alpha^E_{A',B',C'}}
 & E((A'\otimes B')\otimes C')\ar[r]_{E(\alpha_{A',B',C'})}& E(A'\otimes(B'\otimes C'))
 \ar@{-->}[r]^{E(1_{A'}\otimes e_{B'\otimes C'})}
& A'\overset{E}{\otimes}(B'\overset{E}{\otimes}C')
}
$$
If we throw away the inner vertices, we obtain the diagram of naturality for $\alpha^E$:
$$
\xymatrix @R=2.pc @C=5.0pc % @M=14pt
{
(A\overset{E}{\otimes} B)\overset{E}{\otimes} C\ar[d]_{(\ph\overset{E}{\otimes} \chi)\overset{E}{\otimes}\psi}
\ar[r]^{\alpha^E_{A,B,C}} & A\overset{E}{\otimes} (B\overset{E}{\otimes} C)\ar[d]_{\ph\overset{E}{\otimes}
(\chi\overset{E}{\otimes}\psi)}
\\
(A'\overset{E}{\otimes} B')\overset{E}{\otimes} C'\ar[r]_{\alpha^E_{A',B',C'}}& A'\overset{E}{\otimes}(B'\overset{E}{\otimes} C')
}
$$

6. Let us show that $\alpha^E$ satisfies the associativity conditions. For $\alpha$ they look as the pentagon
 \beq\label{5-ugolnik-dlya-otimes}
\xymatrix  @R=1pc @C=-1pc % @M=2pt
{
 & \text{\scriptsize $(A\otimes(B\otimes C))\otimes D$}\ar[rr]^{\alpha_{A,B\otimes C, D}}
 & & \text{\scriptsize $A\otimes((B\otimes C)\otimes D)$}\ar[ddr]!U|{1_A\otimes \alpha_{B,C,D}} &
 \\
 & & \\
\text{\scriptsize $((A\otimes B)\otimes C)\otimes D$}
 \ar[uur]!U|{\alpha_{A,B,C}\otimes 1_D}\ar[drr]_{\alpha_{A\otimes B,C,D}}
 &&&& \text{\scriptsize $A\otimes(B\otimes (C\otimes D))$} \\
  && \text{\scriptsize $(A\otimes B)\otimes (C\otimes D)$}\ar[urr]_{\alpha_{A,B,C\otimes D}}&&
}
 \eeq
Let us apply $E$ and add the diagram to the following prism:
 \beq
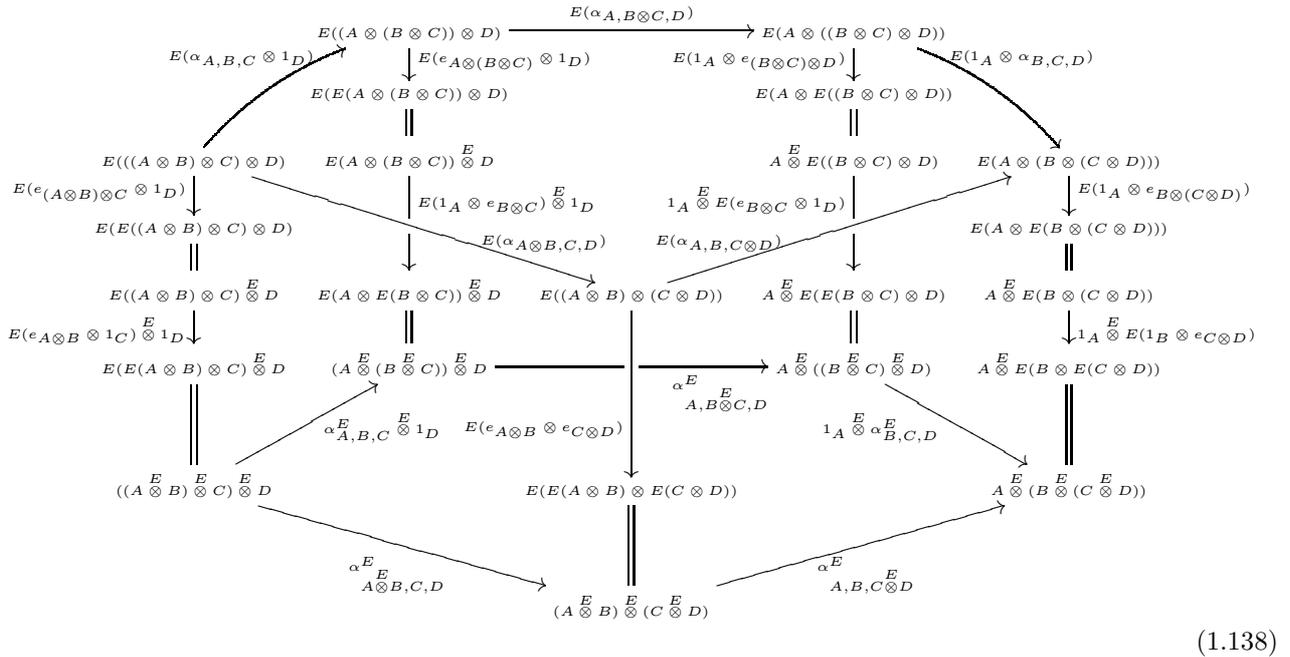
\label{prizma-dlya-otimes^E}
\xymatrix  @R=1pc @C=0pc  %@M=-1pt
{
 & \text{\tiny $E((A\otimes(B\otimes C))\otimes D)$}\ar[d]^{\text{\tiny $E(e_{A\otimes(B\otimes C)}\otimes 1_D)$}} \ar[rr]^{\text{\tiny $E(\alpha_{A,B\otimes C, D})$}}
 & & \text{\tiny $E(A\otimes((B\otimes C)\otimes D))$}\ar@/^2ex/[ddr]^{\text{\tiny $E(1_A\otimes \alpha_{B,C,D})$}}
 \ar[d]_{\text{\tiny $E(1_A\otimes e_{(B\otimes C)\otimes D})$}} &
 \\
& \text{\tiny $E(E(A\otimes(B\otimes C))\otimes D)$}\ar@{=}[d] & & \text{\tiny $E(A\otimes E((B\otimes C)\otimes D))$}\ar@{=}[d] & \\
\text{\tiny $E(((A\otimes B)\otimes C)\otimes D)$}
\ar[d]_{\text{\tiny $E(e_{(A\otimes B)\otimes C}\otimes 1_D)$}}
 \ar@/^2ex/[uur]^{\text{\tiny $E(\alpha_{A,B,C}\otimes 1_D)$}}\ar[ddrr]^(.75){\text{\tiny $E(\alpha_{A\otimes B,C,D})$}}
 &
 \text{\tiny $E(A\otimes(B\otimes C))\overset{E}{\otimes}D$}\ar[dd]^(.3){\text{\tiny $E(1_A\otimes e_{B\otimes C})\overset{E}{\otimes}1_D$}}|!{[d]}{\hole}
 &&
 \text{\tiny $A\overset{E}{\otimes}E((B\otimes C)\otimes D)$}\ar[dd]_(.3){\text{\tiny $1_A\overset{E}{\otimes}E(e_{B\otimes C}\otimes 1_D)$}}|!{[d]}{\hole}
 & \text{\tiny $E(A\otimes(B\otimes (C\otimes D)))$}\ar[d]^{\text{\tiny $E(1_A\otimes e_{B\otimes (C\otimes D)})$}} \\
 \text{\tiny $E(E((A\otimes B)\otimes C)\otimes D)$}\ar@{=}[d]
 &
 &
 &
 &
 \text{\tiny $E(A\otimes E(B\otimes (C\otimes D)))$}\ar@{=}[d]
 \\
  \text{\tiny $E((A\otimes B)\otimes C)\overset{E}{\otimes}D$}
   \ar[d]_{\text{\tiny $E(e_{A\otimes B}\otimes 1_C)\overset{E}{\otimes} 1_D$}}
    & \text{\tiny $E(A\otimes E(B\otimes C))\overset{E}{\otimes}D$}\ar@{=}[d]
  & \text{\tiny $E((A\otimes B)\otimes (C\otimes D))$}
  \ar[ddd]_(.7){\text{\tiny $E(e_{A\otimes B}\otimes e_{C\otimes D})$}}
  \ar[uurr]^(.25){\text{\tiny $E(\alpha_{A,B,C\otimes D})$}}&
  \text{\tiny $A\overset{E}{\otimes}E(E(B\otimes C)\otimes D)$}\ar@{=}[d]
  &
   \text{\tiny $A\overset{E}{\otimes}E(B\otimes (C\otimes D))$}
   \ar[d]^{\text{\tiny $1_A\overset{E}{\otimes} E(1_B\otimes e_{C\otimes D})$}}
    \\
 \text{\tiny $E(E(A\otimes B)\otimes C)\overset{E}{\otimes}D$}\ar@{=}[dd]
  & \text{\tiny $(A\overset{E}{\otimes}(B\overset{E}{\otimes}C))\overset{E}{\otimes}D$}\ar[rr]_(.7){\text{\tiny $\alpha^E_{A,B\overset{E}{\otimes} C, D}$}}|!{[r]}{\hole}
 & & \text{\tiny $A\overset{E}{\otimes}((B\overset{E}{\otimes} C)\overset{E}{\otimes} D)$}\ar[ddr]_(.3){\text{\tiny $1_A\overset{E}{\otimes} \alpha^E_{B,C,D}$}} &
  \text{\tiny $A\overset{E}{\otimes}E(B\otimes E(C\otimes D))$}\ar@{=}[dd]
  \\
 & &
     & & \\
\text{\tiny $((A\overset{E}{\otimes} B)\overset{E}{\otimes} C)\overset{E}{\otimes} D$}
 \ar[uur]_(.7){\text{\tiny $\alpha^E_{A,B,C}\overset{E}{\otimes} 1_D$}}\ar[ddrr]_{\text{\tiny $\alpha^E_{A\overset{E}{\otimes} B,C,D}$}}
 &&
   \text{\tiny $E(E(A\otimes B)\otimes E(C\otimes D))$}\ar@{=}[dd]
 && \text{\tiny $A\overset{E}{\otimes}(B\overset{E}{\otimes} (C\overset{E}{\otimes} D))$}
  \\
 &&  &&
 \\
  && \text{\tiny $(A\overset{E}{\otimes} B)\overset{E}{\otimes} (C\overset{E}{\otimes} D)$}\ar[uurr]_{\text{\tiny $\alpha^E_{A,B,C\overset{E}{\otimes} D}$}}&&
}
 \eeq
The upper base of this prism is commutative, since this is the action of the functor $E$ to the diagram \eqref{5-ugolnik-dlya-otimes}, and the commutativity of the lateral sides can be verified by changing in equivalent way the vertical arrows.

For example, the commutativity of the left side closest to the viewer becomes obvious, if we represent it as a perimeter of the following diagram:
$$
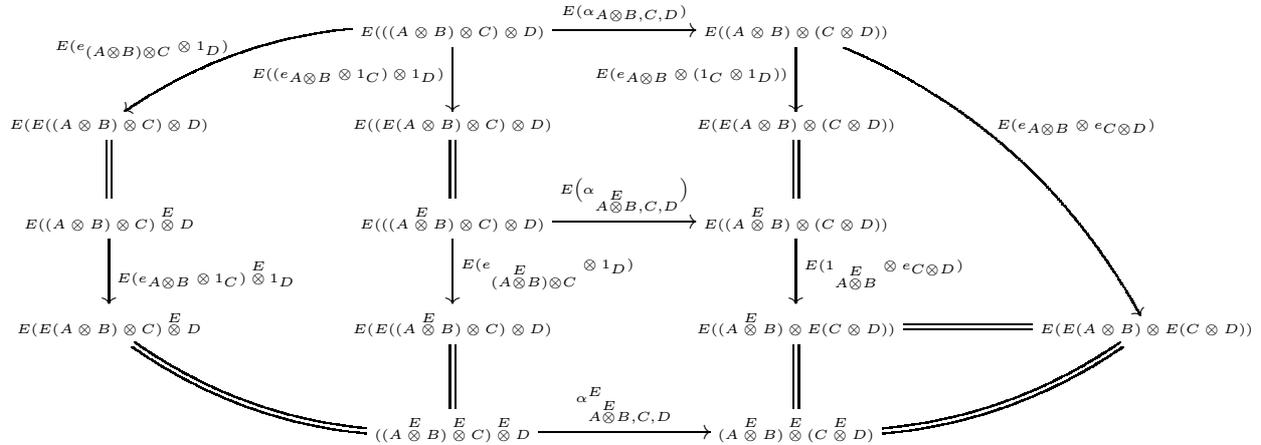
\xymatrix  @R=2pc @C=4pc % @M=2pt
{
 &
 \text{\tiny $E(((A\otimes B)\otimes C)\otimes D)$}\ar[r]^{\text{\tiny $E(\alpha_{A\otimes B,C,D})$}}
 \ar[d]_{\text{\tiny $E((e_{A\otimes B}\otimes 1_C)\otimes 1_D)$}}
 \ar@/_3ex/[dl]_(.7){\text{\tiny $E(e_{(A\otimes B)\otimes C}\otimes 1_D)$}\quad}
 & \text{\tiny $E((A\otimes B)\otimes (C\otimes D))$}
 \ar[d]_{\text{\tiny $E(e_{A\otimes B}\otimes (1_C\otimes 1_D))$}}  \ar@/^4ex/[dddr]^{\text{\tiny $E(e_{A\otimes B}\otimes e_{C\otimes D})$}}
 &
 \\
 \text{\tiny $E(E((A\otimes B)\otimes C)\otimes D)$}
 \ar@{=}[d] & \text{\tiny $E((E(A\otimes B)\otimes C)\otimes D)$}\ar@{=}[d] &
 \text{\tiny $E(E(A\otimes B)\otimes (C\otimes D))$}\ar@{=}[d]
 &
 \\
 \text{\tiny $E((A\otimes B)\otimes C)\overset{E}{\otimes}D$}
 \ar[d]^{\text{\tiny $E(e_{A\otimes B}\otimes 1_C)\overset{E}{\otimes}1_D$}}
 & \text{\tiny $E(((A\overset{E}{\otimes}B)\otimes C)\otimes D)$}
 \ar[r]^{\text{\tiny $E\Big(\alpha_{A\overset{E}{\otimes}B,C,D}\Big)$}}
 \ar[d]^{\text{\tiny $E(e_{(A\overset{E}{\otimes}B)\otimes C}\otimes 1_D)$}}
 & \text{\tiny $E((A\overset{E}{\otimes}B)\otimes (C\otimes D))$}
 \ar[d]^{\text{\tiny $E(1_{A\overset{E}{\otimes}B}\otimes e_{C\otimes D})$}}
 &
 \\
 \text{\tiny $E(E(A\otimes B)\otimes C)\overset{E}{\otimes}D$}
 \ar@/_3ex/@{=}[dr]
 & \text{\tiny $E(E((A\overset{E}{\otimes}B)\otimes C)\otimes D)$}\ar@{=}[d]
 & \text{\tiny $E((A\overset{E}{\otimes}B)\otimes E(C\otimes D))$}\ar@{=}[d]
 &
 \text{\tiny $E(E(A\otimes B)\otimes E(C\otimes D))$}\ar@{=}[l]  \ar@/^3ex/@{=}[dl]
 \\
 & \text{\tiny $((A\overset{E}{\otimes}B)\overset{E}{\otimes}C)\overset{E}{\otimes}D$}
\ar[r]^{\text{\tiny $\alpha^E_{A\overset{E}{\otimes}B,C,D}$}}
 & \text{\tiny $(A\overset{E}{\otimes}B)\overset{E}{\otimes}(C\overset{E}{\otimes} D)$}
 &
}
$$
Here the upper inner hexagon (or it can be called quadrangle) is the result of applying $E$ to the diagram
 $$
\xymatrix  @R=2pc @C=4pc % @M=2pt
{
 ((A\otimes B)\otimes C)\otimes D \ar[r]^{\alpha_{A\otimes B,C,D}} \ar[d]_{(e_{A\otimes B}\otimes 1_C)\otimes 1_D} & (A\otimes B)\otimes (C\otimes D)\ar[d]^{e_{A\otimes B}\otimes (1_C\otimes 1_D)}
\\
 (E(A\otimes B)\otimes C)\otimes D\ar[r]^{\alpha_{E(A\otimes B),C,D}} &  E(A\otimes B)\otimes (C\otimes D)
}
$$
(this is a corollary of \eqref{DIAGR:estestvennost-associativnosti}). The lower inner hexagon is diagram \eqref{DIAGR:opred-alpha-v-L} for $\alpha^E$ on the components $A\overset{E}{\otimes}B$, $C$, $D$. The big octagon can be represented as a rhomb
 $$
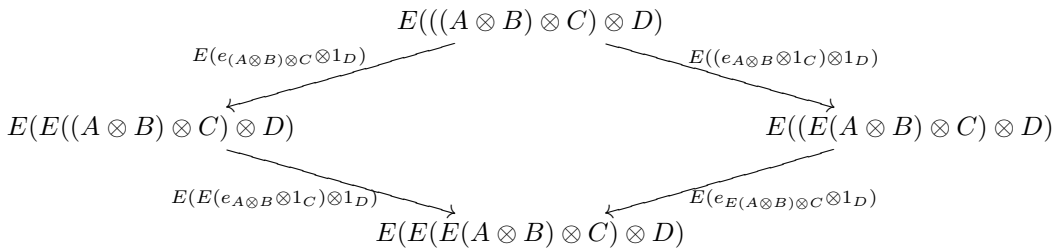
\xymatrix  @R=2pc @C=2pc % @M=2pt
{
  & E(((A\otimes B)\otimes C)\otimes D) \ar[dl]_{E(e_{(A\otimes B)\otimes C}\otimes 1_D)\qquad}
  \ar[dr]^{\qquad E((e_{A\otimes B}\otimes 1_C)\otimes 1_D)}
  &
\\
E(E((A\otimes B)\otimes C)\otimes D)\ar[dr]_{E(E(e_{A\otimes B}\otimes 1_C)\otimes 1_D)\qquad}
 & & E((E(A\otimes B)\otimes C)\otimes D)\ar[dl]^{\qquad E(e_{E(A\otimes B)\otimes C}\otimes 1_D)}
\\
& E(E(E(A\otimes B)\otimes C)\otimes D) &
}
$$
which is a result of allying $E$ to the rhomb
 $$
\xymatrix  @R=2pc @C=2pc % @M=2pt
{
  &((A\otimes B)\otimes C)\otimes D \ar[dl]_{e_{(A\otimes B)\otimes C}\otimes 1_D\quad}
  \ar[dr]^{\quad (e_{A\otimes B}\otimes 1_C)\otimes 1_D}
  &
\\
E((A\otimes B)\otimes C)\otimes D\ar[dr]_{E(e_{A\otimes B}\otimes 1_C)\otimes 1_D\quad}
 & & (E(A\otimes B)\otimes C)\otimes D\ar[dl]^{\quad e_{E(A\otimes B)\otimes C}\otimes 1_D}
\\
& E(E(A\otimes B)\otimes C)\otimes D &
}
$$
which in its turn is a result of multiplication by $D$ to the right the diagram
 $$
\xymatrix  @R=2pc @C=2pc % @M=2pt
{
  &(A\otimes B)\otimes C \ar[dl]_{e_{(A\otimes B)\otimes C}\quad}
  \ar[dr]^{\quad e_{A\otimes B}\otimes 1_C}
  &
\\
E((A\otimes B)\otimes C)\ar[dr]_{E(e_{A\otimes B}\otimes 1_C)\quad}
 & & E(A\otimes B)\otimes C\ar[dl]^{\quad e_{E(A\otimes B)\otimes C}}
\\
& E(E(A\otimes B)\otimes C) &
}
$$
and it can be perceived as diagram \eqref{DIAGR:funktorialnost-env-e-E} where $\alpha$ is replaced by $e_{A\otimes B}\otimes 1_C$. Finally, the upper right pentagon is the triangle
 $$
\xymatrix  @R=2pc @C=5pc % @M=2pt
{
  E((A\otimes B)\otimes(C\otimes D))\ar[d]_{E(e_{A\otimes B}\otimes (1_C\otimes 1_D))}
  \ar[dr]^{\qquad E(e_{A\otimes B}\otimes e_{C\otimes D})}
  &
\\
  E(E(A\otimes B)\otimes(C\otimes D))\ar[r]_{E(1_{A\otimes B}\otimes e_{C\otimes D})} &  E(E(A\otimes B)\otimes E(C\otimes D))
}
$$
which is a result of applying $E$ to the triangle
 $$
\xymatrix  @R=2pc @C=5pc % @M=2pt
{
  (A\otimes B)\otimes(C\otimes D)\ar[d]_{e_{A\otimes B}\otimes (1_C\otimes 1_D)}
  \ar[dr]^{\qquad e_{A\otimes B}\otimes e_{C\otimes D}}
  &
\\
  E(A\otimes B)\otimes(C\otimes D)\ar[r]_{1_{A\otimes B}\otimes e_{C\otimes D}} &  E(A\otimes B)\otimes E(C\otimes D)
}
$$

And the same with the other vertical sides of \eqref{prizma-dlya-otimes^E}. In addition, the vertical arrows are isomorphisms (by Lemma \ref{PROP:E(e_X-otimes-e_Y)-in-Iso} and property \eqref{Iso-E/otimes-Iso-subseteq-Iso}), so we obtain that the lower base of this prism is commutative as well, and this is the diagram that we need for $\alpha^E$.

7. Let $\lambda_X:I\otimes X\to X$ be the left identity in the monoidal category $\tt K$, and $\rho_X:X\otimes I\to X$ the righta identity. For any $A\in\Ob({\tt L})$ we can put
\beq\label{DEF:lambda_A^E,rho_A^E}
\lambda_A^E=E(\lambda_A):I\overset{E}{\otimes}A=E(I\otimes A)\to E(A)=A,\qquad \rho_A^E=E(\rho_A):A\overset{E}{\otimes}I=E(A\otimes I)\to E(A)=A,
\eeq
and this will be the left and the righta identities for $\tt L$. Indeed, for any morphism
$\ph:A\to A'$ in ${\tt L}$ the diagrams
 \beq\label{DIAGR:lambda,rho}
\xymatrix  @R=2pc @C=4pc % @M=2pt
{
I\otimes A \ar[d]_{1_I\otimes \ph} \ar[r]^{\lambda_A} & A\ar[d]^{\ph}
\\
I\otimes A'\ar[r]^{\lambda_{A'}} & A'
}
\qquad\qquad
\xymatrix  @R=2pc @C=4pc % @M=2pt
{
A\otimes I \ar[d]_{\ph\otimes 1_I} \ar[r]^{\rho_A} & A\ar[d]^{\ph}
\\
A'\otimes I\ar[r]^{\rho_{A'}} & A'
}
 \eeq
give
 $$
\xymatrix  @R=2pc @C=2.5pc % @M=2pt
{
I\overset{E}{\otimes}A\ar@{=}[r] & E(I\otimes A) \ar[d]_{1_I\overset{E}{\otimes}\ph=E(1_I\otimes \ph)} \ar[rr]^{\lambda_A^E=E(\lambda_A)} & & A\ar[d]^{\ph}
\\
I\overset{E}{\otimes}A'\ar@{=}[r] & E(I\otimes A')\ar[rr]^{\lambda_{A'}^E=E(\lambda_{A'})} & & A'
}
\qquad\qquad
\xymatrix  @R=2pc @C=2.5pc % @M=2pt
{
A\overset{E}{\otimes}I \ar@{=}[r] & E(A\otimes I)\ar[d]_{\ph\overset{E}{\otimes} 1_I=E(\ph\otimes 1_I)} \ar[rr]^{\rho_A^E=E(\rho_A)} & & A\ar[d]^{\ph}
\\
A'\overset{E}{\otimes}I \ar@{=}[r] & E(A'\otimes I)\ar[rr]^{\rho_{A'}^E=E(\rho_{A'})} & & A'
}
 $$
Besides this, the identity $\lambda_I=\rho_I$ implies the identity $\lambda_I^E=E(\lambda_I)=E(\rho_I)=\rho_I^E$, and the diagram
$$
\xymatrix  @R=2pc @C=2.5pc % @M=2pt
{
(A\otimes I)\otimes B
 \ar[dr]_{\rho_A\otimes 1_B}
 \ar[rr]^{\alpha_{A,I,B}}
 & &
A\otimes (I\otimes B)
 \ar[dl]^{1_A\otimes\lambda_B}
\\
& A\otimes B &
}
$$
gives the upper base of the prism
$$
\xymatrix  @R=2pc @C=5pc % @M=2pt
{
& E((A\otimes I)\otimes B)\ar[d]_{E(e_{A\otimes I}\otimes 1_B)}
 \ar[dr]^{E(\rho_A\otimes 1_B)}
 \ar[rr]^{E(\alpha_{A,I,B})}
 & &
E(A\otimes (I\otimes B))\ar[d]^{E(1_A\otimes e_{I\otimes B})}
 \ar[dl]_{E(1_A\otimes\lambda_B)}
 &
\\
& E(E(A\otimes I)\otimes B)\ar@{=}[d]\ar[r]_{E(\rho_A^E\otimes 1_B)}  & E(A\otimes B)\ar@{=}[dd] & E(A\otimes E(I\otimes B))\ar@{=}[d] \ar[l]^{E(1_A\otimes \lambda_B^E)} &
\\
& (A\overset{E}{\otimes}I)\overset{E}{\otimes}B\ar[dr]^{\rho_A^E\overset{E}{\otimes} 1_B}
\ar[rr]^(.7){\alpha^E_{A,I,B}}|!{[r]}{\hole} && A\overset{E}{\otimes}(I\overset{E}{\otimes}B)\ar[dl]_{1_A\overset{E}{\otimes}\lambda^E_B} &
\\
&& A\overset{E}{\otimes}B &&
}
$$
The commutativity of its lateral sides is obvious, and the vertical arrows are isomorphisms, so the lower base must also be commutative.
\epr

\paragraph{Envelope as a monoidal functor.}

\btm\label{TH:funktor-E-monoidalen}
Let $\Env_\varPhi^\varOmega$ be a regular envelope, coherent with the tensor product in $\tt K$.
Then the functor of the envelope $E:{\tt K}\to {\tt L}$, built in Theorem \ref{TH:regulyarnaya-obolochka},
is monoidal.
\etm
\bpr
To be monoidal the functor $E:{\tt K}\to {\tt L}$ must define a morphism of bifunctors
 $$
\Big((X,Y)\mapsto E(X)\overset{E}{\otimes} E(Y)\Big)\overset{E^{\otimes}}{\rightarrowtail} \Big((X,Y)\mapsto E(X\otimes Y)\Big)
 $$
in this case this is a family of morphisms
$$
E^\otimes_{X,Y}=E(e_X\otimes e_Y)^{-1}:E(X)\overset{E}{\otimes} E(Y)=E(E(X)\otimes E(Y))\to E(X\otimes Y),
$$
(by Lemma \ref{PROP:E(e_X-otimes-e_Y)-in-Iso} all morphisms $E(e_X\otimes e_Y)$ are isomorphisms, so there exist $E(e_X\otimes e_Y)^{-1}$) and a morphism $E^I$ in $\tt L$, that turns the identity object $I$ of $\tt L$ into the image $E(I)$ of the identity object $I$ in $\tt K$, and in this situation this will be the local identity:
$$
E^I=1_{I}:I\to I=\eqref{E(I)=I}=E(I).
$$
Let us check the axioms of monoidal functor for these components. The diagram of coherence with the associativity
\beq\label{DIAGR:assoc-v-L}
\xymatrix @R=2.5pc @C=6.0pc
{
E\big((X\otimes Y)\otimes Z\big)\ar[r]^{E(\alpha_{X,Y,Z})}& E\big(X\otimes (Y\otimes Z)\big)
\\
E(X\otimes Y)\overset{E}{\otimes}E(Z)\ar[u]^{E^{\otimes}_{X\otimes Y,Z}}
& E(X)\overset{E}{\otimes}E(Y\otimes Z)\ar[u]_{E^{\otimes}_{X,Y\otimes Z}}
 \\
\big(E(X)\overset{E}{\otimes}E(Y)\big)\overset{E}{\otimes}E(Z)
\ar[u]^{E^{\otimes}_{X,Y}\overset{E}{\otimes}1_{E(X)}}
\ar[r]^{\alpha^E_{E(X),E(Y),E(Z)}}& E(X)\overset{E}{\otimes}\big(E(Y)\overset{E}{\otimes}E(Z)\big)
\ar[u]_{1_{E(X)}\overset{E}{\otimes}E^{\otimes}_{Y,Z}}
}
\eeq
is translated here as follows:
$$
\xymatrix @R=2.5pc @C=6.0pc
{
E\big((X\otimes Y)\otimes Z\big)\ar[d]_{E(e_{X\otimes Y}\otimes e_Z)}\ar[r]^{E(\alpha_{X,Y,Z})}& E\big(X\otimes (Y\otimes Z)\big)
\ar[d]^{E(e_X\otimes e_{Y\otimes Z})}
\\
E(E(X\otimes Y)\otimes E(Z))\ar[d]_{E(E(e_X\otimes e_Y)\otimes 1_{E(Z)})}
& E(E(X)\otimes E(Y\otimes Z))\ar[d]^{E(1_{E(X)}\otimes E(e_Y\otimes e_Z))}
 \\
E(E(E(X)\otimes E(Y))\otimes E(Z))
\ar[r]^{\alpha^E_{E(X),E(Y),E(Z)}}&
E(E(X)\otimes E(E(Y)\otimes E(Z)))
}
$$
To see that it is commutative, let us represent it as the perimeter of the following diagram:
\beq\label{DIAGR:funktor-E-monoidalen}
\xymatrix @R=3pc @C=4pc
{
&
\text{\tiny $E((X\otimes Y)\otimes Z)$}
\ar@/_3ex/[dl]_{\text{\tiny $E(e_{X\otimes Y}\otimes e_Z)$}}
\ar[d]!U|{\text{\tiny $E((e_X\otimes e_Y)\otimes e_Z)$}}
\ar[r]^{\text{\tiny $E(\alpha_{X,Y,Z})$}}&
\text{\tiny $E(X\otimes (Y\otimes Z))$}
\ar[d]!U|{\text{\tiny $E(e_X\otimes (e_Y\otimes e_Z))$}}
\ar@/^3ex/[dr]^{\text{\tiny $E(e_X\otimes e_{Y\otimes Z})$}}
&
\\
\text{\tiny $E(E(X\otimes Y)\otimes E(Z))$}
\ar@/_3ex/[dr]_(.3){\text{\tiny $E(E(e_X\otimes e_Y)\otimes 1_{E(Z)})$}\quad}
&
\text{\tiny $E((E(X)\otimes E(Y))\otimes E(Z))$}\ar[r]^{\text{\tiny $E(\alpha_{E(X),E(Y),E(Z)}$)}}
\ar[d]!U|{\text{\tiny $E(e_{E(X)\otimes E(Y)}\otimes 1_{E(Z)})$}}
&
\text{\tiny $E(E(X)\otimes (E(Y)\otimes E(Z)))$}
\ar[d]!U|{\text{\tiny $E(1_{E(X)}\otimes e_{E(Y)\otimes E(Z)})$}}
& \text{\tiny $E(E(X)\otimes E(Y\otimes Z))$}
\ar@/^3ex/[dl]^(.3){\quad\text{\tiny $E(1_{E(X)}\otimes E(e_Y\otimes e_Z))$}}
 \\
&
\text{\tiny $E(E(E(X)\otimes E(Y))\otimes E(Z))$}
\ar[r]_{\text{\tiny $\alpha^E_{E(X),E(Y),E(Z)}$}}&
\text{\tiny $E(E(X)\otimes E(E(Y)\otimes E(Z)))$} &
}
\eeq
Here the left inner triangle can be represented in the form
$$
\xymatrix @R=3pc @C=8pc
{
E((X\otimes Y)\otimes Z)
\ar[d]_{E(e_{X\otimes Y}\otimes e_Z)}
\ar[r]^{E((e_X\otimes e_Y)\otimes e_Z)}
&
E((E(X)\otimes E(Y))\otimes E(Z))
\ar[d]^{E(e_{E(X)\otimes E(Y)}\otimes 1_{E(Z)})}
 \\
E(E(X\otimes Y)\otimes E(Z))
\ar[r]_{E(E(e_X\otimes e_Y)\otimes 1_{E(Z)})}
&
E(E(E(X)\otimes E(Y))\otimes E(Z))
}
$$
This is the result of applying $E$ to the diagram
$$
\xymatrix @R=3pc @C=8pc
{
(X\otimes Y)\otimes Z
\ar[d]_{e_{X\otimes Y}\otimes e_Z}
\ar[r]^{(e_X\otimes e_Y)\otimes e_Z}
&
(E(X)\otimes E(Y))\otimes E(Z)
\ar[d]^{e_{E(X)\otimes E(Y)}\otimes 1_{E(Z)}}
 \\
E(X\otimes Y)\otimes E(Z)
\ar[r]_{E(e_X\otimes e_Y)\otimes 1_{E(Z)}}
&
E(E(X)\otimes E(Y))\otimes E(Z)
}
$$
which in its turn is the product of two diagrams
$$
\xymatrix @R=3pc @C=5pc
{
X\otimes Y
\ar[d]_{e_{X\otimes Y}}
\ar[r]^{(e_X\otimes e_Y)}
&
E(X)\otimes E(Y)
\ar[d]^{e_{E(X)\otimes E(Y)}}
 \\
E(X\otimes Y)
\ar[r]_{E(e_X\otimes e_Y)}
&
E(E(X)\otimes E(Y))
}
\qquad
\xymatrix @R=3pc @C=5pc
{
Z
\ar[d]_{e_Z}
\ar[r]^{e_Z}
&
E(Z)
\ar[d]^{1_{E(Z)}}
 \\
E(Z)
\ar[r]_{1_{E(Z)}}
&
E(Z)
}
$$
The left one of them is trivial, and the right one is the transposed diagram \eqref{DIAGR:svyaz-otimes-s-env}.

Further, the upper inner quadrangle in \eqref{DIAGR:funktor-E-monoidalen}
$$
\xymatrix @R=3pc @C=6pc
{
E((X\otimes Y)\otimes Z)
\ar[d]_{E((e_X\otimes e_Y)\otimes e_Z)}
\ar[r]^{E(\alpha_{X,Y,Z})}&
E(X\otimes (Y\otimes Z))
\ar[d]^{E(e_X\otimes (e_Y\otimes e_Z))}
\\
E((E(X)\otimes E(Y))\otimes E(Z))
\ar[r]^{E(\alpha_{E(X),E(Y),E(Z)})}
&
E(E(X)\otimes (E(Y)\otimes E(Z)))
}
$$
is the result of applying $E$ to the diagram
$$
\xymatrix @R=3pc @C=6pc
{
(X\otimes Y)\otimes Z
\ar[d]_{(e_X\otimes e_Y)\otimes e_Z}
\ar[r]^{\alpha_{X,Y,Z}}&
X\otimes (Y\otimes Z)
\ar[d]^{e_X\otimes (e_Y\otimes e_Z)}
\\
(E(X)\otimes E(Y))\otimes E(Z)
\ar[r]^{\alpha_{E(X),E(Y),E(Z)}}
&
E(X)\otimes (E(Y)\otimes E(Z))
}
$$
and this is the special case of \eqref{DIAGR:estestvennost-associativnosti}.

Then, the lower inner quadrangle in \eqref{DIAGR:funktor-E-monoidalen}
$$
\xymatrix @R=3pc @C=6pc
{
E((E(X)\otimes E(Y))\otimes E(Z))
\ar[r]^{E(\alpha_{E(X),E(Y),E(Z)})}
\ar[d]!U|{E(e_{E(X)\otimes E(Y)}\otimes 1_{E(Z)})}
&
E(E(X)\otimes (E(Y)\otimes E(Z)))
\ar[d]!U|{E(1_{E(X)}\otimes e_{E(Y)\otimes E(Z)})}
 \\
E(E(E(X)\otimes E(Y))\otimes E(Z))
\ar[r]_{\alpha^E_{E(X),E(Y),E(Z)}}&
E(E(X)\otimes E(E(Y)\otimes E(Z)))
}
$$
-- is a special case of diagram \eqref{DIAGR:opred-alpha-v-L}.

Finally, it is useful to represent the right inner quadrangle in \eqref{DIAGR:funktor-E-monoidalen} in the form
$$
\xymatrix @R=3pc @C=8pc
{
E(X\otimes (Y\otimes Z))
\ar[d]_{E(e_X\otimes (e_Y\otimes e_Z))}
\ar[r]^{E(e_X\otimes e_{Y\otimes Z})}
&
E(E(X)\otimes E(Y\otimes Z))
\ar[d]^{E(1_{E(X)}\otimes E(e_Y\otimes e_Z))}
\\
E(E(X)\otimes (E(Y)\otimes E(Z)))
\ar[r]_{E(1_{E(X)}\otimes e_{E(Y)\otimes E(Z)})}
&
E(E(X)\otimes E(E(Y)\otimes E(Z)))
}
$$
This is the result of applying $E$ to the diagram
$$
\xymatrix @R=3pc @C=8pc
{
X\otimes (Y\otimes Z)
\ar[d]_{e_X\otimes (e_Y\otimes e_Z)}
\ar[r]^{e_X\otimes e_{Y\otimes Z}}
&
E(X)\otimes E(Y\otimes Z)
\ar[d]^{1_{E(X)}\otimes E(e_Y\otimes e_Z)}
\\
E(X)\otimes (E(Y)\otimes E(Z))
\ar[r]_{1_{E(X)}\otimes e_{E(Y)\otimes E(Z)}}
&
E(X)\otimes E(E(Y)\otimes E(Z))
}
$$
which in its turn is product of two diagrams
$$
\xymatrix @R=3pc @C=6pc
{
X
\ar[d]_{e_X}
\ar[r]^{e_X}
&
E(X)
\ar[d]^{1_{E(X)}}
\\
E(X)
\ar[r]_{1_{E(X)}}
&
E(X)
}
\qquad
\xymatrix @R=3pc @C=6pc
{
Y\otimes Z
\ar[d]_{e_Y\otimes e_Z}
\ar[r]^{e_{Y\otimes Z}}
&
E(Y\otimes Z)
\ar[d]^{E(e_Y\otimes e_Z)}
\\
E(Y)\otimes E(Z)
\ar[r]_{e_{E(Y)\otimes E(Z)}}
&
E(E(Y)\otimes E(Z))
}
$$
The left one here is trivial, and the right one is a little bit changed \eqref{DIAGR:svyaz-otimes-s-env}.

Apart from this we need to verify the commutativity of the diagrams for the left and for the right identities:
$$
\xymatrix @R=2.5pc @C=4.0pc
{
E(I\otimes X)\ar[r]^{E(\lambda_X)} & E(X)
\\
I\overset{E}{\otimes} E(X)\ar[u]^{E^{\otimes}_{I,X}}
& I'\overset{E}{\otimes}E(X)\ar[u]_{\lambda'_{E(X)}}\ar[l]^{E^I\overset{E}{\otimes}1_{E(X)}}
}
\qquad
\xymatrix @R=2.5pc @C=4.0pc
{
E(X\otimes I)\ar[r]^{E(\rho_X)} & E(X)
\\
E(X)\overset{E}{\otimes}I \ar[u]^{E^{\otimes}_{X,I}}
 & E(X)\overset{E}{\otimes} I'\ar[u]_{\rho'_{E(X)}}\ar[l]^{1_{E(X)}\overset{E}{\otimes} E^I}
}
$$
In our situation they have the following form
$$
\xymatrix @R=2.5pc @C=6.0pc
{
E(I\otimes X)\ar[r]^{E(\lambda_X)}\ar[d]_{E(e_I\otimes e_X)} & E(X)
\\
E(I\otimes E(X))\ar[r]_{E(1_{E(X)}\otimes 1_{E(X)})} & E(I\otimes E(X))\ar[u]_{E(\lambda_{E(X)})}
}
\qquad
\xymatrix @R=2.5pc @C=6.0pc
{
E(X\otimes I)\ar[r]^{E(\rho_X)}\ar[d]_{E(e_X\otimes I)} & E(X)
\\
E(E(X)\otimes I)\ar[r]_{E(1_{E(X)}\otimes 1_{E(X)})} & E(E(X)\otimes I)\ar[u]_{E(\rho_{E(X)})}
}
$$
and this is the result of applying $E$ to \eqref{DIAGR:lambda,rho} with $X'=E(X)$ and $\ph=e_X$.
\epr

\bcor\label{COR:Env-sohranyaet-algebry}
Suppose $\Env_\varPhi^\varOmega$ is a regular envelope coherent with the tensor product in $\tt K$.
The operation $\Env_\varPhi^\varOmega$ turns each algebra (respectively, coalgebra, bialgebra, Hopf algebra) $A$ in $\tt K$ into an algebra (respectively, coalgebra, bialgebra, Hopf algebra) $\Env_\varPhi^\varOmega A$ in $\tt L$.
\ecor
\bpr
For the case of algebras and general monoidal functors this fact is pointed out in \cite{Saavedra-Rivano}.
\epr

\section{The category of stereotype spaces ${\tt Ste}$}\label{SEC:kateg-Ste}

In this section we speak about applications of the above results to the theory of stereotype spaces. With the aim to make the exposition more self-contained we give a little summary of the simplest facts of the theory
(see details in the author's works \cite{Akbarov} and \cite{Akbarov-stein-groups}).

\subsection{Pseudocomplete and pseudosaturated spaces}

\paragraph{Totally bounded and capacious sets.}

A set $S$ in a locally convex space $X$ is said to be {\it totally bounded} (or {\it precompact}) \cite{Schaefer}, if for each neighbourhood of zero $U$ in $X$ there is a finite set $A$ such that the shifts of $U$ by elements of  $A$ cover $S$: $S\subseteq U+A$. This is equivalent to the fact that $S$ is totally bounded in the sense of the uniform structure \cite{Engelking} induced from $X$ (i.e. {\it $A$ can be chosen as a subset in $S$}).

A set $D\subseteq X$ is said to be {\it capacious}, if for any totally bounded set $S\subseteq X$ there is a finite set $A\subseteq X$ such that the shifts of $D$ by elements of $A$ cover $S$: $S\subseteq D+A$. (If $D$ is convex, then $A$ can be chosen as a subset in $S$.)

Let $X$ be a locally convex space over the field of complex numbers $\C$. Denote by $X^\star$\label{DEF:X^star} the set of lineat continuous functionals $f:X\to \C$ endowed by the {\it topology of uniform convergence on totally bounded sets} in $X$. We call $X^\star$ the {\it dual space} for the space $X$.

If $B\subseteq X$ and $F\subseteq X^\star$ are arbitrary sets, then by $B^\circ$ and ${^\circ F}$ we denote their  (direct and reverse) polars (in $X^\star$ and in $X$):
$$
B^\circ= \{f\in X^\star: |f|_B:= \sup_{x\in B}|f(x)|\le 1\}, \qquad {^\circ F}=
\{x\in X: |x|_F:= \sup_{f\in F}|f(x)|\le 1\}
$$
Similarly, {\it annihilators} of $B$ and $F$ are the sets
$$
B^\perp= \{f\in X^\star:\quad \forall x\in B\quad f(x)=0\}, \qquad
{^\perp F}=\{x\in X: \quad \forall f\in F\quad f(x)=0\}.
$$

\blm\label{th-2.1} For each locally convex space $X$
\begin{itemize}
\item[(a)] if $B\subseteq X$ is totally bounded, then $B^\circ \subseteq X^\star$ is capacious;

\item[(b)] if $B\subseteq X$ is capacious, then $B^\circ \subseteq X^\star$ is totally bounded;

\item[(c)] if $F\subseteq X^\star$ is totally bounded, then ${^\circ F}\subseteq X$ is capacious;

\item[(d)] if $F\subseteq X^\star$ is capacious, then ${^\circ F}\subseteq X$ is totally bounded.
\end{itemize}
\elm

\blm \label{lm-3.3}
For each LCS $X$, every set $A\subseteq X$ and every subspace $E\subseteq X$
\beq\label{A^circ-cap-E^perp=(A+E)^circ}
A^\circ\cap E^\perp=(A+E)^\circ
\eeq
\elm
\bpr In the trivial situation when $A=\varnothing$ or $E=0$ there is nothing to prove, so we assume that  $A\ne\varnothing$ and $E\ne 0$. Then
\begin{multline*}
f\in A^\circ\cap E^\perp\quad\Longrightarrow\quad
\sup_{a\in A}\abs{f(a)}\le 1\quad \& \quad \forall x\in E\quad f(x)=0\quad\Longrightarrow\\
\Longrightarrow\quad
\sup_{a\in A,x\in E}|f(a+x)|=\sup_{a\in A,x\in E}|f(a)+\underbrace{f(x)}_{\scriptsize\begin{matrix}\| \\ 0\end{matrix}}|\le 1 \quad\Longrightarrow\quad
f\in (A+E)^\circ.
\end{multline*}
and
\begin{multline*}
f\in (A+E)^\circ \quad\Longrightarrow\quad
\sup_{a\in A,x\in E}|f(a+x)|\le 1  \quad\Longrightarrow\\
\Longrightarrow\quad
\sup_{a\in A}|f(a)|=\sup_{a\in A,x=0}|f(a+x)|\le 1 \quad \& \quad
\exists a\in A\quad \forall x\in E\quad |f(a+x)|\le 1 \quad\Longrightarrow\\
\Longrightarrow\quad
\sup_{a\in A}|f(a)|\le 1 \quad \& \quad
\forall x\in E\quad f(x)=0 \quad\Longrightarrow\quad
f\in A^\circ\cap E^\perp
\end{multline*}
\epr

\paragraph{Pseudocomplete and pseudosaturated spaces.}

 \bit{
 \item
A locally convex space $X$ is said to be {\it pseudocomplete}, if every totally bounded Cauchy net in $X$ converges. This is equivalent to the claim that every closed totally bounded set in $X$ is compact.
}\eit
This is connected with the usual completeness and quasicompleteness\footnote{A locally convex space $X$ is said to be {\it quasicomplete}, if every bounded Cauchy net in $X$ converges.} by the implications
$$
\text{$X$ is complete $\Longrightarrow$ $X$ is quasicomplete $\Longrightarrow$ $X$ is pseudocomplete}
$$
In the metrizable case these properties are equivalent.

\bit{
\item A locally convex space $X$ is said to be {\it pseudosaturated}, if each closed convex balanced capacious set $D$ in $X$ is a neighbourhood of zero.
}\eit

\begin{ex}\label{ex-1.13} Every barreled space is pseudosaturated.
\end{ex}

\begin{ex}\label{ex-1.14} Every metrizable (not necessarily complete) space is pseudosaturated.
\end{ex}

\btm[criterion of being pseudosaturated]\label{th-2.5}
For a locally convex space $X$ the following conditions are equivalent:
\begin{itemize}
\item[(i)] $X$ is pseudosaturated,

\item[(ii)] if a set of linear continuous functionals $F\subseteq X'$ is equicontinuous on each totally bounded set $S\subseteq X$, then $F$ is equicontinuous on $X$.

\item[(iii)] if $Y$ is a locally convex space and $\varPhi$ is a set of linear continuous maps $\ph:X\to Y$, equicontinuous on each totally bounded set $S\subseteq X$, then $\varPhi$ is equicontinuous on $X$.
\end{itemize}
\etm

\btm \label{th-2.16} For an arbitrary locally convex space $X$
\begin{itemize}
\item[--] if $X$ is pseudocomplete, then $X^\star$ is pseudosaturated;

\item[--] if $X$ is pseudosaturated, then $X^\star$ is pseudocomplete.
\end{itemize}
\etm

\blm \label{lm-3.3} Let $\ph :X\to Y$ be a morphism of LCS. Then
\begin{align}\label{eq3.1}
& \forall A\subseteq X \quad \ph(A)^\circ_{Y^\star}=(\ph^\star)^{-1}(A^\circ_{X^\star})
&& (\ph^\star)^{-1}(0)=\Big(\overline{\ph(X)}\Big)^\perp,
&& (\ph^\star)^{-1}(0)^\perp=\overline{\ph(X)},
\end{align}
and if $X$ is pseudocomplete, then
\begin{align}\label{eq3.2}
& \forall B\subseteq Y \quad \ph^{-1}(B)^\circ_{X^\star}=\overline{\ph^\star(B^\circ_{Y^\star})}
&& \ph^{-1}(0)=\Big(\overline{\ph^\star(Y^\star)}\Big)^\perp
&& \ph^{-1}(0)^\perp=\overline{\ph^\star(Y^\star)}
\end{align}
\elm

\paragraph{The map $i_X:X\to X^{\star\star}$.}\label{SUBSEC:i_X}

The {\it second dual space} $X^{\star\star}$ for a locally convex space $X$ is the space dual to the first dual:
$$
  X^{\star\star}=(X^\star)^\star
$$
(each star $\star$ means that we take the topology of uniform convergence on totally bounded sets). The formula
$$
     i_X(x)(f)=f(x)
$$
defines a natural map $i_X:X\to X^{\star\star}$.

\bit{
 \item[$\bullet$]\label{DEF:otkr-otobr}
Let us say that a linear map of locally convex spaces $\ph:X\to Y$ is {\it open}\footnote{We use the notion of {\it open map} in the sense different from the one used in General topology \cite{Engelking}.}, if the image $\ph(U)$ of any neighborhood of zero  $U\subseteq X$ is a neighborhood of zero in the subspace $\ph(X)$ of $Y$ (with the topology inherited from $Y$):
$$
\forall
U\in \mathcal{U}(X) \quad \exists V\in \mathcal{U}(Y) \quad \ph(U)\supseteq
\ph(X)\cap V.
$$
Certainly, it is sufficient here to assume that $U$ is open and absolutely convex. By the obvious formula
\beq\label{ph(X)-cap-V=ph(ph^(-1)(V))}
\ph(X)\cap V=\ph\Big(\ph^{-1}(V)\Big),\qquad V\subseteq Y,
\eeq
(valid for any map of sets $\ph:X\to Y$ and for any subset $V\subseteq Y$), this condition can be rewritten as follows:
$$
\forall
U\in \mathcal{U}(X) \quad \exists V\in \mathcal{U}(Y) \quad \ph(U)\supseteq
\ph\Big(\ph^{-1}(V)\Big).
$$
 }\eit

\btm \label{th-2.8} For each LCS $X$ the map $i_X:X\to X^{\star\star}$ is injective, open and has dense set of values in $X^{\star\star}$.
\etm

\btm \label{th-2.12} For an arbitrary LCS $X$ the following conditions are equivalent:
\begin{itemize}
\item[(i)] the space $X$ is pseudocomplete;

\item[(ii)] the map $i_X:X\to X^{\star\star}$ is surjective (and hence, bijective).
\end{itemize}
\etm

\btm \label{th-2.14} For an arbitrary LCS $X$ the following conditions are equivalent:
\begin{itemize}
\item[(i)] the space $X$ is pseudosaturated;

\item[(ii)] the map $i_X:X\to X^{\star\star}$ is continuous.
\end{itemize}
\etm

\btm\label{TH:X-psevdopolno=>X^star-psevdonasysheno} For an arbitrary LCS $X$
\bit{
\item[---] if $X$ is pseudocomplete, then $X^\star$ is pseudosaturated,

\item[---] if $X$ is pseudosaturated, then $X^\star$ is pseudocomplete.
}\eit
\etm

\subsection{Variations of openness and closure}

\paragraph{Open and closed morphisms.}

In the stereotype theory the condition dual to the property of openness defined on page \pageref{DEF:otkr-otobr} is the following.

 \bit{
 \item[$\bullet$]\label{DEF:zamkn-otobr}
Let us say that a linear continuous map of locally convex spaces $\ph:X\to Y$ is {\it closed}, if for any totally bounded set $T\subseteq \overline{\ph(X)}\subseteq Y$ there is a totally bounded set $S\subseteq X$ such that $T\subseteq\ph(S)$. Certainly, this means in particular that the set of values $\ph(X)$ of $\ph$ is closed in $Y$.
 }\eit

\btm\label{TH:otkrytost<->zamknutost} For a linear continuous map $\ph:X\to Y$ of locally convex spaces
 \bit{
\item[(a)] if $X$ pseudosaturated, $Y$ is pseudocomplete and $\ph:X\to Y$ is open, then $\ph^\star:Y^\star\to X^\star$ is closed;

 \item[(b)] if $Y$ is pseudocomplete and $\ph:X\to Y$ is closed, then $\ph^\star:Y^\star\to X^\star$ is open.
 }\eit
\etm

For proof we need the following

\blm\label{LM:sup_(x-in-T-cap-X)|f(x)|<1}
Let $X$ be a closed subspace in a LCS $Y$, $T$ an absolutely convex compact set in $Y$, and $f:X\to\C$ a linear continuous functional such that
\beq\label{sup_(x-in-T-cap-X)|f(x)|<1}
\sup_{x\in T\cap X}\abs{f(x)}<1.
\eeq
Then there exists a linear continuous extension $g:Y\to\C$ of $f$ such that
\beq\label{sup_(y-in-T)|g(y)|<1}
\sup_{y\in T}\abs{g(y)}<1.
\eeq
\elm
\bpr
Take $\e>0$ such that
\beq\label{sup_(x-in-T-cap-X)|f(x)|<1-e}
\sup_{x\in T\cap X}\abs{f(x)}<1-\e.
\eeq
Since $f$ is continuous on $X$, the set $Z=\{x\in X:\ \abs{f(x)}\ge 1-\e\}$ is closed in $X$, and in $Y$ as well. On the other hand, by \eqref{sup_(x-in-T-cap-X)|f(x)|<1-e}, $Z$ is disjoint with $T$. As a corollary, there is an absolutely convex closed neighbourhood of zero $V$ in $Y$ such that
$$
Z\cap(T+V)=\varnothing.
$$
This means, in particular, that
$$
\sup_{x\in (T+V)\cap X}\abs{f(x)}<1-\e.
$$
If we denote by $p$ the Minkowski functional of the set $T+V$ (which is a closed absolutely convex neighbourhood of zero in $Y$), we obtain
$$
\abs{f(x)}\le (1-\e)\cdot p(x),\qquad x\in X.
$$
By the Hahn-Banach theorem there is a linear continuous extension $g:Y\to\C$ of $f$ such that
$$
\abs{g(y)}\le (1-\e)\cdot p(y),\qquad y\in Y.
$$
On the set $T+V$ we have
$$
\sup_{y\in T+V}\abs{g(y)}\le (1-\e)\cdot \sup_{y\in T+V}p(y)=1-\e<1.
$$
\epr

\bpr[Proof of Theorem \ref{TH:otkrytost<->zamknutost}.]
1. Let $\ph:X\to Y$ be open. Take a totally bounded set $F\in\mathcal{BS}(\overline{\ph^\star(Y^\star)})$, i.e.  $F\in\mathcal{BS}(X^\star)$ and $F\subseteq \overline{\ph^\star(Y^\star)}$.
By Lemma \ref{th-2.1}(d), the polar $U={^\circ F}$ is a capacious set in $X$, and since $X$ is pseudosaturated,  $U={^\circ F}$ is a neighbourhood of zero in $X$. Therefore, since $\ph$ is open, there must exist a neighbourhood of zero $V\in\mathcal{BU}(Y)$ such that $\ph(U)\supseteq\ph(X)\cap V$. By Lemma \ref{th-2.1}(b), the polar $G=V^\circ$ is a totally bounded set in $Y^\star$. Let us show that $F\subseteq\ph^\star(G)$.

Take an arbitrary functional $f\in F$ and show that there exists $g\in G$ such that $f=\ph^\star(g)$. Since $Y$ is pseudocomplete, we have $F\subseteq \overline{\ph^\star(Y^\star)}=\eqref{eq3.2}=\ph^{-1}(0)^\perp$, so  $\ph^{-1}(0)\subseteq f^{-1}(0)$. Therefore $f$ can be represented as a composition
$$
f=h\circ\ph,
$$
where $h$ is a (uniquely defined) functional on $\ph(X)$ (and we need to prove that $h$ is continuous). We have:
$$
1\ge \sup_{x\in U}\abs{f(x)}=\sup_{x\in U}\abs{h\big(\ph(x)\big)}=\sup_{y\in\ph(U)}\abs{h(y)}\ge \sup_{y\in\ph(X)\cap V}\abs{h(y)},
$$
i.e. $h$ is bounded by identity on the intersection of the unit ball $V$ of the seminorm $p(y)=\inf\{\lambda>0:\ y\in\lambda\cdot V\}=\sup_{g\in G}\abs{g(y)}$ with the subspace $\ph(X)$, where $h$ is defined. In other words, $h$ is subordinated to the seminorm $p$ on the subspace $\ph(X)$. By the Hahn-Banach theorem $h$ can be extended to some linear continuous functional $g\in Y^\star$, also subordinated to $p$, and as a corollary, lying in $V^\circ=G$:
$$
h=g\Big|_{\ph(X)},\qquad g\in G.
$$
Since on the set $\ph(X)$ the functionals $h$ and $g$ coincide, we have
$$
f=h\circ\ph=g\circ\ph=\ph^\star(g),\qquad g\in G.
$$

2. Suppose $Y$ is pseudocomplete and $\ph:X\to Y$ is closed. Consider a basis open neighbourhood of zero $V$ in $Y^\star$, i.e. the set of the form
$$
V=\{g\in Y^\star:\ \sup_{y\in T}\abs{g(y)}<1\},
$$
where $T$ is a convex balanced compact set in $Y$ (since $Y$ is pseudocomplete, each closed totally bounded set in $Y$ is compact). The map $\ph$ is closed, hence there is a totally bounded set $S\in\mathcal{BS}(X)$ such that  $\ph(S)\supseteq T\cap\overline{\ph(X)}$. Put
$$
U=\{f\in X^\star:\ \sup_{x\in S}\abs{f(x)}<1\}.
$$
If $f\in U\cap\ph^\star(Y^\star)$, then $\sup_{x\in S}\abs{f(x)}<1$ and $f=\ph^\star(g)=g\circ\ph$ for some $g\in Y^\star$. Put $h=g|_{\overline{\ph(X)}}$. Then
$$
\sup_{y\in T\cap \overline{\ph(X)}}\abs{h(y)}\le \sup_{y\in \ph(S)}\abs{h(y)}=\sup_{x\in S}\abs{h(\ph(x))}
=\sup_{x\in S}\abs{f(x)}<1
$$
By Lemma \ref{LM:sup_(x-in-T-cap-X)|f(x)|<1}, there is an extension $h'\in Y^\star$ of $h$ such that
$$
\sup_{y\in T}\abs{h'(y)}<1.
$$
This means that $h'\in V$, and we obtain that $f=\ph^\star(h')\in\ph^\star(V)$. So we have proved the inclusion
$$
U\cap\ph^\star(Y^\star)\subseteq\ph^\star(V).
$$
\epr

\paragraph{Weakly open and weakly closed morphisms.} Here we consider weakenings of the properties defined above.
 \bit{
 \item[$\bullet$]\label{DEF:slabo-otkr-morf}
Let us say that a linear continuous map of locally convex spaces $\ph:X\to Y$ is  {\it weakly open}, if it satisfies the following equivalent conditions:
 \bit{
\item[(i)] each functional $f\in X^\star$ vanishing on the kernel of the map $\ph$,
 $$
 f|_{\Ker\ph}=0,
 $$
can be extended along $\ph$ to a functional $g\in Y^\star$:
 $$
 f=g\circ\ph,
 $$

 \item[(ii)] the image $\ph(U)$ of any $X^\star$-weak neighbourhood of zero $U\subseteq X$ is a $Y^\star$-weak neighbourhood of zero in the subspace $\ph(X)$ of $Y$ (with the topology induced from $Y$):
\beq\label{DEF:slabo-otkrytye-morfizmy}
\forall
U\in \mathcal{U}(X_w) \quad \exists V\in \mathcal{U}(Y_w) \quad \ph(U)\supseteq
\ph(X)\cap V
\eeq
(here $X_w$ denotes the space $X$ with the $X^\star$-weak topology, and, similarly, $Y_w$),

 \item[(iii)] the image $\ph(U)$ of any $X^\star$-weak neighbourhood of zero $U\subseteq X$ is a neighbourhood of zero (not necessarily $Y^\star$-weak) in the subspace$\ph(X)$ of $Y$:
\beq\label{DEF:slabo-otkrytye-morfizmy-0}
\forall
U\in \mathcal{U}(X_w) \quad \exists V\in \mathcal{U}(Y) \quad \ph(U)\supseteq
\ph(X)\cap V.
\eeq
}\eit
}\eit
\bpr[Proof of equivalence.]
1. (i)$\Rightarrow$(ii). Suppose (i) holds, and let $U$ be an $X^\star$-weak neighbourhood of zero in $X$. The set  $\widetilde{U}=U+\ph^{-1}(0)$ is also an $X^\star$-weak neighbourhood of zero in $X$, and in addition it has the following property:
$$
\ph(U)=\ph(\widetilde{U}),\qquad \widetilde{U}+\ph^{-1}(0)=\widetilde{U}.
$$
From the second equality it follows that $\widetilde{U}$ contains the polar $^\circ\{f_1,...,f_k\}$ of some finite sequence of functionals $f_i\in X^\star$ such that $\ph^{-1}(0)\subseteq f_i^{-1}(0)$. By (i), each $f_i$ can be extended to some functional $g_i\in Y^\star$:
$$
f_i=g_i\circ\ph.
$$
Put $V={^\circ\{g_1,...,g_k\}}$, then we have \eqref{DEF:slabo-otkrytye-morfizmy-0}:
\begin{multline*}
y\in \ph(U)=\ph(\widetilde{U})\quad\Longleftarrow\quad y\in \ph({^\circ\{f_1,...,f_k\}})
\quad\Longleftrightarrow\quad \exists x\in {^\circ\{f_1,...,f_k\}}\quad y=\ph(x) \quad\Longleftrightarrow\\ \Longleftrightarrow\quad
\exists x\in X\quad y=\ph(x) \ \& \ \sup_i\abs{f_i(x)}\le 1
\quad\Longleftrightarrow\quad
\exists x\in X\quad y=\ph(x) \ \& \ \sup_i\abs{g_i(\ph(x))}\le 1
\quad\Longleftrightarrow\\ \Longleftrightarrow\quad
\exists x\in X\quad y=\ph(x) \ \& \ \sup_i\ \abs{g_i(y)}\le 1
\quad\Longleftrightarrow\quad
y\in\ph(X) \ \& \ y\in V
\quad\Longleftrightarrow\quad
y\in \ph(X)\cap V.
\end{multline*}

2. The implication (ii)$\Rightarrow$(iii) is obvious.

3. Let us prove (iii)$\Rightarrow$(i). Let $f\in X^\star$ be a functional such that $\Ker\ph\subseteq\Ker f$. Its polar $U={^\circ f}$ is an $X^\star$-weak neighvourhood of zero in $X$, so $\ph(U)\supseteq\ph(X)\cap V$ for some neighbourhood of zero $V$ in $Y$. This means that $f$ can be extended to a functional $h$ on $\ph(X)$, which is bounded on the set $\ph(X)\cap V$:
$$
f=h\circ\ph,\qquad \sup_{y\in \ph(X)\cap V}\abs{h(y)}\le 1.
$$
Hence, $h$ is a continuous functional on $\ph(X)$ (with respect to the topology induced from $Y$). By the Hahn-Banach theorem it can be extended to a functional $g\in Y^\star$, and we have $f=h\circ\ph=g\circ\ph$.
\epr

 \bit{
 \item[$\bullet$]\label{DEF:slabo-zamkn-morf}
Let us say that a linear continuous map $\ph:X\to Y$ is {\it weakly closed}, if its set of values $\ph(X)$ is closed in $Y$:
    $$
    \overline{\ph(X)}^Y=\ph(X).
    $$
 }\eit

 \bprop\label{PROP:otkryt=>slabo-otkryt} For a linear continuous map of locally convex spaces $\ph:X\to Y$
 \bit{
 \item[---] if $\ph$ is open, then $\ph$ is weakly open,

 \item[---] if $\ph$ is closed, then $\ph$ is weakly closed.
  }\eit
 \eprop
\bpr The first part of this proposition follows from condition (iii) in the definition of weak openness on p.\pageref{DEF:slabo-otkr-morf}, and the second part is obvious, and we already noticed this when we defined  closure on p.\pageref{DEF:zamkn-otobr}.
\epr

\btm\label{TH:slabaya-otkrytost<->slabaya-zamknutost} For a linear continuous map $\ph:X\to Y$ of locally convex spaces
 \bit{
 \item[(a)] $\ph:X\to Y$ is weakly open $\Longleftrightarrow$ $\ph^\star:Y^\star\to X^\star$ is weakly closed,

 \item[(b)] if $Y$ pseudosaturated and $\ph:X\to Y$ is weakly closed, then $\ph^\star:Y^\star\to X^\star$ is weakly open.
 }\eit
\etm
\bpr
Here the first propisition is exactly the equivalence of conditions (i) and (ii) in the definition of openness on p.\pageref{DEF:slabo-otkr-morf}. Let us prove the second one. Suppose $Y$ is pseudosaturated and $\ph:X\to Y$ is weakly closed. By proposition (a) (that we have already duscussed), for proving that $\ph^\star:Y^\star\to X^\star$ is weakly open, it is sufficient to verify that the second dual map $\ph^{\star\star}:X^{\star\star}\to Y^{\star\star}$ is closed. Take $h\in\overline{\ph^{\star\star}(X^{\star\star})}$. Since $Y$ is pseudosaturated, by Theorem  \ref{TH:X-psevdopolno=>X^star-psevdonasysheno} $Y^\star$ is pseudocomplete. Therefore,
$$
h\in \overline{\ph^{\star\star}(X^{\star\star})}=\eqref{eq3.2}=(\ph^\star)^{-1}(0)^\perp=
\eqref{eq3.1}=\big(\overline{\ph(X)}\big)^{\perp\perp}=i_Y\big(^\perp(\overline{\ph(X)}^\perp)\big)
$$
(the last equality means that the map $i_Y:Y\to Y^{\star\star}$, bijective by Theorem \ref{th-2.12}, turns the annihilator of the space $\overline{\ph(X)}^\perp$, meant as a subspace in $Y$ to its annihilator, meant as a subspace in $Y^{\star\star}$). This in its turn means that there is $y\in\overline{\ph(X)}$ such that $h=i_Y(y)$. Since $\ph$ is weakly closed, there exists $x\in X$ such that $y=\ph(x)$. If we denote $g=i_X(x)$, then
$$
h=i_Y(y)=i_Y(\ph(x))=\ph^{\star\star}(i_x(x))=\ph^{\star\star}(g).
$$
\epr

\paragraph{Relatively open and relatively closed morphisms.} Another weakening of openness and closure of morphisms is the following.
 \bit{
 \item[$\bullet$]
We say that a linear continuous map of locally convex spaces $\ph:X\to Y$ is
 \bit{
\item[---] {\it relatively open}, if for each neighborhood of zero $U$ in $X$ (without loss of generality we may assume that $U$ is closed and absolutely convex) such that every functional $f\in X^\star$ bounded on $U$ can be extended along the map $\ph$ to some functional $g\in Y^\star$,
    \beq\label{uslovno-otkr-otobrazh}
    \forall f\in X^\star \qquad \Big(\sup_{x\in U}\abs{f(x)}<\infty\quad\Longrightarrow\quad\exists g\in Y^\star\quad f=g\circ\ph\Big),
    \eeq
    its image $\ph(U)$ is a neighborhood of zero in the subspace $\ph(X)$ of the locally convex space $Y$ (with the topology inherited from $Y$):
$$
\ph(U)\supseteq V\cap\ph(X)
$$
for some neighborhood of zero $V$ in $Y$;

\item[---] {\it relatively closed}, if for each absolutely convex compact set $T\subseteq Y$, if $T$ contains in $\ph(X)$, then there is a compact set $S\subseteq X$ such that $T\subseteq \ph(S)$.

 }\eit
 }\eit

The following is obvious:

 \bprop\label{PROP:otkryt=>uslov-otkryt} For a morphism of locally convex spaces $\ph:X\to Y$
 \bit{
 \item[---] if $\ph$ is open, then $\ph$ is relatively open,

 \item[---] if $\ph$ is closed, then $\ph$ is relatively closed.
  }\eit
 \eprop

\btm\label{TH:uslov-otkrytost<->uslov-zamknutost} For a linear continuous map $\ph:X\to Y$ of locally convex spaces
 \bit{
 \item[(a)] $\ph:X\to Y$ is relatively open $\Longleftrightarrow$ $\ph^\star:Y^\star\to X^\star$ is relatively closed;

 \item[(b)] if $X$ is pseudocomplete and $\ph:X\to Y$ is relatively closed, then $\ph^\star:Y^\star\to X^\star$ is relatively open.
 }\eit
\etm
\bpr
a. Suppose $\ph$ is relatively open, and $T$ is a closed absolutely convex totally bounded set in $X^\star$, contained in $\ph^\star(Y^\star)$:
\beq\label{uslovno-otkr-otobrazh-1}
  \forall f\in T \quad\exists g\in Y^\star\quad f=\ph^\star(g)=g\circ\ph.
\eeq
For the polar $U={^\circ T}$ this means condition \eqref{uslovno-otkr-otobrazh}, and, since $U$ is a neighbourhood of zero in $X$, we obtain that the image $\ph(U)$ must be a neighbourhood of zero in the subspace $\ph(X)$ of $Y$ (with the topology induced from $Y$). I.e., there exists a neighbourhood of zero $V$ in $Y$ such that
$$
\ph(U)\supseteq V\cap\ph(X).
$$
Clearly, $V$ can be chosen as closed and absolutely convex in $Y$. Put $S=V^\circ$ and show that $T\subseteq \ph^\star(S)$, i.e.
\beq\label{uslovno-otkr-otobrazh-2}
\forall f\in U^\circ\qquad \exists h\in V^\circ\qquad f=\ph^\star(h)=h\circ\ph.
\eeq
Indeed, take $f\in T=U^\circ$, then by \eqref{uslovno-otkr-otobrazh-1} one can choose $g\in Y^\star$ such that $f=g\circ\ph$. The restriction $g\big|_{\ph(X)}$ of this functional $g$ on $\ph(X)$ is bounded by the isentity on the neighbourhood of zero $V\cap\ph(X)$:
$$
\sup_{y\in V\cap\ph(X)}\abs{g(y)}\le\sup_{y\in \ph(U)}\abs{g(y)}\le \sup_{x\in U}\abs{g(\ph(x))}=\sup_{x\in U}\abs{f(x)}\le 1.
$$
In other words, the functional $g\big|_{\ph(X)}$ on the subspace $\ph(X)$ is subordinated to the seminorm
$$
p(y)=\inf\{\lambda>0:\ y\in\lambda\cdot V\}.
$$
By the Hahn-Banach theorem, the functional $g\big|_{\ph(X)}$ can be extended to some functional $h$ on $Y$, subordinated to $p$:
$$
\abs{h(y)}\le p(y)\quad (y\in Y),\qquad h|_{\ph(X)}=g.
$$
Here from the first condition it follows that $\sup_{y\in V}\abs{h(y)}\le \sup_{y\in V}p(y)\le 1$,
i.e. $h\in V^\circ=S$. And from the second one, that $h(\ph(x))=g(\ph(x))=f(x)$. Together this means  \eqref{uslovno-otkr-otobrazh-2}.

The reverse implication. Let $\ph^\star:Y^\star\to X^\star$ be relatively closed and $U$ be an absolutely convex neighbourhood of zero in $X$ satisfying \eqref{uslovno-otkr-otobrazh}. Consider the polar $T=U^\circ$. This is a closed absolutely convex totally bounded set in $X^\star$, and for it the condition \eqref{uslovno-otkr-otobrazh} is equivalent to \eqref{uslovno-otkr-otobrazh-1}. This in its turn means $T\subseteq \ph^\star(Y^\star)$,
and since $\ph^\star$ is relatively closed, there must exist an absolutely convex totally bounded set $S\subseteq Y^\star$ such that
$$
T\subseteq\ph^\star(S)
$$
Hence
$$
{^\circ T}=\{x\in X: \ \forall f\in T\quad \abs{f(x)}\le 1\}\supseteq \{x\in X: \ \forall g\in S\quad \abs{g(\ph(x))}\le 1\}=\{x\in X: \ \ph(x)\in{^\circ S}\}=\ph^{-1}\big({^\circ S}\big).
$$
Now if we put $V={^\circ S}$ (this is a neighbourhood of zero in $Y$), then we obtain
$$
U\supseteq\ph^{-1}(V)\quad\Longrightarrow\quad
\ph(U)\supseteq\ph\big(\ph^{-1}(V)\big)=\eqref{ph(X)-cap-V=ph(ph^(-1)(V))}=\ph(X)\cap V.
$$

b. Suppose $X$ is pseudocomplete, $\ph:X\to Y$ is relatively closed and $U$ an absolutely convex neighbourhood of zero in $Y^\star$, satisfying \eqref{uslovno-otkr-otobrazh}, i.e. in this case
$$
    \forall \upsilon\in Y^{\star\star} \qquad \Big(\sup_{g\in U}\abs{\upsilon(g)}<\infty\quad\Longrightarrow\quad\exists \xi\in X^{\star\star}\quad \upsilon=\xi\circ\ph^\star\Big).
$$
In particular, for any $y\in T={^\circ U}$ there must exist a functional $\xi\in X^{\star\star}$ such that
$$
i_Y(y)=\xi\circ\ph^\star.
$$
Since $X$ is pseudocomplete, by Theorem \ref{th-2.12} there exists a point $x\in X$ such that $i_X(x)=\xi$. Then
$$
\forall g\in Y^\star\qquad g(y)=i_Y(y)(g)=(i_X(x)\circ\ph^\star)(g)=i_X(x)(\ph^\star(g))=\ph^\star(g)(x)=g(\ph(x)),
$$
and therefore $y=\ph(x)$. We have proved that $T\subseteq \ph(X)$, and since $\ph$ is relatively closed, there must exist an absolutely convex totally bounded set $S\subseteq X$ such that
$$
T\subseteq\ph(S)
$$
We have for it
$$
T^\circ\supseteq\Big(\ph(S)\Big)^\circ=\eqref{eq3.1}=(\ph^\star)^{-1}\big(S^\circ\big)
$$
Now if we put $V=S^\circ$ (this is a neighbourhood of zero in $X^\star$), then
$$
U\supseteq(\ph^\star)^{-1}(V)\quad\Longrightarrow\quad
\ph^\star(U)\supseteq\ph^\star\big((\ph^\star)^{-1}(V)\big)=\eqref{ph(X)-cap-V=ph(ph^(-1)(V))}=\ph^\star(Y^\star)\cap V.
$$
This is what we need.
\epr

\paragraph{Connections between the three variations of openness and closure.}
Propositions \ref{PROP:otkryt=>slabo-otkryt} and \ref{PROP:otkryt=>uslov-otkryt} can be strengthened as follows.

\btm For a morphism of locally convex spaces $\ph:X\to Y$
 \bit{
 \item[(a)] $\ph$ is open $\Longleftrightarrow$ $\ph$ is weakly open and relatively open;

 \item[(b)] $\ph$ is closed $\Longleftrightarrow$ $\ph$ is weakly closed and relatively closed.
 }\eit
\etm
\bpr
To the right direction this was already noticed in Propositions \ref{PROP:otkryt=>slabo-otkryt} and \ref{PROP:otkryt=>uslov-otkryt}, so we must check the reverse implications.

a. Let $\ph$ be weakly open and relatively open. For each neighbourhood of zero $U$ in $X$ the set $U+\ph^{-1}(0)$ is also a neighbourhood of zero in $X$. If a functional $f\in X^\star$ is bounded on $U+\ph^{-1}(0)$, then on the subspace $\ph^{-1}(0)$ it must vanish, $f|_{\ph^{-1}(0)}=0$, so by the weak openness of $\ph$, $f$ can be extended to a functional $g\in Y^\star$. This means that the neighbourhood of zero $U+\ph^{-1}(0)$ satisfies the condition  \eqref{uslovno-otkr-otobrazh}. Since $\ph$ is relatively open, we have
$$
\ph(U)=\ph(U+\ph^{-1}(0))\supseteq \ph(X)\cap V
$$
for some neighbourhood of zero $V$ in $Y$.

b. If $\ph$ is weakly closed and relatively closed, then, first, $\overline{\ph(X)}=\ph(X)$, and, second, each closed absolutely convex totally bounded set $T\subseteq\ph(X)$ is an image of some totally bounded set $S\subseteq X$ under tha map $\ph$. Together this means that $T$ can be chosen as a subset in $\overline{\ph(X)}$, and the same will be true. This means that $\ph$ is closed.
\epr

\paragraph{Embeddings and coverings.}

\begin{itemize}
\item
A linear continuous map $\ph:X\to Y$ of locally convex spaces will be called
\begin{itemize}
\item[--] an {\it embedding} (respectively, a {\it weak embedding}, a {\it relative embedding}), if it is injective and open (respectively, weakly open, relatively open),

\item[--] a {\it dense embeddding} (respectively, a {\it dense weak embedding}, a {\it dense relative embedding}), if in addition the set of values $\ph(X)$ is dense in $Y$.

\item[--] a {\it covering} (respectively, a {\it weak covering}, a {\it relative covering}), if it is surjective and closed (respectively, weakly closed, relatively closed),

\item[--] an {\it exact covering} (respectively, an {\it exact weak covering}, an {\it exact relative covering}), if in addition it is injective.
\end{itemize}
\end{itemize}

\brem
If a LCS $X$ is pseudocomplete and $\ph:X\to Y$ is an exact covering, then for any totally bounded set $S\subseteq X$ the restriction $\ph|_S:S\to\ph(S)$ is a homeomorphism of topological spaces.
\erem

\bex\label{cor-2.13} If a locally convex space $X$ is pseudocomplete, then a
(continuous and bijective) map $i_X^{-1}:X^{\star\star} \to X$ is defined, and it is an exact covering.
\eex

\bex\label{cor-2.15} If a locally convex space $X$ is pseudosaturated, then the map $i_X:X\to X^{\star\star}$ is a dense embedding.
\eex

The following proposition is proved in \cite[Theorems 3.2, 3.1]{Akbarov}.

\btm For a linear continuous map $\ph:X\to Y$ of locally convex spaces
 \bit{

\item[--] if $X$ is pseudosaturated and $\ph:X\to Y$ is a dense embedding, then $\ph^\star:Y^\star\to X^\star$ is an exact covering,

\item[--] if $X$ is pseudocomplete and $\ph:X\to Y$ is an exact covering, then $\ph^\star:Y^\star\to X^\star$ is a dense embedding.
 }\eit
\etm

\subsection{Pseudocompletion and pseudosaturation.}

\paragraph{Pseudocompletion.}

Like in the case of completeness, each locally convex space $X$ has pseudocompletion, i.e. the ``outside-nearest'' pseudocomplete space. Formally this construction is described in the following

\btm \label{th-1.5} There exists a map $X\mapsto\triangledown_X$ that assigns to each locally convex space $X$ a linear continuous map ${\triangledown}_X :X \to X^{\triangledown}$ into a pseudocomplete locally convex space
$X^{\triangledown}$ in such a way that the following conditions are fulfilled:
\begin{itemize}
\item[(i)] $X$ is pseudocomplete if an only if ${\triangledown}_X :X\to X^{\triangledown}$ is an isomorphism;

\item[(ii)] for any linear continuous map $\ph :X\to Y$ of locally convex spaces there is a unique linear continuous map $\ph^{\triangledown} :X^{\triangledown} \to Y^{\triangledown}$ such that the following diagram is commutative:
\beq
\begin{diagram}
\node{X} \arrow{e,t}{\triangledown_X} \arrow{s,l}{\varphi}
\node{X^\triangledown} \arrow{s,r,--}{\varphi^\triangledown}
\\
\node{Y} \arrow{e,t}{\triangledown_Y} \node{Y^\triangledown}
\end{diagram}.
\label{eqI.12}
\eeq \end{itemize}
\etm

From $(i),(ii)$ it follows that for any linear continuous map $\ph :X\to Y$ into a pseudocomplete space $Y$ there exists a unique linear continuous map $X^{\triangledown} \to Y$ such that the following diagram is commutative:
\beq\label{eqI.13}
\begin{diagram}
\node{X} \arrow[2]{e,t}{\triangledown_X} \arrow{se,b}{\varphi}
\node[2]{X^\triangledown} \arrow{sw,--}{}
\\
\node[2]{Y}
\end{diagram}.
\eeq
This means by the way that, the morphism $\triangledown_X:X\to X^\triangledown$ is an extension of $X$ in $\Ob({\tt LCS})$ with respect to the object $\C$. Since $\C$ differs morphisms on the outside in ${\tt LCS}$, by Theorem \ref{TH:M-razdel-moprfizmy}, $\triangledown_X:X\to X^\triangledown$ is a bimorphism. This implies in its turn that the morphism  $\triangledown_X:X\to X^\triangledown$ is unique up to an isomorphism in $\Epi^X$.

 \bit{
 \item
The space $X^{\triangledown}$ is called the {\it pseudocompletion}, and the map ${\triangledown}_X :X \to
X^{\triangledown}$ the {\it pseudocompletion map} of the locally convex space $X$. From $(ii)$, it follows also, that the map $\ph \mapsto \ph ^{\triangledown}$ is a covariant functor of the category ${\tt LCS}$ into itself: $(\psi \circ \ph)^{\triangledown}=\psi^{\triangledown} \circ \ph^{\triangledown}$. We call it the {\it pseudocompletion functor}.
 }\eit

\btm \label{th-1.6} For any locally convex space $X$ the pseudocompletion map ${\triangledown}_X :X \to X^{\triangledown}$ is a dense embedding.
\etm

Like usual completion, the operation of pseudocompletion $X\mapsto X^\triangledown$ adds new elements to $X$, but does not change the topology of $X$.

\paragraph{Pseudosaturation.}
It is remarkable that there exists a dual construction, that assigns to each locally convex space $X$ an ``inside-nearest'' pseudosaturated locally convex space $X^\vartriangle$:

\btm \label{th-1.16} There exists a map $X\mapsto\vartriangle_X$, that assigns to each locally convex space $X$ a linear continuous map $\vartriangle _X :X^\vartriangle \to X$ from a pseudosaturated locally convex space $X^\vartriangle$ in such a way that the following conditions are fulfilled:
\begin{itemize}
\item[(i)] $X$ is pseudosaturated if and only if $\vartriangle _X:X^\vartriangle \to X$ is an isomorphism;

\item[(ii)] for any linear continuous map $\ph:Y\to X$ of locally convex spaces there is a unique linear continuous map $\ph^\vartriangle :Y^\vartriangle\to X^\vartriangle$ such that the following diagram is commutative:
\beq
\begin{diagram}
\node{X} \node{X^\vartriangle} \arrow{w,t}{\vartriangle_X}
\\
\node{Y} \arrow{n,l}{\varphi} \node{Y^\vartriangle} \arrow{w,t}{\vartriangle_Y}
\arrow{n,r,--}{\varphi^\vartriangle}
\end{diagram}.
\label{eqI.25}
\eeq \end{itemize}
\etm

From $(i),(ii)$ it follows that for any linear continuous map $\ph :Y\to X$ from a pseudosaturated locally convex space $Y$ there is a unique linear continuous map $Y\to X^{\vartriangle}$ such that the following diagram is commutative:
\beq\label{eqI.26}
\begin{diagram}
\node{X} \node[2]{X^\vartriangle} \arrow[2]{w,t}{\vartriangle_X}
\\
\node[2]{Y} \arrow{nw,b}{\varphi} \arrow{ne,b,--}{}
\end{diagram}.
\eeq
This means by the way that the morphism $\vartriangle_X:X^\vartriangle\to X$ is an enrichment of $X$ in the class $\Ob({\tt LCS})$ by means of the object $\C$. Since $\C$ differs morphisms on the inside in ${\tt LCS}$, by Theorem \ref{TH:M-razdel-moprfizmy-iznutri}, $\vartriangle_X:X^\vartriangle\to X$ is a bimorphism. This implies in its turn that the morphism $\vartriangle_X:X^\vartriangle\to X$ is unique up to an isomorphism in $\Mono_X$.

\bit{\label{DEF:pseudosaturation}
\item  The space $X^\vartriangle$ is called the {\it pseudosaturation}, and the map $\vartriangle _X :X^\vartriangle \to X$ the {\it pseudosaturation map} of the space $X$. From $(ii)$ it follows that the map $\ph \mapsto \ph ^{\vartriangle}$ is a covariant functor of the category ${\tt LCS}$ into itself: $(\psi\circ\ph)^{\vartriangle}=\psi^{\vartriangle} \circ \ph^{\vartriangle}$. We call it {\it pseudosaturation functor}.
}\eit

\btm\label{th-1.17} For any locally convex space $X$ the pseudosaturation map $\vartriangle _X :X^\vartriangle \to X$ is an exact covering.
\etm

The pseudosaturation $X^\vartriangle$ can be imagined as a new, stronger topologization of the space $X$, which preserves the system of totally bounded sets and the topology on each totally bounded set in $X$.

Each of the operations $X\mapsto X^\triangledown$ and $X\mapsto X^\vartriangle$ preserves the properties of being pseudocomplete and pseudosaturated:

\btm\label{TH:X-psevdopolno=>X^vartriangle-psevdopolno} For a locally convex space $X$
\bit{
\item[---] if $X$ is pseudocomplete, then its pseudosaturation $X^\vartriangle$ is also pseudocomplete,

\item[---] if $X$ is pseudosaturated, then its pseudocompletion $X^\triangledown$ is also pseudosaturated.
}\eit
\etm

The following examples show that pseudocompletion and psudosaturation are independent.

\begin{ex}\label{ex-4.9} Let $X$ be an infinite-dimensional Banach space, and $Y=X'_\sigma$ its dual space with the $X$-weak topology. The space $Y$ is pseudocomplete, but not pseudosaturated.
\end{ex}

\begin{ex}\label{ex-4.10} An arbitrary non-complete metrizable locally convex space is pseudosaturated, but not pseudocomplete space.
\end{ex}

\paragraph{Duality between pseudocompletion and pseudosaturation.}

The passage to a dual space $X\mapsto X^\star$ reshuffles pseudocompleteness and pseudosaturateness:

\btm \label{th-3.14} Let $X$ be a pseudocomplete LCS. Then
\begin{itemize}
\item[(a)] there is a unique isomorphism of locally convex spaces
\beq\label{(X^vartriangle)^star=(X^star)^triangledown}
 \xymatrix @C=5.0pc{  (X^\vartriangle)^\star\ar@{~>}[r] & (X^\star)^\triangledown}
\eeq
such that the following diagram is commutative:
\beq \label{eqIII.18}
 \xymatrix  %  @R=2.5pc @C=1.0pc
 {
 (X^\vartriangle)^\star \ar@{~>}[rr] & & (X^\star)^\triangledown
\\
 & X^\star \ar[ul]^{(\vartriangle_X)^\star} \ar[ur]_{\triangledown_{X^\star}} &
 };
\eeq

\item[(b)] for any linear continuous map $\ph :X\to Y$ of locally convex spaces the following diagram is commutative:
   \beq\label{eqIII.19}
 \xymatrix  %  @R=2.5pc @C=1.0pc
 {
 (X^\vartriangle)^\star \ar@{~>}[r] & (X^\star)^\triangledown
 \\
 (Y^\vartriangle)^\star \ar@{~>}[r] \ar[u]^{(\varphi^\vartriangle)^\star}
 & (Y^\star)^\triangledown \ar[u]_{(\varphi^\star)^\triangledown}
 }
\eeq

\end{itemize}
\etm

\btm\label{th-3.15} Let $X$ be a pseudosaturated locally convex space. Then
\begin{itemize}
\item[(a)] there is a unique isomorphism of locally convex spaces
 \beq\label{(X^triangledown)^star=(X^star)^vartriangle}
 \xymatrix @C=5.0pc
 {  (X^\triangledown)^\star\ar@{~>}[r] & (X^\star)^\vartriangle}
 \eeq
such that the following dioagram is commutative
\beq \label{eqIII.20}
 \xymatrix  %  @R=2.5pc @C=1.0pc
 {
 (X^\triangledown)^\star \ar@{~>}[rr]
 \ar[dr]_{(\triangledown_X)^\star} && (X^\star)^\vartriangle
\ar[dl]^{\vartriangle_{X^\star}}
\\
 & X^\star &
 };
\eeq

\item[(b)] for any linear continuous map $\ph :X\to Y$  the following diagram is commutative:
\beq \label{eqIII.21}
 \xymatrix  %  @R=2.5pc @C=1.0pc
 {
 (X^\triangledown)^\star \ar@{~>}[r] & (X^\star)^\vartriangle
\\
 (Y^\triangledown)^\star \ar@{~>}[r]
\ar[u]^{(\varphi^\triangledown)^\star} & (Y^\star)^\vartriangle
\ar[u]_{(\varphi^\star)^\vartriangle}
 }.
\eeq \end{itemize}
\etm

\subsection{Stereotype spaces}

 \bit{
 \item[$\bullet$] A locally convex space $X$ is said to be {\it stereotype}, if its natural map to the second dual space
$$
\i_X:X\to (X^\star)^\star \quad | \quad \i_X(x)(f)=f(x), \quad x\in X, f\in
X^\star
$$
is an isomorphism of locally convex spaces (both times $\star$ means the dual space in the sense of definition on p.\pageref{DEF:X^star}).
}\eit
Certanly, if $X$ is a stereotype space, then its dual space $X^\star$ is also stereotype.
Theorems \ref{th-2.12} and \ref{th-2.14} imply the following criterion:

\btm\label{TH:kriterij-ster}
A locally convex space $X$ is stereotype if and only if it is pseudocomplete and pseudosaturated.
\etm

This mean in particular, that there are non-stereotype locally convex spaces (since there are non-pesudocomplet and non-pseudosaturated spaces, see Examples \ref{ex-4.9} and \ref{ex-4.10}). Nevertheless, the class of stereotype spaces ${\tt Ste}$ turns out to be amazingly wide. This is seen from the following series of examples, generalizing each other.

\bex All {\sl Banach spaces} are stereotype.
\eex

\bex All {\sl Fr\'{e}chet spaces} are stereotype.
\eex

\bex\label{ex-4.3} All {\sl quasicomplete barreled spaces} are stereotype.
\eex

As a corollary, the place of stereotype spaces among other frequently used classes of spaces can be illustrated by the following diagram:

 \beq\label{DIAGR:ster-prostranstva}
 {\small
 \begin{picture}(380,220)
\put(195,105){\oval(330,200)} \put(110,180){\text{\sf\large STEREOTYPE SPACES}} \put(195,90){\oval(270,140)} \put(110,130){\text{\sf\large quasicomplete barreled spaces}}
 \put(160,80){\oval(140,60)\put(-40,15){\text{\sf Fr\'{e}chet spaces}}}
 \put(160,70){\oval(100,25)\put(-40,-2){\text{\sf Banach spaces}}}
 \put(243,65){\oval(110,60)} \put(245,60){\text{\sf reflexive}}
 \put(249,50){\text{\sf spaces}}
\end{picture}
 }
 \eeq
This picture is supplemented by the examples of spaces, dual to the already mentioned, and having quite unwonted\footnote{Because of the non-standard notion of dual space.} properties:

\bex A locally convex space $X$ is called a {\it Smith space}\footnote{After M.F.Smith \cite{Smith}.}, if it is a complete $k$-space\footnote{\label{DEF:k-space}A topological space $X$ is called {\it $k$-space} or {\it Kelley space}, if every set $M\subseteq X$ having closed trace $M\cap K$ on each compact set $K\subseteq X$ is closed in $X$.} and has a {\it universal compact set}, i.e. a compact set $K\subset X$ that absorbs any other compact set $T\subset X$:  $T\subseteq \lambda K$ for some $\lambda\in \C$. It is known that {\it a locally convex space $X$ is a Smith space if and only if it is stereotype and its dual space $X^\star$ is a Banach space}.
\eex

\bex A locally convex space $X$ is called a {\it Brauner space}\footnote{After K.Brauner \cite{Brauner}.}, if it is a complete $k$-space\footnote{See footnote \ref{DEF:k-space}.} and has a {\it countable fundamental system of compact sets}, i.e. a sequence of compact sets $K_n\subseteq X$ such that every compact set $T\subseteq X$ is contained in some $K_n$. {\it A locally convex space $X$ is a Brauner space if and only if it is stereotype and its dual space $X^\star$ is a Fr\'{e}chet space}.
\eex

The connections between the spaces of Fr\'{e}chet, Brauner, Banach, and Smith are illustrated in the
following diagram (where the 180 degree rotation corresponds to the passage to the dual class):
$$
 
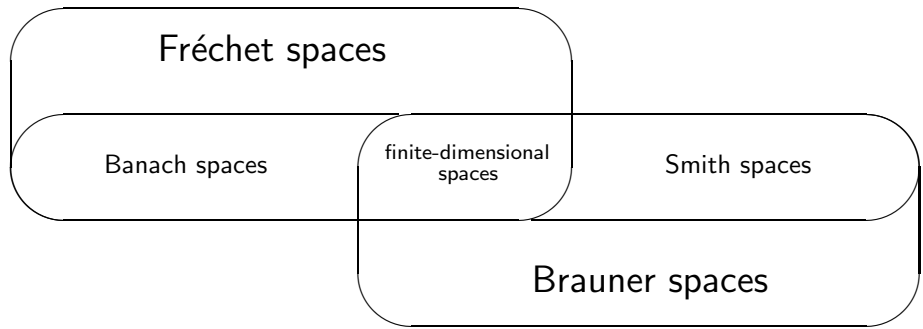
\begin{picture}(400,140)
 \put(130,90){\oval(210,80)} \put(80,110){\text{\sf\Large Fr\'{e}chet spaces}}
 \put(170,70){\oval(290,40)[l]} \put(60,68){\text{\sf Banach spaces}}
 \put(260,50){\oval(210,80)} \put(220,23){\text{\sf\Large Brauner spaces}}
 \put(220,70){\oval(290,40)[r]} \put(270,68){\text{\sf Smith spaces}}
 \put(165,73){\text{\sf\footnotesize finite-dimensional}}
 \put(185,65){\text{\sf\footnotesize spaces}}
\end{picture}
$$

It is clear from the definition that each stereotype space $X$ can be restored from its dual space $X^\star$. So different properties of $X$ have their dual analogs in $X^\star$. The most obvious facts of that type are listed in the following

\btm \label{th-4.11} Let $X$ be a stereotype space. Then
\begin{itemize}
\item[(a)] $X$ is normable $\Longleftrightarrow$  $X$ is a Banach space
$\Longleftrightarrow$ $X^\star$ is a Smith space;

\item[(b)] $X$ is metrizable $\Longleftrightarrow$  $X$ is a Fr\'{e}chet space
$\Longleftrightarrow$ $X^\star$ is a Brauner space;

\item[(c)] $X$ is barreled $\Longleftrightarrow$ $X^\star$ has the Heine-Borel psoperty;

\item[(d)] $X$ is quasibarreled $\Longleftrightarrow$ in $X^\star$ each subset $T$ absorbed by any barrel, is totally bounded;

\item[(e)] $X$ is a Mackey space $\Longleftrightarrow$ in $X^\star$ each $(X^\star)^\star$-weak compact set is compact;

\item[(f)] $X$ is a Montel space $\Longleftrightarrow$ $X$ is barreled and has the Heine-Borel property $\Longleftrightarrow$ $X^\star$ is a Montel space;

\item[(g)] $X$ is a space with a weak topology $\Longleftrightarrow$ in $X^\star$ every compact set $T$ is finite-dimensional;

\item[(h)] $X$ is separable (i.e. has a countable everywhere dense set) $\Longleftrightarrow$ in $X^\star$ there is a sequence of closed subspaces $L_n$ of finite co-dimension with trivial intersection: $\bigcap_{n=1}^{\infty}
L_n=\{0\}$.

\item[(i)] $X$ has the (classical) approximation property $\Longleftrightarrow$ $X^\star$ has the (classical) approximation property;

\item[(j)] $X$ is complete  $\Longleftrightarrow$ $X^\star$ co-complete\footnote{A locally convex space $X$ is said to be {\it co-complete}\cite{Akbarov}, if each linear functional $f:X\to \C$ continuous on each totally bounded set $S\subseteq X$, is continuous on $X$.} $\Longleftrightarrow$ $X^\star$ is saturated\footnote{\label{DEF:saturated-space} A locally convex space  $X$ is said to be {\it saturated}\cite{Akbarov}, if for an absolutely convex set $B$ being a neighbourhood of zero in $X$ is equivalent to the following: for any totally bounded set $S\subseteq X$ there is a closed neighbourhood of zero $U$ in $X$ such that $B\cap S=U$.} ;

\item[(k)] $X$ is a Pt\'{a}k space\footnote{A locally convex space $X$ is called {\it Pt\'{a}k space}
\cite{Schaefer} or {\it fully complete} \cite{Robertson-Robertson}, in the dual space $X^\star$ every subspace $Q\subseteq X^\star$ is $X$-weakly closed, when it leaves the $X$-weakly closed trace $Q\cap U^\circ$ on the polar $U^\circ$ of each neighbourhood of zero $U\subseteq X$.}
$\Longleftrightarrow$ in $X^\star$ a subspace $L$ is closed, if it leaves a closed trace $L\cap K$ on each compact set  $K\subseteq X^\star$;

\item[(l)] $X$ is hypercomplete\footnote{A locally convex space $X$ is said to be {\it hypercomplete}
 \cite{Robertson-Robertson}, if in the dual space $X^\star$ an absolutely convex set $Q\subseteq X^\star$ is $X$-weakly closed, when it leaves an $X$-weakly closed trace $Q\cap U^\circ$ on the polar $U^\circ$ of each neighbourhood of zero $U\subseteq X$.}
$\Longleftrightarrow$ in $X^\star$ an absolutely convex set $B$ is closed if it leaves a closed trace $B\cap K$ on each compact set $K\subseteq X^\star$.
\end{itemize}
\etm

\bprop \label{lm-4.12}
Let $E$ be a closed subspace in a locally convex space $X$, considered as a locally convex space with the topology induced from $X$, and let the annihilator $E^\perp$ be also endowed with the topology induced from $X^\star$. Then
\begin{itemize}
\item[(a)] there is a natural isomorphism of locally convex spaces
\beq\label{eq4.2}
  E^\star \cong X^\star/E^\perp,
\eeq
and if in addition $E$ is pseudocomplete (for example, if $X$ is pseudocomplete), then the isomorphism  \eqref{eq4.2} generates isomorphisms of stereotype spaces
\beq
  (E^\vartriangle)^\star \cong (X^\star/E^\perp)^\triangledown  \quad\quad\quad
  E^\vartriangle \cong [(X^\star/E^\perp)^\triangledown]^\star
\label{eq4.3}
\eeq
\item[(b)] if $X$ is stereotype, then there is a natural isomorphism of locally convex spaces
\beq  \label{eq4.4}
  (E^\perp)^\star \cong X/E
\eeq
generating isomorphisms of stereotype spaces
\beq
  ((E^\perp)^\vartriangle)^\star \cong (X/E)^\triangledown  \quad\quad\quad
  (E^\perp)^\vartriangle \cong [(X/E)^\triangledown]^\star
\label{eq4.5}
\eeq
\end{itemize}
\eprop

The following example is due to O.G.Smolyanov \cite{Smolyanov} and it was mentioned in \cite{Akbarov} (as Example  3.22). We will use it later as an important technical result:

\begin{ex}\label{ex-3.22} There is a stereotype space $Z$ with the following properties:
\begin{itemize}
\item[$(i)$] $Z$ and $Z^\star$ are complete and saturated\footnote{See footnote \ref{DEF:saturated-space}.};

\item[$(ii)$] $Z$ has a closed subspace $Y$ such that
\begin{itemize}
\item[$(a)$] the quotient space $Z/Y$ is metrizable, but not complete;

\item[$(b)$] the annihilator $Y^\perp$ (with the topology induced from $Z^\star$) is not a pseudosaturated space.
\end{itemize}
\end{itemize}
\end{ex}
\bpr This is the space $Z=\mathcal  D(\Bbb R)$ of smooth functions with compact support on $\Bbb R$. It is complete  (as a strong unductive limit of a sequence of complete spaces \cite{Robertson-Robertson}) and saturated (as an inductive limit of a system of saturated spaces). By Theorem \ref{TH:X-psevdopolno=>X^star-psevdonasysheno}, the dual space $Z^\star=\mathcal D^\star(\Bbb R)$ is also complete and saturated. In \cite{Smolyanov} O.G.Smolyanov showed that $Z$ contains a closed subspace $Y$, such that the quotient space $Z/Y$ is metrizable, but not complete. Hence, $Z/Y$ is not pseudocomplete.

Put $X=Z^\star, \, E=Y^\perp$. By Proposition \ref{lm-4.12}(a), $Z/Y=X^\star/E^\perp=E^\star$. So if $E$ was pseudosaturated, then $Z/Y$ would be pseudocomplete by Theorem \ref{TH:X-psevdopolno=>X^star-psevdonasysheno}. \epr

\begin{ex}\label{ex-3.23} There exists a complete locally convex space $E$ (and thus, $E$ can be represented as a projective limit of Banach spaces in the category ${\tt LCS}$), such that the dual space $E^\star$ is metrizable, but not complete. As a corollary, $E$ is not pseudosaturated, and there is a discontinuous linear functional $f:E\to\C$, which is continuous with respect to the topology of pseudosaturation $E^\vartriangle$.
\end{ex}
\bpr
This si the space $E=Y^\perp$ from Example \ref{ex-3.22}. It is complete, since it is closed subspace in the complete space $Z^\star=\mathcal D^\star(\Bbb R)$. On the other hand, by Proposition \ref{lm-4.12}(a), $E^\star\cong X^\star / E^\perp \cong Z/Y$, and the last space in this chain is metrizable, but not complete. That is
$$
E^\star\ne(E^\star)^\triangledown,
$$
and this can be extended to the chain
$$
E^\star\ne(E^\star)^\triangledown\cong\eqref{(X^vartriangle)^star=(X^star)^triangledown}\cong (E^\vartriangle)^\star,
$$
which means that there exists a functional $f\in (E^\vartriangle)^\star\setminus E^\star$. (It is important here that $E$ is pseudocomplete, while $E^\star$ is not pseudocomplete.)
\epr

\paragraph{Spaces of operators and continuous bilinear maps.}

\bit{
\item Let $X$ and $Y$ be stereotype spaces. Let us denote
 \bit{
\item[---] by $Y:X$ the space of linear continuous maps $\ph:X\to Y$, endowed with the topology of uniform convergence on totally bounded sets in $X$,

\item[---] by $Y\oslash X$ the pseudosaturation of the space $Y:X$,
 \beq\label{DEF:Y-oslash-X}
  Y\oslash X = (Y:X)^\vartriangle
 \eeq
 }\eit
The space $Y\oslash X$ is stereotype, and we call it the {\it inner space of operators} from $X$ into $Y$. Again, it   consists of all linear continuous maps $\ph:X\to Y$, but its topology is formally stronger than the topology of uniform convergence on totally bounded sets in $X$.\footnote{\label{FOOTNOTE:Y-oslash-X=Y:X?} Thus, $Y:X$ and $Y\oslash X$ coincide as linear spaces, but may have different topologies. So far, however, it is not clear, whether $Y:X$ and $Y\oslash X$ indeed are different, since the examples of non-pesudosaturated spaces of the form $Y:X$ (with stereotype $X$ and $Y$) are not constructed yet.}
}\eit

\btm \label{th-5.1} Let $X$ and $Y$ be arbitrary locally convex spaces. A set of morphisms
$\varPhi\subseteq Y:X$ is totally bounded in $Y:X$ if and only  if it satisfies the following two conditions:
\begin{itemize}
\item[(a)] equicontinuity on totally bounded sets:
$$
\forall S\in \mathcal  S(X) \quad \forall V\in \mathcal U(Y) \quad \exists U\in
\mathcal  U(X) \quad \forall a,b\in S \quad a-b\in U \Rightarrow \forall \ph\in
\varPhi \quad \ph(a)-\ph(b)\in V
$$
\item[(b)] uniform total boundedness on totally bounded sets:
$$
\forall S\in \mathcal  S(X) \quad \varPhi(S)=\{\ph(x),x\in S,\ph\in\varPhi\}\in
\mathcal  S(Y)
$$
\end{itemize}
\noindent The condition $(b)$ in this pair can be replaces by the weakened condition
\begin{itemize}
\item[(c)] pointwice total boundedness:
$$
\forall x\in X  \quad \varPhi(x)=\{\ph(x),\ph\in\varPhi\}\in \mathcal  S(Y)
$$
\end{itemize}
Apart from this,
\begin{itemize}
\item[--] if $Y$ is a Heine-Borel space, then $(a) \Rightarrow (b) \&
(c)$;

\item[--] if $X$ is barreled, then $(c)\Rightarrow (a)\& (b)$.
\end{itemize}
\etm

For any linear continuous map $\ph:X\to Y$ its {\it dual map} is the map $\ph^\star:Y^\star\to X^\star$ acting by formula
$$
\ph^\star(f)=f\circ\ph,\qquad f\in Y^\star.
$$

\btm\label{th-6.2} The map $\ph \mapsto \ph^\star$ is an isomorphism of stereotype spaces
$$
X^\star \oslash Y^\star \cong Y\oslash X
$$
\etm

\begin{ex}\label{ex-6.3} If $X$ is a Smith space, and $Y$ a Banach space, then $Y\oslash X =
Y:X$ is a Banach space.
\end{ex}

\begin{ex}\label{ex-6.4}
If $X$ is a Banach space, and $Y$ a Smith space, then $Y\oslash X =Y:X$ is a Smith space.
\end{ex}

\begin{ex}\label{ex-6.5}
If $X$ is a Brauner space, and $Y$ a Fr\'{e}chet space, then $Y\oslash X= Y:X$ is a Fr\'{e}chet space.
\end{ex}

\begin{ex}\label{ex-6.6}
If $X$ is a Fr\'{e}chet space, and $Y$ a Brauner space, then $Y\oslash X= Y:X$ is a Brauner space.
\end{ex}

\bit{
\item Let $X,Y,Z$ be stereotype spaces. Then
 \bit{
\item[---]\label{DEF:bilin-otobrazhenie} we say that  a bilinear map $\beta:X\times Y\to Z$ is {\it continuous}\footnote{This type of continuity is called sometimes  $(\mathcal{K}(X),\mathcal{K}(Y))$-hypocontinuity (cf.\cite{Schaefer}), where $\mathcal{K}(X)$ and $\mathcal{K}(Y)$ are systems of compact sets in $X$ and $Y$ respectively.}, if
 \bit{
\item[(1)] for each compact set $K$ in $X$ and for each neighbourhood of zero $W$ in $Z$ there is a neighbourhood of zero $V$ in $Y$ such that
$$
\beta(K,V)\subseteq W,
$$
\item[(2)] for each compact set $L$ in $Y$ and for each neighbourhood of zero $W$ in $Z$ there is a neighbourhood of zero $U$ in $X$ such that
$$
\beta(U,L)\subseteq W,
$$
 }\eit

\item[---] we denote by $Z:(X,Y)$ the space of continuous bilinear maps $\beta:X\times Y\to Z$ endowed with the topology of uniform convergence on compact sets in $X\times Y$,

\item[---] we denote by $Z\oslash(X,Y)$ the pseudosaturation of the space $Z:(X,Y)$,
 \beq\label{DEF:Z-oslash-(X,Y)}
  Z\oslash(X,Y) = (Z:(X,Y))^\vartriangle
 \eeq
 }\eit
The space $Z\oslash(X,Y)$ is stereotype, and we call it the {\it inner space of bilinear maps} from $X\times Y$ into $Z$. Like $Z:(X,Y)$, it consists of continuous bilinear maps $\beta:X\times Y\to Z$, but the topologies of $Z:(X,Y)$ and $Z\oslash(X,Y)$ may be different.\footnote{Cf. footnote \ref{FOOTNOTE:Y-oslash-X=Y:X?}, the situation with $Z:(X,Y)$ and $Z\oslash(X,Y)$ is the same.}
 }\eit

\begin{ex}\label{ex-6.8} If $X$ and $Y$ are Smith spaces, and $Z$ a Banach space, then $Z\oslash (X,Y)=Z:(X,Y)$ is a Banach space.
\end{ex}

\begin{ex}\label{ex-6.9} If $X$ and $Y$ are Banach spaces, and $Z$ is a Smith space, then $Z\oslash (X,Y)=Z:(X,Y)$ is a Smith space.
\end{ex}

\begin{ex}\label{ex-6.10} If $X$ and $Y$ are Brauner spaces, and $Z$ a Fr\'{e}chet space, then $Z\oslash (X,Y)=Z:(X,Y)$ is a Fr\'{e}chet space.
\end{ex}

\begin{ex}\label{ex-6.11} If $X$ and $Y$ are Fr\'{e}chet spaces, and $Z$ a Brauner space, then $Z\oslash (X,Y)=Z:(X,Y)$ is a Brauner space.
\end{ex}

\btm \label{th-6.12} If $X,Y,Z$ are stereotype spaces, then the formula
\beq
\beta (x,y)= \ph(y)(x) \label{eq6.1}
\eeq
defines an isomorphism of stereotype spaces
\beq
Z\oslash (X,Y) = (Z\oslash X)\oslash Y \label{eq6.2}
\eeq
\etm
\begin{rem} In the special case when $Z=\C$ we have
\beq
\C\oslash (X,Y)= X^\star \oslash Y \label{eq6.3}
\eeq
\beq
Y\oslash X = \C\oslash (Y^\star,X) \label{eq6.4}
\eeq
\end{rem}

\btm \label{th-6.13} For all stereotype spaces $X,Y,Z$ the composition map
$$
(\beta,\alpha) \in (Z\oslash Y)\times (Y\oslash X) \mapsto  \beta \circ \alpha
\in (Z\oslash X)
$$
is a continuous bilinear form.
\etm

 \bit{
\item Let $\alpha :E\to F$ and $\beta :G\to H$ be linear continuous maps of stereotype spaces. Denote by $\beta
\oslash \alpha$ the map
$$
(\beta \oslash \alpha): (G\oslash F)\to (H\oslash E)
$$
acting by formula
 \beq\label{DEF:beta-oslash-alpha}
(\beta \oslash \alpha)(\psi)= \beta \circ \psi \circ \alpha
 \eeq
 }\eit

\btm \label{th-6.14} For all stereotype spaces $X,Y,Z$ the bilinear map
\beq
(\beta,\alpha)\in (H\oslash G)\times (F\oslash E) \mapsto \beta \oslash \alpha
\in (H\oslash E)\oslash (G\oslash F) \label{eq6.5}
\eeq
is continuous.
\etm

\paragraph{Tensor products.}

{\it A projective (stereotype) tensor product} $X\circledast Y$ of stereotype spaces $X$ and $Y$ is defined by the equality
\beq\label{eq7.1}
X\circledast Y=(X^\star \oslash Y)^\star
\eeq
or, equivalently, due to \eqref{eq6.3},
\beq
X\circledast Y=(\C\oslash (X,Y))^\star \label{eq7.2}
\eeq
For $x\in X$ and $y\in Y$ the elementary tensor $x\circledast y\in X\circledast Y$ is defined by the equality
\beq
(x\circledast y)(\ph)= \ph(y)(x) \label{eq7.3}
\eeq
(where $\ph \in X^\star \oslash Y$, and $x\circledast y$ is considered as the element of $(X^\star \oslash Y)^\star$), or, equivalently,
\beq
(x\circledast y)(\beta)= \beta(x,y) \label{eq7.4}
\eeq
(where $\beta \in \C\oslash(X,Y)$, and $x\circledast y$ is considered as an element of $\C\oslash (X,Y)^\star$).

\begin{prop}\label{prop-7.1} The map
$\iota :(x,y)\in X\times Y \mapsto x\circledast y \in X\circledast Y$ is a continuous bilinear form.
\end{prop}

\begin{prop}\label{prop-7.2} The algebraic tensor product $X\otimes Y$ is injectively and denseley embedded into the projective tensor product$X\circledast Y$ by the formula
$$
x\otimes y \mapsto x\circledast y
$$
\end{prop}

\btm [universality of projective tensor product]\label{th-7.3} For any stereotype spaces $X,Y,Z$ and for any continuous bilinear form $\beta :X\times Y\to Z$ there is a unique linear continuous map of stereotype spaces $\widetilde{\beta}:X\circledast Y\to Z$ such that the following diagram is commutative:
$$
\begin{diagram}
\node{X\times Y} \arrow[2]{e,t}{\iota} \arrow{se,b}{\beta}
\node[2]{X\circledast Y} \arrow{sw,b}{\widetilde\beta}
\\
\node[2]{Z}
\end{diagram},
$$
where $\iota$ is defined in Proposition \ref{prop-7.1}. Moreover, the maps $\beta \mapsto \widetilde{\beta}$ is an isomorphism of stereotype spaces
\beq
Z\oslash (X,Y)=Z\oslash (X\circledast Y) \label{eq7.7}
\eeq
\etm

An {\it injective (stereotype) tensor product $X\odot Y$} of stereotype spaces $X$ and $Y$ is defined by the formula
\beq
X\odot Y= Y\oslash X^\star \label{eq7.8}
\eeq
or, equivalently, due to \eqref{eq6.4}, by the formula
\beq
X\odot Y= \C \oslash (X^\star,Y^\star) \label{eq7.9}
\eeq
For $x\in X$ and $y\in Y$ the elementary operator $x\odot y\in X\odot Y$ is defined by
\beq
(x\odot y)(f)=f(x)y, \quad f\in X^\star \label{eq7.10}
\eeq
(if $x\odot y$ is considered as an element of $Y\oslash X^\star$), or by
\beq
(x\odot y)(f,g)=f(x)g(y), \quad f\in X^\star, g\in Y^\star \label{eq7.11}
\eeq
(if $x\odot y$ is considered as an element of $\C\oslash (X^\star,Y^\star)$).

\begin{prop}\label{prop-7.4} The map $\iota :(x,y)\in X\times Y \mapsto x\odot y \in X\odot Y$ is a continuous bilinear form.
\end{prop}

\begin{prop}\label{prop-7.5} The algebraic tensor product $X\otimes Y$ is injectively (but not necessarily dense) embedded into the injective tensor product $X\odot Y$ by the formula
$$
x\otimes y \mapsto x\odot y
$$
\end{prop}

\bex \label{th-7.24} If $X$ and $Y$ are Banach spaces, then $X\circledast Y$ are $X\odot Y$ Banach spaces.
\eex

\bex \label{th-7.25} If $X$ and $Y$ are Smith spaces, then $X\circledast Y$ and $X\odot Y$ are Smith spaces.
\eex

\bex \label{th-7.22} If $X$ and $Y$ are Fr\'{e}chet spaces, then $X\circledast Y$ are $X\odot Y$ Fr\'{e}chet spaces.
\eex

\bex \label{th-7.23} If $X$ and $Y$ are Brauner spaces, then $X\circledast Y$ and $X\odot Y$ are Brauner spaces.
\eex

\paragraph{Category of stereotype spaces.}

The class of stereotype spaces ${\tt Ste}$ forms a category with the linear continuous maps as morphisms.

\medskip
\centerline{\bf Properties of the category ${\tt Ste}$ of stereotype spaces:}

\bit{\it

\item[$1^\circ$.] ${\tt Ste}$ is pre-Abelian.

\item[$2^\circ$.]\label{TH:STE-polnaya-kategoriya} ${\tt Ste}$ is complete: each covariant (and each contravariant) system has injective and projective limits. In the case of the direct coproducts and direct products these constructions coincide with the standard ones in the category ${\tt LCS}$ of locally convex spaces, while in general case the difference is that the injective limits in ${\tt LCS}$ must be pseudocompleted, while the projective limits pseudosaturated:
\begin{align}
& \label{summa-v-Ste}
\overset{\tt Ste}{\bigoplus_{i\in I}} X_i=\overset{\tt LCS}{\bigoplus_{i\in I}} X_i  && \overset{\tt Ste}{\prod_{i\in I}} X_i=\overset{\tt LCS}{\prod_{i\in I}} X_i \\
& \label{rightlim-v-Ste}
\overset{\tt Ste}{\rightlim_{i\to \infty}}X_i=\l \overset{\tt LCS}{\rightlim_{i\to \infty}}X_i\r^\triangledown, &&
\overset{\tt Ste}{\leftlim_{i\to \infty}}X_i=\l \overset{\tt LCS}{\leftlim_{i\to \infty}}X_i\r^\vartriangle.
 \end{align}

\item[$3^\circ$.] The tensor products $\circledast$ and $\odot$ the fraction $\oslash$ are connected with each other through the following isomorphisms of functors:
    \begin{align}
& \label{eq7.16}
(X\circledast Y)^\star  \cong   Y^\star \odot X^\star  &&
(X\odot Y)^\star
\cong   Y^\star \circledast X^\star  \\
& \label{eq7.17} Z\oslash (X\circledast Y)  \cong  (Z\oslash X)\oslash Y &&
(X\odot Y)\oslash Z \cong X\odot (Y\oslash Z)
 \end{align}

\item[$4^\circ$.] ${\tt Ste}$ is a symmetric monoidal category with respect to each of the two tensor products $\circledast$ and $\odot$:
    \begin{align}
& \label{eq7.18}
\C\circledast X \cong  X  \cong  X\circledast \C &&
\C\odot X \cong  X  \cong  X\odot \C \\
& \label{eq7.19} X\circledast Y \cong  Y\circledast X &&  X\odot Y \cong Y\odot
X \\
& \label{eq7.20} (X\circledast Y)\circledast Z  \cong  X\circledast
(Y\circledast Z) && (X\odot Y)\odot Z  \cong  X\odot (Y\odot Z)
 \end{align}

\item[$5^\circ$.] The projective tensor product in ${\tt Ste}$ commutes with injective limits, and the injective product with the projective limits:
\begin{align}
& \label{perestanovochnost-circledast-s-summami}
\left(\bigoplus_{i\in I} X_i \right)\circledast \l\bigoplus_{j\in J}Y_j\r \cong
\bigoplus_{i\in I, j\in J} (X_i \circledast Y_j) && \left(\prod_{i\in I} X_i \right)\odot
\l\prod_{j\in J}Y_j\r \cong \prod_{i\in I, j\in J} (X_i \odot Y_j) \\
& \label{perestanovochnost-circledast-i-oplus-s-predelami-2}
 \Big(\rightlim_{i\to \infty}X_i\Big)\circledast \Big(\rightlim_{j\to
\infty}Y_j\Big)\cong \rightlim_{i,j\to\infty}\Big(X_i\circledast Y_j\Big), &&
\Big(\leftlim_{i\to\infty}X_i\Big)\odot \Big(\leftlim_{j\to\infty}Y_j\Big)\cong
\leftlim_{i,j\to\infty}\Big(X_i\odot Y_j\Big),
 \end{align}
}\eit

\subsection{Subspaces}

\label{SUBSEC:neposr-podpr}

\bit{
\item[$\bullet$]
Let $Y$ be a subset in a stereotype space $X$ endowed with the structure of stereotype space in such a way that the set-theoretic enclosure  $Y\subseteq X$ becomes a morphism of stereotype spaces (i.e. a linear continuous map). Then the stereotype space $Y$ is called a {\it subspace} of the stereotype space $X$, and the set-theoretic enclosure $\sigma:Y\subseteq X$ its {\it representing monomorphism}. The record
$$
Y\subarr X
$$
or
$$
X\suparr Y
$$
will mean that $Y$ is a subspace of the stereotype space $X$. If in addition we write
$$
Y=X
$$
then this means that the stereotype spaces $Y$ and $X$ coincide not only as sets but also with their algebraic and topological structure.

\item[$\bullet$] The system of subspaces of a stereotype space $X$ will be denoted by the symbol $\Sub(X)$.
}\eit

\bprop\label{PROP:mu-cong-podpr-v-Ste} For a morphism $\mu:Z\to X$ in the category ${\tt Ste}$ of stereotype spaces the following conditions are equivalent:
 \bit{
\item[(i)] $\mu$ is a monomorphism,

\item[(ii)] there exists a subspace $Y$ in $X$ with the representing monomorphism $\sigma:Y\subarr X$ and an isomorphism $\theta:Z\to Y$ of stereotype spaces such that the following diagram is commutative:
$$
 \xymatrix @R=1.0pc @C=2.5pc
 {
 Z\ar[rd]^{\mu}\ar@{~>}[dd]_{\theta} & \\
  & X\\
 Y\ar@{-->}[ru]_{\sigma} &
 }
$$
 }\eit
\eprop

\bcor\label{COR:Sub(V)-v-Ste} For any stereotype space $X$ the system $\Sub(X)$ of its subspaces is a system of subobjects in $X$ (in the sense of definition on page \pageref{DEF:toch-sist-podobj-a}).
\ecor

Certainly, for a stereotype space $P$ the relation $\subarr$ is a partial order on the set $\Sub(P)$ of subspaces of $P$.

\paragraph{Immediate subspaces.}
 \bit{
 \item[$\bullet$]
Suppose we have a sequence of two subspaces
$$
Z\subarr Y\subarr X,
$$
and the enclosure $Z\subarr Y$ is a bimorphism of stereotype spaces, i.e. apart from the other requirements, $Z$ is dense in $Y$ (with respect to the topology of $Y$):
$$
\overline{Z}^Y=Y.
$$
Then we will say that the subspace $Y$ is a {\it mediator} for the subspace $Z$ in the space $X$.

 \item[$\bullet$]
We call a subspace $Z$ of a stereotype space $X$ an {\it immediate subspace} in $X$, if it has no non-isomorphic mediators, i.e. for any mediator $Y$ in $X$ the corresponding enclosure $Z\subarr Y$ is an isomorphism. In this case we use the record $Z\osubarr X$:
 $$
 Z\osubarr X\qquad\Longleftrightarrow\qquad \forall Y\quad \bigg( \Big(Z\subarr Y\subarr X \quad\&\quad \overline{Z}^Y=Y\Big)\quad\Longrightarrow\quad Z=Y\bigg).
 $$

 }\eit

\brem In the category of locally convex spaces ${\tt LCS}$ the same construction gives a widely used object: immediate subspaces in a locally convex space $X$ are exactly closed subspaces in $X$ with the topology inherited from $X$. Below in Examples \ref{EX:zamk-neposr-podpr-so-strogo-mazhor-topol} and \ref{EX:nezamk-neposr-podpr} we will see that in the category ${\tt Ste}$ of stereotype spaces the situation becomes sufficiently more complicated.
\erem

Recall that immediate monomorphisms were defined on page \pageref{DEF:immediate-mono}.

\bprop\label{PROP:mu-cong-neposr-podpr-v-Ste} For a morphism $\mu:Z\to X$ in the category ${\tt Ste}$ the following conditions are equivalent:
 \bit{
\item[(i)] $\mu$ is an immediate monomorphism,

\item[(ii)] there exists an immediate subspace $Y$ of $X$ with a representing monomorphism $\sigma:Y\subarr X$ and an isomorphism $\theta:Z\to Y$ such that the following diagram is commutative
\beq\label{predstavlenie-monomorfizma}
 \xymatrix @R=1.0pc @C=2.5pc
 {
 Z\ar[rd]^{\mu}\ar@{~>}[dd]_{\theta} & \\
  & X\\
 Y\ar@{-->}[ru]_{\sigma} &
 }
\eeq
 }\eit\noindent
 The subspaces $Y$ and the morphism $\theta$ here are uniquely defined by $Z$ and $\mu$.
\eprop
\bpr The implication (i)$\Longleftarrow$(ii) is obvious, so we need to prove only (i)$\Longrightarrow$(ii).
 Put $Y=\mu(Z)$, and denote by $\theta:Z\to Y$ the co-restriction of $\mu$ on $Y$, i.e. $\theta$ is the same map as $\mu$ but it is assumed that $\theta$ acts from $Z$ into $Y$. Since $\mu$ is injective, $\theta$ is bijective. Let us endow $Y$ by the topology under which $\theta$ is an isomorphism of locally convex spaces. Then $Y$ becomes a subspace of $X$, since for any neighborhood of zero $U$ in $X$ its inverse image  $\mu^{-1}(U)$ must be a neighborhood of zero in $Z$, and thus the set $Y\cap U=\theta(\mu^{-1}(U))$ must be a neighborhood of zero in $Y$.
\epr

\bprop\footnote{In \cite{Akbarov} Theorem 4.14, which is equivalent to Proposition \ref{PROP:stroenie-zamk-neposr-popdpr} here, and the more general Theorem 11.7, contain an inaccuracy: the requirement of closure of $\sigma$ is omitted there.}\label{PROP:stroenie-zamk-neposr-popdpr} For an immediate subspace $Y$ of a stereotype space $X$ with a representing monomorphism $\sigma:Y\subseteq X$ the following conditions are equivalent:
 \bit{
 \item[(i)] $\sigma$ is a closed map,

 \item[(ii)] $\sigma$ is a weakly closed map,

 \item[(iii)] $Y$ as a set is a closed subspace in the locally convex space $X$, and the topology of $Y$ is a pseudosaturation of the topology inherited from $X$.
  }\eit\noindent
\eprop
 \bit{
 \item[$\bullet$] If the conditions (i)-(iii) of this proposition are fulfilled, then we say that the immediate subspace $Y$ of the space $X$ is {\it closed}.
  }\eit

\bpr
1. The implication (i)$\Longrightarrow$(ii) is a special case of the common situation stated in Proposition \ref{PROP:otkryt=>slabo-otkryt}.

2. Let us prove (ii)$\Longrightarrow$(iii). Let $\sigma:Y\subseteq X$ be a weakly closed map, i.e. $Y$ as a set is closed in $X$. Denote by $E$ the space $Y$ with the topology inherited from $X$. Clearly, $Y$ is continuously embedded into $E$, and, since $Y$ is pseudosaturated, this enclosure preserves its continuity after passage from $E$ to its pseudosaturation $E^\vartriangle$ (we use here the reasoning stated in Diagram  \cite[(1.26)]{Akbarov}). Thus, we obtain a sequence of subspaces
$$
Y\subarr E^\vartriangle\subarr X,
$$
and, since $Y$ and $E^\vartriangle$ coincide as sets, the first of these monomorphisms is a bimorphism. Hence, $E^\vartriangle$ is a mediator for  $Y$, and we obtain that $Y=E^\vartriangle$.

3. The implication (iii)$\Longrightarrow$(i) follows from the fact the pseudosaturation does not change the system of totally bounded subsets.
\epr

\bex\label{EX:zamk-neposr-podpr-so-strogo-mazhor-topol} There exists a stereotype space $P$ with a closed immediate subspace $Q$, which topology is not inherited from $P$, and, moreover, some continuous functionals $g\in Q^\star$ cannot be continuously extended on $P$ (in the formal language this means that the representing monomorphism $Q\osubarr P$ is closed, but not a weakly open map).
\eex
\bpr
Consider the space $E$ from Example \ref{ex-3.23}. It is complete, so it can be represented as a complete subspace in some stereotype space $P$ with the topology inherited from $P$ (for example, one can take as $P$ the direct product of all Banach quotient spaces $E/F$ of $E$). The space  $Q=E^\vartriangle$ is the one with the required properties. Indeed, it is closed in $P$, since $E$ is closed in $P$. On the other hand, the functional $f:Q\to\C$, described in Example \ref{ex-3.23}, is continuous on $Q=E^\vartriangle$, but it cannot be continuously extended to $P$, since otherwise it would be continuous on $E$.
\epr

\bex\label{EX:nezamk-neposr-podpr} There exists a stereotype space $X$ with an immediate subspace $Z$, which is not closed as a subset in $X$. Hence the enclosure $Z\subseteq X$ is not a weakly closed morphism in the sense of definition on page \pageref{DEF:slabo-zamkn-morf} (in particular, the enclosure $Z\subseteq X$ is not isomorphic in $\Mono_X$ to a kernel of some other morphism $\ph:X\to A$ in ${\tt Ste}$).
\eex
\bpr
Let $E$ and $f$ be the space and the functional from Example \ref{ex-3.23}. Consider the kernel $F=\{x\in E^\vartriangle:\ f(x)=0\}$ of $f$, and endow $F$ with the topology inherited from $E^\vartriangle$ (as a locally convex space $F$ is a closed subspace in $E^\vartriangle$). By \cite[Proposition 3.19]{Akbarov}, $E^\vartriangle$ is complete, hence $F$ is also complete, and again by \cite[Proposition 3.19]{Akbarov}, its pseudosaturation
$Z=F^\vartriangle$ must be complete. In addition, $Z=F^\vartriangle$ is pseudosaturated, and thus, stereotype. Note then, that since $E$ is complete, it can be represented (as a locally convex space) as a closed subspace in a direct product $X$ of some Banach space (in such a way that the topology of $E$ is inherited from $X$). We will show that $Z$ is an immediate subspace, but not a closed set in $X$.

First let us show that $Z$ is not closed in $X$. As a set $Z$ coincides with $F$, which is dense in $E$ (in the topology of $E$, which is inherited from  $X$). Hence,
$$
\overline{Z}^X=\overline{F}^X=E\ne F=Z
$$
(here $\overline{\phantom{Z}}^X$ means closure in $X$, as we settled on page \pageref{DEF:overline^X}). Now let us show that $Z$ is an immediate subspace in $X$. Let $Y$ be a mediator of $Z$ in $X$:
$$
Z\subseteq Y\subseteq X
$$
From the fact that $Z$ is dense in $Y$ we obtain the following chain
$$
\overline{Z}^Y=Y
$$
$$
\Downarrow
$$
$$
\overline{Y}^X=\overline{\overline{Z}^Y}^X=\overline{Z}^X=E
$$
$$
\Downarrow
$$
$$
Y\subseteq E.
$$
This is an enclosure of sets. Note now that since $Y$ is a subspace in $X$, the topology of $Y$ must majorate the topology inherited from $X$, or, what is the same, the topology inherited from $E$. That is why the enclosure $Y\subseteq E$ is continuous, and therefore $Y$ is a subspace in $E$. This implies that the pseudosaturation of $Y$ must be a subspace in the pseudosaturation of $E$, and, since $Y$ is pseudosaturated, we obtain a continuous enclosure:
$$
Y=Y^\triangledown\subseteq E^\triangledown.
$$
Thus, $Y$ is a subspace in $E^\triangledown$.

Let us now forget about $X$ and consider the following chain of subspaces:
$$
Z\subseteq Y\subseteq E^\triangledown.
$$
From the fact that $Z$ is a dense subspace in $Y$ we obtain a new logical chain:
$$
\overline{Z}^Y=Y
$$
$$
\Downarrow
$$
$$
\overline{Y}^{E^\triangledown}=\overline{\overline{Z}^Y}^{E^\triangledown}=\overline{Z}^{E^\triangledown}=F
$$
$$
\Downarrow
$$
$$
Y\subseteq F.
$$
Again this is an enclosure of sets. Then we note that since $Y$ is a subspace in $E^\triangledown$, the topology of $Y$ must majorate the topology inherited from $E^\triangledown$, or, what is the same, the topology inherited from $F$. Thus the enclosure $Y\subseteq F$ is continuous, and, as a corollary, $Y$ is a subspace in $F$. This implies that the pseudosaturation of $Y$ must be a  subspace in the pseudosaturation of $F$, and, since $Y$ is pseudosaturated, we obtain a continuous enclosure:
$$
Y=Y^\triangledown\subseteq F^\triangledown=Z.
$$
Thus, $Y$ is a subspace in $F^\triangledown=Z$. On the other hand, from the very beginning $Z$ was a subspace in $Y$. Hence, $Z=Y$.
\epr

\paragraph{Envelope $\Env^X M$ of a set $M$ of elements in a space $X$.}

Theorem \ref{TH:Env^X(M)-obolochka-v-kateg-smysle} which we will prove later, justifies the following definition.

\bit{

\item[$\bullet$] The {\it envelope} of a set $M\subseteq X$ in a stereotype space $X$ is a subspace in $X$, denoted by $\Env^X M$ or by $\Env M$, and defined as the projective limit in the category $\tt Ste$
        \beq\label{DEF:Span_infty}
    \Env^X M=\Env M={\tt Ste}\text{-}\leftlim E_i
    \eeq
    of a contravariant system $\{E_i;\ i\in{\tt Ord}\}$ of subspaces in $X$, indexed by ordinal numbers and defined by the following inductive rules:
    \bit{

    \item[0)] the space $E_0$ is defined as the pseudosaturation of the closure of linear span $\Sp M$ of the set $M$ in the space $X$:
    $$
    E_0=\left(\overline{\Sp M}^X\right)^\vartriangle
    $$

    \item[1)] suppose that for some ordinal number $j\in{\tt Ord}$ all the spaces $\{E_i;\ i<j\}$ are already defined, then the space $E_j$ is defined as follows:
    \bit{

    \item[---] if $j$ is an isolated ordinal number, i.e. $j=i+1$ for some $i$, then $E_j=E_{i+1}$ is defined as the pseudosaturation of the closure of linear span $\Sp M$ of the set $M$ in the space $E_i$:
    $$
    E_j=E_{i+1}=\left(\overline{\Sp M}^{E_i}\right)^\vartriangle
    $$

    \item[---] if $j$ is a limit ordinal number, i.e. $j\ne i+1$ for any $i$, then $E_j$ is defined as the projective limit in the category ${\tt Ste}$ of the net $\{E_i;\ i\to j\}$:
    $$
    E_j=\lim_{j\gets i}E_i,
    $$
    -- this means that as a set $E_j$ is the intersection of the spaces $\{E_i;\ i\to j\}$,
    $$
    E_j=\bigcap_{i<j}E_i,
    $$
    and the topology in $E_j$ is the weakest stereotype locally convex topology, under which all the enclosures $E_j\subseteq E_i$ are continuous.

    }\eit

    }\eit
    Since the transfinite sequence $\{E_i;\ i\in{\tt Ord}\}$ cannot be an injective map from the class of all ordinal numbers ${\tt Ord}$ to the set  $\Sub(X)$ of all subspaces of a stereotype space $X$, it must stabilize, i.e. after some number $k\in{\tt Ord}$ all the spaces $E_i$ must coincide (with their topologies):
\beq\label{E_(k+1)=E_k}
\forall l\ge k\qquad E_l=E_k.
\eeq
This implies that the contravariant system $\{E_i;\ i\in{\tt Ord}\}$ indeed has a projective limit, and this is exactly the subspace $E_k$ in $X$.
    }\eit

\bex\label{EX:M-poln-v-X=>Sp_infty(M)=X}
If a set $M$ is total in $X$, then its envelope coincides with $X$:
    $$
\overline{\Sp M}^X=X\quad\Longrightarrow\quad \Env M=X
    $$
\eex
\bpr
The equality $\overline{\Sp M}^X=X$ implies $E_0=\l\overline{\Sp M}^X\r^\vartriangle=X$, and after that all the spaces $E_i$ become equal to $X$
$$
X=E_0=E_1=...
$$
Hence, $\Env M=X$.
\epr

\bex
If a set $M$ forms a closed subspace in $X$ (as in a locally convex space), then its envelope coincides with the pseudosaturation of $M$ with respect to the topology inherited from $X$:
    $$
\overline{\Sp M}^X=M \quad\Longrightarrow\quad    \Env M=M^{\vartriangle}
    $$
\eex
\bpr
From $\overline{\Sp M}^X=M$ we have $E_0=\l\overline{\Sp M}^X\r^\vartriangle=M^\vartriangle$, then $E_1=\l\overline{\Sp M}^{E_0}\r^\vartriangle= M^\vartriangle=E_0$, and all the other spaces $E_i$ coincide with $E_0$. Thus, $\Env M=E_0=M^\vartriangle$.
\epr

\btm\label{TH:Span_infty^X(M)-osubarr-X} The envelope $\Env^X M$ of each set $M\subseteq X$ is an immediate subspace in $X$, containing $M$ as a total subset:
\beq\label{M-subseteq-Sp_infty-M^X-osubarr-X}
M\subseteq \Env^X M\osubarr X,
\qquad
\overline{\Sp M}^{\Env^X M}=\Env^X M.
\eeq
\etm
\bpr
1. First let us verify that $M$ is total in $\Env^X M$. Suppose $k$ is an ordinal number  after which the sequence $\{E_i;\ i\in{\tt Ord}\}$ is stabilized, i.e. \eqref{E_(k+1)=E_k} holds. Then $\Env M=E_k$, and if it turned out that $M$ is not total in $E_k$, then we would have a contradiction with \eqref{E_(k+1)=E_k}:
$$
E_{k+1}=\overline{\Sp M}^{E_k}\ne E_k.
$$

2. Let us show that $\Env M$ is an immediate subspace in $X$. Suppose $Y$ is a subspace in $X$ such that
$$
\Env M\subseteq Y\subseteq X,
$$
and $\Env M$ is dense in $Y$. Since, as we already understood, $\Sp M$ is dense in $\Env M$, we have
 \beq\label{Y=overline(Sp_M)^Y}
Y=\overline{\Sp M}^Y.
 \eeq
Now by induction we have that $Y$ is continuously embedded into each $E_i$:
 \bit{
\item[0)] for $i=0$ we have a chain
$$
Y\subarr X
\quad\Longrightarrow\quad
Y=\eqref{Y=overline(Sp_M)^Y}=\overline{\Sp M}^Y\subarr \overline{\Sp M}^X
\quad\Longrightarrow\quad
Y=Y^\vartriangle\subarr  \left(\overline{\Sp M}^X\right)^\vartriangle=E_0.
$$

\item[1)] suppose that we proved $Y\subarr E_i$ for all $i$ less that some $j$, then
 \bit{
\item[---] if $j$ is an isolated ordinal number, i.e. $j=i+1$ for some $i$, then
$$
\kern-40pt
Y\subarr E_i
\quad\Longrightarrow\quad
Y=\eqref{Y=overline(Sp_M)^Y}=\overline{\Sp M}^Y\subarr \overline{\Sp M}^{E_i}
\quad\Longrightarrow\quad
Y=Y^\vartriangle\subseteq  \left(\overline{\Sp M}^{E_i}\right)^\vartriangle=E_{i+1}=E_j,
$$

\item[---] if $j$ is a limit ordinal number, then from the continuous enclosures $Y\subarr E_i$ for $i<j$ we obtain a continuous enclosure of locally convex spaces
$$
Y\subarr \text{\tt LCS-}\kern-2pt\lim_{j\gets i}E_i,
$$
and this implies a continuous enclosure of stereotype spaces
$$
Y=Y^\vartriangle\subarr\l \text{\tt LCS-}\kern-2pt\lim_{j\gets i}E_i\r^\vartriangle=\text{\tt Ste-}\kern-2pt\lim_{j\gets i}E_i=E_j.
$$
 }\eit
 }\eit\noindent
From the fact that $Y$ is continuously embedded into each $E_i$ we obtain a continuous enclosure $Y\subarr \Env M$. Together with the initial enclosure $\Env M\subarr Y$ this means the equality $\Env M=Y$ (with topologies).
\epr

The following theorem shows that in an immediate subspace the topology is automatically defined by the set of its elements:

\btm\label{Y-subarr-X=>Y=Sp_infty(Y)} Every subspace $Y$ in a stereotype space $X$ is a subspace in its envelope $\Env^X Y$
\beq\label{Y-subarr-X=>Y-subarr-Sp_infty^X(Y)}
Y\subarr X\quad\Longrightarrow\quad Y\subarr \Env^X Y,
\eeq
and $Y$ is an immediate subspace in $X$ iff it coincide (with the topologies) with its envelope in $X$:
\beq\label{Y-osubarr-X=>Y-=-Sp_infty^X(Y)}
Y\osubarr X\quad\Longleftrightarrow\quad Y=\Env^X Y.
\eeq
\etm
\bpr
The continuity of the enclosure $Y\subarr\Env^X Y$ is proved by induction:
 \bit{
\item[0)] at the zero step we have a continuous enclosure of locally convex spaces
$$
Y\subarr \overline{\Sp Y}^X=\overline{Y}^X,
$$
which implies a continuous enclosure of stereotype spaces
$$
Y=Y^\vartriangle\subarr  \left(\overline{Y}^X\right)^\vartriangle=E_0,
$$

\item[1)] suppose that the continuous enclosure $Y\subseteq E_i$ is proved for all $i$ less than some $j$, then
 \bit{
\item[---] if $j$ is an isolated ordinal number, i.e. $j=i+1$ for some $i$, then we obtain a continuous enclosure of locally convex spaces
$$
Y\subseteq \overline{\Sp Y}^{E_i}=\overline{Y}^{E_i},
$$
which implies a continuous enclosure of stereotype spaces
$$
Y=Y^\vartriangle\subseteq  \left(\overline{Y}^{E_i}\right)^\vartriangle=E_{i+1}=E_j,
$$

\item[---] if $j$ is a limit ordinal number, then from the continuous enclosures $Y\subseteq E_i$ for all $i<j$ we obtain a continuous enclosure of locally convex spaces
$$
Y\subseteq \text{\tt LCS-}\kern-2pt\lim_{j\gets i}E_i,
$$
which implies a continuous enclosure of stereotype spaces
$$
Y=Y^\vartriangle\subseteq \l \text{\tt LCS-}\kern-2pt\lim_{j\gets i}E_i\r^\vartriangle=\text{\tt Ste-}\kern-2pt\lim_{j\gets i}E_i=E_j.
$$
 }\eit
 }\eit
Let us now consider a special case when $Y$ is an immediate subspace in $X$. Then by Theorem \ref{TH:Span_infty^X(M)-osubarr-X}, $Y$ is dense in  $\Env Y$, hence in the chain of enclosures
$$
Y\subseteq \Env Y\subseteq  X
$$
the second space is a mediator. Therefore, it coincides with the first one: $Y=\Env Y$.
\epr

\bcor The representing monomorphism $\sigma:Y\osubarr X$ of an immediate subspace $Y$ in a stereotype space $X$ is always relatively closed.
\ecor
\bpr
By Theorem \ref{Y-subarr-X=>Y=Sp_infty(Y)}, $Y$ is the projective limit in $\tt Ste$ of the chain of subspaces  $\{E_i\}$ (which defines $\Env^X Y$):
$$
Y=\Env^X Y=\leftlim_{i\in{\tt Ord}}E_i=\bigcap_{i\in{\tt Ord}}E_i.
$$
Let $T$ be an absolutely convex compact set in $X$, lying in $Y$ as a set. Then $T$ lies in  $E_0=\big(\overline{Y}^X\big)^\vartriangle$, and since on passing from the topology of $X$ to the topology of $E_0$ the system of compact sets (as well as the topology on each compact set) is inherited from $X$ (this is one of the fundamental properties of the pseudosaturation $\vartriangle$, \cite[Theorem 1.17]{Akbarov}), we get that $T$ is a compact set in $E_0$. After that with the same technique we show that $T$ is compact in $E_1$, and more generally, on passing from each ordinal $i$ to its successor $i+1$. When we need to pass to a limit ordinal $j$, we come to the situation where $T$ is a compact set in each $E_i$ with the index $i<j$. As a corollary, $T$ is compact in the projective limit  $\leftlim_{i<j} E_i=\bigcap_{i<j} E_i$. When we come to enough big ordinal, we obtain that $T$ is compact in $Y$.
\epr

\btm\label{ph(M)-sub-N=>ph:Sp8(M)->Sp8(N)} If $\ph:Y\to X$ is a morphism of stereotype spaces, turning a set $N\subseteq Y$ into a set $M\subseteq X$,
$$
\ph(N)\subseteq M,
$$
then $\ph$ continuously maps $\Env^Y N$ into $\Env^X M$:
 $$
 \xymatrix{
 Y\ar[r]^{\ph} & X\\
 \Env^Y N\ar[u]\ar@{-->}[r] & \Env^X M \ar[u]
 }
 $$
In the special cases:
 \begin{align}
& \begin{Bmatrix}Y  \subarr  X \\ \text{\rotatebox{90}{$\subseteq$}} \quad \ \text{\rotatebox{90}{$\subseteq$}}\\
N  \subseteq  M \end{Bmatrix}\quad\Longrightarrow\quad \Env^Y N\subarr \Env^X M,
\label{N-subseteq-M,Y=X}
\\
& \begin{Bmatrix}Y  \osubarr  X \\ \text{\rotatebox{90}{$\subseteq$}} \quad \ \text{\rotatebox{90}{$\subseteq$}}\\
N  \subseteq  M \end{Bmatrix}\quad\Longrightarrow\quad \Env^Y N\osubarr \Env^X M,
\label{N-subseteq-M,Y-osubarr-X}
\\
& \begin{Bmatrix}Y  \osubarr  X \\ \text{\rotatebox{90}{$\subseteq$}} \quad \ \text{\rotatebox{90}{$\subseteq$}}\\
N  =  M \end{Bmatrix}\quad\Longrightarrow\quad \Env^Y M=\Env^X M.
\label{N=M,Y-osubarr-X}
 \end{align}
\etm
\bpr
Take a morphism $\ph:Y\to X$ of stereotype spaces turning a set $N\subseteq Y$ into a set $M\subseteq X$, $\ph(N)\subseteq M$.
If by $\{F_i;\ i\in{\tt Ord}\}$ and $\{E_i;\ i\in{\tt Ord}\}$ we denote the sequences of subspaces in $Y$ and $X$, which define
$\Env N$ and $\Env M$ respectively,
$$
\Env N=\leftlim  F_i,\qquad \Env M=\leftlim  E_i,
$$
then we can prove by induction, that for each $i$ the map $\ph$ continuously maps $F_i$ into $E_i$,
$$
\ph:F_i\to E_i,
$$
and this implies that $\ph$ continuously maps $\Env N$ into $\Env M$,
$$
\ph:\Env N\to \Env M.
$$
Let us now observe the special cases.

1. If $N\subseteq M$ and $Y\subarr X$, then we consider the sequences $\{F_i;\ i\in{\tt Ord}\}$ and $\{E_i;\ i\in{\tt Ord}\}$ of subspaces in $X$, which define $\Env^X N$ and $\Env^X M$. By induction we obtain an enclosure of subspaces $F_i\subarr E_i$ for each $i$, and this gives the enclosure  $\Env^X N\subarr \Env^X M$.

2. Suppose that $N\subseteq M$ and $Y\osubarr X$. Then, by implication \eqref{N-subseteq-M,Y=X} we already proved, $\Env^Y N\subarr\Env^X M$. Let us show that in this enclosure $\Env^Y N$ is an immediate subspace in $\Env^X M$. Let $Z$ be a mediator for $\Env^Y N$ in $\Env^X M$:
$$
\Env^Y N\subarr Z\subarr \Env^X M,\qquad \overline{\Env^Y N}^Z=Z.
$$
Consider the envelope $\Env^X (Y\cup Z)$ of the set $Y\cup Z$ in the space $X$. We can include it into a diagram (where all the arrows are theoretic set enclosures, which are continuous maps):
 $$
 \xymatrix{
 Y\ar[r] & \Env^X (Y\cup Z)\ar[r] & X\\
 \Env^Y N\ar[u]\ar[r] & Z\ar[r]\ar[u] & \Env^X M \ar[u]
 }
 $$
By Theorem \ref{TH:Span_infty^X(M)-osubarr-X}, $N$ is total in $\Env^Y N$, which in its turn is total in $Z$ (since $Z$ is a mediator). Hence,
$N$ is total in $Z$. On the other hand, $N\subseteq Y$, hence $Y$ is dense in $Z$ (in the topology of $Z$, and thus in the topology of $X$ as well). From this we have that $Y$ is dense in the subset $Y\cup Z$ of the space $X$, and again by Theorem \ref{TH:Span_infty^X(M)-osubarr-X}, $Y$ is dense in $\Env^X (Y\cup Z)$.

This means that $\Env^X (Y\cup Z)$ is a mediator for $Y$ in the space $X$:
$$
Y\subarr \Env^X (Y\cup Z)\subarr X,\qquad \overline{Y}^{\Env^X (Y\cup Z)}=\Env^X (Y\cup Z).
$$
The condition $Y\osubarr X$ implies the equality of stereotype spaces $Y=\Env^X (Y\cup Z)$. This in its turn implies that $Z\subarr Y$, i.e. $Z$ is a mediator for $\Env^Y N$ in $Y$:
$$
\Env^Y N\subarr Z\subarr Y,\qquad \overline{\Env^Y N}^Z=Z.
$$
By Theorem \ref{TH:Span_infty^X(M)-osubarr-X}, $\Env^Y N$ is an immediate subspace in $Y$, so we obtain an equality of stereotype spaces $\Env^Y N=Z$.

3. Suppose that $N=M\subseteq Y\osubarr X$. Then by property \eqref{N-subseteq-M,Y-osubarr-X} which we already proved,
$$
\Env^Y M\osubarr \Env^X M.
$$
On the other hand by property \eqref{N-subseteq-M,Y=X} which has already been proved as well, the chain $M\subseteq Y\osubarr X$ implies
$$
\Env^X M\osubarr \Env^X Y=\eqref{Y-osubarr-X=>Y-=-Sp_infty^X(Y)}=Y
$$
Together this gives a chain
$$
\Env^Y M\osubarr \Env^X M\osubarr Y.
$$
By Theorem \ref{TH:Span_infty^X(M)-osubarr-X}, the set $M$ is total in $\Env^X M$, hence the space $\Env^Y M$ is total in $\Env^X M$. Thus, $\Env^X M$ is a mediator in this chain, and we obtain the equality $\Env^Y M=\Env^X M$.
\epr

\btm The envelope $\Env^X M$ of any set $M\subseteq X$ is a minimal subspace among all the immediate subspaces in $X$, which contain $M$, and
in each of those immediate subspaces $Y\osubarr X$ the space $\Env^X M$ is an immediate subspace:
\beq\label{Span_infty(M)-osubarr-Y-osubarr-X}
\forall Y\qquad \Big(M\subseteq Y\osubarr X\quad\Longrightarrow\quad \Env^X M\osubarr Y\Big).
\eeq
\etm
\bpr
$$
\Env^X M\overset{\tiny\eqref{N=M,Y-osubarr-X}}{=}\Env^Y M\overset{\tiny\eqref{M-subseteq-Sp_infty-M^X-osubarr-X}}{\osubarr} Y.
$$
\epr

\bprop\label{PROP:Y-osubarr-X,Z-subarr-X,Z-subseteq-Y=>Z-subarr-Y}
If $Y\osubarr X$ and $Z\subarr X$, then the condition $Z\subseteq Y$ implies $Z\subarr Y$. In a special case, when $Y\osubarr X$ and $Z\osubarr X$, the condition  $Z\subseteq Y$ implies $Z\osubarr Y$.
\eprop
\bpr
If $Y\osubarr X$, $Z\subarr X$, $Z\subseteq Y$, then
$$
Z\overset{\eqref{Y-subarr-X=>Y-subarr-Sp_infty^X(Y)}}{\subarr}\Env^XZ\overset{\tiny\eqref{N=M,Y-osubarr-X}}{=}\Env^YZ
\overset{\eqref{M-subseteq-Sp_infty-M^X-osubarr-X}}{\osubarr} Y.
$$
If $Y\osubarr X$, $Z\osubarr X$, $Z\subseteq Y$, then
$$
Z\overset{\eqref{Y-osubarr-X=>Y-=-Sp_infty^X(Y)}}{=}\Env^XZ\overset{\tiny\eqref{N=M,Y-osubarr-X}}{=}\Env^YZ
\overset{\eqref{M-subseteq-Sp_infty-M^X-osubarr-X}}{\osubarr} Y.
$$
\epr

\subsection{Quotient spaces}

 \bit{
\item[$\bullet$] Let $X$ be a stereotype space, and
 \bit{

 \item[1)] in $X$ as in a locally convex space we take a closed subspace $E$,

 \item[2)] on the quotient space $X/E$ we consider an arbitrary locally convex topology $\tau$, which is majorated by the natural quotient topology of $X/E$,

 \item[3)] in the completion $(X/E)^\blacktriangledown$ of the locally convex space $X/E$ with the topology $\tau$ we take a subspace $Y$, which contains $X/E$ and is a stereotype space with respect to the topology inherited from $(X/E)^\blacktriangledown$.
 }\eit
Then we call the stereotype space $Y$ a {\it quotient space of the stereotype space} $X$, and the composition $\upsilon=\sigma\circ\pi$ of the quotient map $\pi:X\to X/E$ and the natural enclosure $\sigma:X/E\to Y$ is called the {\it representing epimorphism} of the quotient space $Y$. The record
$$
Y\quarr X
$$
or the record
$$
X\qparr Y
$$
will mean that $Y$ is a quotient space of the stereotype space $X$. The class of all quotient spaces of $X$ will be denoted by $\Quot(X)$. From its construction it is clear that $\Quot(X)$ is a set.

 }\eit

The following is evident:

\bprop\label{PROP:e-cong-neposr-faktor-pr-v-Ste} For a morphism $\e:Z\gets X$ in the category ${\tt Ste}$ the following conditions are equivalent:
 \bit{
\item[(i)] $\e$ is an epimorphism,

\item[(ii)] there is a quotient space $Y$ of $X$ with the representing epimorphism $\upsilon:Y\quarr X$, and an isomorphism $\theta:Z\gets Y$ such that the following diagram is commutative:
\beq\label{predstavlenie-epimorfizma}
 \xymatrix @R=1.0pc @C=2.5pc
 {
 Z & \\
  & X\ar[lu]_{\e} \ar@{-->}[ld]^{\upsilon} \\
 Y\ar@{~>}[uu]^{\theta} &
 }
\eeq
 }\eit
\eprop

\bcor\label{COR:Quot(V)-v-Ste} For a stereotype space $X$ the system $\Quot(X)$ of all its quotient spaces is a system of quotient objects for $X$.
\ecor

The formalization of the idea of quotient object we have presented here has a qualitative shortcoming in comparison with the notion of subspace which we considered above: the problem is that the relation $\quarr$ does not establish a partial order in the system $\Quot(P)$ of quotient spaces of a stereotype space $P$. By the set-theoretic reasons no one of axioms of partial order (reflexivity, antisymmetry and transitivity) holds for $\quarr$. In particular, the first two axioms do not hold since the situation when $Y\quarr X$ and at the same time $Y=X$ is impossible. To explain this, let us agree for simplicity that we do not take into account the necessity to pass to a subspace in the completion which was stated in the step 3 of our definition -- then $Y\quarr X$ (and $Y\ne\varnothing$) implies by the axiom of regularity \cite[Appendix, Axiom VII]{38} that there exists an element  $y\in Y$ such that $y\cap Y=\varnothing$. But if in addition $Y=X$, then the element $y$, being a coset of $X$, i.e. a non-empty subset in $X$, must have non-empty intersection $y\cap Y=y\cap X=y\ne\varnothing$ with $X=Y$. As to the transitivity, in the situation when $Z\quarr Y$ and $Y\quarr X$ the elements of $Z$ are non-empty sets of elements of $Y$, and each such element is a non-empty set of elements of $X$. From the point of view of set theory this is not the same as if elements of $Z$ were sets of elements of $X$, so in this situation the relation $Z\quarr X$ is also impossible. This forces us to introduce a new binary relation.
 \bit{
 \item[$\bullet$] Suppose $Y\quarr X$ and $Z\quarr X$. We will say that the quotient space $Y$ {\it subordinates} the quotient space $Z$, and we write in this situation $Z\le Y$, if there exists a morphism $\varkappa:Y\to Z$ such that the following diagram is commutative:
 \beq\label{DEF:le-dlya-faktor-prostranstv}
\xymatrix @R=1pc @C=2pc
{
Y\ar[dd]_{\varkappa} &   \\
  & X\ar[lu]_{\upsilon_Y}\ar[ld]^{\upsilon_Z} \\
Z &
}
\eeq
(here $\upsilon_Y$ and $\upsilon_Z$ are representing epimorphisms for $Y$ and $Z$). The morphism $\varkappa$, if exists, must be, first, unique, and, second, an epimorphism.
 }\eit

For any stereotype space $P$ the relation $\le$ is a partial order on the set $\Quot(P)$ of quotient spaces of $P$.

\paragraph{Immediate quotient spaces.}
 \bit{
 \item[$\bullet$]
Let $Y$ and $Z$ be two quotient spaces of $X$ such that
$$
Z\le Y,
$$
and the epimorphism $\varkappa:Z\gets Y$ in diagram \eqref{DEF:le-dlya-faktor-prostranstv} is a monomorphism (and hence, a bimorphism) of stereotype spaces. Then we will say that the quotient space $Y$ is a {\it mediator} for the quotient space $Z$ of the space $X$.

 \item[$\bullet$]
We call a quotient space $Z$ of a stereotype space $X$ an {\it immediate quotient space} in $X$, if it has no non-isomorphic mediators, i.e. for any its mediator $Y$ in $X$ the corresponding epimorphism $Z\quarr Y$ is an isomorphism. We write in this case $Z\oquarr X$:
 $$
 Z\oquarr X\qquad\Longleftrightarrow\qquad \forall Y\quad \bigg( \Big(Y\quarr X \quad\&\quad Z\le Y\quad \&\quad  \varkappa\in\Mono\Big)\quad\Longrightarrow\quad Z=Y\bigg).
 $$

 \item[$\bullet$] Let us say that an immediate quotient space $Y\oquarr X$ {\it strongly subordinates} an immediate quotient space $Z\oquarr X$, and write $Z\ole Y$, if there exists a strong epimorphism $\varkappa:Y\to Z$ such that diagram \eqref{DEF:le-dlya-faktor-prostranstv} is commutative.
 }\eit

\brem In the category of locally convex spaces ${\tt LCS}$ the immediate quotient spaces of a locally convex space $X$ are exactly quotient space of  $X$ by closed subspaces with the usual quotient topologies. Like in the case of subspaces, in the category ${\tt Ste}$ of stereotype spaces the situation becomes more complicated (see below Examples \ref{EX:otkr-neposr-faktor-pr-ne-X/E} and \ref{EX:nezamk-neposr-faktor-pr}).
\erem

Recall that the notion of immediate epimorphism was defined on page \pageref{DEF:immediate-epi}. The following statement is dual to Proposition  \ref{PROP:mu-cong-neposr-podpr-v-Ste}, and can be proved by the dual reasoning:

\bprop\label{PROP:mu-cong-neposr-faktor-pr-v-Ste} For a morphism $\e:Z\gets X$ in the category ${\tt Ste}$ the following conditions are equivalent:
 \bit{
\item[(i)]  $\e$ is an immediate epimorphism,

\item[(ii)] there exists an immediate quotient space $Y$ of the stereotype space $X$ with the representing morphism $\upsilon:Y\quarr X$ and an isomorphism $\theta:Z\gets Y$ such that the following diagram is commutative:
\beq\label{predstavlenie-neposr-epimorfizma}
 \xymatrix @R=1.0pc @C=2.5pc
 {
 Z & \\
  & X\ar[lu]_{\e} \ar@{-->}[ld]^{\upsilon} \\
 Y\ar@{~>}[uu]^{\theta} &
 }
\eeq
 }\eit\noindent
The quotient space $Y$ and the morphism $\theta$ are uniquely defined by $Z$ and $\e$.
\eprop

\bprop\footnote{In author's paper \cite{Akbarov} Theorem 4.16, which is equivalent to Proposition \ref{PROP:stroenie-zamk-neposr-faktor-pr} here, as well as the more general proposition, Theorem 11.9, contain an inaccuracy: the requirement of openness of $\upsilon$ is omitted there.}\label{PROP:stroenie-zamk-neposr-faktor-pr} For an immediate quotient space $Y$ of a stereotype space $X$ with the representing epimorphism $\upsilon:Y\gets X$ the following conditions are equivalent:
 \bit{
 \item[(i)] $\upsilon$ is an open map,

 \item[(ii)] $\upsilon$ is a weakly open map,

 \item[(iii)] $Y$ is a pseudocompletion $(X/E)^\triangledown$ of the quotient space $X/E$ of the locally convex space $X$ (with the usual quotient topology) by some closed locally convex subspace $E$.
  }\eit
\eprop
 \bit{
 \item[$\bullet$] If the conditions (i)-(iii) of this proposition are fulfilled, then we say that the immediate quotient space $Y$ of $X$ is {\it open}.
  }\eit

\bpr
1. The implication (i)$\Longrightarrow$(ii) is a special case of the common situation described in Proposition \ref{PROP:otkryt=>slabo-otkryt}.

2. Let us prove (ii)$\Longrightarrow$(iii). Suppose the representing epimorphism $\upsilon:Y\gets X$ is a weakly open map. Denote by $E$ its kernel. By definition of stereotype quotient space, $Y$ is a pseudocomplete locally convex subspace in the completion $(X/E)^\blacktriangledown$ of the locally convex space $X/E$ under some topology $\tau$ which is majorated by the quotient topology $X/E$, and $X/E$ lies in $Y$ as set.
Thus, we can represent $\upsilon$ as a diagram
$$
\xymatrix  @R=2.5pc @C=4.0pc
{
X/E\ar@{-->}[d]_{\sigma} &  X\ar[l]_{\pi}\ar[ld]^{\upsilon} \\
  Y &
}
$$
where $\pi:X\to X/E$ is the usual quotient map of locally convex spaces, and $\sigma:X/E\to Y$ is a natural bimorphism. Since $Y$ is pseudocomplete, $\sigma$ can be extended to some morphism $\sigma^\triangledown$ on pseudocompletion $(X/E)^\triangledown$ of the space $X/E$ (we use here the reasoning stated in diagram \cite[(1.13)]{Akbarov}):
$$
\xymatrix @R=2.5pc @C=4.0pc
{
(X/E)^\triangledown\ar@{-->}[dr]_{\sigma^\triangledown} & X/E\ar@{-->}[d]_{\sigma}\ar[l]_{\triangledown_{X/E}} &  X\ar[l]_{\pi}\ar[ld]^{\upsilon} \\
  & Y &
}
$$
Note that $\sigma^\triangledown$ is not only epimorphism (this follows from the property of epimorphisms $3^\circ$ on page \pageref{PROP:e-ph=epi-=>-e=epi}, since the composition $\upsilon=\sigma^\triangledown\circ\triangledown_{X/E}\circ\pi$ is an epimorphism), but also a monomorphism. This is proved as follows. The fact that $\upsilon$ is weakly open, implies that $\sigma$ is weakly open as well. This means that every linear continuous functional on $X/E$ can be extended along the map $\sigma$ to a linear continuous functional on $Y$. In other words, the dual map $\sigma':Y'\to X'$ is a surjection. This implies that the pseudosaturation $\sigma^\triangledown$ must be an injection\footnote{We use here the following obvious property of pseudocompletion: {\it if $\ph:X\to Y$ is a monomorphism of locally convex space, such that the dual map $\ph':X'\gets Y'$ is a surjection, then its pseudocompletion $\ph^\triangledown:X^\triangledown\to Y^\triangledown$ is also a monomorphism of locally convex spaces}.}.

As a result, we have a chain of epimorphisms
$$
Y\overset{\sigma^\triangledown}{\gets} (X/E)^\triangledown\overset{\triangledown_{X/E}\circ\pi}{\quarr} X,
$$
where the first morphism $\sigma^\triangledown$ is a bimorphism. Thus, $(X/E)^\triangledown$ is a mediator for $Y$, and we obtain the equality $Y=(X/E)^\triangledown$.

3. The implication (iii)$\Longrightarrow$(i) follows from the fact that pseudocompletion does not change the topology.
\epr

The following example is dual to Example \ref{EX:zamk-neposr-podpr-so-strogo-mazhor-topol}:

\bex\label{EX:otkr-neposr-faktor-pr-ne-X/E} There exists a stereotype space $P$ with an immediate quotient space of the form $Y=(P/E)^\triangledown$, which cannot be represented in the form $Y=P/F$ for a subspace $F\subseteq P$ (in formal language this means that the representing epimorphism  $Y\oquarr P$ is open, but not closed).
\eex
\bpr
The space $Z$ from example \cite[3.22]{Akbarov} is such a space. It contains a closed subspace $E$ such that the locally convex quotient space $Z/E$ is metrizable, but not complete. As a corollary, in the stereotype sense the space $(Z/E)^\triangledown$ is an immediate quotient space, but it cannot be represented in the form $Z/F$, since $F$ is uniquely defined as the kernel of the map $Z\to Y$, and hence must coincide with $E$.
\epr

From Example \ref{EX:nezamk-neposr-podpr} we have

\bex\label{EX:nezamk-neposr-faktor-pr} There exists a stereotype space $P$ with an immediate quotient space $Y$ such that the representing epimorphism $Y\oquarr P$ is not weakly open (in the sense of definition on page \pageref{DEF:slabo-otkr-morf}). As a corollary, $Y$ is not representable in the form $Y=(P/E)^\triangledown$ for a subspace $E\subseteq P$ (and hence is not isomorphic in $\Epi^P$ to a cokernel of some morphism $\ph:A\to P$ in ${\tt Ste}$).
\eex

\paragraph{Refinement $\Rf^X F$ of a set $F$ of functionals on a space $X$.}

Theorem \ref{TH:Imp^X(F)-otpechatok-v-kateg-smysle} which we will prove later justifies the following definition.

\bit{

\item[$\bullet$] Let $F$ be a set of linear continuous functionals on a stereotype space $X$. The
{\it refinement} of the set of functionals $F$ on $X$ is a quotient space of $X$, denoted by $\Rf^X F$, or by $\Rf F$, and defined as the injective limit in the category ${\tt Ste}$
    \beq\label{DEF:Cospan_infty}
    \Rf^X F=\Rf F={\tt Ste}\text{-}\rightlim E_i
    \eeq
of the covariant system $\{E_i;\ i\in{\tt Ord}\}$ of quotient spaces of $X$ indexed by ordinal numbers and defined by the following inductive rules:
    \bit{

    \item[0)] the space $E_0$ is the pseudocompletion of the quotient space $X/\Ker F$ (with the usual quotient topology) of $X$ by the common kernel $\Ker F=\bigcap_{f\in F}\Ker f$ of functionals from $F$:
    $$
    E_0=\left(X/\Ker F\right)^\triangledown,
    $$
    after that a set $F_0$ of linear contunous functionals on $E_0$ is defined as the set of extensions to $E_0$ of functionals from $F$ (every functional $f\in F$ vanishes on the common kernel $\Ker F$, so it can be uniquely extended to a linear continuous functional on the quotient space $X/\Ker F$, and then to its pseudocompletion  $E_0=\left(X/\Ker F\right)^\triangledown$),

    \item[1)] if for an ordinal number $j\in{\tt Ord}$ all the spaces $\{E_i;\ i<j\}$ are already defined, then the space $E_j$ is defined as follows:
    \bit{

    \item[---] if $j$ is an isolated ordinal, i.e. $j=i+1$ for some $i$, then $E_j=E_{i+1}$ is defined as the pseudocompletion of the quotient space  $E_i/\Ker F$ (with the usual quotient topology):
    $$
    E_j=E_{i+1}=\left(E_i/\Ker F\right)^\triangledown,
    $$
    after that a set $F_{i+1}$ of linear continuous functionals on $E_{i+1}$ is defined as the set of extensions of functionals from $F_i$,

    \item[---] if $j$ is a limit ordinal, i.e. $j\ne i+1$ for all $i$, then $E_j$ is defined as the injective limit in the category  ${\tt Ste}$ of stereotype spaces of the net $\{E_i;\ i\to j\}$:
    $$
    E_j={\tt Ste}\text{-}\lim_{i\to j}E_i=\l {\tt LCS}\text{-}\lim_{i\to j}E_i\r^\triangledown.
    $$
    after that a set $F_j$ of linear continuous functionals on $E_j$ is defined as the system of functionals which being restricted to every $E_i$ coincide with $F_i$.
     }\eit
    }\eit
Since the transfinite sequence $\{E_i;\ i\in{\tt Ord}\}$ cannot be an injective map from the class ${\tt Ord}$ of all ordinal numbers into the set $\Quot(X)$ of quotient spaces of $X$, it must stabilize, i.e. after some number $i$ all the spaces $E_i$ must coincide together with the topology. As a corollary, the formula \eqref{DEF:Cospan_infty} uniquely defines some quotient space $\Rf F$ of $X$.
}\eit

\bex\label{EX:M-razd-X=>Imp(M)=X}
If a set of functionals $F$ separates elements of $X$ (in other words, the common kernel $\Ker F$ of functionals from $F$ is zero), then the refinement of $F$ on $X$ coincide with $X$:
    $$
  \Ker F=0\qquad\Longrightarrow\quad  \Rf^X F=X.
    $$
    \eex
\bpr
From $\Ker F=\{x\in X:\ \forall f\in F\ f(x)=0\}=0$ we have $E_0=\left(X/\Ker F\right)^\triangledown=X$. As a corollary, all the other spaces $E_i$ coincide with $X$
$$
X=E_0=E_1=...
$$
Thus, $\Rf F=X$.
\epr

\bex
If a set of functionals $F$ is a closed subspace in $X^\star$ (as in a locally convex space), then the refinement of $F$ on $X$ is the open immediate quotient space of $X$ by the common kernel $\Ker F$, i.e. coincides with the pseudocompletion of the locally convex quotient space $X/\Ker F$ with the usual quotient topology:
    $$
\overline{\Sp F}^{X^\star}=F \quad\Longrightarrow\quad    \Rf^X F=(X/\Ker F)^\triangledown
    $$
\eex
\bpr
In this case
$$
E_0=(X/\Ker F)^\triangledown=(X/F^\perp)^\triangledown=\cite[(4.3)]{Akbarov}=(F^\vartriangle)^\star,
$$
hence $\Ker F_0=\{y\in (F^\vartriangle)^\star:\ \forall f\in F\quad f(y)=0\}=0$, and $E_1=E_0/0=E_0$. And further all the spaces $E_i$ coincide with $E_0$.
\epr

The following two theorems are dual to Theorems \ref{TH:Span_infty^X(M)-osubarr-X} and \ref{Y-subarr-X=>Y=Sp_infty(Y)}, and therefore do not require proof.

\btm\label{TH:Cospan_infty^X(M)-oquarr-X} The refinement $\Rf^X F$ of any set of functionals $F\subseteq X^\star$ on a stereotype space $X$ is an immediate quotient space of $X$, to which functionals from $F$ can be continuously extended:
\beq\label{F-subseteq-Cosp_infty-F^X-osubarr-X}
\Rf^X F\overset{\upsilon}{\oquarr} X,\qquad \forall f\in F\quad \exists g\in \l\Rf^X F\r^\star:\quad f=g\circ\upsilon.
\eeq
\etm

\btm\label{Y-oquarr-X=>Y=Cosp_infty(Y)} Every quotient space $Y$ of a stereotype space $X$ is subordinated to the refinement $\Rf^X (Y^\star\circ\upsilon)$ of the system of functionals $Y^\star\circ\upsilon=\{g\circ\upsilon;\ g\in Y^\star\}$ on the space $X$, where  $\upsilon:Y\quarr X$ is the representing epimorphism of $Y$:
\beq\label{Y-quarr-X=>Y-le-Cosp_infty^X(Y^star)}
\upsilon:Y\quarr X\quad\Longrightarrow\quad Y\le \Rf^X (Y^\star\circ\upsilon),
\eeq
and $Y$ is an immediate quotient subspace of $X$, iff $Y$ coincides (as a locally convex space) with this refinement:
\beq\label{Y-oquarr-X=>Y-=-Cosp_infty^X(Y^star)}
\upsilon:Y\oquarr X\quad\Longleftrightarrow\quad Y=\Rf^X (Y^\star\circ\upsilon).
\eeq
\etm

\bcor The representing epimorphism $\upsilon:Y\oquarr X$ of any continuous quotient space $Y$ of a stereotype space $X$ is always relatively open.
\ecor

The following theorem is dual to Theorem \ref{ph(M)-sub-N=>ph:Sp8(M)->Sp8(N)}.

\btm\label{ph(F)-sub-G=>ph:Cosp8(F)->Cosp8(G)} If $\ph:Y\gets X$ is a morphism of stereotype spaces, turning a set of functionals $G\subseteq Y^\star$ into a set of functionals $F\subseteq X^\star$,
$$
G\circ\ph\subseteq F,
$$
then there exists a unique morphism $\e:\Rf^Y G\gets \Rf^X F$ such that the following diagram is commutative:
 $$
 \xymatrix{
 Y\ar[d] & X\ar[l]_{\ph}\ar[d]\\
 \Rf^Y G & \Rf^X F \ar@{-->}[l]_{\e}
 }
 $$
In the special cases:
 \begin{align}
& \begin{Bmatrix}\ph: Y  \quarr  X \\
G\circ\ph\subseteq F \end{Bmatrix}\quad\Longrightarrow\quad \text{$\e$ is an epimorphism},
\label{G-subseteq-F,Y-quarr-X}
\\
& \begin{Bmatrix}\ph:Y  \oquarr  X \\
G\circ\ph\subseteq F \end{Bmatrix}\quad\Longrightarrow\quad \text{$\e$ is an immediate epimorphism},
\label{G-subseteq-F,Y-oquarr-X}
\\
& \begin{Bmatrix}\ph: Y  \oquarr  X \\
G\circ\ph= F \end{Bmatrix}\quad\Longrightarrow\quad \text{$\e$ is an isomorphism}.
\label{G=F,Y-osubarr-X}
 \end{align}
\etm

\btm The refinement $\Rf^X F$ of a set $F\subseteq X^\star$ of functionals on a stereotype space $X$ is a minimal quotient space among immediate quotient spaces of $X$ to which functionals $F$ can be extended. Moreover, every such quotient space $Y$ strongly subordinates $\Rf^X F$:
\beq\label{minimalnost-Cospan_infty(F)}
\forall Y\qquad \Big(F\subseteq Y^\star\ \&\ Y\oquarr X\quad\Longrightarrow\quad \Rf^X F\ole Y\Big).
\eeq
\etm

\bprop\label{PROP:Y-oquarr-X,Z-quarr-X,...=>Z-le-Y}
If $\alpha:Y\oquarr X$ and $\beta:Z\quarr X$, then the condition $Z^\star\circ\alpha\subseteq Y^\star\circ\beta$ implies $Z\le Y$. In a special case, when $Y\osubarr X$ and $Z\osubarr X$, the condition $Z^\star\circ\alpha\subseteq Y^\star\circ\beta$ implies $Z\ole Y$.
\eprop

\subsection{Decompositions, factorizations, envelope and refinement in ${\tt Ste}$.}

\paragraph{Pre-abelian property and basic decomposition in ${\tt Ste}$.}\label{SUBSEC:predabelevost}

Since any two parallel morphisms $\xymatrix{X\ar@/^1ex/[r]^{\ph}\ar@/_1ex/[r]_{\psi} & Y}$ in the category ${\tt Ste}$ of stereotype spaces can be added and subtracted one from another, it is clear that ${\tt Ste}$ is an additive category. In \cite{Akbarov} it was noticed that this category is pre-abelian:

\btm \label{th-4.17} In the category ${\tt Ste}$ of stereotype spaces for each morphism $\ph:X\to Y$ the formulas
\begin{align}\label{DEF:Ker,Coker,...-v-Ste}
& \Ker\ph=\Big(\ph^{-1}(0)\Big)^\vartriangle,
&& \Coker\ph=\Big(Y/ \overline{\ph(X)}\Big)^\triangledown,
&& \Coim\ph=\Big(X/\ph^{-1}(0)\Big)^\triangledown,
&& \Im\ph=\Big(\overline{\ph(X)}\Big)^\vartriangle
\end{align}
define respectively kernel, cokernel, coimage and image. The operation $\ph\mapsto \ph^\star$ of taking dual map establishes the following connections between these objects:
\begin{align}
&  (\ker\ph)^\star=\coker \ph^\star
&&  (\coker\ph)^\star=\ker \ph^\star
&&  (\im\ph)^\star=\coim \ph^\star
&&  (\coim\ph)^\star=\im \ph^\star
\label{eq4.10}
\\
& (\Ker \ph)^\perp{}^\vartriangle=\Im \ph^\star
&&
(\Im\ph)^\perp{}^\vartriangle=\Ker \ph^\star
&&
\Ker \ph=(\Im\ph^\star)^\perp{}^\vartriangle
&&
\Im \ph=(\Ker\ph^\star)^\perp{}^\vartriangle
\label{eq4.11}
\end{align}
\etm

The pre-abelian property of ${\tt Ste}$ implies

\btm Each morphism $\ph:X\to Y$ in ${\tt Ste}$ has basic decomposition \eqref{EX:bazis-razlozh}. The operation
$\ph\mapsto \ph^\star$ of taking dual map establishes the following identities:
\begin{align}
&  (\im\ph)^\star=\coim \ph^\star
&&  (\coim\ph)^\star=\im\ph^\star
\label{svyaz-im-i-coim}
\\
&  (\Im\ph)^\star=\Coim\ph^\star
&&  (\Coim\ph)^\star=\Im\ph^\star
\label{svyaz-Im-i-Coim}
\end{align}
\etm

Formulas \eqref{DEF:Ker,Coker,...-v-Ste} imply

\btm For any morphism of stereotype spaces $\ph:X\to Y$
 \bit{
\item[---] its kernel $\Ker\ph$ and image $\Im\ph$ are closed immediate subspaces (in $X$ and $Y$ respectively),

\item[---] its coimage $\Coim\ph$ and cokernel $\Coker\ph$ are open immediate quotient spaces (of $X$ and $Y$ respectively).
}\eit
\etm

\bex\label{EX:reg(ph)-ne-Bim} There exists a morphism of stereotype spaces $\ph$ such that the reduced morphism $\red\ph$ is not a bimorphism.
\eex
\bpr
Let $E$ be a space from Example \ref{ex-3.23}, i.e. a complete locally convex space with a discontinuous linear functional $f:E\to\C$ which is continuous in the topology of pseudosaturation $E^\vartriangle$ of space $E$. The kernel $F=\Ker f$ of this functional is a closed subspace in the pseudosaturation $E^\vartriangle$ of the space $E$, different from $E^\vartriangle$, but in the space $E$ the subspace $F$ is dense. Since $E$ is complete, we can embed it as a closed subspace into a direct product of Banach spaces, let us denote it by $Y$. Let $\ph:F^\vartriangle\to Y$ be the composition of the injections
$$
F^\vartriangle\subset F\subset E^\vartriangle\subset E\subset Y.
$$
Since $F$ is a closed subspace in the pseudocomplete space $E^\vartriangle$, it is pseudocomplete. Hence, its pseudosaturation $F^\vartriangle$ is a stereotype space. On the other hand, $Y$ is a direct product of Banach spaces, therefore it is stereotype as well. Finally, since $\ph$ is an injection, its kernel is zero, hence its coimage coincides with $F^\vartriangle$:
$$
\Coim\ph=F^\vartriangle,
$$
On the other hand, the image of $\ph$ is the pseudosaturation of the space $\ph(F^\vartriangle)=F$ in $Y$, i.e. pseudosaturation of the space $E$:
$$
\Im\ph=\Big(\overline{\ph(F^\vartriangle)}^Y\Big)^\vartriangle=E^\vartriangle.
$$
Thus, the reduced morphism $\red\ph$ is just the enclosure
$$
F^\vartriangle\subset E^\vartriangle,
$$
and this cannot be a bimorphism, since $F^\vartriangle$ is closed in $E^\vartriangle$, but not equal to $E^\vartriangle$. Diagram  \eqref{EX:bazis-razlozh} for $\ph$ takes the following form:
$$
\begin{diagram}
\node{F^\vartriangle} \arrow{e,t}{\varphi} \arrow{s,l}{\coim\varphi} \node{Y}
\\
\node{F^\vartriangle} \arrow{e,t}{\red\varphi} \node{E^\vartriangle}
\arrow{n,r}{\im\varphi}
\end{diagram}\quad.
$$
 \epr

\bcor
The category ${\tt Ste}$ of stereotype spaces is not quasi-abelian in the sense of J.-P.~Schneiders \cite{Schneiders}.
\ecor
\bpr
Example \ref{EX:reg(ph)-ne-Bim} contradicts to \cite[Corollary 1.1.5]{Schneiders}.
\epr

\paragraph{Nodal decomposition in ${\tt Ste}$.} In \cite[Theorem 4.21]{Akbarov} it was noticed that the category ${\tt Ste}$ is complete. On the other hand, from Corollaries \ref{COR:Sub(V)-v-Ste} and \ref{COR:Quot(V)-v-Ste} it follows that ${\tt Ste}$ is well-powered and co-well-powered. Together with the existence of basic decomposition, this by Theorem \ref{TH:faktorizatsija-v-lok-maloi-kategorii} means that ${\tt Ste}$ is a category with nodal decomposition:

\btm\label{TH:uzlovoe-razlozhenie-v-Ste} In the category ${\tt Ste}$ of stereotype spaces each morphism $\ph:X\to Y$ has nodal decomposition
\eqref{DEF:oboznacheniya-dlya-uzlov-razlozh}. The operation $\ph\mapsto \ph^\star$ of taking dual map establishes the following identities:
\begin{align}
&  (\im_\infty\ph)^\star=\coim_\infty \ph^\star
&&  (\coim_\infty\ph)^\star=\im_\infty \ph^\star
\label{eq4.10-infty}
\\
&  (\Im_\infty\ph)^\star=\Coim_\infty \ph^\star
&&  (\Coim_\infty\ph)^\star=\Im_\infty \ph^\star
\label{eq4.10-infty-BIG}
\end{align}
\etm

As we noticed above, the basic and the nodal decomposition are connected with each other through diagram \eqref{svayz-bazisnogo-i-uzlovogo-razlozheniya}:
$$
\xymatrix  @R=2.5pc @C=4.0pc
{
X\ar[rrr]^{\ph}\ar[rd]^{\coim_\infty\ph}\ar[dd]_{\coim\ph} & & & Y \\
& \Coim_\infty\ph\ar[r]^{\red_\infty\ph} & \Im_\infty\ph\ar[ru]^{\im_\infty\ph}\ar@{-->}[rd]_{\tau} & \\
\Coim\ph\ar[rrr]^{\red\ph}\ar@{-->}[ru]_{\sigma} & & & \Im\ph \ar[uu]_{\im\ph}
}
$$
where morphisms $\sigma$ and $\tau$ are uniquely defined (by $\ph$).

\bex For a morphism described in Example \ref{EX:reg(ph)-ne-Bim} diagram \eqref{svayz-bazisnogo-i-uzlovogo-razlozheniya} has the from
$$
\xymatrix  @R=2.5pc @C=4.0pc
{
X^\vartriangle\ar[rrr]^{\ph}\ar[rd]^{\coim_\infty\ph}\ar[dd]_{\coim\ph} & & & Y \\
& X^\vartriangle\ar[r]^{\red_\infty\ph} & X^\vartriangle\ar[ru]^{\im_\infty\ph}\ar@{-->}[rd]_{\tau} & \\
X^\vartriangle\ar[rrr]^{\red\ph}\ar@{-->}[ru]_{1_{X^\vartriangle}} & & & E^\vartriangle \ar[uu]_{\im\ph}
}
$$
This shows that $\tau$ is not necessary an isomorphism. If we consider the dual map $\ph^\star$, we can conclude that $\sigma$ is not necessarily an isomorphism as well.
\eex

\btm\label{TH:Im_infty-ph=Sp_infty-ph(X)} For any morphism of stereotype spaces $\ph:X\to Y$
\bit{
\item[---] its nodal image $\Im_\infty\ph$ coincides with the envelope in $Y$ of its set of values $\ph(X)$:
    \beq\label{Im_infty-ph=Sp_infty-ph(X)}
    \Im_\infty\ph=\Env^Y\ph(X)
    \eeq

\item[---] its nodal coimage $\Coim_\infty\ph$ coincides with the refinement on $X$ of a set of functionals $\ph^\star(Y^\star)$:
    \beq\label{Coim_infty-ph=Cosp_infty-ph*(Y*)}
    \Coim_\infty\ph=\Rf^X\ph^\star(Y^\star)
    \eeq
    }\eit
\etm
\bpr
By Remark  \ref{REM:struktura-uzlov-razlozh-v-kateg-s-nulem}, the nodal image $\Im_\infty\ph$ is a projective limit of a sequence of ``usual'' images $\Im\ph^i$ of transfinite system of morphisms, defined by the inductive rule $\ph^{i+1}=\red\ph^i$. And each space $\Im\ph^i$ exactly coincides with the space $E_i$ from the definition of the envelope $\Env^YM$ in $Y$ of the set $M=\ph(X)$.

Similarly, the nodal coimage $\Coim_\infty\ph$ is an injective limit of transfinite system of ``usual'' coimages $\Coim\ph^i$, and each such space coincide with the space $E_i$ from the definition of the refinement $\Rf^X F$ on $X$ of the set of functionals $F=\ph^\star(Y^\star)$.
\epr

\paragraph{Factorizations in ${\tt Ste}$.}

Recall that by definition on page \pageref{paragraph:faktorizatsija}, a {\it factorization} of a morphism $X\overset{\ph}{\longrightarrow}Y$ is its representation as a composition $\ph=\mu\circ\e$ of an epimorphism $\e$ and a monomorphism $\mu$. Theorem \ref{TH:faktorizatsija-v-lok-maloi-kategorii} implies

\btm\label{TH:faktorizatsija-v-Ste} In the category ${\tt Ste}$ of stereotype spaces
 \bit{
\item[(i)] each morphism $\ph$ has a factorization,

\item[(ii)] among all factorizations of $\ph$ there is a minimal one $(\e_{\min},\mu_{\min})$ and a maximal one $(\e_{\max},\mu_{\max})$, i.e. each factorization $(\e,\mu)$ lies between them:
    $$
    (\e_{\min},\mu_{\min})\le (\e,\mu)\le (\e_{\max},\mu_{\max})
    $$
 }\eit
\etm

\paragraph{Characterization of strong morphisms in ${\tt Ste}$.}

\btm\label{TH:strogie-monomorfizmy-v-Ste} In the category ${\tt Ste}$ for a morphism $\mu:Z\to X$ the following conditions are equivalent:
 \bit{

 \item[$(i)$] $\mu$ is an immediate monomorphism,

 \item[$(i)'$] in diagram \eqref{predstavlenie-monomorfizma} the space $Y$ is an immediate subspace in $X$,

 \item[$(ii)$] $\mu$ is a strong monomorphism,

 \item[$(ii)'$] in diagram \eqref{predstavlenie-monomorfizma} the morphism $\sigma$ is a strong monomorphism,

 \item[$(iii)$] $\mu\cong\im_\infty\mu$,

 \item[$(iv)$] $\coim_\infty\mu$ and $\red_\infty\mu$ are isomorphisms.

 }\eit
\etm
\bpr
The equivalences $(i)\Longleftrightarrow(ii)\Longleftrightarrow(iii)\Longleftrightarrow(iv)$ follow from Theorem \ref{TH:harakter-strogih-morf-v-kat-s-uzlov-razlozh}, since ${\tt Ste}$ is a category with nodal decomposition. In addition, Proposition  \ref{PROP:mu-cong-neposr-podpr-v-Ste} imply equivalences $(i)\Longleftrightarrow(i)'$ and $(ii)\Longleftrightarrow(ii)'$.
\epr

The dual proposition is proved by analogy:

\btm\label{TH:strogie-epimorfizmy-v-Ste} In the category ${\tt Ste}$ for a morphism $\e:Z\to X$ the following conditions are equivalent:
 \bit{

 \item[$(i)$] $\e$ is an immediate epimorphism,

 \item[$(i)'$] in diagram \eqref{predstavlenie-epimorfizma} the space $Y$ is an immediate quotient space for $X$,

 \item[$(ii)$] $\e$ is a strong epimorphism,

 \item[$(ii)'$] in diagram \eqref{predstavlenie-epimorfizma} the morphism $\pi$ is a strong epimorphism,

 \item[$(iii)$] $\e\cong\coim_\infty\e$,

 \item[$(iv)$] $\im_\infty\mu$ and $\red_\infty\mu$ are isomorphisms.
 }\eit
\etm

\paragraph{Envelope and refinement in ${\tt Ste}$.}

Since the category ${\tt Ste}$ is complete, well-powered, co-well-powered and has nodal decomposition, this implies the existence of some envelopes and refinements in ${\tt Ste}$.

\btm\label{TH:obolochki-i-nachinki-otn-klassa-morfizmov-v-Ste} In the category ${\tt Ste}$ of stereotype spaces

 \bit{
\item[(a)] each space $X$ has envelopes in the classes $\Epi$ of all epimorphisms and $\SEpi$ of all strong epimorphisms with respect to arbitrary class of morphisms $\varPhi$, among which there is at least one going from $X$; in addition,
     \bit{
\item[(i)] if $\varPhi$ differs morphisms on the outside in ${\tt Ste}$, then the envelope in $\Epi$ is also an envelope in the class $\Bim$ of all bimorphisms:
    $$
    \env_\varPhi^{\Epi} X=\env_\varPhi^{\Bim} X,
    $$

\item[(ii)] if $\varPhi$ differs morphisms on the outside and is an ideal in ${\tt Ste}$, then the envelope in $\Epi$ is also an envelope in the class $\Bim$ of all bimorphisms, and in any other class $\varOmega$ which contains $\Bim$ (for example, in the class $\Mor$ of all morphisms):
    $$
    \env_\varPhi^{\Epi} X=\env_\varPhi^{\Bim} X=\env_\varPhi^\varOmega X=\env_\varPhi X,\qquad \varOmega\supseteq\Bim.
    $$
}\eit

\item[(b)] each space $X$ has refinements in the classes $\Mono$ of all monomorphisms and $\SMono$ of all strong monomorphisms by means of arbitrary class of morphisms $\varPhi$, among which there is at least one coming to $X$; in addition,
      \bit{
\item[(i)] if $\varPhi$ differs morphisms on the inside in ${\tt Ste}$, then the refinement in $\Mono$ is also a refinement in the class $\Bim$ of all bimorphisms:
    $$
    \rf_\varPhi^{\Mono} X=\rf_\varPhi^{\Bim} X.
    $$

\item[(ii)] if $\varPhi$ differs morphisms on the inside and is a left ideal in ${\tt Ste}$, then the refinement in $\Mono$ is a refinement in the class  $\Bim$ of all bimorphisms, and of any other class $\varOmega$ which contains $\Bim$ (for example, in the class $\Mor$ of all morphisms):
    $$
    \rf_\varPhi^{\Mono} X=\rf_\varPhi^{\Bim} X=\rf_\varPhi^\varOmega X=\rf_\varPhi X,\qquad \varOmega\supseteq\Bim.
    $$

 }\eit
 }\eit
\etm
\bpr Due to duality it is sufficient to prove (a). Let $X$ be a stereotype space, and $\varPhi$ a class of morphisms, which contains at least one going from $X$. Then the envelopes $\env_\varPhi^{\Epi} X$ and $\env_\varPhi^{\SEpi} X$ exist by $5^\circ$ on p.\pageref{5^0:obolochka-env_Phi^Epi-otn-klassa-morphizmov}. Suppose now that $\varPhi$ differs morphisms on the outside in ${\tt Ste}$. Then by Theorem \ref{TH:Phi-razdel-moprfizmy} the existence of envelope $\env_\varPhi^{\Epi} X$ automatically implies the existence of envelope $\env_\varPhi^{\Bim} X$ and their equality: $\env_\varPhi^{\Epi} X=\env_\varPhi^{\Bim} X$. Finally, suppose that $\varPhi$ differs morphisms on the outside in ${\tt Ste}$ and in addition is a right ideal. Then by Theorem \ref{TH:Phi-razdel-moprfizmy-*} the existence of envelope $\env_\varPhi^{\Bim} X$ (which is already proved) implies that for any class $\varOmega\supseteq\Bim$ the envelope $\env_\varPhi^\varOmega X$ also exists, and these envelopes coincide: $\env_\varPhi^{\Bim} X=\env_\varPhi^\varOmega X$.
\epr

\btm\label{TH:Env^X(M)-obolochka-v-kateg-smysle} The envelope $\Env^X M$ of a set $M$ in a stereotype space $X$ coincides with the envelope of the space\footnote{We use here the notations of \cite[p.478]{Akbarov-stein-groups}.} $\C_M$ in the class $\Epi$ of all epimorphisms of the category ${\tt Ste}$ with respect to the morphism $\ph:\C_M\to X$, $\ph(\alpha)=\sum_{x\in M}\alpha_x\cdot x$,
$$
\Env^X M=\Env_{\ph}^{\Epi}\C_M.
$$
\etm
\bpr This follows from $1^\circ$ on p.\pageref{1^0:obolochka-env_ph^Epi} and from Theorem  \ref{TH:Im_infty-ph=Sp_infty-ph(X)}:
$$
\Env_{\ph}^{\Epi}\C_M=\eqref{obolochka--env_ph^Epi}=
\Im_\infty\ph=\eqref{Im_infty-ph=Sp_infty-ph(X)}=\Env^X\ph(\C_M)=\Env^X \Sp M=\Env^X M.
$$
\epr

\btm\label{TH:Imp^X(F)-otpechatok-v-kateg-smysle} The refinement $\Rf^X F$ of a set $F$ of functionals on a stereotype space $X$ coinsides with the refinement of the space\footnote{The notations of \cite[p.477]{Akbarov-stein-groups} are used here.} $\C^F$ in the class $\Mono$ of all monomorphisms in the category ${\tt Ste}$ by means of the morphism $\ph:X\to\C^F$, $\ph(x)^f=f(x)$, $f\in F$
$$
\Rf^X F=\Rf_{\ph}^{\Mono}\C^F.
$$
\etm
\bpr This follows from $1^\circ$ on p.\pageref{1^0:otpechatok-imp_ph^Mono} and from Theorem  \ref{TH:Im_infty-ph=Sp_infty-ph(X)}:
$$
\Rf_{\ph}^{\Mono}\C^F=\eqref{otpechatok-imp_ph^Mono}=\Coim_\infty\ph=
\eqref{Coim_infty-ph=Cosp_infty-ph*(Y*)}=\Rf^X\ph^\star(Y^\star)=
\Rf^X\Sp F=\Rf^XF.
$$
\epr

\subsection{On homologies in ${\tt Ste}$}

As is known, in the homology theory, in opposition to the well-established methods of Abelian categories, there have always been attempts to find alternative approaches, where it is considered desirable to get rid of the Abelian property and even of the additivity with the aim to cover the widest spectrum of situations (one can make an impression of this by the works \cite{Raikov}, \cite{Palamodov}, \cite{Yakovlev}, \cite{Helemskii-homology}, \cite{Grandis}, \cite{Generalov}, \cite{Schneiders}, \cite{Kuzminov-Shvedov}, \cite{Rump}, \cite{Borceux-Bourn}, \cite{Buhler}, \cite{Kopylov-Wegner}, \cite{Sieg-Wegner}, \cite{Kopylov}). We hope that the following effect will be interesting in this connection: in the (non-Abelian, but pre-Abelian) category ${\tt Ste}$ of stereotype spaces the standard definition of homology breaks up into two non-equivalent notions. Let us start with the following definition (taken from \cite{Kopylov}):

\bit{
\item[$\bullet$]
Suppose in a pre-Abelian category ${\tt K}$ we have a pair of morphisms $X\stackrel{\ph}{\to}Y\stackrel{\psi}{\to}Z$ which form a complex:
$$
\psi\circ\ph=0.
$$
By the definitions of kernel and cokernel, this equality defines two natural morphisms $X\stackrel{\ph^{\Ker\psi}}{\longrightarrow} \Ker\psi$ and $\Coker\ph\stackrel{\psi_{\Coker\ph}}{\longrightarrow}Z$ such that the following diagram is commutative
$$
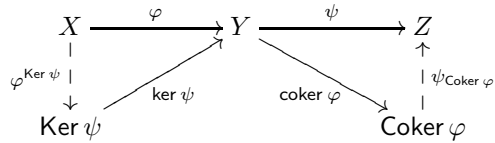
\xymatrix  @R=2pc @C=3.5pc
{
 X\ar[r]^{\ph}\ar@{-->}[d]_{\ph^{\Ker\psi}} & Y\ar[r]^{\psi}\ar[rd]_{\coker\ph} & Z \\
 \Ker\psi\ar[ru]_{\ker\psi} &  & \Coker\ph \ar@{-->}[u]_{\psi_{\Coker\ph}}
}
$$
The cokernel of the morphism $\ph^{\Ker\psi}$ is called the {\it left homology} of the pair $(\ph,\psi)$ and is denoted by
     \beq\label{H_-(psi:ph)}
     \H_-(\psi:\ph)=\Coker(\ph^{\Ker\psi}).
     \eeq
and the kernel of the morphism $\psi_{\Coker\ph}$ is called the {\it right homology} of the pair $(\ph,\psi)$ and is denoted by
     \beq\label{H_+(psi:ph)}
     \H_+(\psi:\ph)=\Ker(\psi_{\Coker\ph}).
     \eeq
}\eit

The following observation belongs to folklore:

\bprop\label{PROP:h(psi:ph)} In a pre-Abelian category ${\tt K}$ for any pair of morphisms $X\stackrel{\ph}{\to}Y\stackrel{\psi}{\to}Z$ forming a complex, $\psi\circ\ph=0$, there exists a unique morphism $\h(\psi:\ph):\H_-(\psi:\ph)\to \H_+(\psi:\ph)$ such that the following diagram is commutative:
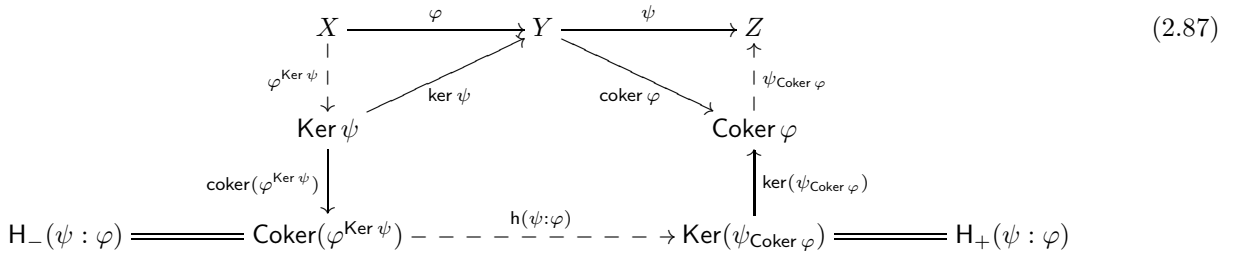
\beq\label{DIAGR:h(psi:ph)}
\xymatrix  @R=2pc @C=3.5pc
{
 & X\ar[r]^{\ph}\ar@{-->}[d]_{\ph^{\Ker\psi}} & Y\ar[r]^{\psi}\ar[rd]_{\coker\ph} & Z & \\
 & \Ker\psi\ar[ru]_{\ker\psi}\ar[d]_{\coker(\ph^{\Ker\psi})} &  & \Coker\ph \ar@{-->}[u]_{\psi_{\Coker\ph}} & \\
\H_-(\psi:\ph)\ar@{=}[r] & \Coker(\ph^{\Ker\psi})\ar@{-->}[rr]^{\h(\psi:\ph)} & & \Ker(\psi_{\Coker\ph}) \ar[u]_{\ker(\psi_{\Coker\ph})}\ar@{=}[r] & \H_+(\psi:\ph)
}
\eeq
\eprop

In each autodual category (for instance, in ${\tt Ste}$) the purely categorial duality reasoning gives the following identities:
 \beq
 \H_+(\psi:\ph)^\star\cong \H_-(\ph^\star:\psi^\star),\qquad
  \H_-(\psi:\ph)^\star\cong \H_+(\ph^\star:\psi^\star)
 \eeq

\bex\label{EX:hat(H)-ne-check(H)} In the category ${\tt Ste}$ of stereotype spaces the morphism $\xymatrix{\H_-(\psi:\ph)\ar[r]^{\h(\psi:\ph)} & \H_+(\psi:\ph)}$ is not always an epimorphism.
\eex
\bpr
Let $E$ be a space from Example \ref{ex-3.23}, i.e. a complete locally convex space with a discontinuous linear functional $f:E\to\C$, which is continuous in the topology of pseudosaturation $E^\vartriangle$. The kernel $F=\Ker f$ of this functional is a dense subspace in $E$, but in the space $E^\vartriangle$ it is a closed subspace, different from $E^\vartriangle$ (since $f\ne 0$). As a corollary, the natural enclosure $\sigma:F\to E$ is dense (i.e. has a dense image in $E$), but its pseudosaturation $\sigma^\vartriangle:F^\vartriangle\to E^\vartriangle$ does not have this property.

Let us represent $E$ as a closed subspace in a stereotype space $Y$ (with the topology inherited from $Y$; for example, we can consider the system of Banach quotient spaces of $E$ and say that $Y$ is the direct product of these spaces). Let
$$
\ph:F^\vartriangle\to E^\vartriangle\to Y
$$
be the corresponding composition of monomorphisms, and
$$
\psi:Y\to (Y/E^\vartriangle)^\triangledown
$$
the corresponding epimorphism. Then, first,
$$
\Ker\psi=E^\vartriangle
$$
$$
\Downarrow
$$
$$
\Im\ph^{\Ker\psi}=\Big(\overline{\ph(F)}^{E^\vartriangle}\Big)^\vartriangle=\Big(\overline{F}^{E^\vartriangle}\Big)^\vartriangle=F^\vartriangle
$$
$$
\Downarrow
$$
$$
\Coker(\ph^{\Ker\psi})=(E^\vartriangle/F^\vartriangle)^\triangledown\cong\C^\triangledown=\C.
$$
And, second,
$$
\Im\ph=\Big(\overline{\ph(F)}^Y\Big)^\vartriangle=\Big(\overline{F}^Y\Big)^\vartriangle=E^\vartriangle,
$$
$$
\Downarrow
$$
$$
\Coker\ph=(Y/E^\vartriangle)^\triangledown
$$
$$
\Downarrow
$$
$$
{\psi_{\Coker\ph}}=1_{(Y/E^\vartriangle)^\triangledown}
$$
$$
\Downarrow
$$
$$
\Ker({\psi_{\Coker\ph}})=0
$$
As a result diagram \eqref{DIAGR:h(psi:ph)} takes the form
$$
\xymatrix  @R=2pc @C=3.5pc
{
 & F^\vartriangle\ar[r]^{\ph}\ar@{-->}[d]_{\ph^{\Ker\psi}} & Y\ar[r]^{\psi}\ar[rd]_{\coker\ph} &  (Y/E^\vartriangle)^\triangledown & \\
 &  E^\vartriangle\ar[ru]_{\ker\psi}\ar[d]_{\coker(\ph^{\Ker\psi})} &  & (Y/E^\vartriangle)^\triangledown\ar@{-->}[u]_{\psi_{\Coker\ph}} & \\
\H_-(\psi:\ph)\ar@{=}[r] & (E^\vartriangle/F^\vartriangle)^\triangledown\cong\C\ar@{-->}[rr]^{\h(\psi:\ph)} & & 0 \ar[u]_{\ker(\psi_{\Coker\ph})}\ar@{=}[r] & \H_+(\psi:\ph)
}
$$
and clearly, $\h(\psi:\ph)$ cannot be an isomorphism.
 \epr

\section{The category of stereotype algebras ${\tt Ste}^\circledast$}
\label{SEC:proj-ster-algebry}

\subsection{Stereotype algebras and stereotype modules}

\paragraph{Stereotype algebras.}

A stereotype space $A$ over $\C$ is called a {\it stereotype algebra}, if $A$ is endowed with a structure of unital associative algebra over $\C$, and the multiplication is a continuous bilinear form in the  sense of the definition on p.\pageref{DEF:bilin-otobrazhenie}: for any compact set $K$ in $A$ and for any neighborhood of zero $U$ in $A$ there exists a neighborhood of zero $V$ in $A$ such that
$$
K\cdot V\subseteq U\quad \& \quad V\cdot K\subseteq U.
$$
This is equivalent to the fact that $A$ is a monoid in the category ${\tt Ste}$ of stereotype spaces with respect to  the tensor product $\circledast$ (defined in \eqref{eq7.1}). Certainly, each stereotype algebra $A$ is a topological algebra (but not vice versa). The class of all stereotype algebras is denoted by ${\tt Ste}^\circledast$. It is a category, where morphisms are linear, continuous, multiplicative and preserving identity maps $\ph :A\to B$.

In contrast to the category ${\tt Ste}$ of stereotype spaces, the category ${\tt Ste}^\circledast$ of stereotype algebras is not additive. In addition, in ${\tt Ste}^\circledast$ there arise an asymmetry between monomorphisms and epimorphisms, since epimorphisms are not inherited from ${\tt Ste}$:
\bit{
 \item[---]\label{EX:mono-v-Ste^circledast} a morphism $\ph :A\to B$ of stereotype algebras is a monomorphism, iff $\ph$ is an injective map (i.e. a monomorphism of stereotype spaces).

 \item[---]\label{EX:epi-v-Ste^circledast} on the other hand, an epimorphism $\ph :A\to B$ of stereotype algebras not necessarily have dense image in $B$ (i.e., not necessarily is an epimorphism of stereotype spaces). A counterexample is the enclosure of the algebra $\mathcal{P}(\C)$ of polynomials on $\C$ into the algebra  $\mathcal{P}(\C^\times)$ of Laurent polynomials on $\C^\times=\C\setminus\{0\}$ (both algebras are endowed with the strongest locally convex topology).
}\eit

The following lemma will be useful in the further considerations:

\blm\label{LM:o-plotnoi-podalgebre} Let $A$ and $B$ be topological algebras (with the separately continuous multipication), and $\ph:A\to B$ -- a linear continuous map, which is multiplicative on some dense subalgebra $A_0$ in $A$:
$$
\ph(x\cdot y)=\ph(x)\cdot\ph(y),\qquad x,y\in A_0.
$$
Then $\ph$ is multiplicative on $A$:
$$
\ph(x\cdot y)=\ph(x)\cdot\ph(y),\qquad x,y\in A.
$$
\elm
\bpr
For any $x,y\in A$ we find nets $x_i,y_j\in A_0$ such that
$$
x_i\underset{i\to\infty}{\longrightarrow}x,\qquad y_j\underset{j\to\infty}{\longrightarrow}y
$$
and then we have:
$$
\ph(x\cdot y)\underset{\infty\gets j}{\longleftarrow}\ph(x\cdot y_j)\underset{\infty\gets i}{\longleftarrow}
\ph(x_i\cdot y_j)=\ph(x_i)\cdot\ph(y_j)\underset{i\to\infty}{\longrightarrow}\ph(x_i)\cdot\ph(y)
\underset{j\to\infty}{\longrightarrow}\ph(x)\cdot\ph(y).
$$
\epr

Let us give some examples of stereotype algebras. First, two abstract example.

\begin{ex}\label{ex-10.1} {\sl Fr\'{e}chet algebras.} For a Fr\'{e}chet space $A$ being a stereotype algebra is equivalent to the joint continuity of multiplication. Hence, {\it each unital Fr\'{e}chet algebra is a stereotype algebra}.
\end{ex}

\begin{ex}\label{ex-10.2} {\sl Operator algebra ${\mathcal L}(X)$.}
Theorem \ref{th-6.13} implies that for any stereotype space $X$ the space ${\mathcal L}(X)=X\oslash X$ of linear continuous maps $\ph :X\to X$ is a stereotype algebra with respect to the composition $\ph \circ \psi$.
\end{ex}

After that the four functional algebras.

\begin{ex}\label{ex-10.3} {\sl Algebra of continuous functions $\mathcal{C}(M)$ on a paracompact locally compact space $M$.} Let us recall that a topological space $M$ is said to be {\it $\sigma$-compact}, if it is a union of a countable system of compact sets:
$$
M=\bigcup_{n=1}^\infty S_n
$$
($S_n$ are compact sets in $M$). For locally compact spaces $M$ this condition is equivalent to the {\it Lindl\"{o}f property}: every open covering of $M$ has a countable subcovering (cf. \cite[3.8.C(b)]{Engelking}). As a corollary, if $M$ is a {\it Lindl\"{o}f space} (i.e. has the Lindl\"{o}f property) and is locally compact, then the space $\mathcal{C}(M)$ of continuous functions $u:M\to \C$ will be a Fr\'{e}chet space with respect to the topology of uniform convergence on compact sets $S\subseteq M$.

Consider a more general class of topological spaces. Let $M$ be a paracompact locally compact topological space. Then it can be decomposed into a direct sum
$$
M=\coprod_{i\in I} M_i
$$
of Lindl\"{o}f locally compact spaces $M_i$ (\cite[Theorem 5.1.27]{Engelking}). Therefore the space $\mathcal C(M)$\label{C(M)} of continuous functions $u:M\to\C$ (with the topology of uniform convergence on compact sets $S\subseteq M$) is a stereotype space, as a direct product of Fr\'{e}chet spaces:
$$
\mathcal{C}(M)=\prod_{i\in I} \mathcal{C}(M_i)
$$
Certainly, $\mathcal C(M)$ is an algebra with respect to the pointwise multiplication
$$
  (u\cdot v)(t)=u(t)\cdot v(t)
$$
Easy to check that this is a continuous bilinear form. Hence, $\mathcal{C}(M)$ is a stereotype algebra.
\end{ex}

\begin{ex}\label{ex-10.4} {\sl Algebra of smooth functions $\mathcal{E}(M)$ on a smooth manifold $M$} (with the usual topology of uniform convergence on compact sets with any derivative) is a stereotype algebra (with the usual pointwise multiplication).
\end{ex}

\begin{ex}\label{ex-10.5} {\sl Algebra of holomorphic functions $\mathcal{O}(M)$ on a Stein manifold $M$} (with the topology of uniform convergence on compact sets in $M$) is a stereotype algebra (with pointwise multiplication).
\end{ex}

\begin{ex}\label{ex-10.6} {\sl Algebra of polynomials (i.e. regular functions) $\mathcal{P}(M)$ on an affine algebraic variety $M$} (with the strongest locally convex topology) is a stereotype algebra (with pointwise multiplication).
\end{ex}

Finally, the four group algebras,

\begin{ex}\label{ex-10.7} {\sl Algebra of measures $\mathcal{C}^\star (G)$ on a locally compact group $G$.} As is known, each locally compact group $G$ is paracompact \cite[2.8.13]{48}, hence the space $\mathcal{C}(G)$ of continuous functions on $G$ (with the topology of uniform convergence on compact sets in $G$) can be considered as a special case of Example \ref{ex-10.3}, and what is important for us, $\mathcal{C}(G)$ is stereotype. Its dual space  $\mathcal{C}^\star (G)$ consists of measures with compact support on $G$. The {\it convolution of measures} $\alpha, \beta \in \mathcal{C}^\star (G)$ is defined by the formula
\beq
\alpha * \beta (u)=(\alpha\otimes \beta)(w)\Big|_{w(s,t)=u(s\cdot t)}=
\int_G \l \int_G u(s\cdot t) \d \alpha(s) \r \d \beta(t)= \int_G \l \int_G
u(s\cdot t) \d \beta(t)\r \d \alpha(s) \label{eq10.1}
\eeq
This operation is associative and has unit (this is the delta-functional $\delta^{1_G}$ of the unit in $G$). In addition it is continuous as a bilinear map, so {\it the space  of measures $\mathcal{C}^\star (G)$ on a locally compacty group $G$ is a stereotype algebra with convolution $(\alpha, \beta) \mapsto \alpha * \beta$ as multiplication (and with $\delta^{1_G}$ as unit)}.
\end{ex}

\begin{ex}\label{ex-10.8} {\sl Algebra of distributions $\mathcal{E}^\star (G)$ on a Lie group $G$.} Let $G$ be a real Lie group \cite{Warner,Vinberg}. Consider the space $\mathcal{E}^\star (G)$ of distributions with compact support on
$G$ (i.e. the space dual to $\mathcal{E}(G)$ from Example \ref{ex-10.4}).

The {\it convolution of distributions} $\alpha, \beta \in \mathcal E^\star (G)$
is defined by formula \eqref{eq10.1}. {\it The space of distributions $\mathcal{E}^\star (G)$ on
$G$ is a stereotype algebra with the convolution $(\alpha, \beta) \mapsto \alpha * \beta$ as multiplication (and with $\delta^{1_G}$ as unit)}.
\end{ex}

\begin{ex}\label{ex-10.9} {\sl Algebra of analytic functionals $\mathcal{O}^\star (G)$ on a Stein group $G$.} Let $G$ be a {\it Stein group}, i.e. a complex Lie group \cite{Bourbaki-LieAlg}, which is a Stein manifold \cite{Shabat}. Consider the space $\mathcal{O}^\star(G)$ of analytic functionals on $G$ (i.e. the space dual to the space  $\mathcal{O}(G)$ from Example \ref{ex-10.5}).

The {\it convolution of analytic functionals} $\alpha, \beta \in \mathcal{O}^\star (G)$ is defined by formula \eqref{eq10.1}. {\it The space of analytic functionals $\mathcal{O}^\star (G)$ on $G$ is a stereotype algebra with convolution $(\alpha, \beta) \mapsto \alpha * \beta$ as multiplication (and with $\delta^{1_G}$ as unit)}.
\end{ex}

\begin{ex}\label{ex-10.10} {\sl Algebra of currents $\mathcal{P}^\star (G)$ on an affine algebraic group $G$.} Recall some facts from the theory of algebraic groups \cite{Vinberg}. The general linear group $GL(n,\C)$ is a basic open subset in the vector space $L(n,\C)$, therefore it can be represented as a closed in the Zariski topology subset in some affine algebraic space $\C^m$. This means that $GL(n,\C)$ is an affine algebraic variety. Its polynomials (regular functions) have the form
\beq
u(g)=P(g)/D(g)^k \label{eq10.3}
\eeq
where $D(g)$ is the determinant of the matrix $g\in L(n,\C)$, $k$ belongs to $\Bbb N$, and $P$ is a polynomial on $L(n,\C)$ \cite{Vinberg}.

Let now $G$ be an affine algebraic group, i.e. a closed in the Zariski topology subgroup in $GL(n,\C)$ \cite{Vinberg}, or, equivalently, a set of common zeroes of a system of functions $u:GL(n,\C)\to \C$ of the form \eqref{eq10.3}, which is closed under the group operation in $GL(n,\C)$. Since $G$ is a closed subset in  $GL(n,\C)$, it is an affine variety.

Therefore the space $\mathcal{P}(G)$ of polynomials on $G$ is a special case of the general construction from Example \ref{ex-10.6}. In this case $\mathcal P(G)$ consists of functions $v:G\to \C$ which can be extended to functions $u:GL(n,\C)\to \C$ of the form \eqref{eq10.3}. The dual space $\mathcal{P}^\star(G)$ consists of linear (and automatically continuous) functionals $f:\mathcal{P}(G) \to \C$, called {\it currents} on $G$.

The {\it convolution of currents} $\alpha, \beta \in \mathcal{P}^\star (G)$ is defined by the formula \eqref{eq10.1}. {\it The space of currents $\mathcal{P}^\star(G)$ on an affine algebraic group $G$ is a stereotype algebra with convolution $(\alpha, \beta) \mapsto \alpha* \beta$ as multiplication (and $\delta^{1_G}$ as unit)}.
\end{ex}

\paragraph{Stereotype modules.}

A stereotype space $X$ over $\C$ with a given structure of left (or right) $A$-module is called a {\it stereotype $A$-module}, if the multiplication by elements of $A$ is a continuous bilinear form in the sense of definition on page \pageref{DEF:bilin-otobrazhenie}. Theorem \ref{th-7.3} implies that $X$ is a stereotype (left) module over $A$ if and only if the multiplication $\mu$ by elements of $A$ can be continuously factored through the projective stereotype tensor product
$$
\begin{diagram}
\node{A\times X} \arrow[2]{e,t}{} \arrow{se,b}{\mu} \node[2]{A\circledast X}
\arrow{sw,b,--}{}
\\
\node[2]{X}
\end{diagram}.
$$

\begin{ex}\label{ex-11.1} Each stereotype space $X$ is a stereotype left module over the stereotype algebra $\mathcal  L(X)$ (from Example \ref{ex-10.2}).
\end{ex}

\btm [on representation]\label{th-11.2}
Let $A$ be a stereotype algebra. A stereotype space $X$ with the structure of left (right) $A$-module is a stereotype $A$-module if and only if the operation of multiplication by elements of $A$ defines a continuous homomorphism (respectively, antihomomorphism) of $A$ into $\mathcal L(X)$.
\etm

The classes $_A{\tt Ste}$ and ${\tt Ste}_A$ of left and right stereotype modules over a stereotype algebra $A$ form categories with continuous $A$-linear maps as morphisms.

\medskip
\centerline{\bf Properties of the categories $_A{\tt Ste}$ and ${\tt Ste}_A$ of stereotype modules:}

\bit{\it

\item[$1^\circ$.] $_A{\tt Ste}$ and ${\tt Ste}_A$ are pre-Abelian categories.

\item[$2^\circ$.]\label{TH:STE-polnaya-kategoriya} $_A{\tt Ste}$ and ${\tt Ste}_A$ are complete:
 each covariant (and each contravariant) system has an injective and a projective limit.

\item[$3^\circ$.] $_A{\tt Ste}$ and ${\tt Ste}_A$ are enriched categories over the monoidal category ${\tt Ste}$.

}\eit

\subsection{Subalgebras, quotient algebras, limits and completeness of ${\tt Ste}^\circledast$}\label{SUBSEC:polnota-Ste^circledast}

\paragraph{Subalgebras, products and projective limits.}

\bit{
\item[$\bullet$]
Suppose $B$ is a subset in a stereotype algebra $A$ endowed with a structure of stereotype algebra in such a way that the set-theoretic enclosure $B\subseteq A$ is a morphism of stereotype algebras (i.e. a linear, multiplicative and preserving identity continuous map). Then the stereotype algebra  $B$ is called a {\it subalgebra} of the stereotype algebra $A$, and the set-theoretic enclosure $\sigma:B\subseteq A$ its {\it representing monomorphism}.

\item[$\bullet$] We say that a subalgebra $B$ of a stereotype algebra $A$ is {\it closed}\label{DEF:zamknutaya-podalgebra}, if its representing monomorphism $\sigma:B\to A$ is a closed map in the sense of definition on page \pageref{DEF:zamkn-otobr}.

}\eit

The following fact was stated in \cite{Akbarov} (Theorem 10.13):

\btm \label{th-10.13} Let $A$ be a stereotype algebra and $B$ its subalgebra (in the purely algebraic sense), and at the same time a closed subspace of the locally convex space $A$. Then the pseudosaturation $B^\vartriangle$ is a (stereotype algebra and a) closed subalgebra in $A$.
\etm

\btm \label{TH:prod-proj-ster-algebr} Each family $\{A_i;i\in I\}$ of stereotype algebras has a direct product in the category ${\tt Ste}^{\circledast}$ of stereotype algebras, and as a stereotype space this product is exactly the direct product of the family of stereotype spaces $\{A_i;\ i\in I\}$:
$$
{\tt Ste}^{\circledast}\text{-}\prod_{i\in I}A_i={\tt Ste}\text{-}\prod_{i\in I}A_i
$$.
\etm
\bpr
We have to verify that the direct product is the usual direct product of locally convex spaces $A=\prod_{i\in I}A_i$ with the coordinate-wise multiplication:
$$
(x\cdot y)_i=x_i\cdot y_i,\qquad i\in I.
$$
By \cite[Theorem 4.20]{Akbarov}, this is a stereotype space, so we only need to prove that the multiplication is continuous. Let $U$ be a neighborhood of zero and $K$ a compact set in $A$. We must find a neighborhood of zero $V$ in $A$ such that
$$
V\cdot K\subseteq U,\qquad K\cdot V\subseteq U.
$$
It is sufficient to consider a base neighborhood of zero $U$, i.e.
$$
U=\{x\in A:\ \forall i\in J\quad x_i\in U_i\}
$$
where $J\subseteq I$ is a finite subset in $I$, and for any $i\in J$ the set $U_i$ is a neighborhood of zero in $A_i$, and $x_i$ is the projection of $x\in A$ onto $A_i$. If $U$ has this form, then for any $i\in J$ we can consider the neighborhood of zero $U_i$ in $A_i$, and (since $A_i$ is a stereotype algebra) we can choose a neighborhood of zero $V_i$ such that
$$
V_i\cdot K_i\subseteq U_i,\qquad K_i\cdot V_i\subseteq U_i
$$
(where $K_i$ is the projection of the compact set $K\subseteq A$ onto $A_i$). Then we put
$$
V=\{x\in A:\ \forall i\in J\quad x_i\in V_i\}
$$
and for each $x\in V$ and $y\in K$ we get:
$$
\Big(\forall i\in J\qquad (x\cdot y)_i=x_i\cdot y_i\in V_i\cdot K_i\subseteq U_i\Big)\quad\Longrightarrow\quad x\cdot y\in U,
$$
This means that $V\cdot K\subseteq U$. Similarly,
$$
\Big(\forall i\in J\qquad (y\cdot x)_i=y_i\cdot x_i\in K_i\cdot V_i\subseteq U_i\Big)\quad\Longrightarrow\quad y\cdot x\in U,
$$
and this means that $K\cdot V\subseteq U$.
\epr

\btm \label{TH:lim-proj-ster-algebr} Each covariant system $\{A_i;\ \pi_i^j\}$ of stereotype algebras has a projective limit in the category ${\tt Ste}^{\circledast}$ of stereotype algebras, and as a stereotype space this limit is exactly the projective limit of the covariant system of stereotype spaces $\{A_i;\ \pi_i^j\}$:
$$
{\tt Ste}^{\circledast}\text{-}\leftlim A_i={\tt Ste}\text{-}\leftlim A_i
$$
\etm
\bpr By Theorem \ref{TH:prod-proj-ster-algebr}, the direct product $A=\prod_{i\in I} A_i$ with the coordinate-wise multiplication is a direct product of the family of algebras $\{A_i\}$ in ${\tt Ste}^{\circledast}$, and by Theorem \ref{th-10.13}, the subalgebra $B$ in $A$, consisting of families  $\{x_i;\ i\in I\}$ with the properties
$$
x_i=\pi^j_i(x_j),\qquad i\le j\in I,
$$
and endowed with the topology of pseudosaturation of the topology inherited from $A$, is a stereotype algebra. The same mode as in the case of stereotype spaces, prove that $B$ is the projective limit in ${\tt Ste}^{\circledast}$.
\epr

\paragraph{Quotient algebras, coproducts and injective limits.}

\bit{
\item[$\bullet$] Let $A$ be a stereotype algebra, and let
 \bit{

 \item[1)] $I$ be a two-sided ideal in $A$ (as in algebra), and at the same time a closed set in $A$ (as in a topological space), we will further call such ideals {\it closed ideals} in $A$,

 \item[2)] $\tau$ be a locally convex topology on the quotient algebra $A/I$, such that $\tau$ is majorated by the usual quotient topology,

 \item[3)] $B$ be a subspace in the completion $(A/I)^\blacktriangledown$ of the locally convex space $A/I$ with respect to $\tau$, such that $B$ contains $A/I$ and is a stereotype algebra with respect to the algebraic operations and the topology inherited from $(A/I)^\blacktriangledown$.
 }\eit
Then we call the stereotype algebra $B$ the {\it quotient algebra of the stereotype algebra} $A$, and the composition $\upsilon=\sigma\circ\pi$ of the quotient map $\pi:A\to A/I$ and the natural embedding $\sigma:A/I\to B$ is called the {\it representing epimorphism} of the quotient algebra $B$.

\item[$\bullet$] A quotient algebra $B$ of a stereotype algebra $A$ is said to be {\it open}\label{DEF:otkrytaya-faktor-algebra}, if its representing epimorphism $\upsilon:B\gets A$ is an open map in the sense of definition on page \pageref{DEF:otkr-otobr}.

}\eit

The symmetry between projective and injective constructions which was obvious for stereotype space (see \cite{Akbarov}), is preserved in some sense for stereotype algebras, but the difference is that the injective constructions in ${\tt Ste}^\circledast$ become more complicated and as a corollary, the proofs become more difficult (however, the situation here is the same as for algebras in purely algebraic sense). For example, the analog of Theorem \ref{th-10.13} uses the theory of modules over algebras (see proof of Theorem 10.14 in \cite{Akbarov}):

\btm \label{th-10.14} Let  $A$ be a stereotype algebra and $I$ a closed ideal in $A$. Then the pseudocompletion $(A/I)^\triangledown$ is a stereotype algebra (and is called an open quotient algebra of $A$ by the ideal $I$).
\etm

\brem\label{REM:nesushestvennost-unialnosti-v-th-10.14} In Theorem \ref{th-10.14} the unitality requirement (i.e. the existence of identity) for the algebra $A$ is unessential.
\erem

Suppose $\{A_i;i\in I\}$ is a family of stereotype algebras.
Let us construct an algebra $\coprod_{i\in I}A_i$ in the following way. First let us say that a sequence of indices $i=\{i_1,...,i_n\}\in I$ {\it alternates}, if any two neighboring elements there do not coincide:
    $$
    \forall k=1,...,n-1\qquad i_k\ne i_{k+1}.
    $$
The set of all alternating seuqences in $I$ of (various) finite length will be denoted as $I_{\N}^{\alt}$. Let us introduce the operation of multiplication on $I_{\N}^{\alt}$ as follows: if $\iota,\varkappa\in I_{\N}^{\alt}$ have lengths $m$ and $n$ respectively, then their producet is
$$
\iota * \varkappa=
\begin{cases}
(\iota_1,...,\iota_m,\varkappa_1,...,\varkappa_n),& \text{for $\iota_m\ne\varkappa_1$}\\
(\iota_1,...,\iota_m,\varkappa_2,...,\varkappa_n),& \text{for $\iota_m=\varkappa_1$}
\end{cases}
$$
(the length of $\iota*\varkappa$ is $m+n$ if $\iota_m\ne \varkappa_1$, and $m+n-1$ if $\iota_m=\varkappa_1$).
For each sequence $\iota\in I_{\N}^{\alt}$ we put
$$
A_\iota=A_{\iota_1}\circledast A_{\iota_2}\circledast...\circledast A_{\iota_m}.
$$
(where $\circledast$ is the projective tensor product from \eqref{eq7.1}). Let us note that for all seuqences   $\iota,\varkappa\in I_{\N}^{\alt}$ the spaces $A_\iota\circledast A_\varkappa$ and $A_{\iota*\varkappa}$ are naturally related through the following linear continuous map:
\beq\label{mu_(iota,varkappa)}
\mu_{\iota,\varkappa}:A_\iota\circledast A_\varkappa\to A_{\iota*\varkappa}
\quad \Big|\quad \mu_{\iota,\varkappa}=
\begin{cases}
1_{A_\iota\circledast A_\varkappa}, & \iota_m\ne\varkappa_1\\
1_{A_{\iota_1}}\circledast...\circledast 1_{A_{\iota_{m-1}}}
\circledast \mu_{\iota_m} \circledast 1_{A_{\varkappa_2}}\circledast...\circledast 1_{A_{\varkappa_n}}
, & \iota_m=\varkappa_1
\end{cases}
\eeq
where $\mu_i:A_i\circledast A_i\to A_i$ is the morphism of multiplication in the algebra $A_i$.

Consider the stereotype space
    $$
    A_*=\bigoplus_{\iota\in I_{\N}^{\alt}} A_\iota
    $$
and note that the formula
    $$
    \Big(a_{i_1}\circledast a_{i_2}\circledast...\circledast a_{i_m}\Big)\cdot\Big( b_{j_1}\circledast b_{j_2}\circledast...\circledast b_{j_n}\Big)=
    \begin{cases}
    a_{i_1}\circledast a_{i_2}\circledast...\circledast a_{i_m}\circledast b_{j_1}\circledast b_{j_2}\circledast...\circledast b_{j_n},& i_m\ne j_1
    \\
    a_{i_1}\circledast a_{i_2}\circledast...\circledast \Big( a_{i_m}\cdot b_{j_1}\Big)\circledast b_{j_2}\circledast...\circledast b_{j_n}
    ,& i_m=j_1.
    \end{cases}
    $$
defines a multiplication in $A_*$, which is a continuous bilinear map. This becomes obvious if we represent this operation as the composition of the following maps:
$$
A_*\times A_*\to A_*\circledast A_*=\Big(\bigoplus_{\iota\in I_{\N}^{\alt}} A_\iota\Big)\circledast \Big(\bigoplus_{\varkappa\in I_{\N}^{\alt}} A_\varkappa\Big)\to \bigoplus_{\iota,\varkappa\in I_{\N}^{\alt}} A_\iota\circledast A_\varkappa
\to \bigoplus_{\iota,\varkappa\in I_{\N}^{\alt}} A_{\iota*\varkappa}
\to \bigoplus_{\lambda\in I_{\N}^{\alt}} A_\lambda= A_*.
$$
Here the first arrow $A_*\times A_*\to A_*\circledast A_*$ is the mapping described in Proposition \ref{prop-7.1}, the second arrow $\Big(\bigoplus_{\iota\in I_{\N}^{\alt}} A_\iota\Big)\circledast \Big(\bigoplus_{\varkappa\in I_{\N}^{\alt}} A_\varkappa\Big)\to \bigoplus_{\iota,\varkappa\in I_{\N}^{\alt}} A_\iota\circledast A_\varkappa$ is the natural isomorphism \eqref{perestanovochnost-circledast-s-summami} that connects the direct sum and the projective tensor product, the third arrow $\bigoplus_{\iota,\varkappa\in I_{\N}^{\alt}} A_\iota\circledast A_\varkappa \to \bigoplus_{\iota,\varkappa\in I_{\N}^{\alt}} A_{\iota*\varkappa}$ is the direct sum
$\bigoplus_{\iota,\varkappa\in I_{\N}^{\alt}}\mu_{\iota,\varkappa}$ of the morphisms \eqref{mu_(iota,varkappa)},
and the final arrow $\bigoplus_{\iota,\varkappa\in I_{\N}^{\alt}} A_{\iota*\varkappa}
\to \bigoplus_{\lambda\in I_{\N}^{\alt}} A_\lambda$ is the result of identification of each summand of the form  $A_{\iota*\varkappa}$ (there can be many of those) with the space $A_\lambda$ in the sum $\bigoplus_{\lambda\in I_{\N}^{\alt}} A_\lambda$  (which is unique).

Obviuously, this operation of multiplication in $A_*$ is associative. If we take the quotient algebra of the  (non-unital) algebra $A_*$ by the closed ideal $M$ (here we use Remark \ref{REM:nesushestvennost-unialnosti-v-th-10.14} to Theorem \ref{th-10.14}) generated by the elements of the form
    $$
    1_{A_i}-1_{A_j},\qquad i,j\in I,
    $$
then the quotient algebra $(A_*/M)^\triangledown$ will be the stereotype algebra with the identity
   $$
   1_{(A_*/M)^\triangledown}=\pi(1_{A_i})
   $$
   (here at the right side we mean the image of the identity $1_{A_i}$ in the arbitrary algebra $A_i$ under the quotient map $\pi:A_*\to (A_*/M)^\triangledown$).

\bit{
\item[$\bullet$] Following \cite{Nica-Speicher} we call the algebra $(A_*/M)^\triangledown$ the {\it free product} of the algebras  $\{A_i;i\in I\}$ and we denote it by ${\tt Ste}^{\circledast}\text{-}\coprod_{i\in I}A_i$, or by
   $$
{\tt Ste}^{\circledast}\text{-}\coprod_{i\in I}A_i=(A_*/M)^\triangledown.
   $$
   This is justified by the following theorem.
}\eit

\btm \label{TH:sum-proj-ster-algebr} For each family $\{A_i;i\in I\}$ of stereotype algebras its free product ${\tt Ste}^{\circledast}\text{-}\coprod_{i\in I}A_i$ is a co-product in the category of stereotype algebras ${\tt Ste}^{\circledast}$.
\etm

\btm \label{TH:lim-inj-ster-algebr} Each covariant system $\{A_i;\ \iota_i^j\}$ of stereotype algebras has an injective limit in the category ${\tt Ste}^{\circledast}$.
\etm
\bpr This is the open quotient algebra $\Big(\big(\coprod_{i\in I}A_i\big)/N\Big)^\triangledown$ of the free product $\coprod_{i\in I}A_i$ by the closed ideal $N$ generated by elements of the form
$$
\iota_i(x)-\iota_j(\iota_i^j(x)),\qquad x\in A_i,
$$
where $\iota_k:A_k\to \coprod_{i\in I}A_i$ are natural embeddings.
\epr

As one can notice (this is an illustration to the difference between the projective and the injective constructions in ${\tt Ste}^{\circledast}$) the injective limits in ${\tt Ste}^{\circledast}$ do not necessarily coincide as stereotype spaces with the injective limits in ${\tt Ste}$. For instance for co-products we have inequality:
$$
{\tt Ste}^{\circledast}\text{-}\coprod_{i\in I}A_i\ne{\tt Ste}\text{-}\coprod_{i\in I}A_i
$$
(although there is a natural map from the right side to the left side). This asymmetry however diasppears in the case, when the index set $I$ is directed:

\btm \label{TH:lim-inj-ster-algebr-I-napravleno} If $\{A_i;\ \iota_i^j\}$ is a covariant system of stereotype algebras over a directed set $I$, then the natural map
$$
{\tt Ste}\text{-}\rightlim A_i\to {\tt Ste}^{\circledast}\text{-}\rightlim A_i
$$
between its injective limit in the category ${\tt Ste}$ and the injective limit in the category ${\tt Ste}^{\circledast}$ is an isomorphism of stereotype spaces:
$$
{\tt Ste}\text{-}\rightlim A_i\cong {\tt Ste}^{\circledast}\text{-}\rightlim A_i.
$$
\etm
\bpr
Denote by $A$ the injective limit of the system $\{A_i;\ \iota_i^j\}$ in ${\tt Ste}$:
$$
A={\tt Ste}\text{-}\rightlim A_i
$$
and let $\rho_i:A_i\to A$ be the corresponding morphisms of stereotype spaces:
\beq\label{DIAGR:inj-lim-siln-algebr}
\begin{diagram}
\node[2]{A}  \\
\node{A_i}\arrow[2]{e,b}{\iota_i^j}\arrow{ne,t}{\rho_i}
\node[2]{A_j}\arrow{nw,t}{\rho_j}
\end{diagram}
\eeq
We will show that $A$ has a natural structure of stereotype algebra, and with this structure $A$ is an injective limit of the covariant system of stereotype algebras $\{A_i;\ \iota_i^j\}$.

1. Take $i\in I$ and note that for any $j\ge i$ the homomorphism $\iota^j_i:A_i\to A_j$ induces on $A_j$ a structure of left $A_i$-module by formula
\beq\label{a-underset(cdot)(i)-b}
a\underset{i}{\cdot}b=\iota_i^j(a)\underset{A_j}{\cdot}b,\qquad a\in A_i,\ b\in A_j.
\eeq
(here $\underset{i}{\cdot}$ means the left multiplication by elements of $A_i$, and $\underset{A_j}{\cdot}$ the multiplication in $A_j$). Besides this, for $i\le j\le k$ the maps $\iota_j^k:A_j\to A_k$ turn out to be morphisms of left $A_i$-modules:
$$
\iota_j^k(a\underset{i}{\cdot}b)=\eqref{a-underset(cdot)(i)-b}=\iota_j^k\big(\iota_i^j(a)\underset{A_j}{\cdot}b\big)=
\iota_j^k\big(\iota_i^j(a)\big)\underset{A_k}{\cdot}\iota_j^k(b)=\iota_i^k(a)\underset{A_k}{\cdot}\iota_j^k(b)=
\eqref{a-underset(cdot)(i)-b}=a\underset{i}{\cdot}\iota_j^k(b),\quad a\in A_i,\ b\in A_j.
$$
This means that $\{A_j;\ j\ge i\}$ can be considered as a covariant system of left stereotype $A_i$-modules. By \cite[Theorem 11.17]{Akbarov}, it has an injective limit, which as a stereotype space coincide with the injective limit of the system of stereotype spaces $\{A_i; j\ge i\}$. And the latter one coincides with the injective limit of all the system of stereotype spaces $\{A_i;\ \iota_i^j\}$, since $I$ is directed:
$$
{_{A_i}{\tt Ste}}\text{-}\lim_{i\le j\to\infty}A_j={\tt Ste}\text{-}\lim_{i\le j\to\infty}A_j={\tt Ste}\text{-}\lim_{j\to\infty}A_j=A.
$$
An important conclusion for us is that for any $i\in I$ the space $A$ has a structure of stereotype $A_i$-module, and under this structure the maps in diagram \eqref{DIAGR:inj-lim-siln-algebr} become morphisms of $A_i$-modules, in particular,
\beq\label{rho_i(a-cdot-b)=a-cdot-rho_i(b)}
\rho_j(a\underset{i}{\cdot}b)=a\underset{i}{\cdot}\rho_j(b),\qquad i\le j,\ a\in A_i,\ b\in A_j.
\eeq

2. Note then that for $i\le j$ the structures of left $A_i$-module and of left $A_j$-module on $A$ are coherent with each other by the identity
\beq\label{a-underset(i)(cdot)-x=iota_i^j(a)-underset(j)(cdot)x}
\iota_i^j(a)\underset{j}{\cdot}x=a\underset{i}{\cdot}x,\qquad a\in A_i,\ x\in A.
\eeq
To prove this we should first consider a special case when $x=\rho_k(b)$, $b\in A_k$, $k\ge j$. We have in this situation:
\begin{multline*}
\iota_i^j(a)\underset{j}{\cdot}x=\iota_i^j(a)\underset{j}{\cdot}\rho_k(b)=\eqref{rho_i(a-cdot-b)=a-cdot-rho_i(b)}=
\rho_k\big(\iota_i^j(a)\underset{j}{\cdot}b\big)=\eqref{a-underset(cdot)(i)-b}=
\rho_k\Big(\iota_j^k\big(\iota_i^j(a))\underset{A_k}{\cdot}b\Big)=\\=
\rho_k\big(\iota_i^k(a)\underset{A_k}{\cdot}b\big)=\eqref{a-underset(cdot)(i)-b}=
\rho_k(a\underset{i}{\cdot}b)=\eqref{rho_i(a-cdot-b)=a-cdot-rho_i(b)}=
a\underset{i}{\cdot}\rho_k(b)=a\underset{i}{\cdot}x.
\end{multline*}
After that let us recall that the family of spaces $A_k$ is dense in its injective limit $A$ (we use here the left formula of \cite[(4.15)]{Akbarov} and the fact that $I$ is directed). This means that for any $x\in A$ there is a net $x_k\in \rho_k(A_k)$ tending to $x$ in $A$:
$$
x_k\overset{A}{\underset{k\to\infty}{\longrightarrow}}x.
$$
Since for any $x_k$ the equality \eqref{a-underset(i)(cdot)-x=iota_i^j(a)-underset(j)(cdot)x} is already proved, we obtain a relation which proves \eqref{a-underset(i)(cdot)-x=iota_i^j(a)-underset(j)(cdot)x} for this $x$:
$$
\iota_i^j(a)\underset{j}{\cdot}x\overset{A}{\underset{\infty\gets k}{\longleftarrow}}
\iota_i^j(a)\underset{j}{\cdot}x_k=a\underset{i}{\cdot}x_k
\overset{A}{\underset{k\to\infty}{\longrightarrow}}a\underset{i}{\cdot}x
$$
(the possibility to take limits follows from the continuity of the multiplication in a stereotype module).

3. From the fact that $A$ is a left $A_i$-module we obtain by \cite[Theorem 11.2]{Akbarov}, that the formula
$$
\ph_i(a)(x)=a\underset{i}{\cdot} x,\quad a\in A_i,\ x\in A,
$$
defines a homomorphism od stereotype algebras
$$
\ph_i:A_i\to{\mathcal L}(A).
$$
The fact that this is a homomorphism means that we have the identity
\beq\label{inj-limit-algebr:ph_i-gomomorfizm}
\ph_i(a\cdot b)=\ph_i(a)\circ\ph_i(b),\qquad a,b\in A_i,
\eeq
and equality
\beq\label{ph_i(1_(A_i))=id_A}
\ph_i(1_{A_i})=\id_A.
\eeq
formula \eqref{rho_i(a-cdot-b)=a-cdot-rho_i(b)} in this style of writing turns into the identity
\beq\label{ph_i(a)(rho_i(b))=a-cdot-rho_i(b)=rho_i(a-cdot-b)}
\ph_i(a)\big(\rho_i(b)\big)=a\underset{i}{\cdot}\rho_i(b)=\rho_i(a\underset{A_i}{\cdot} b),\qquad a,b\in A_i,
\eeq
and formula \eqref{a-underset(i)(cdot)-x=iota_i^j(a)-underset(j)(cdot)x} into the identity
$$
\ph_j\big(\iota_i^j(a)\big)(x)=\ph_i(a)(x),\qquad a\in A_i,\ x\in A,
$$
which is equivalet to the equality
\beq\label{a-underset(i)(cdot)-x=iota_i^j(a)-underset(j)(cdot)x-new}
\ph_j\circ \iota_i^j=\ph_i,\qquad i\le j.
\eeq
The latter one means that the following diagram in ${\tt Ste}^\circledast$ is commutative:
$$
\begin{diagram}
\node[2]{{\mathcal L}(A)}  \\
\node{A_i}\arrow[2]{e,b}{\iota_i^j}\arrow{ne,t}{\ph_i}
\node[2]{A_j}\arrow{nw,t}{\ph_j}
\end{diagram}
$$
One can interpret this as an injective cone of the covariant system $\{A_i;\ \iota_i^j\}$ in the category ${\tt Ste}$ of stereotype spaces. Then we can conclude that there exists a linear continuous map $\ph$ from the injective limit $A=\rightlim A_i$ of this system into the space  ${\mathcal L}(A)$ such that for any $i$ the following diagram is commutative
\beq\label{DIAGR:inj-lim-siln-algebr-L(A)}
\begin{diagram}
\node{A}\arrow[2]{e,t}{\ph}\node[2]{{\mathcal L}(A)}  \\
\node[2]{A_i}\arrow{nw,b}{\rho_i}\arrow{ne,b}{\ph_i}
\end{diagram}
\eeq
Let us put
\beq\label{DEF:umnozhenie-v-inj-lim-algebr}
x\cdot y=\ph(x)(y),\qquad x,y\in A,
\eeq
and verify that this multiplication turns $A$ into a stereotype algebra.

4. Let us note that the bilinear form $(x,y)\mapsto x\cdot y$ is continuous. Indeed, if $K$ is a compact set in $A$, then its image $\ph(K)$ is a compact set in ${\mathcal L}(A)$. Hence, $\ph(K)$ is a compact set in the space of operators $A:A$. By \cite[Theorems 5.1 and 2.5]{Akbarov}, this means that $\ph(K)$ is equicontinuous on $A$. Hence for a neighborhood of zero $W$ in $A$ there is a neighborhood of zero $V$ in $A$ such that
$$
K\cdot V=\ph(K)(V)\subseteq W.
$$
On the other hand, for any compact set $K$ and for any neighborhood of zero $W$ in $A$ the set $W\oslash K$ is a neighborhood of zero in ${\mathcal L}(A)$, hence from the continuity of $\ph$ it follows that there is a neighborhood of zero $V$ in $A$ such that
$$
\ph(V)\subseteq W\oslash K,
$$
and this is equivalent to the inclusion
$$
V\cdot K=\ph(V)(K)\subseteq V.
$$

5. Besides this, the formula
\beq\label{1_A=rho_i(1_(A_i))}
1_A=\rho_i(1_{A_i})
\eeq
defines some element of the space $A$, so if $i\le j$, then
$$
\rho_j(1_{A_j})=\rho_j\big(\iota_i^j(1_{A_i})\big)=\rho_i(1_{A_i}).
$$
At the same time the chain
\beq\label{ph(1_A)=id_A}
\ph(1_A)=\ph\big(\rho_i(1_{A_i})\big)=\ph_i(1_{A_i})=\eqref{ph_i(1_(A_i))=id_A}=\id_A
\eeq
implies that this element is the identity for the multiplication \eqref{DEF:umnozhenie-v-inj-lim-algebr}: first, for any $y\in A$ we have
$$
1_A\cdot y=\ph(1_A)(y)=\id_A(y)=y.
$$
And, second, for any $x\in A$ we can find a net $a_k\in A_k$ such that
$$
\rho_k(a_k)\overset{A}{\underset{k\to\infty}{\longrightarrow}}x,
$$
and by the already proven continuity of the multiplication in $A$, we have:
\begin{multline*}
x\cdot 1_A\overset{A}{\underset{\infty\gets k}{\longleftarrow}}\rho_k(a_k)\cdot 1_A=
\rho_k(a_k)\cdot\rho_k(1_{A_k})=\eqref{DEF:umnozhenie-v-inj-lim-algebr}=
\ph\Big(\rho_k(a_k)\Big)\Big(\rho_k(1_{A_k})\Big)=\\=
\ph_k(a_k)\Big(\rho_k(1_{A_k})\Big)=\eqref{ph_i(a)(rho_i(b))=a-cdot-rho_i(b)=rho_i(a-cdot-b)}=
\rho_k(a_k\underset{A_k}{\cdot}1_{A_k})=\rho_k(a_k)\overset{A}{\underset{k\to\infty}{\longrightarrow}}x,
\end{multline*}
Thus,
$$
x\cdot 1_A=x.
$$

6. Now we notice that the map $\rho_i$ in \eqref{DIAGR:inj-lim-siln-algebr-L(A)} must be a homomorphism of algebras. Indeed its turns identity into identity just by the definition of $1_A$ in \eqref{1_A=rho_i(1_(A_i))}. On the other hand, it preserves multiplication since for all $a,b\in A_i$
\beq\label{inj-lim-alg:rho_i-gomomorfizm}
\rho_i(a\underset{A_i}{\cdot} b)=\eqref{ph_i(a)(rho_i(b))=a-cdot-rho_i(b)=rho_i(a-cdot-b)}=\ph_i(a)\big(\rho_i(b)\big)=
 \eqref{DIAGR:inj-lim-siln-algebr-L(A)}=\ph\big(\rho_i(a)\big)\big(\rho_i(b)\big)=\eqref{DEF:umnozhenie-v-inj-lim-algebr}=
 \rho_i(a)\cdot \rho_i(b).
\eeq

7. The same for the map $\ph$. The preserving of identities were already stated in the chain \eqref{1_A=rho_i(1_(A_i))}. And to prove multiplicativity we first have to note the formula
\beq\label{ph(rho_i(a)-cdot-rho_j(b))=ph(rho_i(a))-circ-ph(rho_j(b))}
\ph\Big(\rho_i(a)\cdot\rho_j(b)\Big)=\ph\Big(\rho_i(a)\Big)\circ\ph\Big(\rho_j(b)\Big),\qquad i,j\in I,\ a\in A_i,\ b\in A_j
\eeq
Indeed, for $k\in I$ such that $k\ge i$ and $k\ge j$ we have:
\begin{multline*}
\ph\Big(\rho_i(a)\cdot\rho_j(b)\Big)=\eqref{DIAGR:inj-lim-siln-algebr}=\ph\Big(\rho_k\big(\iota_i^k(a)\big)\cdot\rho_k\big(\iota_j^k(b)\big)\Big)= \eqref{inj-lim-alg:rho_i-gomomorfizm}=
\ph\Big(\rho_k\big(\iota_i^k(a)\cdot\iota_j^k(b)\big)\Big)=\\=\eqref{DIAGR:inj-lim-siln-algebr-L(A)}=
\ph_k\big(\iota_i^k(a)\cdot\iota_j^k(b)\big)=\eqref{inj-limit-algebr:ph_i-gomomorfizm}=
\ph_k\big(\iota_i^k(a)\big)\circ\ph_k\big(\iota_j^k(b)\big)=\eqref{DIAGR:inj-lim-siln-algebr-L(A)}=\\=
\ph\Big(\rho_k\big(\iota_i^k(a)\big)\Big)\circ\ph\Big(\rho_k\big(\iota_j^k(b)\big)\Big)=\eqref{DIAGR:inj-lim-siln-algebr}=
\ph\Big(\rho_i(a)\Big)\circ\ph\Big(\rho_j(b)\Big)
\end{multline*}
Then we take $x,y\in A$ and find $a_i\in A_i$ and $b_j\in A_j$ such that
$$
\rho_i(a_i)\overset{A}{\underset{i\to\infty}{\longrightarrow}}x,\qquad
\rho_j(b_j)\overset{A}{\underset{j\to\infty}{\longrightarrow}}y.
$$
We obtain:
\begin{multline*}
\ph(x\cdot y)\overset{{\mathcal L}(A)}{\underset{\infty\gets i}{\longleftarrow}}
\ph\Big(\rho_i(a_i)\cdot y\Big)\overset{{\mathcal L}(A)}{\underset{\infty\gets j}{\longleftarrow}}
\ph\Big(\rho_i(a_i)\cdot \rho_j(b_j)\Big)=\eqref{ph(rho_i(a)-cdot-rho_j(b))=ph(rho_i(a))-circ-ph(rho_j(b))}=\\=
\ph\Big(\rho_i(a_i)\Big)\circ\ph\Big(\rho_j(b_j)\Big)
\overset{{\mathcal L}(A)}{\underset{j\to\infty}{\longrightarrow}}
\ph\Big(\rho_i(a_i)\Big)\circ\ph(y)
\overset{{\mathcal L}(A)}{\underset{i\to\infty}{\longrightarrow}}
\ph(x)\circ\ph(y),
\end{multline*}
hence,
$$
\ph(x\cdot y)=\ph(x)\circ\ph(y).
$$
This formula proves in addition the associativity of the multiplication in $A$,
$$
x\cdot(y\cdot z)=\ph(x)\Big(y\cdot z\Big)=\ph(x)\Big(\ph(y)(z)\Big)=\Big(\ph(x)\circ\ph(y)\Big)(z)=
\ph(x\cdot y)(z)=(x\cdot y)\cdot z,
$$
and this was the last what we needed to understand that $A$ is a stereotype algebra.

8. We only have to verify that the cone of algebras $\{A_i;\ \rho_i\}$ is an injective limit of the covariant system  of algebras $\{A_i;\ \iota_i^j\}$. Let $\{B_i;\ \sigma_i\}$ be another cone of algebras. Since it is also a cone of stereotype spaces, there exists a unique linear continuous map $\sigma:A\to B$ such that the following diagram is commutative:
\beq\label{inj-lim-algebr:A-sigma-B}
\begin{diagram}
\node{A}\arrow[2]{e,t,--}{\sigma}\node[2]{B}  \\
\node[2]{A_i}\arrow{nw,b}{\rho_i}\arrow{ne,b}{\sigma_i}
\end{diagram}
\eeq
We must check that $\sigma$ is a homomorphism of algebras. The preserving of identities follows from the fact that all $\sigma_i$ preserve identity:
$$
\sigma(1_A)=\sigma\Big(\rho_i(1_{A_i})\Big)=\sigma_i(1_{A_i})=1_B.
$$
For proving the multiplicativity we note first the following identity:
\beq\label{sigma(rho_i(a)-cdot-rho_j(b))=sigma(rho_i(a))-cdot-sigma(rho_j(b))}
\sigma\Big(\rho_i(a)\cdot\rho_j(b)\Big)=\sigma\Big(\rho_i(a)\Big)\cdot\sigma\Big(\rho_j(b)\Big),\qquad i,j\in I,\ a\in A_i,\ b\in A_j,
\eeq
This can be proved by the same reasoning as \eqref{ph(rho_i(a)-cdot-rho_j(b))=ph(rho_i(a))-circ-ph(rho_j(b))} above: take $k\in I$ such that $k\ge i$ and $k\ge j$, then
\begin{multline*}
\sigma\Big(\rho_i(a)\cdot\rho_j(b)\Big)=\eqref{DIAGR:inj-lim-siln-algebr}=\sigma\Big(\rho_k\big(\iota_i^k(a)\big)\cdot\rho_k\big(\iota_j^k(b)\big)\Big)= \eqref{inj-lim-alg:rho_i-gomomorfizm}=
\sigma\Big(\rho_k\big(\iota_i^k(a)\cdot\iota_j^k(b)\big)\Big)=\\=\eqref{inj-lim-algebr:A-sigma-B}=
\sigma_k\big(\iota_i^k(a)\cdot\iota_j^k(b)\big)=
\sigma_k\big(\iota_i^k(a)\big)\cdot\sigma_k\big(\iota_j^k(b)\big)=\eqref{inj-lim-algebr:A-sigma-B}=\\=
\sigma\Big(\rho_k\big(\iota_i^k(a)\big)\Big)\cdot\sigma\Big(\rho_k\big(\iota_j^k(b)\big)\Big)=\eqref{DIAGR:inj-lim-siln-algebr}=
\sigma\Big(\rho_i(a)\Big)\cdot\sigma\Big(\rho_j(b)\Big)
\end{multline*}
After that we take $x,y\in A$ and choose $a_i\in A_i$ and $b_j\in A_j$ such that
$$
\rho_i(a_i)\overset{A}{\underset{i\to\infty}{\longrightarrow}}x,\qquad
\rho_j(a_j)\overset{A}{\underset{j\to\infty}{\longrightarrow}}y.
$$
We obtain:
\begin{multline*}
\sigma(x\cdot y)\overset{B}{\underset{\infty\gets i}{\longleftarrow}}
\sigma\Big(\rho_i(a_i)\cdot y\Big)\overset{B}{\underset{\infty\gets j}{\longleftarrow}}
\sigma\Big(\rho_i(a_i)\cdot \rho_j(b_j)\Big)=\eqref{sigma(rho_i(a)-cdot-rho_j(b))=sigma(rho_i(a))-cdot-sigma(rho_j(b))}=\\=
\sigma\Big(\rho_i(a_i)\Big)\cdot\sigma\Big(\rho_j(b_j)\Big)
\overset{B}{\underset{j\to\infty}{\longrightarrow}}
\sigma\Big(\rho_i(a_i)\Big)\cdot\sigma(y)
\overset{B}{\underset{i\to\infty}{\longrightarrow}}
\sigma(x)\cdot\sigma(y),
\end{multline*}
and thus,
$$
\sigma(x\cdot y)=\sigma(x)\cdot\sigma(y).
$$
\epr

\paragraph{Completeness of the category ${\tt Ste}^\circledast$.}

Theorems \ref{TH:lim-proj-ster-algebr} and \ref{TH:lim-inj-ster-algebr} imply

\btm\label{TH:polnota-Ste^circledast} The category ${\tt Ste}^{\circledast}$ of stereotype algebras is complete.
\etm

\subsection{Nodal decomposition, envelope and refinement in ${\tt Ste}^\circledast$}\label{SUBSEC:uzlov-razlozh-v-Ste^circledast}

\paragraph{Discerning properties of strong epimorphisms in ${\tt Ste}^\circledast$.}

\btm\label{TH:strogie-epi-v-Ste^circledast=strogie-epi-v-Ste} For a morphism of stereotype algebras $\e:A\to B$ the following conditions are equivalent:
 \bit{
 \item[(i)] $\e$ is an immediate epimorphism in ${\tt Ste}^\circledast$,

 \item[(ii)] $\e$ is a strong epimorphism in  ${\tt Ste}^\circledast$,

 \item[(iii)] $\e$ is an immediate epimorphism in  ${\tt Ste}$,

 \item[(iv)] $\e$ is a strong epimorphism in  ${\tt Ste}$.
 }\eit
\etm
\bpr
Let us note that the connections (i)$\Leftarrow$(ii) and (iii)$\Leftrightarrow$(iv) are already known. So it is sufficient to prove  (i)$\Rightarrow$(iii) and (ii)$\Leftarrow$(iv).

1. Let us start with (i)$\Rightarrow$(iii). Let $\e:A\to B$ be an immediate epimorphism in ${\tt Ste}^\circledast$. Consider its minimal factorization in ${\tt Ste}$, i.e. a diagram with linear continuous maps
$$
\xymatrix
{
A\ar[rr]^{\e}\ar[rd]_{\coim_\infty\e} & & B \\
& \Coim_\infty\e\ar[ru]_{\mu} &
}
$$
where $\Coim_\infty\e$ is the nodal coimage in ${\tt Ste}$. Our aim is to show that $\Coim_\infty\e$ has a structure of stereotype algebra, under which the morphisms $\coim_\infty\e$ and $\mu$ become morphisms in ${\tt Ste}^\circledast$ -- this will mean that the epimorphism $\coim_\infty\e$ is a mediator for  $\e$ in the category ${\tt Ste}^\circledast$, and, since $\e$ is an immediate epimorphism, $\mu$ must be an isomorphism in ${\tt Ste}^\circledast$, and hence in ${\tt Ste}$ as well. This allows to conclude that the epimorphism $\e$ is isomorphic in ${\tt Ste}$ to the epimorphism $\Coim_\infty\e$, which is an immediate epimorphism in ${\tt Ste}$, and thus, $\e$ is also an immediate epimorphism in ${\tt Ste}$.

The existence of the structure of stereotype algebra on $\Coim_\infty\e$ follows from Theorems \ref{th-10.14} and \ref{TH:lim-inj-ster-algebr-I-napravleno}: on the one hand, any operation of the form $A'\mapsto (A'/I)^\triangledown$ (where $I$ is a closed two-sided ideal in $A'$) turns each stereotype algebra $A'$ into a stereotype algebra, and on the other hand, the injective limit in ${\tt Ste}$ of the system of stereotype algebras that one can form from $A$ in this way, is a stereotype algebra. Theorem  \ref{TH:lim-inj-ster-algebr-I-napravleno} implies also that the natural map of $A$ into this injective limit $\Coim_\infty\e$ is a morphism of stereotype algebras.

It remains to check that $\mu$ is a morphism of stereotype algebras as well, i.e. it is multiplicative and it preserves identity. Preserving identity follows from the same property for $\e$ and $\coim_\infty\e(1_A)$:
$$
\mu(1_C)=\mu\big(\coim_\infty\e(1_A)\big)=\e(1_A)=1_B.
$$
The multiplicativity of $\mu$ on the subalgebra $\coim_\infty\e(A)$ follows from the multiplicativity of $\e$ and $\coim_\infty\e(1_A)$: for any  $a,b\in A$ we have
$$
\mu(\coim_\infty\e(a)\cdot\coim_\infty\e(b))=\mu(\coim_\infty\e(a\cdot b))=\e(a\cdot b)=\e(a)\cdot\e(b)=\mu(\coim_\infty\e(a))\cdot\mu(\coim_\infty\e(b))
$$
After that we should recall that $\coim_\infty\e$ is an epimorphism in ${\tt Ste}$, so the algebra $\coim_\infty\e(A)$ is dense in $\Coim_\infty\e$. Hence, by Lemma \ref{LM:o-plotnoi-podalgebre}, $\mu$ must be multiplicative on $\Coim_\infty\e$.

2. Let us now prove (ii)$\Leftarrow$(iv). Suppose $\e:A\to B$ is a strong epimorphism in ${\tt Ste}$. Consider a diagram in ${\tt Ste}^\circledast$
$$
\begin{diagram}
\node{A}\arrow{s,l}{\alpha}\arrow{e,t}{\e}\node{B}\arrow{s,r}{\beta}
\\
\node{C}\arrow{e,b}{\mu}\node{D}
\end{diagram}
$$
where $\mu$ is a monomorphism. It  can be considered as a diagram in ${\tt Ste}$, and since $\mu$ is a monomorphism in ${\tt Ste}$ (by Example \ref{EX:mono-v-Ste^circledast}), and $\e$ a strong epimorphism in ${\tt Ste}$, there must exist a morphism $\delta$ in ${\tt Ste}$ (i.e. a linear continuous map) such that the following diagram is commutative:
$$
\begin{diagram}
\node{A}\arrow{s,l}{\alpha}\arrow{e,t}{\e}\node{B}\arrow{s,r}{\beta}\arrow{sw,r,--}{\delta}
\\
\node{C}\arrow{e,b}{\mu}\node{D}
\end{diagram}
$$
It remains to check that the map $\delta$ is a homomorphism of algebras. Preserving identity follows from monomorphity of $\mu$:
$$
\mu(1_C)=1_D=\beta\big(\e(1_A)\big)=\mu\Big(\delta\big(\e(1_A)\big)\Big)=\mu\Big(\delta\big(1_{B}\big)\Big)
\quad\Longrightarrow\quad 1_C=\delta\big(1_{B}\big).
$$
By the same reason $\delta$ is multiplicative on the subalgebra $\e(A)$: for each $a,b\in A$
$$
\mu(\delta(\e(a\cdot b)))=\beta(\e(a\cdot b)) =\beta(\e(a)\cdot\beta(\e(b))=\mu(\delta(\e(a)))\cdot\mu(\delta(\e(b))) =\mu(\delta(\e(a))\cdot\delta(\e(b)))
$$
$$
\Downarrow
$$
$$
\delta(\e(a\cdot b))=\delta(\e(a))\cdot\delta(\e(b)).
$$
After that the multiplicativity of $\delta$ on $B$ follows from Lemma \ref{LM:o-plotnoi-podalgebre}.
\epr

\btm\label{TH:strogie-epi-razlich-mono-v-Ste^circledast} If a morphism of stereotype algebras $\ph:A\to B$ is not a monomorphism, then there exists a decomposition $\ph=\ph'\circ\e$, where $\e$ is a strong epimorphism, but not an isomorphism.
\etm
\bpr If $\ph$ is not a monomorphism, then its kernel $I=\Ker\ph$ is a nonzero closed ideal in $A$. By Theorem \ref{th-10.14} the quotient space  $(A/I)^\triangledown$ is a stereotype algebra. The homomorphism of algebras $\ph$ can be lifted to some homomorphism of algebras $\psi:A/I\to B$, which by definition of usual quotient topology is a continuous map:
$$
\xymatrix
{
A\ar[dr]_{\ph}\ar[r]^{\pi} & A/I\ar@{-->}[d]^{\psi} \\
 & B
}
$$
Since the space $B$ is pseudocomplete, the map $\psi$ can be extended to a continuous map $\ph':(A/I)^\triangledown\to B$
$$
\xymatrix
{
A\ar[dr]_{\ph}\ar[r]^{\pi} & A/I\ar[d]^{\psi}\ar[r]^{\triangledown_{A/I}} & (A/I)^\triangledown \ar@{-->}[dl]^{\ph'} \\
 & B &
}
$$
By Theorem \ref{TH:strogie-epi-v-Ste^circledast=strogie-epi-v-Ste} the map $\upsilon=\triangledown_{A/I}\circ\pi:A\to (A/I)^\triangledown$ is a strong epimorphism of stereotype algebras, so we only have to verify that $\ph'$ is a homomorphism of algebras. It preserves identity since $1_{(A/I)^\triangledown}=1_{A/I}$:
$$
\ph'(1_{(A/I)^\triangledown})=\psi(1_{A/I})=1_B
$$
And its multiplicativity follows from Lemma \ref{LM:o-plotnoi-podalgebre}, since $\psi$ is multiplicative.
\epr

\paragraph{Discerning properties of strong monomorphisms in ${\tt Ste}^\circledast$.}

\blm\label{LM:neposr-podpr-porozhd-podalgebroi=proj-ster-algebra} Let $A$ be a stereotype algebra and $B$ a subalgebra in $A$ (in the purely algebraic sense). Then the envelope $\Env^A B$ of the set $B$ in the stereotype space $A$ is a stereotype algebra.
\elm
\bpr This follows from the completeness of the category ${\tt Ste}^\circledast$ (Theorem \ref{TH:polnota-Ste^circledast}) and from the fact that the pseudosaturation of closure $\overline{C}^\vartriangle$ of any subalgebra $C$ in $A$ is always a stereotype algebra by Theorem \ref{th-10.13}. \epr

\blm\label{TH:neposr-mono=strogie-mono-v-Ste^circledast} In the category ${\tt Ste}^\circledast$ of stereotype algebras the immediate monomorphisms coincide with the strong monomorphisms.
\elm
\bpr We already noticed (the property $2^\circ$ on page \pageref{mu-Smono-=>-mu-Nmono}) that each strong monomorphism is an immediate monomorphism, so we have to verify that in ${\tt Ste}^\circledast$ the inverse is also true. Let $\mu:C\to D$ be an immediate monomorphism of stereotype algebras. Consider a diagram
$$
\begin{diagram}
\node{A}\arrow{s,l}{\alpha}\arrow{e,t}{\e}\node{B}\arrow{s,r}{\beta}
\\
\node{C}\arrow{e,b}{\mu}\node{D}
\end{diagram}
$$
where $\e$ is an epimorphism. Consider the subset in $\mu(C)\cup\beta(B)$ in $D$. Let $\alg(\mu(C)\cup\beta(B))$ be the subalgebra (in the purely algebraic sense) in $D$ generated by $\mu(C)\cup\beta(B)$, and $R=\Env^D \Big(\alg(\mu(C)\cup\beta(B))\Big)$ the envelope of the set $\alg(\mu(C)\cup\beta(B))$ in $D$  (in the sense of the definition on page \pageref{DEF:Span_infty}). By Lemma  \ref{LM:neposr-podpr-porozhd-podalgebroi=proj-ster-algebra}, $R$ is a stereotype algebra. Let $\sigma:R\to D$ denote its natural enclosure in $D$. Since $\mu(C)\subseteq R$, and $R$ is an immediate subspace in $D$, the morphism of stereotype spaces $\mu$ can be factored through the morphism of stereotype spaces $\sigma:R\to D$,
$$
\mu=\sigma\circ\pi
$$
Here $\pi$ must be multiplicative, since from the identities
$$
\sigma(\pi(x\cdot y))=\mu(x\cdot y)=\mu(x)\cdot\mu(y)=\sigma(\pi(x))\cdot\sigma(\pi(y))=\sigma(\pi(x)\cdot\pi(y))
$$
imply by monomorphity of $\sigma$ the identity
$$
\pi(x\cdot y)=\pi(x)\cdot\pi(y).
$$
So we can conclude that $\pi$ is a morphism of stereotype algebras. Similarly, the enclosure $\beta(B)\subseteq R$ implies that the morphism of stereotype spaces $\beta$ can be factored through the morphism of stereotype spaces $\sigma:R\to D$,
$$
\beta=\sigma\circ\rho
$$
and again the monomorphity of $\sigma$ implies that $\rho$ is a morphism of stereotype algebras.

So we obtain a diagram in the category ${\tt Ste}^\circledast$:
$$
\begin{diagram}
\node{A}\arrow[2]{s,l}{\alpha}\arrow[2]{e,t}{\e}\node[2]{B}\arrow[2]{s,r}{\beta} \arrow{sw,r,--}{\rho}
\\
\node[2]{R}\arrow{se,r,--}{\sigma} \\
\node{C}\arrow[2]{e,b}{\mu}\arrow{ne,r,--}{\pi}\node[2]{D}
\end{diagram}
$$
Let us show that $\pi$ is an epimorphism (in ${\tt Ste}^\circledast$). Let $\zeta,\eta:R\rightrightarrows T$ be two parallel morphisms of stereotype algebras. Then the equality
$$
\zeta\circ\pi=\eta\circ\pi
$$
implies, on the one hand, the identity
 $$
\zeta\Big|_{\pi(C)}=\eta\Big|_{\pi(C)},
 $$
and, on the other hand, it implies the chain
$$
\zeta\circ\rho\circ\e=\zeta\circ\pi\circ\alpha=\eta\circ\pi\circ\alpha=\eta\circ\rho\circ \underset{\scriptsize\begin{matrix}\text{\rotatebox{90}{$\owns$}}\\ \Epi\end{matrix}}{\e}\quad\Longrightarrow\quad \zeta\circ\rho=\eta\circ\rho
\quad\Longrightarrow\quad \zeta\Big|_{\rho(B)}=\eta\Big|_{\rho(B)}.
$$
Together they give
 $$
\zeta\Big|_{\pi(C)\cup\rho(B)}=\eta\Big|_{\pi(C)\cup\rho(B)}\quad\Longrightarrow\quad
\zeta\Big|_{\alg(\pi(C)\cup\rho(B))}=\eta\Big|_{\alg(\pi(C)\cup\rho(B))}.
 $$
Let us recall that formally $R$ is a subset in $B$, so the set $\alg(\pi(C)\cup\rho(B))$ formally coincides with the set $\alg(\mu(C)\cup\beta(B))$. As a corollary, $\alg(\pi(C)\cup\rho(B))=\alg(\mu(C)\cup\beta(B))$ is dense in $R$, and we obtain that $\zeta=\eta$.

This proves that $\pi$ is an epimorphism of stereotype algebras. Thus, $\mu$ is decomposed into a composition of an epimorphism $\pi$ and a monomorphism $\sigma$. Since $\mu$ is an immediate monomorphism, $\pi$, being a mediator, must be an isomorphism. Now we can put  $\delta=\pi^{-1}\circ\rho$, and we obtain the required diagram
$$
\begin{diagram}
\node{A}\arrow{s,l}{\alpha}\arrow{e,t}{\e}\node{B}\arrow{s,r}{\beta}\arrow{sw,b,--}{\delta}
\\
\node{C}\arrow{e,b}{\mu}\node{D}
\end{diagram}
$$
\epr

\btm\label{TH:strogie-mono-razlich-epi-v-Ste^circledast} If a morphism of stereotype algebras $\ph:A\to B$ is not an epimorphism, then there exists a decomposition $\ph=\lambda\circ\ph'$ (in ${\tt Ste}^\circledast$), where $\lambda$ is a strong monomorphism, but not an ismorphism.
\etm
\bpr 1. Denote by $P$ the envelope in $B$ of the set $\ph(A)$:
$$
P=\Env^B\ph(A).
$$
By Lemma \ref{LM:neposr-podpr-porozhd-podalgebroi=proj-ster-algebra}, $P$ is a stereotype algebra, and the set-theoretic enclosure $\iota:P\to B$ is a monomorphism of stereotype algebras (and an immediate monomorphism of stereotype spaces). Let $\varPhi$ be the class of all factorizations of the morphism $\iota$ in ${\tt Ste}^\circledast$,
 \beq\label{PROOF:strogie-mono-razlich-epi-v-Ste^circledast-1}
\begin{diagram}
\node{P}\arrow[2]{e,t}{\iota}\arrow{se,b}{\Epi\owns\pi}\node[2]{B}\\
\node[2]{X}\arrow{ne,b}{\mu\in\Mono}
\end{diagram}
 \eeq
where the algebra $X$ as a set lies between $P$ and $B$:
 \beq\label{P<=X<=B}
P\subseteq X\subseteq B
 \eeq
This class is not empty, since it contains the factorization $\iota=\iota\circ 1$, and it is full in the class of all factorizations (i.e. each factorization of $\iota$ is isomorphic to some factorization from $\varPhi$). Every factorization from $\varPhi$ is uniquely defined by the set $X$ in $B$ and a topology on $X$, i.e. by a subspace $X$ in the topological space $B$. Since all subspaces of a given topological spaces form a set, we obtain that $\varPhi$ must be a set (not just a class). For simplicity we can conceive $\varPhi$ as just a set of subalgebras $X$ in $B$ satisfying  \eqref{P<=X<=B} and endowed a topology that turns $X$ into stereotype algebras in such a way that the enclosures \eqref{P<=X<=B} are continuous maps (this will mean that they are morphisms of stereotype algebras). For any $X\in\varPhi$ the set-theoretic enclosures $P\subseteq X$ and $X\subseteq B$ will be denoted by $\pi_X$ and $\mu_X$. Thus, diagram \eqref{PROOF:strogie-mono-razlich-epi-v-Ste^circledast-1} turns into diagram
 \beq\label{PROOF:strogie-mono-razlich-epi-v-Ste^circledast-2}
\begin{diagram}
\node{P}\arrow[2]{e,t}{\iota}\arrow{se,b}{\pi_X}\node[2]{B}\\
\node[2]{X}\arrow{ne,b}{\mu_X}
\end{diagram}
 \eeq

Put
$$
Y=\bigcup_{X\in\varPhi} X,
$$
then
$$
Q=\Env^B \alg Y
$$
and $\varkappa$ and $\lambda$ are the enclosures $P\subseteq Q$ and $Q\subseteq B$ respectively:
$$
\begin{diagram}
\node{P}\arrow[2]{e,t}{\iota}\arrow{se,b}{\varkappa}\node[2]{B}\\
\node[2]{Q}\arrow{ne,b}{\lambda}
\end{diagram}
$$
By Lemma \ref{LM:neposr-podpr-porozhd-podalgebroi=proj-ster-algebra}, $Q$ is a stereotype algebra, and this means that $\varkappa$ and $\lambda$ are  (mono)morphisms of stereotype algebras. For any $X\in\varPhi$ we denote by $\sigma_X$ the enclosure $X\subseteq Q$. The topology of $X$ majorizes the topology of $Q$, hence $\sigma_X$ is a continuous map, and we obtain a diagram in the category ${\tt Ste}^\circledast$:
 \beq\label{PROOF:strogie-mono-razlich-epi-v-Ste^circledast-3}
 \xymatrix % @R=2.5pc @C=2.5pc
 {
 P\ar[dr]_{\pi_X}\ar[rr]^{\iota}\ar@/_4ex/[ddr]_{\varkappa} & & B\\
  & X\ar[ur]_{\mu_X} \ar[d]_{\sigma_X} & \\
   & Q\ar@/_4ex/[uur]_{\lambda}  &
 }
 \eeq

2. Let us show now that $\varkappa$ is (not only a monomorphism, but also) an epimorphism of stereotype algebras. Indeed, for any two morphisms $\zeta,\eta:Q\rightrightarrows T$ we have the following chain:
$$
\zeta\circ\varkappa=\eta\circ\varkappa
$$
$$
\Downarrow
$$
\begin{multline*}
\forall X\in\varPhi\qquad \zeta\circ\sigma_X\circ\pi_X=\eta\circ\sigma_X\circ\underset{\scriptsize\begin{matrix}\text{\rotatebox{90}{$\owns$}}\\ \Epi\end{matrix}}{\pi_X}
\quad\Longrightarrow\quad
\forall X\in\varPhi\qquad \zeta\circ\sigma_X=\eta\circ\sigma_X
\quad\Longrightarrow\quad
\forall X\in\varPhi\qquad \zeta\Big|_X=\eta\Big|_X
\quad\Longrightarrow \\ \Longrightarrow\quad
\zeta\Big|_Y=\zeta\Big|_{\bigcup_{X\in\varPhi} X}=\eta\Big|_{\bigcup_{X\in\varPhi} X}=\eta\Big|_Y
\quad\Longrightarrow\quad
\zeta\Big|_{\alg Y}=\eta\Big|_{\alg Y}
\quad\Longrightarrow\quad
\zeta=\zeta\Big|_Q=\eta\Big|_Q=\eta
\end{multline*}
(the last implication follows from  the fact the the vector space $\alg Y$ is dense in its envelope).

3. Let us show that $\lambda:Q\to B$ is an immediate monomorphism (in ${\tt Ste}^\circledast$). Suppose $\lambda=\lambda'\circ\e$ is its arbitrary factorization. Denote by $R$ the range of $\e$ (and the domain of $\lambda'$), then we have a diagram:
\beq\label{diagr:dok-neposr-monom-lambda}
 \xymatrix % @R=2.5pc @C=2.5pc
 {
 P\ar[dr]_{\varkappa}\ar[rr]^{\iota}\ar@/_4ex/[ddr]_{\e\circ\varkappa} & & B\\
  & Q\ar[ur]_{\lambda} \ar[d]_{\e} & \\
   & R\ar@/_4ex/[uur]_{\lambda'}  &
 }
\eeq
The morphism $\e\circ\varkappa$ is an epimorphism (as a composition of two epimorphisms), so the decomposition  $\iota=\lambda'\circ(\e\circ\varkappa)$ is a factorization of $\iota$. As a corollary, it is isomorphic to some standard factorization $\iota=\mu_X\circ\pi_X$ for some $X\in\varPhi$:
$$
 \xymatrix % @R=2.5pc @C=2.5pc
 {
 P\ar[dr]_{\varkappa}\ar[rr]^{\iota}\ar@/_4ex/[ddr]_{\e\circ\varkappa}\ar@/_10ex/[dddr]_{\pi_X} & & B\\
  & Q\ar[ur]_{\lambda} \ar[d]_{\e} & \\
   & R\ar@/_4ex/[uur]_{\lambda'}\ar@{-->}[d]_{}  & \\
   & X\ar@/_10ex/[uuur]_{\mu_X} &
 }
$$
(here the dashed arrow is some isomorphism of stereotype algebras). So from the very beginning we can think that in \eqref{diagr:dok-neposr-monom-lambda} some $X\in\varPhi$ stands instead of $R$:
$$
 \xymatrix % @R=2.5pc @C=2.5pc
 {
 P\ar[dr]_{\varkappa}\ar[rr]^{\iota}\ar@/_4ex/[ddr]_{\pi_X} & & B\\
  & Q\ar[ur]_{\lambda} \ar[d]_{\e} & \\
   & X\ar@/_4ex/[uur]_{\mu_X}  &
 }
$$
Here every arrow is a set-theoretic enclosure, and the topology on the beginning of the arrow majorizes the topology on its end. In particular, the arrow $\e$ means that $Q$ is a subset of $X$, and the topology of $Q$ majorizes the topology of $X$. But on the other hand the arrow $\sigma_X$ in diagram \eqref{PROOF:strogie-mono-razlich-epi-v-Ste^circledast-3} means that on the contrary $X$ is a subset in $Q$, and the topology of $X$ majorizes the topology of $Q$. Together this means that $X$ and $Q$ coincide with the topologies:
$$
X\cong Q.
$$
In particular, $\e$ is an isomorphism, and this is what we had to verify.

4. Since $\lambda$ is an immediate monomorphism, by Lemma \ref{TH:neposr-mono=strogie-mono-v-Ste^circledast}, we obtain that $\lambda$ is a strong monomorphism.

5. Note that since $\ph(A)\subseteq P$, the morphism $\ph$ is factored through $P$:
$$
\ph=\iota\circ\theta,
$$
for some morphism $\theta:A\to P$. We obtain a diagram in ${\tt Ste}^\circledast$:
$$
\begin{diagram}
\node{A}\arrow{e,t}{\ph}\arrow{s,l}{\theta}\node{B}\\
\node{P}\arrow{ne,b}{\iota}\arrow{e,b}{\varkappa}\node{Q}\arrow{n,r}{\lambda}
\end{diagram}
$$
We see now that $\lambda$ cannot be an isomorphism, since otherwise $\ph$ would be an epimorphism, as a composition of two epimorphisms $\theta$ and $\varkappa$, and an isomorphism $\lambda$. So if we put $\ph'=\varkappa\circ\theta$, we obtain a decomposition $\ph=\lambda\circ\ph'$, where  $\lambda$ is a strong monomorphism, but not an isomorphism.
\epr

\paragraph{Nodal decomposition in ${\tt Ste}^\circledast$.}
Let us notice the following two properties of the category ${\tt Ste}^\circledast$.

\btm\label{TH:Ste^circledast-lok-mala-v-pod-objektah} The category ${\tt Ste}^\circledast$ of stereotype algebras is well-powered.
\etm
\bpr
A morphism $\mu:A\to B$ in ${\tt Ste}^\circledast$ is a monomorphism in ${\tt Ste}^\circledast$ iff it is a monomorphism in ${\tt Ste}$, and the latter category is well-powered.
\epr

\btm\label{TH:Ste^circledast-strogo-lok-mala-v-faktor-objektah} The category ${\tt Ste}^\circledast$ of stereotype algebras is co-well-powered in strong epimorphisms.
\etm
\bpr
By Theorem \ref{TH:strogie-epi-v-Ste^circledast=strogie-epi-v-Ste}, a morphism $\e:A\to B$ in ${\tt Ste}^\circledast$ is a strong epimorphism in ${\tt Ste}^\circledast$ iff it is a strong  epimorphism in ${\tt Ste}$, and the latter category is co-well-powered.
\epr

On the other hand, as we already know the category ${\tt Ste}^\circledast$ is complete (by Theorem \ref{TH:polnota-Ste^circledast}), and in ${\tt Ste}^\circledast$ the strong epimorphisms discern monomorphisms, and the strong monomorphisms discern epimorphisms (Theorems  \ref{TH:strogie-epi-razlich-mono-v-Ste^circledast} and \ref{TH:strogie-mono-razlich-epi-v-Ste^circledast}). Thus, we can apply Theorem  \ref{TH:sush-uzlov-razlozh-v-polnoi-lok-maloi-kateg}, and we get

\btm\label{TH:uzlovoe-razlozhenie-v-Ste^circledast} In the category ${\tt Ste}^\circledast$ of stereotype algebras each morphism $\ph:X\to Y$ has a nodal decomposition \eqref{DEF:oboznacheniya-dlya-uzlov-razlozh}.
\etm

\brem Theorem \ref{TH:strogie-epi-v-Ste^circledast=strogie-epi-v-Ste} implies in addition that the nodal coimage $\Coim_\infty\ph$ in ${\tt Ste}^\circledast$ coincides with the nodal coimage in ${\tt Ste}$, and as a corollary with the refinement (as a quotient space of a stereotype space) on $X$ of a set of functionals $\ph^\star(Y^\star)$:
    \beq\label{Coim_infty-ph=Cosp_infty-ph*(Y*)-v-Ste^circledast}
    \Coim_\infty\ph=\Rf^X\ph^\star(Y^\star)
    \eeq
For the nodal image $\Im_\infty\ph$ the analogous proposition is not true.
\erem

\btm\label{TH:uzlovoe-razlozhenie-v-Ste-kak-v-Ste^circledast} For each morphism $\ph:A\to B$ in the category ${\tt Ste}^\circledast$ of stereotype algebras its nodal decomposition $\ph=\im_\infty\ph\circ\red_\infty\ph\circ\coim_\infty\ph$ in the category ${\tt Ste}$ of stereotype spaces is a decomposition (not necessarily, nodal) in the category ${\tt Ste}^\circledast$.
\etm
\bpr
We need to see here that the stereotype spaces $\Coim_\infty\ph$ and $\Im_\infty\ph$ have natural structure of stereotype algebras, and that the morphisms of stereotype spaces $\coim_\infty\ph:A\to\Coim_\infty\ph$, $\red_\infty\ph:\Coim_\infty\ph\to\Im_\infty\ph$, $\im_\infty\ph:\Im_\infty\ph\to B$, are morphisms of stereotype algebras (i.e., homomorphisms of algebras). This follows from the way of constructing $\Coim_\infty\ph$ and $\Im_\infty\ph$: since $\ph:A\to B$ is a morphism of stereotype algebras, its reduced morphism  $\ph^1=\red\ph:\Coim\ph\to\Im\ph$ is also a morphism of stereotype algebras (together with the morphisms $\coim\ph:A\to\Coim\ph$ and $\im\ph:\Im\ph\to B$). By the same reason the second reduced morphism $\ph^2=\red\ph^1$ must be a morphism of stereotype algebras, and so on. We have to organize a transfinite induction by the degree of this operation, and we will obtain that the nodal coimage $\Coim_\infty\ph$ is a stereotype algebra (as an injective limit of stereotype algebras $\Coim\ph^i$), the nodal image $\Im_\infty\ph$ is a stereotype algebra (as a projective limit of stereotype algebras $\Im\ph^i$), and the morphisms $\coim_\infty\ph:A\to\Coim_\infty\ph$, $\red_\infty\ph:\Coim_\infty\ph\to\Im_\infty\ph$, $\im_\infty\ph:\Im_\infty\ph\to B$ are homomorphisms of algebras.
\epr

\paragraph{Envelopes and refinements in ${\tt Ste}^\circledast$.}

Since it is not clear whether the category ${\tt Ste}^\circledast$ is co-well-powered in the class $\Epi$, in the analogue of Theorem \ref{TH:obolochki-i-nachinki-otn-klassa-morfizmov-v-Ste} for the class of envelopes in $\Epi$ one should assume that the class of test morphisms $\varPhi$ is a set (so that in the proof the property $5^\circ$ on p.\pageref{5^0:obolochka-otn-klassa-morphizmov} could be replaced by $3^\circ$):

\btm\label{TH:obolochki-i-nachinki-otn-klassa-morfizmov-v-Ste^circledast} In the category ${\tt Ste}^\circledast$ of stereotype algebras
 \bit{
\item[(a)] each algebra $A$ has an envelope in the class $\Epi$ of all epimorphisms (respectively, in the class $\SEpi$ of all strong epimorphisms) of the category ${\tt Ste}^\circledast$ with respect to arbitrary {\rm set} (respectively, {\rm class}) of morphisms $\varPhi$, going from $A$; in addition,
     \bit{
\item[(i)] if $\varPhi$ differs morphisms on the outside in ${\tt Ste}^\circledast$, then the envelope in $\Epi$ is an envelope in the class $\Bim$ of all bimorphisms:
    $$
    \env_\varPhi^{\Epi} A=\env_\varPhi^{\Bim} A,
    $$

\item[(ii)] if $\varPhi$ differs morphisms on the outside and is a right ideal in ${\tt Ste}^\circledast$, then the envelope in the class $\Epi$ is also an envelope in any other class $\varOmega$, containing $\Bim$ (for example, in the class $\Mor$ of all morphisms):
    $$
\env_\varPhi^{\Epi} A=\env_\varPhi^{\Bim} A=\env_\varPhi^\varOmega A=\env_\varPhi A,\qquad \varOmega\supseteq\Bim.
    $$
 }\eit

\item[(b)] each algebra $A$ has a refinement in the class $\Mono$ of all monomorphisms (respectively, the class $\SMono$ of all strong monomorphisms) in the category ${\tt Ste}^\circledast$ by means of an arbitrary class $\varPhi$ of morphisms, going to $A$; in addition,
     \bit{
\item[(i)] if $\varPhi$ differs morphisms on the inside in ${\tt Ste}^\circledast$, then the refinement in the class $\Mono$ is also a refinement in the class $\Bim$ of all bimorphisms:
    $$
    \rf_\varPhi^{\Mono} A=\rf_\varPhi^{\Bim} A.
    $$

\item[(ii)] if $\varPhi$ differs morphisms on the inside and is a left ideal in ${\tt Ste}^\circledast$, then the refinement in $\Mono$ is also a refinement in any other class $\varGamma$, containing $\Bim$ (for example, the class  $\Mor$ of all morphisms):
    $$
    \rf_\varPhi^{\Mono} A=\rf_\varPhi^{\Bim} A=\rf_\varPhi^\varGamma A=\rf_\varPhi A,\qquad \varGamma\supseteq\Bim.
    $$
 }\eit
 }\eit
\etm
\bpr Consider the case of envelopes. If $\varPhi$ is a set, then the existence of $\env_\varPhi^{\Epi({\tt Ste}^\circledast)} A$ follows from $3^\circ$ on p.\pageref{3^0:obolochka-otn-mnozhestva-morphizmov}.
If $\varPhi$ differs morphisms on the outside, then by Theorem \ref{TH:Phi-razdel-moprfizmy} the existence of  $\env_\varPhi^{\Epi} A$ implies the existence of $\env_\varPhi^{\Epi\cap\Mono} A=\env_\varPhi^{\Bim} A$, the the coincidence of these envelopes: $\env_\varPhi^{\Epi} A=\env_\varPhi^{\Bim} A$. If $\varPhi$ differs morphisms on the outside and is a right ideal, then by Theorem \ref{TH:Phi-razdel-moprfizmy-*} the existence of $\env_\varPhi^{\Bim} A$ implies the existence of $\env_\varPhi^{\varOmega}A$ for any $\varOmega\supseteq\Bim$, and the coincidence of these envelopes: $\env_\varPhi^{\Bim}A=\env_\varPhi^{\varOmega} A$.
\epr

From Theorems \ref{TH:sushestvovanie-seti-pri-faktorizatsii} and \ref{TH:regulyarnaya-obolochka} (with $\varOmega=\Epi$) we have

\btm\label{TH:sushestvovanie-seti-v-Ste^circledast-dlya-Epi} Let $\varPhi$ be a class of morphisms in ${\tt Ste}^\circledast$, which goes from ${\tt Ste}^\circledast$ and is a right ideal:
$$
\varPhi\circ\Mor({\tt Ste}^\circledast)\subseteq\varPhi.
$$
Then the classes of morphisms $\Epi$ and $\varPhi$ define in ${\tt Ste}^\circledast$ a semiregular envelope $\Env_\varPhi^{\Epi}$, and for each object $A$ in ${\tt Ste}^\circledast$ the envelope is described by the formula
\beq\label{sushestvovanie-seti-v-Ste^circledast-dlya-Epi-1}
\red_\infty \leftlim{\mathcal N}^A\circ\coim_\infty \leftlim{\mathcal N}^A=\env_\varPhi^{\Epi} A,
\eeq
where ${\mathcal N}$ is the net of epimorphisms, generated by the classes $\Epi$ and $\varPhi$, and
$\red_\infty \leftlim{\mathcal N}^A$ and $\coim_\infty \leftlim{\mathcal N}^A$ are elements of nodal decomposition
\eqref{DEF:oboznacheniya-dlya-uzlov-razlozh} of the morphism
$\leftlim{\mathcal N}^A:A\to A_{\mathcal N}$ in the category ${\tt Ste}^\circledast$. If in addition the class $\Epi$ pushes $\varPhi$, then the envelope $\Env_\varPhi^{\Epi}$ is regular (and thus, it can be defined as an idempotent functor).
 \etm

\paragraph{Dense epimorphisms.}

 \bit{
\item[$\bullet$]\label{DEF:DEpi} Let us say that a morphism of stereotype (or, in general case, topological) algebras $\ph:A\to B$ is {\it dense}, if its set of values $\ph(A)$ is dense in $B$:
$$
\overline{\ph(A)}=B.
$$
Certainly, dense morphisms are epimorphisms, so we also call them dense epimorphisms. The class of all dense epimorphisms in ${\tt Ste}^\circledast$ (or in ${\tt TopAlg}$) will be denoted by $\DEpi$. It is connected with the classes $\Epi$ of epimorphisms and $\SEpi$ of strong epimorphisms by the inclusions
\beq\label{SEpi-subset-DEpi-subset-Epi}
\SEpi\subset\DEpi\subset\Epi.
\eeq
 }\eit

\brem
Inclusions \eqref{SEpi-subset-DEpi-subset-Epi} are not equalities. An example of a dense epimorphism, which is not strong, is the set-theoretic inclusion of the algebra ${\mathcal C}^\infty(M)$ of smooth functions into the algebra ${\mathcal C}(M)$ of continuous functions on a smooth manifold $M$ (this inclusion is a bimorphism of stereotype algebras, so if it was a strong epimorphism, this automatically would mean that this is an isomorphism, but this is not true). An example of a non-dense epimorphism is the standard inclusion of the algebra ${\mathcal P}(\C)$ of polynomials on $\C$ into the algebra ${\mathcal P}(\C^\times)$ of Laurent polynomials on $\C^\times$ (we already mentioned this example on page \pageref{EX:epi-v-Ste^circledast}).
\erem

\btm\label{TH:SMono_Ste-circledcirc-DEpi=Ste^circledast}
The class $\DEpi$ of dense epimorphisms is monomorphically complementable in ${\tt Ste}^\circledast$.
\etm
\bpr
The monomorphic complement for $\DEpi$ is the class $\SMono_{\tt Ste}$ of strong monomorphisms in ${\tt Ste}^\circledast$, which are strong monomorphisms in ${\tt Ste}$
\beq\label{SMono_Ste-circledcirc-DEpi=Ste^circledast}
\SMono_{\tt Ste}\circledcirc\DEpi={\tt Ste}^\circledast
\eeq
\epr

For dense epimorphisms the first part of Theorem \ref{TH:obolochki-i-nachinki-otn-klassa-morfizmov-v-Ste^circledast} can be strengthened as follows:

\btm\label{TH:obolochki-i-nachinki-otn-plotnyh-epimorfizmov-v-Ste^circledast} In the category ${\tt Ste}^\circledast$ every algebra $A$ has an envelope in the class $\DEpi$ of dense epimorphisms with respect to the arbitrary {\rm class} of morphisms $\varPhi$ going from $A$. If in addition $\varPhi$ differs morphisms on the outside in ${\tt Ste}^\circledast$, then the envelope in $\DEpi$ is also an envelope in the class $\DBim$ of all dense bimorphisms:
    $$
    \env_\varPhi^{\DEpi} A=\env_\varPhi^{\DBim} A.
    $$
\etm
\bpr
The existence of $\env_\varPhi^{\DEpi({\tt Ste}^\circledast)} A$ follows from $5^\circ$ on p.\pageref{5^0:obolochka-otn-klassa-morphizmov}. If $\varPhi$ differs morphisms on the outside, then by Theorem \ref{TH:Phi-razdel-moprfizmy} the existence of $\env_\varPhi^{\DEpi} A$ implies the existence of
$\env_\varPhi^{\DEpi\cap\Mono} A=\env_\varPhi^{\DBim} A$, and their coincidence,
$\env_\varPhi^{\DEpi} A=\env_\varPhi^{\DBim} A$. If in addition $\varPhi$ is a right ideal, then by Theorem  \ref{TH:Phi-razdel-moprfizmy-*} the existence of $\env_\varPhi^{\DEpi} A$ implies the existence of
$\env_\varPhi^{\DEpi\cap\Mono} A=\env_\varPhi^{\DBim} A$ and their coincidence,
$\env_\varPhi^{\DEpi} A=\env_\varPhi^{\DBim} A$.
\epr

From Theorems \ref{TH:sushestvovanie-seti-pri-faktorizatsii} and \ref{TH:regulyarnaya-obolochka} (with $\varOmega=\DEpi$) it follows

\btm\label{TH:sushestvovanie-seti-v-Ste^circledast-dlya-DEpi}
Let $\varPhi$ be a class of morphisms in ${\tt Ste}^\circledast$, which goes from ${\tt Ste}^\circledast$ and is a right ideal:
$$
\varPhi\circ\Mor({\tt Ste}^\circledast)\subseteq\varPhi.
$$
Then the classes $\DEpi$ and $\varPhi$ define a semiregular envelope
$\Env_\varPhi^{\DEpi}$ in ${\tt Ste}^\circledast$, and for any object $A$ in ${\tt Ste}^\circledast$ it is described by the formula
\beq\label{sushestvovanie-seti-v-Ste^circledast-dlya-DEpi-1}
\env_\varPhi^{\DEpi} A=\red_\infty \leftlim{\mathcal N}^A\circ\coim_\infty \leftlim{\mathcal N}^A,
\eeq
where ${\mathcal N}$ is the net of epimorphisms, generated by classes $\DEpi$ and $\varPhi$, and
$\red_\infty \leftlim{\mathcal N}^A$ and $\coim_\infty \leftlim{\mathcal N}^A$ are elements of nodal decomposition
\eqref{DEF:oboznacheniya-dlya-uzlov-razlozh} of the morphism $\leftlim{\mathcal N}^A:A\to A_{\mathcal N}$ in the category ${\tt Ste}$ of stereotype spaces (not algebras!). If in addition $\DEpi$ pushes $\varPhi$, then the envelope $\Env_\varPhi^{\DEpi}$ is regular (and thus, it can be defined as an idempotent functor).
 \etm

\subsection{Holomorphic envelope}\label{SUBSEC;Arens-Michael}

A.~Ya.~Helemskii introduced in \cite{Helemskii} the notion of the Arens-Michael envelope in the category of topological algebras. The properties of this construction used in the duality theory for complex Lie groups \cite{Akbarov-stein-groups} have different formal interpretations (under preserving the essential results) as an envelope in the sense of definition of \ref{obolochka,otpechatok,uzl-razlozh} in the category of stereotype algebras. For one of them, which seems to be most natural, we use the (working) name {\it holomorphic envelope}. The choice of the term has the aim to emphasize the connection of this construction with the complex analysis and the analogy with the continuous envelope (which we define below on p.\pageref{DEF:C*-obolochka}) and smooth envelope from \cite{Akbarov-C-infty}.

\paragraph{Net of Banach quotient maps and the stereotype Arens-Michael envelope.} Here we define the analogue of the Arens-Michael envelope in the category ${\tt Ste}^\circledast$ of stereotype algebras. All the definitions and results can be easily transferred to the category ${\tt TopAlg}$ of topological algebras.

Recall that an absolutely convex closed neighborhood of zero $U$ in a topological algebra $A$ is said to be {\it submultiplicative}, if $U\cdot U\subseteq U$. The set of all submultiplicative absolutely convex closed neighbourhoods of zero in $A$ we denote by $\mathcal{SU}(A)$. To any such neighborhood of zero $U$ in $A$ one can assign a two-sided closed ideal $\Ker U=\bigcap_{\varepsilon>0}\varepsilon\cdot U$ in $A$ and a quotient algebra $A/\Ker U$ endowed with (not the quotient topology as one could expect, but) the topology of normed space with the unit ball $U+\Ker U$. Then the completion $(A/\Ker U)^\blacktriangledown$ is a Banach algebra, and we denote it by $A/U$ and call it the {\it quotient algebra of $A$ by the neighborhood of zero $U$}\label{DEF:A/U}. The natural map from $A$ into $A/U$
$$
 \xymatrix  @R=2.5pc @C=4pc
 {
A \ar@/^5ex/[rr]^{\rho_U}\ar[r]_{\tau_U}  & A/\Ker U \ar[r]_{\blacktriangledown_{A/\Ker U}}
 & (A/\Ker U)^\blacktriangledown=A/U
 }
$$
(where $\tau_U$ is a quotient map, and $\blacktriangledown_{A/\Ker U}$ is the completion map) will be called the {\it Banach quotient map of $A$ by the neighborhood of zero $U$}.

Denote by $\mathcal B$ the class of all Banach quotient maps $\{\rho_U:A\to A/U\}$, where $A$ runs over the class of topological algebras, and $U$ the set of all submultiplicative neighborhoods of zero in $A$.

\bprop\label{PROP:B-set-v-TopAlg} The class $\mathcal B$ of Banach quotient maps is a net of epimorphisms in the category ${\tt Ste}^\circledast$ of stereotype algebras, and the relation of pre-order\footnote{The pre-order $\to$ on the class $\Epi^X$ of all epimorphisms going from a given object $X$ of a category ${\tt K}$ was defined on p.\pageref{DEF:le-in-F_X}.}
$\to$ is equivalent to the embedding of the corresponding neighborhoods of zero up to a positive scalar multiplier:
 \beq\label{rho_V-to-rho_U-<=>-V-subseteq-U}
\rho_V\to\rho_U\quad\Longleftrightarrow\quad\exists\e>0\quad \e\cdot V\subseteq
U.
 \eeq
 \eprop
\bpr 1. Let us first verify \eqref{rho_V-to-rho_U-<=>-V-subseteq-U}. Suppose $U$ and $V$ are submultiplicative closed absolutely convex neighborhoods of zero in $A$, and $\e\cdot V\subseteq U$ for some $\e>0$. Then $\Ker V\subseteq\Ker U$, and the formula
$$
x+\Ker V\mapsto x+\Ker U
$$
defines a linear continuous map $A/\Ker V\to A/\Ker U$ which can be extended by continuity to an operator
$$
\pi_V^U: A/V=(A/\Ker V)^\triangledown\to (A/\Ker U)^\triangledown=A/U.
$$
Obviously, the following diagram is commutative:
 \beq\label{konus-banach-faktor-algebr}
 \xymatrix % @R=2.5pc @C=2.5pc
 {
 & A \ar@/_2ex/[dl]_{\rho_V} \ar@/^2ex/[dr]^{\rho_U} &
\\
A/V \ar[rr]_{\pi_V^U} & & A/U
 },
 \eeq
In particular, $\rho_V\to\rho_U$. On the contrary, if for some morphism $\iota:A/V\to A/U$ we have a commutative diagram
\beq\label{konus-banach-faktor-algebr-*}
 \xymatrix % @R=2.5pc @C=2.5pc
 {
 & A \ar@/_2ex/[dl]_{\rho_V} \ar@/^2ex/[dr]^{\rho_U} &
\\
A/V \ar[rr]_{\iota} & & A/U
 },
\eeq
then we can put $\widetilde{U}=\overline{\rho_U(U)}$ and $\widetilde{V}=\overline{\rho_V(V)}$, and these will be ball centered in zeroes in $A/U$
and $A/V$ respectively, so the continuity of the operator $\iota:A/V\to A/U$ implies that
$$
\e\cdot\widetilde{V}\subseteq\iota^{-1}(\widetilde{U})
$$
for some $\e>0$. And we have
$$
\e\cdot
V=(\rho_V)^{-1}\l\e\cdot\widetilde{V}\r\subseteq(\rho_V)^{-1}\l\iota^{-1}(\widetilde{U})\r=
(\rho_U)^{-1}(\widetilde{U})=U.
$$

2. Let us check now axiom (a) of the net of epimorphisms from the page \pageref{AX:set-Epi-a}. For each topological algebra $A$ the set
${\mathcal B}_A$ of its Banach quotient maps is non-empty, since always there exists at least one submultiplicative neighborhood of zero $U$ in $A$, namely, $U=A$ (and the corresponding quotient map is zero, $\rho_U:A\to 0$). Besides this, if $U$ and $V$ are two submultiplicative closed absolutely convex neighborhoods of zero in $A$, then, clearly, its intersection $U\cap V$ is also a submultiplicative (and closed absolutely convex) neighborhood of zero in $A$. That is, the submultiplicative absolutely convex neighborhoods of zero in $A$ form a system directed in descending order. Together with the rule \eqref{rho_V-to-rho_U-<=>-V-subseteq-U} this means that the system of epimorphisms $\{\rho_U:A\to A/U\}$ is directed to the left with respect to the pre-order $\to$.

3. Then we check axiom (b). For each topological algebra $A$ the system of connecting morphisms $\Bind({\mathcal B}_A)$ has a projective limit, since the category ${\tt Ste}^\circledast$ is complete. This limit can be defined as a map $A\mapsto\leftlim\Bind({\mathcal B}_A)$, since it is directly constructed as a set in the product of algebras $A/U$.

4. It remains to check axiom (c). Let $\alpha:A\to B$ be a morphism of topological algebras and $\rho_V:B\to B/V$ a Banach quotient map. The set  $U=\alpha^{-1}(V)$ is a submultiplicative closed absolutely convex neighborhood of zero in $A$. The map
$$
x+\Ker U\mapsto \alpha(x)+\Ker V,
$$
is extended by continuity to a map $\alpha_U^V:A/U\to B/V$ such that the following diagram is commutative:
 $$ \xymatrix @R=2.5pc @C=4.0pc {
 X\ar[r]^{\alpha}\ar[d]_{\rho_U} & Y\ar[d]^{\rho_V} \\
 A/U\ar@{-->}[r]_{\alpha_\sigma^\tau} & B/V
 }
 $$
 \epr

 \bit{
\item[$\bullet$] The net $\mathcal B$ will be called a {\it net of Banach quotient maps}.

\item[$\bullet$] For each algebra $A$ diagram \eqref{konus-banach-faktor-algebr} means that the family of quotient maps $\rho_U:A\to A/U$ is a projective cone of the contravariant system $\Bind({\mathcal B}_A)=\{\pi_V^U\}$. The projective limit of this cone in the category ${\tt Ste}^\circledast$ of stereotype algebras is called the {\it stereotype Arens-Michael envelope} of the algebra $A$ and is denoted by
 \beq\label{DEF:Arens-Michael=lim-B_A}
\leftlim{\mathcal B}_A:A\to A_{\mathcal B}
 \eeq
(this limit exists since ${\tt Ste}^\circledast$ is projectively complete). The range of this morphism
 \beq\label{DEF:Arens-Michael}
A_{\mathcal B}=\Ran\leftlim{\mathcal B}_A=
{\tt Ste}^\circledast\text{-}\kern-10pt\leftlim_{U\in\mathcal{SU}(A)}\kern-10pt A/U=
\Big({\tt TopAlg}\text{-}\kern-10pt\leftlim_{U\in\mathcal{SU}(A)}\kern-10pt A/U\Big)^\vartriangle.
 \eeq
will also be called the  {\it stereotype Arens-Michael envelope} of the algebra $A$.
 }\eit

If $U$ and $V$ are submultiplicative neighbourhoods of zero such that $\e\cdot V\subseteq U$ for some $\e>0$, then we have a commutative diagram
$$
 \xymatrix  @R=2.pc @C=4.pc
 {
 & A \ar@{-->}[d]^{\leftlim{\mathcal B}_A}\ar@/_5ex/[ddl]_{\rho_V} \ar@/^5ex/[ddr]^{\rho_U} &
\\
& A_{\mathcal B} \ar[dl]^{\pi_V} \ar[dr]_{\pi_U} &
\\
A/V \ar[rr]_{\pi_V^U} & & A/U
 }
$$
Theorem  \ref{TH:funktorialnost-obolochki_F} implies

\btm\label{TH:Arens-Michael=funktor} The Arens-Michael envelope is an envelope in the class of all morphisms in ${\tt Ste}^\circledast$ with respect to the system of Banach quotient maps  ${\mathcal B}$,
 \beq
 A_{\mathcal B}=\Env_{\mathcal B}^{\Mor({\tt Ste}^\circledast)}A,
 \eeq
and to each morphism $\ph:A\to B$ in ${\tt Ste}^\circledast$ the formula
 \beq\label{DEF:ph^heartsuit}
\ph_{\mathcal B}=\leftlim_{\tau\in{\mathcal B}_B}\leftlim_{\sigma\in{\mathcal
B}_A}\ph_\sigma^\tau\circ\sigma_{\mathcal B}
 \eeq
assigns a morphism $\ph_{\mathcal B}:A_{\mathcal B}\to B_{\mathcal B}$ such that the following diagram is commutative:
 \beq\label{DIAGR:funktorialnost-heartsuit}
\xymatrix @R=2.pc @C=5.0pc % @M=14pt
{
A\ar[d]^{\ph}\ar[r]^{\leftlim{\mathcal B}_A} &  A_{\mathcal B}\ar@{-->}[d]^{\ph_{\mathcal B}} \\
B\ar[r]^{\leftlim{\mathcal B}_B} &  B_{\mathcal B} \\
},
 \eeq
and the map $(A,\ph)\mapsto(A_{\mathcal B},\ph_{\mathcal B})$ can be defined as a functor from ${\tt Ste}^\circledast$ into ${\tt Ste}^\circledast$.
 \etm

\paragraph{Holomorphic envelope of stereotype algebra.}
Let us remind that on p.\pageref{DEF:DEpi} we defined the dense epimorphisms of topological algebras $\ph:A\to B$.

 \bit{
 \item
By {\it holomorphic envelope} of a stereotype algebra $A$ we mean its envelope in the class $\DEpi$ of dense epimorphisms of the category ${\tt Ste}^\circledast$ with respect to the class ${\tt BanAlg}$ of Banach algebras. We use the following notation for this construction:
 \beq\label{DEF:A^heartsuit}
 A^\heartsuit=\Env_{\tt BanAlg}^{\DEpi}A, \qquad \heartsuit_A=\env_{\tt BanAlg}^{\DEpi}A.
 \eeq
Thus,
 $$
 \big(\heartsuit_A:A\to A^\heartsuit\big)=\big(\env_{\tt BanAlg}^{\DEpi}A:A\to\Env_{\tt BanAlg}^{\DEpi}A\big).
 $$
 }\eit

 \vglue10pt
\centerline{\bf Properties of holomorphic envelopes:}
 \vglue10pt

\bit{\it

\item[$1^\circ$.]\label{1^0:golomorfnaya-obolochka-sushestvuet}
Each stereotype algebra $A$ has holomorphic envelope $A^\heartsuit$.

\item[$2^\circ$.] The holomorphic envelope $A^\heartsuit$ is connected with the stereotype Arens-Michael envelope $A_{\mathcal B}$ through the formulas
  \beq\label{heartsuit_A=red-circ-coim}
\heartsuit_A=\red_\infty\leftlim{\mathcal B}_A\circ\coim_\infty\leftlim{\mathcal B}_A,\qquad
A^\heartsuit=\Dom\im_\infty\leftlim{\mathcal B}_A
  \eeq
where $\coim_\infty\leftlim{\mathcal B}_A$, $\red_\infty\leftlim{\mathcal B}_A$, $\im_\infty\leftlim{\mathcal B}_A$ are elements of nodal decomposition of the morphism $\leftlim{\mathcal B}_A$ in the category ${\tt Ste}$ of stereotype spaces (not algebras!).

\item[$3^\circ$.] For any morphism $\ph:A\to B$ of stereotype algebras and for each choice of holomorphic envelopes $\heartsuit_A: A\to A^\heartsuit$ and $\heartsuit_B: B\to B^\heartsuit$ there exists a unique morphism $\ph^\heartsuit:A^\heartsuit\to B^\heartsuit$ such that the following diagram is commutative
    \beq\label{DIAGR:existence-of-ph^heartsuit}
\xymatrix @R=2.pc @C=5.0pc % @M=14pt
{
A\ar[d]^{\ph}\ar[r]^{\heartsuit_A} & A^\heartsuit\ar@{-->}[d]^{\ph^\heartsuit} \\
B\ar[r]^{\heartsuit_B} & B^\heartsuit \\
}
\eeq

\item[$4^\circ$.] The correspondence $(X,\alpha)\mapsto(X^\heartsuit,\alpha^\heartsuit)$ can be defined as a covariant functor from ${\tt Ste}^\circledast$ into ${\tt Ste}^\circledast$:
\beq\label{funktorialnost-heartsuit}
(1_A)^\heartsuit=1_{A^\heartsuit},\qquad (\beta\circ\alpha)^\heartsuit=\beta^\heartsuit\circ\alpha^\heartsuit,
\qquad (\alpha^\heartsuit)^\heartsuit=\alpha^\heartsuit.
\eeq

\item[$5^\circ$.] If an algebra $A$ is dense in its stereotype Arens-Michael envelope $A_{\mathcal B}$, i.e.
$$
\leftlim{\mathcal B}_A\in\DEpi({\tt Ste}^\circledast),
$$
then the holomorphic envelope of $A$ coincides with its envelope in the class $\Epi$ of all epimorphisms in ${\tt Ste}^\circledast$ and with the stereotype Arens-Michael envelope:
\beq\label{A^heartsuit=Env_BanAlg^DEpi-A=Env_BanAlg^Epi-A=A_B}
A^\heartsuit=\Env_{\tt BanAlg}^{\DEpi}A=\Env_{\tt BanAlg}^{\Epi}A=A_{\mathcal B}.
\eeq

\item[$6^\circ$.]\label{TH:golom-obolochka-soglasovana-s-circledast}
The holomorphic envelope is coherent with the projective tensor product $\circledast$ in the category
${\tt Ste}^\circledast$.

}\eit

In proof we shall need the following

\blm\label{LM:A->B-propusk-cherez-rho_U} In the category ${\tt Ste}^\circledast$ of stereotype algebras the net  $\mathcal B$ of Banach quotient maps consists of dense epimorphisms and generates on the inside the class $\Mor({\tt Ste}^\circledast,{\tt BanAlg})$ of morphisms with values in Banach algebras:
\beq\label{Ban-porozhdaet-Mor(Ste,Ban)-iznutri}
{\mathcal B}\subseteq \Mor({\tt Ste}^\circledast,{\tt BanAlg})\subseteq \Mor({\tt Ste}^\circledast)\circ{\mathcal B}.
\eeq
\elm
 \bpr
The class $\mathcal B$ consists of dense epimorphisms, since the image $\rho_U(A)$ of any algebra $A$ is always dense in its Banach quotient algebra $A/U=(A/\Ker U)^\blacktriangledown$. Let us show that $\mathcal B$ generates the class of morphisms with values in Banach algebras. It is important here to verify the second embedding in the chain \eqref{DEF:morfizmy-porozhdayutsya-iznutri}. Let $\ph:A\to B$ be a morphism into a Banach algebra  $B$. If $V$ is a unit ball in $B$, then the set $U=\ph^{-1}(V)$ is a neighborhood of zero in $A$, and the condition $V\cdot V\subseteq V$ implies the condition $U\cdot U\subseteq U$:
$$
x,y\in U\quad\Rightarrow\quad \ph(x),\ph(y)\in V\quad\Rightarrow\quad
\ph(x\cdot y)=\ph(x)\cdot\ph(y)\in V\quad\Rightarrow\quad x\cdot y\in
U=\ph^{-1}(V)
$$
Consider the normed algebra $A/\Ker U$ and the quotient map $\tau_U:A\to A/\Ker U$. From the obvious equality $\Ker\ph=\Ker U$ it follows that the morphism $\ph$ can be decomposed in the category ${\tt Alg}$ of algebras as follows:
$$
 \xymatrix  @R=2.5pc @C=2.5pc
 {
A \ar[r]^{\tau_U} \ar[rd]_{\ph} & A/\Ker U
\ar@{-->}[d]^{\chi}
\\
 & B
 }
$$
On the other hand, the equality $\chi^{-1}(V)=U+\Ker\ph=U+\Ker U$ implies continuity of $\chi$. So it will be continuously extended to the completion $(A/\Ker U)^\blacktriangledown=A/U$ of the space $A/\Ker U$:
$$
 \xymatrix  @R=2.5pc @C=4pc
 {
A \ar@/^5ex/[rr]^{\rho_U}\ar[r]_{\tau_U} \ar[rd]_{\ph} & A/\Ker U \ar[r]_{\blacktriangledown_{A/\Ker U}}
\ar[d]^{\chi} & A/U \ar@{-->}[dl]^{\chi^\blacktriangledown}
\\
 & B &
 }
$$
and since $A/\Ker U$ is dense in its completion $(A/\Ker U)^\blacktriangledown=A/U$, the map $\chi^\blacktriangledown$ must be multiplicative by Lemma \ref{LM:o-plotnoi-podalgebre}. At the same time, obviously, $\chi^\blacktriangledown$ preserves the identity. Hence, $\chi^\blacktriangledown$ is a morphism in ${\tt Ste}^\circledast$.
\epr

\bpr[Proof of the properties $1^\circ$-$6^\circ$.]
1. By Lemma \ref{LM:A->B-propusk-cherez-rho_U} the net of Banach quotient maps generates on the inside the class $\Mor({\tt Ste}^\circledast,{\tt BanAlg})$ of morphisms with values in Banach algebras. On the other hand, by Theorem \ref{TH:SMono_Ste-circledcirc-DEpi=Ste^circledast} the class $\DEpi$ of dense epimorphisms is monomorphically complementable in ${\tt Ste}^\circledast$. Therefore by Theorem \ref{TH:funktorialnost-pri-seti-Epi-i-dolonyaemosti} each object $A$ in ${\tt Ste}^\circledast$ there is an envelope in $\DEpi$ with respect to the class $\Mor({\tt Ste}^\circledast,{\tt BanAlg})$, and by definition this is the holomorphic envelope of $A$.

2 and 3. Formulas \eqref{heartsuit_A=red-circ-coim} follow immediately from \eqref{im_infty-lim-F_X=env_M^Epi-X-1}, and the diagram \eqref{DIAGR:existence-of-ph^heartsuit} from diagram \eqref{DIAGR:funktorialnost-env_varPhi^Epi-v-kat-s-uzl-razl-1}.

4. The category ${\tt Ste}^\circledast$ is projectively complete and co-well-powered in the quotient objects of the class $\DEpi$, and the class $\Mor({\tt Ste}^\circledast,{\tt BanAlg})$ goes from ${\tt Ste}^\circledast$ (since each algebra $A$ can be mapped at least into the Banach algebra consisiting of only one element, zero 0) and is the right ideal. Therefore, the holomorphic envelope $\heartsuit$ is semiregular, and by Theorem  \ref{TH:sushestvovanie-seti-pri-faktorizatsii} it can be defined as a functor. Moreover, by Remark \ref{Omega-podtalkivaet-Mor(K,M)}, each class, in particular, $\DEpi$ pushes $\Mor({\tt Ste}^\circledast,{\tt BanAlg})$, hence the holomorphic envelope is regular, and by Theorem \ref{TH:regulyarnaya-obolochka} it can be defined as an idempotent functor.

5. Suppose $\leftlim{\mathcal B}_A$ is a dense epimorphism. By Lemma \ref{LM:A->B-propusk-cherez-rho_U}, the net  $\mathcal B$ generates on the inside the class of morphisms with values in Banach algebras, hence by Theorem \ref{TH:morfizmy-porozhdayutsya-iznutri} (with $\varOmega=\DEpi$),
$$
\heartsuit_A=\env_{\tt BanAlg}^{\DEpi}A=\env_{\Mor({\tt Ste}^\circledast,{\tt BanAlg})}^{\DEpi}A=\eqref{env_Psi=env_Phi}=\env_{\mathcal B}^{\DEpi}A.
$$
Further, the condition $\leftlim{\mathcal B}_A\in\DEpi$ implies by Lemma \ref{LM:obolochka-konusa}
$$
\env_{\mathcal B}^{\DEpi}A=\eqref{obolochka-konusa-Omega}=\leftlim{\mathcal B}_A.
$$
Again, by Lemma \ref{LM:obolochka-konusa} from $\leftlim{\mathcal B}_A\in\DEpi\subseteq\Epi$ we have
$$
\leftlim{\mathcal B}_A=\eqref{obolochka-konusa-Omega}=\env_{\mathcal B}^{\Epi}A.
$$
And again by Theorem \ref{TH:morfizmy-porozhdayutsya-iznutri} (now with $\varOmega=\Epi$),
$$
\env_{\mathcal B}^{\Epi}A=\eqref{env_Psi=env_Phi}=
\env_{\Mor({\tt Ste}^\circledast,{\tt BanAlg})}^{\Epi}A=\env_{\tt BanAlg}^{\Epi}A.
$$

6. We need to verify here that the holomorphic envelope satisfies the conditions T.1 and T.2 on p.\pageref{DEF:obolochka-soglasovana-s-tenz-proizv}. First, let $\rho:A\to A'$ and $\sigma:B\to B'$ be two holomorphic extensions. Then for any Banach algebra $C$ and for any morphism $\ph:A\circledast B\to C$ there is a pair of morphisms $\ph_A:A\to C$ and $\ph_B:B\to C$ such that
$$
\ph(a\circledast b)=\ph_A(a)\cdot\ph_B(b)=\ph_B(b)\cdot\ph_A(a).
$$
Since $\ph_A$ and $\ph_B$ are morphisms into the Banach algebra $C$, they can be extended along $\rho$ and $\sigma$:
$$
\ph_A=\ph'_A\circ\rho,\qquad \ph_B=\ph'_B\circ\sigma.
$$
Put
$$
\ph'(x\circledast y)=\ph'_A(x)\cdot\ph'_B(y)=\ph'_B(y)\cdot\ph'_A(x),\qquad x\in A',\ y\in B',
$$
then
$$
\ph'((\rho\circledast\sigma)(a\circledast b))=\ph'(\rho(a)\circledast\sigma(b))=
\ph'_A(\rho(a))\cdot\ph'_B(\sigma(b))=\ph_A(a)\cdot\ph_B(b)=\ph(a\circledast b).
$$
Second, let $\sigma:\C\to B$ be a holomorphic extension of the algebra $\C$. It must be a dense epimorphism, and since $\C$ is a finite-dimensional space, this means that $\sigma$ is an epimorphism.
\epr

\bit{

\item We say that a stereotype algebra $A$ is {\it holomorphic}, if it is a complete object with respect to the envelope $\heartsuit$, i.e. its holomorphic envelope is an isomorphism: $\heartsuit_A\in\Iso$.

}\eit

From the property $6^\circ$ and Theorems \ref{TH:sushestvovanie-tenz-proizv-v-L} and \ref{TH:funktor-E-monoidalen} it follows

\btm\label{TH:sushestvovanie-tenz-proizv-golomorfnyh-algebr} The formulas
\beq\label{sushestvovanie-tenz-proizv-golomorfnyh-algebr}
A\stackrel{\heartsuit}{\circledast}B=(A\circledast B)^\heartsuit,\qquad
\ph\stackrel{\heartsuit}{\circledast}\psi=(\ph\circledast\psi)^\heartsuit
\eeq
define the monoidal structure on the category of holomorphic algebras, and the functor of taking holomorphic envelope $A\mapsto A^\heartsuit$ is a monoidal functor (from the category ${\tt Ste}^\circledast$ of stereotype algebras with $\circledast$ as tensor product into the category of holomorphic algebras with $\stackrel{\heartsuit}{\circledast}$ as tensor product).
\etm

One can describe the tensor product \eqref{sushestvovanie-tenz-proizv-golomorfnyh-algebr} in terms of the net $\mathcal{B}$ of Banach quotient maps as follows. Let $A$ and $A'$ be two stereotype algebras. If $U,V,U',V'$ are submultiplicative closed absolutely convex neighbourhoods of zero such that
$$
V\subseteq U\subseteq A,\qquad V'\subseteq U'\subseteq A',
$$
then by multiplying the arising pair of diagrams \eqref{konus-banach-faktor-algebr} we have:
 \beq\label{circledast-konus-banach-faktor-algebr}
 \xymatrix % @R=2.5pc @C=2.5pc
 {
 & A\circledast A' \ar@/_2ex/[dl]_{\rho_V\circledast\rho_{V'}} \ar@/^2ex/[dr]^{\rho_U\circledast\rho_{U'}} &
\\
A/V\circledast A'/V' \ar[rr]_{\pi_V^U\circledast \pi_{V'}^{U'}} & & A/U\circledast A'/U'
 }.
 \eeq
This means that the system of morphisms $\rho_U\circledast\rho_{U'}:A\circledast A'\to A/U\circledast A'/U'$, $U\in\mathcal{SU}(A)$, $U'\in\mathcal{SU}(A')$ is a projective cone of the covariant system $\pi_V^U\circledast \pi_{V'}^{U'}$. As a corollary, there exists a unique morphism $\vartheta:A\circledast A'\to\kern-10pt\leftlim_{\scriptsize\begin{matrix}U\in\mathcal{SU}(A),\\ U'\in\mathcal{SU}(A')\end{matrix}}\kern-10pt A/U\circledast A'/U'$ such that the following diagrams are commutative:
$$
\xymatrix @R=2.5pc @C=2.0pc {
   A\circledast A' \ar@/_2ex/[dr]_{\rho_V\circledast\rho_{V'}}\ar@{-->}[rr]^(.4){\vartheta} & & \kern-10pt
   \leftlim_{\scriptsize\begin{matrix}U\in\mathcal{SU}(A),\\ U'\in\mathcal{SU}(A')\end{matrix}}\kern-10pt A/U\circledast A'/U' \ar@/^2ex/[ld]^{\pi_{V,V'}}  \\
   & A/V\circledast A'/V' &
}
$$
where $V\in\mathcal{SU}(A)$, $V'\in\mathcal{SU}(A')$, and $\pi_{V,V'}$ is the cone of morphisms from the projective limit into the covariant system.

\bprop
For any stereotype algebras $A$ and $A'$
\beq\label{(A-circledast-A')^heartsuit}
(A\circledast A')^\heartsuit=\Im_\infty\vartheta.
\eeq
where $\Im_\infty$ is the element of the nodal decomposition in the category $\tt Ste$ of stereotype spaces (not algebras). In particular, if the algebras $A$ and $A'$ are holomorphic, then
\beq\label{A-stackrel-heartsuit-circledast-A'}
A\stackrel{\heartsuit}{\circledast}A'=\Im_\infty\vartheta.
\eeq
\eprop
\bpr
We need to verify that the map
$\red_\infty\vartheta\circ\coim_\infty\vartheta:A\circledast A'\to \Ran\Im_\infty\vartheta$ is a holomorphic envelope of the algebra $A\circledast A'$ (where $\red_\infty$ and $\coim_\infty$ are elements of the nodal decomposition in $\tt Ste$).

1. Let us show first that this is a holomorphic extension. Take a morphism $\ph:A\circledast A'\to B$ into a Banach algebra $B$. Put
$$
\eta(x)=\ph(x\circledast 1_{A'}),\qquad \eta'(y)=\ph(1_A\circledast a'),\qquad x\in A,\quad y\in A',
$$
then
\beq\label{eta(x)-cdot-eta'(y)=eta'(y)-cdot-eta(x)}
\eta(x)\cdot\eta'(y)=\eta'(y)\cdot\eta(x),\qquad x\in A,\quad y\in A'.
\eeq
and
$$
\ph(x\circledast y)=\eta(x)\cdot\eta'(y)=\eta'(y)\cdot\eta(x),\qquad x\in A,\quad y\in A'.
$$
Let $U$ be a unit ball in $B$. Consider its preimages in $A$ and $A'$
$$
V=\eta^{-1}(U),\qquad V'=(\eta')^{-1}(U),
$$
and morphisms $\psi$ and $\psi'$, such that the following diagrams are commutative:
$$
 \xymatrix  @R=2.5pc @C=2.5pc
 {
& A \ar[ld]_{\eta}\ar[rd]^{\rho_V}
\\
B & & A/V\ar[ll]_{\psi}
 }
 \qquad
 \xymatrix  @R=2.5pc @C=2.5pc
 {
& A' \ar[ld]_{\eta'}\ar[rd]^{\rho_{V'}}
\\
B & & A'/V'\ar[ll]_{\psi'}
 }
$$
From \eqref{eta(x)-cdot-eta'(y)=eta'(y)-cdot-eta(x)} we have the identity
\beq\label{psi(x)-cdot-psi'(y)=psi'(y)-cdot-psi(x)}
\psi(s)\cdot\psi'(t)=\psi'(t)\cdot\psi(s),\qquad s\in A/V,\quad t\in A'/V',
\eeq
which means in its turn that a morphism is defined
$$
\ph_{V,V'}:A/V\circledast A'/V'\to B\qquad \Big|\qquad \ph_{V,V'}(x\circledast y)=\psi(x)\cdot\psi'(y).
$$
We have
$$
\ph(x\circledast y)=\eta(x)\cdot\eta'(y)=\psi(\rho_V(x))\cdot\psi'(\rho_{V'}(y))=
\ph_{V,V'}(\rho_V(x)\circledast\rho_{V'}(y))=\ph_{V,V'}((\rho_V\circledast\rho_{V'})(x\circledast y))
$$
hence the following diagram is commutative
$$
 
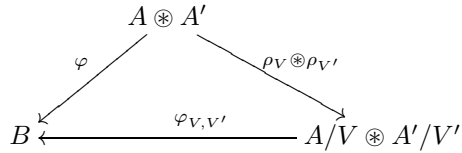
\xymatrix  @R=2.5pc @C=2.5pc
 {
& A\circledast A' \ar[ld]_{\ph}\ar[rd]^{\rho_V\circledast\rho_{V'}}
\\
B & & A/V\circledast A'/V'\ar[ll]_{\ph_{V,V'}}
 }
$$
It can be inserted into the diagram
$$
 
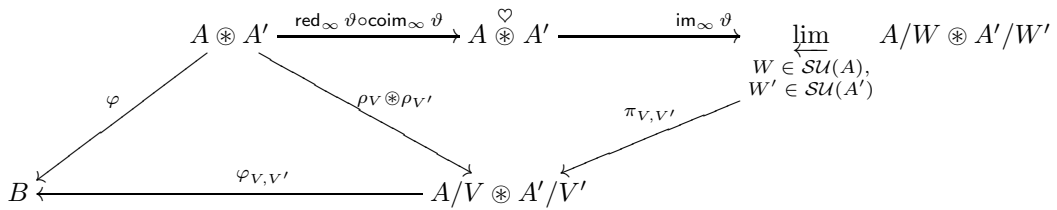
\xymatrix  @R=2.5pc @C=4.5pc
 {
& A\circledast A'
\ar[dl]_{\ph}\ar[rd]^{\rho_V\circledast\rho_{V'}}
\ar[r]^{\red_\infty\vartheta\circ\coim_\infty\vartheta}
&
A\stackrel{\heartsuit}{\circledast}A'
\ar[r]^{\im_\infty\vartheta}
&
\leftlim_{\scriptsize\begin{matrix}W\in\mathcal{SU}(A),\\ W'\in\mathcal{SU}(A')\end{matrix}}
A/W\circledast A'/W'
\ar[dl]_(.6){\pi_{V,V'}}
 \\
B & & A/V\circledast A'/V' \ar[ll]_{\ph_{V,V'}} &
 }
$$
which we can transform into the diagram
$$
 
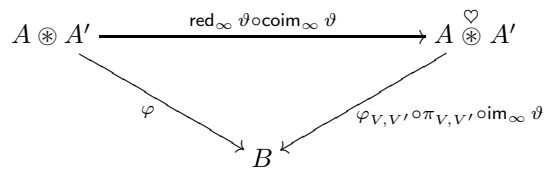
\xymatrix  @R=2.5pc @C=4.5pc
 {
 A\circledast A'
\ar[dr]_{\ph}\ar[rr]^{\red_\infty\vartheta\circ\coim_\infty\vartheta}
& &
A\stackrel{\heartsuit}{\circledast}A'
\ar[dl]^{\quad\ph_{V,V'}\circ\pi_{V,V'}\circ\im_\infty\vartheta}
 \\
 & B &
 }
$$
and it means that $\ph$ is extended along $\red_\infty\vartheta\circ\coim_\infty\vartheta$.

2. Let us show now that $\red_\infty\vartheta\circ\coim_\infty\vartheta$ is a holomorphic envelope. Suppose $\sigma:A\circledast A'\to C$ is another holomorphic extension. Then for any submultiplicative neighbourhoods of zero $V\subseteq A$ and $V'\subseteq A'$ the morphism $\rho_V\circledast\rho_{V'}:A\circledast A'\to A/V\circledast A'/V'$ will be a morphism into a Banach algebra, hence there exists a unique morphism $\widetilde{\rho_V\circledast\rho_{V'}}$ such that the following diagram is commutative:
$$
 \xymatrix  @R=2.5pc @C=2.pc
 {
 A\circledast A'
\ar[dr]_{\rho_V\circledast\rho_{V'}}\ar[rr]^{\sigma}
& &
C
\ar[dl]^{\widetilde{\rho_V\circledast\rho_{V'}}}
 \\
 & A/V\circledast A'/V' &
 }
$$
At the same time for a system of submultiplicative neighbourhoods $W\subseteq V\subseteq A$ and $W'\subseteq V'\subseteq A'$ a diagram arises
$$
 \xymatrix  @R=2.pc @C=4.pc
 {
 & A\circledast A' \ar[d]^{\sigma}\ar@/_5ex/[ddl]_{\rho_W\circledast\rho_{W'}} \ar@/^5ex/[ddr]^{\rho_V\circledast\rho_{V'}} &
\\
& C \ar@{-->}[dl]_{\widetilde{\rho_W\circledast\rho_{W'}}} \ar@{-->}[dr]^{\widetilde{\rho_V\circledast\rho_{V'}}} &
\\
A/W\circledast A'/W' \ar[rr]_{\pi_V^W\circledast\pi_{V'}^{W'}} & & A/V\circledast A'/V'
 }
$$
where the perimeter and the two inner triangles, bordering on the upper vertex, are commutative, and, since $\sigma$ is an epimorphism, this means that the last lower inner triangle is also commutative.

This diagram implies that the system of morphisms $\widetilde{\rho_V\circledast\rho_{V'}}$ forms a projective cone of the contravariant system $\pi_V^W\circledast\pi_{V'}^{W'}$. As a corollary, there exists a unique morphism  $\varkappa$ such that all diagrams
$$
 \xymatrix  @R=2.5pc @C=2.pc
 {
 C
\ar[dr]_{\widetilde{\rho_V\circledast\rho_{V'}}}\ar@{-->}[rr]^{\varkappa}
& &
\leftlim_{W,W'}A/W\circledast A'/W'
\ar[dl]^(.6){\pi_{V,V'}}
 \\
 & A/V\circledast A'/V' &
 }
$$
are cmmutative. In the diagrams of the form
\beq\label{(A-circledast-A')^heartsuit-1}
 \xymatrix  @R=2.5pc @C=2.pc
 {
A/V\circledast A'/V'
\ar@/_5ex/[ddr]_{\rho_V\circledast\rho_{V'}}
\ar[rr]^{\sigma}
\ar[dr]_{\vartheta}
& &
C \ar@/^5ex/[ddl]^{\widetilde{\rho_V\circledast\rho_{V'}}}\ar@{-->}[dl]^{\varkappa}
\\
&
\leftlim_{W,W'}A/W\circledast A'/W'
\ar[d]^{\pi_{V,V'}}
&
 \\
 & A/V\circledast A'/V' &
 }
\eeq
there preimeter and the two lower triangles are commutative. Hence for all $V,V'$
$$
\begin{cases}
\pi_{V,V'}\circ\vartheta=\rho_V\circledast\rho_{V'}\\
\pi_{V,V'}\circ\varkappa\circ\sigma=\rho_V\circledast\rho_{V'}
\end{cases}
$$
and from the uniqueness of $\vartheta$ satisfying these equalities it follows that
$$
\vartheta=\varkappa\circ\sigma
$$
i.e. \eqref{(A-circledast-A')^heartsuit-1} is commutative. Now we obtain a commutative diagram
$$
 \xymatrix  @R=2.5pc @C=6.pc
 {
 A\circledast A'
 \ar@/^5ex/[rr]^{\vartheta}
 \ar@/_2ex/[dr]_{\sigma}\ar[r]_{\red_\infty\vartheta\circ\coim_\infty\vartheta}
&  A\stackrel{\heartsuit}{\circledast}A'
\ar[r]_{\im_\infty\vartheta}
&
\leftlim_{W,W'}A/W\circledast A'/W'
 \\
 & C\ar@/_2ex/[ur]_{\varkappa} &
 }
$$
Here $\sigma\in\DEpi({\tt Ste}^\circledast)=\Epi({\tt Ste})$, and $\im_\infty\vartheta\in\SMono({\tt Ste})$, hence there exists a diagonal $\delta$:
$$
 \xymatrix  @R=2.5pc @C=6.pc
 {
 A\circledast A'
 \ar@/_2ex/[dr]_{\sigma}\ar[r]_{\red_\infty\vartheta\circ\coim_\infty\vartheta}
&  A\stackrel{\heartsuit}{\circledast}A'
\ar[r]_{\im_\infty\vartheta}
&
\leftlim_{W,W'}A/W\circledast A'/W'
 \\
 & C\ar@/_2ex/[ur]_{\varkappa}\ar@{-->}[u]^{\delta} &
 }
$$
Initially $\delta$ is built as a morphism in the category $\tt Ste$, but since $\sigma$ is a dense morphism, $\delta$ must be a homomorphism of algebras, i.e. a morphism in ${\tt Ste}^\circledast$. We have that every extension $\sigma$ are subordinated to $\red_\infty\vartheta\circ\coim_\infty\vartheta$, thus $\red_\infty\vartheta\circ\coim_\infty\vartheta$ is an envelope.
\epr

\paragraph{Fourier transform on a commutative Stein group.}

Let $G$ be a commutative compactly generated Stein group, $\mathcal{O}(G)$ the algebra of holomorphic functions on $G$, and $\mathcal{O}^\star(G)$ the algebra of analytic functionals from Examples \ref{ex-10.5} and \ref{ex-10.9}. Let $G^\bullet$ be the dual group of complex characters on $G$, i.e. continuous homomorphisms $\chi:G\to\C^\times$ into the multiplicative group $\C^\times$ of non-zero complex numbers ($G^\bullet$ is endowed with the pointwise multiplication and the topology of uniform convergence on compact sets in $G$), ${\mathcal F}_G:{\mathcal
O}^\star(G)\to {\mathcal O}(G^\bullet)$ the Fourier transform on $G$, i.e. the homomorphism of algebras acting by formula
$$
\overbrace{{\mathcal F}_G(\alpha)(\chi)}^{\scriptsize \begin{matrix}
\text{value of the function ${\mathcal F}_G(\alpha)\in {\mathcal O}(G^\bullet)$}\\
\text{in the point $\chi\in G^\bullet$} \\ \downarrow \end{matrix}}\kern-35pt=\kern-50pt\underbrace{\alpha(\chi)}_{\scriptsize \begin{matrix}\uparrow \\
\text{action of the functional $\alpha\in{\mathcal O}^\star(G)$}\\
\text{on the function $\chi\in G^\bullet\subseteq {\mathcal O}(G)$ }\end{matrix}}
\kern-50pt \qquad (\chi\in G^\bullet,\quad \alpha\in {\mathcal O}^\star(G))
$$

\btm\label{TH:Fourier-Stein=obolochka} For a compactly generated commutative Stein group $G$ its Fourier transform ${\mathcal F}_G:{\mathcal O}^\star(G)\to {\mathcal O}(G^\bullet)$ is a holomorphic envelope of the algebra ${\mathcal O}^\star(G)$, and coincides with the stereotype Arens-Michael envelope and with the envelope with respect to the class of Banach algebras in the class $\Epi$ of all epimorphisms (in the categories $\tt TopAlg$ and ${\tt
Ste}^\circledast$):
 \beq\label{Fourier-Stein=obolochka}
{\mathcal F}_G=\heartsuit_{{\mathcal O}^\star(G)}=\env_{\tt
BanAlg}^{\DEpi}{\mathcal O}^\star(G)=\env_{\tt BanAlg}^{\Epi}{\mathcal
O}^\star(G)=\leftlim{\mathcal B}_{{\mathcal O}^\star(G)}.
 \eeq
 \etm
\bpr In  \cite{Akbarov-stein-groups} it was shown that in the category $\tt TopAlg$ the local limit of the net of Banach quotient maps on the object ${\mathcal O}^\star(G)$ coincides with ${\mathcal O}(G^\bullet)$:
 \beq\label{O(G^bullet)=leftlim-B_(O^star(G))}
{\mathcal O}(G^\bullet)=\leftlim{\mathcal B}_{{\mathcal O}^\star(G)}.
 \eeq
Here ${\mathcal O}(G^\bullet)$ is a Fr\'echet algebra, so it coincides with its pseudosaturation,
and this implise that \eqref{O(G^bullet)=leftlim-B_(O^star(G))} holds in the category of stereotype algebras. In addition, the morphism ${\mathcal F}_G:{\mathcal O}^\star(G)\to {\mathcal O}(G^\bullet)$, being a local limit in $\tt TopAlg$, is a dense epimorphism, therefore in ${\tt Ste}^\circledast$ it is also a dense epimorphism. Thus, by
\eqref{A^heartsuit=Env_BanAlg^DEpi-A=Env_BanAlg^Epi-A=A_B} we have \eqref{Fourier-Stein=obolochka}.
 \epr

\subsection{Continuous envelope}
\label{SUBSEC:C^*-envelopes}

 \bit{
 \item[$\bullet$]
Let us say that a stereotype algebra $A$ is {\it involutive}, if an operation of involution $x\mapsto\overline{x}$ is defined on $A$ (in the usual sense, see e.g. \cite{Helemskii} or \cite{Murphy}), and this operation is continuous as a map from $A$ into $A$. The involutive stereotype algebras form a category ${\tt InvSte}^\circledast$ where morphisms are continuous involutive unital homomorphisms $\ph:A\to B$:
$$
\ph(\lambda\cdot x+\mu\cdot y)=\lambda\cdot \ph(x)+\mu\cdot \ph(y),\qquad
\ph(x\cdot y)=\ph(x)\cdot\ph(y),\qquad \ph(1)=1,\qquad
\ph(\overline{x})=\overline{\ph(x)}
$$
}\eit
All $C^*$-algebras are obvious examples (\cite{Helemskii}, \cite{Murphy}). Another example is the algebra $\mathcal{C}(M)$ of continuous functions on a paracompact locally compact topological space $M$ from Example \ref{ex-10.3}.

\paragraph{Net of $C^*$-quotient-maps and the Kuznetsova envelope.}

\bit{
\item[$\bullet$]\label{DEF:C^*-seminorm} By {\it $C^*$-seminorm} on an involutive algebra $A$ we mean any seminorm $p:A\to \R_+$ satisfying the following condition:
\beq\label{DEF:C*-polunorma}
p(x\cdot\overline{x})=p(x)^2,\quad x\in A.
\eeq
}\eit
By the Z.~Sebestyen theorem \cite{Sebestyen}, any such seminorm automatically preserves involution and is submultiplicative:
$$
p(\overline{x})=p(x),\qquad p(x\cdot y)\le p(x)\cdot p(y).
$$
The identity \eqref{DEF:C*-polunorma} implies in particular the equality
$$
p(1)=p(1\cdot \overline{1})=p(1)^2,
$$
which mean that $p$ must turn 1 either into 1, or into 0,
$$
p(1)=1\quad\vee\quad p(1)=0,
$$
and the second one means that $p$ vanishes, since in this case
$$
p(x)=p(x\cdot 1)\le p(x)\cdot p(1)=p(x)\cdot 0=0.
$$
Further we will be interested in continuous $C^*$-seminorms on involutive topological algebras.

\bit{
\item[$\bullet$]\label{DEF:C^*-okrestnost} Let us call a {\it $C^*$-neighborhood} in a topological algebra $A$ any closed absolutely convex neighborhood of zero $U$, for which the Minkowski functional
$$
p(x)=\inf\{\lambda>0:\ \lambda\cdot x\in U\}
$$
is a $C^*$-seminorm on $A$. For any such neighborhood of zero $U$ the quotient algebra $A/U$ (defined on p.\pageref{DEF:A/U}) is a
$C^*$-algebra, and we call it the {\it $C^*$-quotient algebra} of $A$, and the natural map $\rho_U:A\to A/U$ will be called a {\it $C^*$-quotient map} of $A$. The symbol ${\mathcal C}^*$ will denote the class of all $C^*$-quotient maps $\{\rho_U:A\to A/U\}$, where $A$ runs over the class of involutive topological algebras, and $U$ over the set of all $C^*$-neighborhoods of zero in $A$.

 }\eit

The following fact is an analog of Proposition \ref{PROP:B-set-v-TopAlg}.

\bprop\label{PROP:C^*-set-v-TopAlg} The class ${\mathcal C}^*$ of all $C^*$-quotient maps is a net of epimorphisms in the category ${\tt InvSte}^\circledast$ of involutive stereotype algebras, and the pre-order $\to$ in ${\mathcal C}^*$ is equivalent to the embedding of the neighborhoods of zero:
 \beq\label{rho_V-to-rho_U-<=>-V-subseteq-U-C^*}
\rho_V\to\rho_U\quad\Longleftrightarrow\quad V\subseteq U.
 \eeq
 \eprop
 \bpr
By definition, the relation $\rho_V\to\rho_U$ means the existence of an involutive continuous homomorphism of $C^*$-algebras $\iota:A/V\to A/U$, such that diagram \eqref{konus-banach-faktor-algebr-*} is commutative. By the well-known property of $C^*$-algebras \cite[Theorem 2.1.7]{Murphy}, the homomorphism $\iota$ cannot increase the $C^*$-norm: $\norm{\iota(x)}\le\norm{x}$. Being applied to $C^*$-seminorms $p_U$ and $p_V$, which correspond to the neighborhoods $U$ and $V$, this means inequality $p_U(x)\le p_V(x)$, which in its turn is equivalent to embedding $V\subseteq U$.
 \epr

 \bit{
\item[$\bullet$] The net ${\mathcal C}^*$ will be called the {\it net of the $C^*$-quotient maps}.

\item[$\bullet$] For each involutive stereotype algebra $A$ the family of the $C^*$-quotient maps $\rho_U:A\to A/U$ is a projective cone of the covariant system $\Bind({\mathcal C}_A)$. The projective limit of this cone in the category  ${\tt InvSte}^\circledast$ of involutive stereotype algebras will be called the {\it Kuznetsova envelope}\footnote{Our terminology and notations differ from those used in \cite{Kuznetsova}.} of $A$ and will be denoted by
 \beq\label{DEF:Kuznetsova=lim-C_A}
\leftlim{\mathcal C}_A:A\to A_{\mathcal C}
 \eeq
(this limit exists, since ${\tt InvSte}^\circledast$ is projectively complete). The range of this morphism
 \beq\label{DEF:Kuznetsova}
A_{\mathcal C}=\Ran\leftlim{\mathcal B}_A=
{\tt InvSte}^\circledast\text{-}\kern-10pt\leftlim_{U\in\mathcal{C^*U}(A)}\kern-10pt A/U=
\Big({\tt InvTopAlg}\text{-}\kern-10pt\leftlim_{U\in\mathcal{C^*U}(A)}\kern-10pt A/U\Big)^\vartriangle.
 \eeq
is also called the {\it Kuznetsova envelope}\label{DEF:C*-obolochka} of $A$.
 }\eit

Theorem \ref{TH:funktorialnost-obolochki_F} implies

\btm\label{TH:Kuznetsova=funktor} The Kuznetsova envelope is an envelope in the class of all morphisms in ${\tt InvSte}^\circledast$ with respect to the system of all $C^*$-quotient maps ${\mathcal C}^*$,
 \beq
 A_{\mathcal C}=\Env_{{\mathcal C}^*}^{\Mor({\tt InvSte}^\circledast)}A,
 \eeq
and to each morphism $\ph:A\to B$ in ${\tt InvSte}^\circledast$ the formula
 \beq\label{DEF:ph^mathcal-C}
\ph_{\mathcal C}=\leftlim_{\tau\in{\mathcal C}_B}\leftlim_{\sigma\in{\mathcal
B}_A}\ph_\sigma^\tau\circ\sigma_{\mathcal C}
 \eeq
assigns a morphism $\ph_{\mathcal C}:A_{\mathcal C}\to B_{\mathcal C}$ such that the following diagram is commutative,
 \beq\label{DIAGR:funktorialnost-mathcal-C}
\xymatrix @R=2.pc @C=5.0pc % @M=14pt
{
A\ar[d]^{\ph}\ar[r]^{\leftlim{\mathcal C}_A} &  A_{\mathcal C}\ar@{-->}[d]^{\ph_{\mathcal C}} \\
B\ar[r]^{\leftlim{\mathcal C}_B} &  B_{\mathcal C} \\
},
 \eeq
and the correspondence $(A,\ph)\mapsto(A_{\mathcal C},\ph_{\mathcal C})$ can be defined as a functor from
${\tt InvSte}^\circledast$ into ${\tt InvSte}^\circledast$.
 \etm

\paragraph{Continuous envelope of an involutive stereotype algebra.}
By dense epimorphism of involutive stereotype algebras we mean the same object as for general (non-involutive) stereotype algebras, i.e. a morphism $\ph:A\to B$ such that the set of values $\ph(A)$ is dense in $B$.

 \bit{
 \item
A {\it continuous envelope} of an involutive stereotype algebra $A$ is its envelope in the class $\DEpi$ of dense epimorphisms in the category ${\tt InvSte}^\circledast$ with respect to the class ${\tt C^*}$ of $C^*$-algebras.
We use the following notation for this construction:
 \beq\label{DEF:A^diamondsuit}
 A^\diamondsuit=\Env_{\tt C^*}^{\DEpi}A, \qquad \diamondsuit_A=\env_{\tt C^*}^{\DEpi}A.
 \eeq
Thus,
 $$
 \big(\diamondsuit_A:A\to A^\diamondsuit\big)=\big(\env_{\tt C^*}^{\DEpi}A:A\to\Env_{\tt C^*}^{\DEpi}A\big).
 $$
 }\eit
The following properties are proved by analogy with the properties of holomorphic envelopes on p.\pageref{1^0:golomorfnaya-obolochka-sushestvuet}.

 \vglue10pt
\centerline{\bf Properties of continuous envelopes:}
 \vglue10pt

\bit{\it

\item[$1^\circ$.] Each involutive stereotype algebra $A$ has a continuous envelope $A^\diamondsuit$.

\item[$2^\circ$.] The continuous envelope $A^\diamondsuit$ is connected with the Kuznetsova envelope $A_{\mathcal C}$ through the formulas
  \beq\label{diamondsuit_A=red-circ-coim}
\diamondsuit_A=\red_\infty\leftlim{\mathcal C}_A\circ\coim_\infty\leftlim{\mathcal C}_A,\qquad
A^\diamondsuit=\Dom\im_\infty\leftlim{\mathcal C}_A
  \eeq
where $\coim_\infty\leftlim{\mathcal C}_A$, $\red_\infty\leftlim{\mathcal C}_A$,
$\im_\infty\leftlim{\mathcal C}_A$ are elements of the nodal decomposition of the morphism $\leftlim{\mathcal C}_A$ in the category ${\tt Ste}$ of stereotype spaces (not algebras!).

\item[$3^\circ$.] For any morphism $\ph:A\to B$ of involutive stereotype algebras and for each choice of continuous envelopes $\diamondsuit_A: A\to A^\diamondsuit$ and $\diamondsuit_B: B\to B^\diamondsuit$ there exists a unique morphism $\ph^\diamondsuit:A^\diamondsuit\to B^\diamondsuit$ such that the following diagram is commutative
    \beq\label{DIAGR:existence-of-ph^diamondsuit}
\xymatrix @R=2.pc @C=5.0pc % @M=14pt
{
A\ar[d]^{\ph}\ar[r]^{\diamondsuit_A} & A^\diamondsuit\ar@{-->}[d]^{\ph^\diamondsuit} \\
B\ar[r]^{\diamondsuit_B} & B^\diamondsuit \\
}
\eeq

\item[$4^\circ$.] The correspondence $(X,\alpha)\mapsto(X^\diamondsuit,\alpha^\diamondsuit)$
can be defined as an idempotent functor from ${\tt InvSte}^\circledast$ into ${\tt InvSte}^\circledast$:
\beq\label{funktorialnost-diamondsuit}
(1_A)^\diamondsuit=1_{A^\diamondsuit},\qquad (\beta\circ\alpha)^\diamondsuit=\beta^\diamondsuit\circ\alpha^\diamondsuit,
\qquad (\alpha^\diamondsuit)^\diamondsuit=\alpha^\diamondsuit.
\eeq

\item[$5^\circ$.] If an algebra $A$ is dense in its Kuznetsova envelope $A_{\mathcal C}$,
$$
\leftlim{\mathcal C}_A\in\DEpi({\tt Ste}^\circledast),
$$
then the continuous envelope of $A$ coincides with its envelope in the class $\Epi$ of all epimorphisms in ${\tt InvSte}^\circledast$ and with the Kuznetsova envelope:
\beq\label{A^diamondsuit=Env_BanAlg^DEpi-A=Env_BanAlg^Epi-A=A_B}
A^\diamondsuit=\Env_{\tt C^*}^{\DEpi}A=\Env_{\tt C^*}^{\Epi}A=A_{\mathcal C}.
\eeq

\item[$6^\circ$.]\label{TH:neprer-obolochka-soglasovana-s-circledast}
The continuous envelope is coherent with the projective tensor product $\circledast$ in
${\tt InvSte}^\circledast$.
}\eit

The following lemma is used in the proof:

\blm\label{LM:A->B-propusk-cherez-rho_U-C^*} In the category ${\tt InvSte}^\circledast$ of involutive stereotype algebras the net $\mathcal C$ of the $C^*$-quotient maps consists of dense epimorphisms and generates on the inside the class $\Mor({\tt InvSte}^\circledast,{\tt C^*})$ of all morphisms with values in $C^*$-algebras:
\beq\label{C^*-porozhdaet-Mor(InvSte,C^*)-iznutri}
{\mathcal C}\subseteq \Mor({\tt InvSte}^\circledast,{\tt C^*})\subseteq
\Mor({\tt InvSte}^\circledast)\circ{\mathcal C}.
\eeq
\elm
\bpr Let $\ph:A\to B$ be amorphism into a $C^*$-algebra $B$, and $V$ be a unit ball in $B$. Consider its pre-image  $U=\ph^{-1}(V)$. This is a neighbourhood of zero in $A$, and its Minkowski functional $p$ is a composition of $\ph$ and the norm on $B$:
$$
p(x)=\inf\{\lambda>0:\ \lambda\cdot x\in \ph^{-1}(V)\}=\inf\{\lambda>0:\
\lambda\cdot\ph(x)\in V\}=\norm{\ph(x)}.
$$
This imples that $p$ is a $C^*$-seminorm on $A$:
$$
p(x\cdot\overline{x})=\norm{\ph(x\cdot\overline{x})}=\norm{\ph(x)\cdot\overline{\ph(x)}}=\norm{\ph(x)}^2=p(x)^2.
$$
That is $U$ is a $C^*$-neighbourhood of zero in $A$. After that the proof of Lemma \ref{LM:A->B-propusk-cherez-rho_U} works.\epr

\bit{

\item An involutive stereotype algebra $A$ is said to be {\it continuous}, if it is a complete object with respect to the envelope $\diamondsuit$, i.e. its continuous envelope is an isomorphism: $\diamondsuit_A\in\Iso$.

}\eit

From $6^\circ$ and Theorems \ref{TH:sushestvovanie-tenz-proizv-v-L} and \ref{TH:funktor-E-monoidalen} it follows

\btm\label{TH:sushestvovanie-tenz-proizv-nepreryvnyh-algebr} Formulas
\beq\label{sushestvovanie-tenz-proizv-nepreryvnyh-algebr}
A\stackrel{\diamondsuit}{\circledast}B=(A\circledast B)^\diamondsuit,\qquad
\ph\stackrel{\diamondsuit}{\circledast}\psi=(\ph\circledast\psi)^\diamondsuit
\eeq
define a monoidal structure on the category of continuous algebras, and the functor of taking continuous envelope $A\mapsto A^\diamondsuit$ is monoidal (from the category os involutive stereotype algebras with the tensor product
$\circledast$ into the category of continuous algebras with the tensor product
$\stackrel{\diamondsuit}{\circledast}$).
\etm

The continuous envelope can be described in terms of the net $\mathcal C$ of the $C^*$-quotient maps by the formula \eqref{(A-circledast-A')^heartsuit} with obvious corrections.

\paragraph{The Gelfand transform as a continuous envelope of a commutative algebra.}
 \bit{

 \item[$\bullet$]
By {\it involutive spectrum} $\Spec(A)$ of an involutive topological (respectively, stereotype) algebra $A$ over $\C$ we mean the set of its involutive characters, i.e. homomorphisms $\chi:A\to \C$ (also continuous, involutive and preserving identity). This set is endowed with the topology of uniform convergence on the totally bounded sets in $A$.

 \item[$\bullet$]
By {\it Gelfand transform} of an involutive stereotype algebra $A$ we mean the natural map ${\mathcal G}_A:A\to C(M)$ of $A$ into the algebra
$C(M)$ of functions on the involutive spectrum $M=\Spec(A)$, continuous on each compact set $K\subseteq M$:
 \beq\label{TOP:vlozh-A-v-C(Spec(A))}
{\mathcal G}_A(x)(t)=t(x),\qquad t\in M=\Spec(A),\ x\in A.
 \eeq
We endow algebra $C(M)$ with the topology which is a pseudosaturation\footnote{The operation of pseudosaturation was defined on p.\pageref{DEF:pseudosaturation}.} of the topology of uniform convergence on compact sets in $M$ -- this turns $C(M)$ into a stereotype algebra. In the special case, when $M$ is a paracompact locally compact space, the topology of uniform convergence on compact sets in $M$ is already a pseudosaturated (and complete) topology on $C(M)$, so $C(M)$ becomes a stereotype algebra already at this step \cite[Sec.8.1]{Akbarov} (and the operation of pseudosaturation do not change this topology anymore).

 \item[$\bullet$]
For each compact set $K\subseteq M$ let us consider the restriction map
$$
\pi_K:C(M)\to C(K),\qquad y\mapsto y\big|_K,
$$
and let ${\mathcal G}_K=\pi_K\circ{\mathcal G}_A$ be the composition
 \beq\label{DEF:preobr-Gelfanda_K}
 \xymatrix @R=2pc @C=1.2pc
 {
 A\ar[rr]^{{\mathcal G}_A}\ar[dr]_{{\mathcal G}_K} & & C(M)\ar[dl]^{\pi_K}\\
  & C(K) &
 }
 \eeq
If $K$ and $L$ are two compact sets in $M$, and $K\subseteq L\subseteq M$, then by symbol $\pi_K^L$ we denote the restriction map
$$
\pi_K^L:C(L)\to C(K),\qquad y\mapsto y\big|_K.
$$
Obviously, the algebra $C(M)$ with the system of projections $\rho_K:C(M)\to C(K)$, $K\subseteq M$, is a projective limit of the system of binding morphisms $\pi_K^L:C(L)\to C(K)$, $K\subseteq L\subseteq M$ (in the category ${\tt InvSte}^\circledast$):
$$
C(M)={\tt InvSte}^\circledast\text{-}\kern-3pt\leftlim_{K\subseteq M}\kern-3pt C(K).
$$
 }\eit

\bprop\label{LM:rho-in-DEpi} For any involutive stereotype algebra $A$ its Gelfand transform ${\mathcal G}_A:A\to C(M)$ is a morphism of stereotype algebras. In the special case when the spectrum $M=\Spec(A)$ of $A$ is a paracompact locally compact space, the morphism ${\mathcal G}_A:A\to C(M)$ is a dense epimorphism.
\eprop
\bpr  In the first part of this proposition only the continuity of the map ${\mathcal G}_A$ is not obvious. Take a base neighborhood of zero $U$ in $C(M)$, i.e. $U=\{f\in C(M):\ \sup_{t\in T}|f(t)|\le \e\}$ for some compact set $T\subseteq M$ and some $\e>0$. Its preimage under the map ${\mathcal G}_A:A\to C(M)$ is the set $\{x\in A:\ \sup_{t\in T}|t(x)|\le \e\}=\e\cdot {^\circ T}$, i.e. the homothety of the polar ${^\circ T}$ of the compact set $T$. Since $A$ is stereotype, ${^\circ T}$ is a neighborhood of zero in it. This proves that the map ${\mathcal G}_A:A\to C(M)$ is continuous if the space $C(M)$ is endowed with the topology of uniform connvergence on compact sets in $M$. Since the space $A$, being stereotype, is pseudosaturated, this means that under the pseudosaturation of the topology in $C(M)$ the map ${\mathcal G}_A:A\to C(M)$ remains continuous (this follows, for example, from \cite[Theorem 1.16]{Akbarov}).

Suppose further that $M=\Spec(A)$ is a paracompact locally compact space. For each compact set $K\subseteq M$ the image ${\mathcal G}_K(A)$ of the algebra $A$ in $C(K)$ under the map ${\mathcal G}_K$ is an involutive subalgebra in $C(K)$, and it contains the identity (and hence, all constant functions) and differs the points $t\in K$. So by the Stone-Weierstrass theorem, ${\mathcal G}_K(A)$ is dense in $C(K)$. This is true for each map  ${\mathcal G}_K=\pi_K\circ\gamma$, where $K$ is a compact set in $M$. Since the topology in $C(M)$ is the projective topology with respect to the maps $\pi_K$, we have that the image ${\mathcal G}_A(A)$ of $A$ in $C(M)$ is dense in $C(M)$. \epr

\btm\label{LM:gamma_A=lim-C^*_A} For each commutative involutive stereotype algebra $A$ the system of morphisms ${\mathcal G}_K:A\to
C(K)$ consists of dense epimorphisms and is isomorphic in the category $\Epi^A$ to the system $\rho_U:A\to A/U$ of all $C^*$-quotient maps of $A$,
 \beq\label{gamma_K-cong-C^*_A}
\{{\mathcal G}_K:A\to C(K), \ K\subseteq\Spec(A)\}\cong {\mathcal C}_A^*.
 \eeq
Under this isomorphism
 \bit{

\item[---] the system of restrictions $\pi_K^L:C(L)\to C(K)$, $K\subseteq L\subseteq M$ turns into the system $\Bind({\mathcal C}_A^*)$ of binding morphisms of the net ${\mathcal C}^*$ on the algebra $A$:
 \beq\label{pi_K^L-cong-Bind(C^*_A)}
\{\pi_K^L:C(L)\to C(K), \ K\subseteq L\subseteq \Spec(A)\}\cong \Bind({\mathcal
C}_A^*).
 \eeq

\item[---] the Gelfand transform ${\mathcal G}_A:A\to C(M)$ is a local limit of the net ${\mathcal C}^*$ on the algebra $A$ (and hence, it coincides with the Kuznetsova envelope of the algebra $A$):
 \beq\label{gamma_A=lim-C^*_A}
 {\mathcal G}_A=\leftlim{\mathcal C}^*_A
 \eeq
 }\eit
 \etm
\bpr On each compact set  $K\subseteq M$ the algebra of functions of the form ${\mathcal G}_A(x)$, where $x\in A$, differs the points, contains constant functions, and is invariant with respect to involution, so it is dense in $C(K)$ by the Stone-Weierstrass theorem. This implies, that the algebra $C(M)$, which contains ${\mathcal G}_A(A)$, is also dense in $C(K)$, so both morphisms ${\mathcal G}_K:A\to C(K)$ and $\pi_K:C(M)\to C(K)$ are dense epimorphisms (in the category ${\tt InvSte}^\circledast$).

The range $A/U$ of each $C^*$-quotient map $\rho_U:A\to A/U$ must be a commutative $C^*$-algebra, hence it is isomorphic to the algebra $C(T_U)$ of continuous functions on its spectrum $T_U$. Under the dual map $\rho_U^\star:\Spec(A)\gets\Spec(A/U)$ this spectrum $T_U$ is homeomorphically turned into a compact set $K_U=\rho_U^\star(T_U)$ in the space $M=\Spec(A)$, and we get the following diagram
$$
 
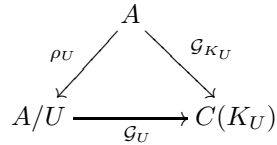
\xymatrix @R=2pc @C=1.2pc
 {
 & A\ar[ld]_{\rho_U}\ar[dr]^{{\mathcal G}_{K_U}} &  \\
 A/U\ar[rr]_{{\mathcal G}_U} & & C(K_U)
 }
$$
where ${\mathcal G}_U$ is the Gelfand transform of the algebra $A/U$ in composition with the map, dual to the homeomorphism $T_U\cong K_U$.

On the contrary, for each compact set $K\subseteq M$ the set
$$
U_K=\{a\in A:\ \sup_{t\in K}\abs{t(a)}\le 1\}
$$
is a $C^*$-neighborhood of zero in $A$. The corresponding quotient algebra $A/U_K$ will be commutative, hence it is isomorphic to the algebra  $C(T_K)$ of continuous functions on its spectrum $T_K$, which is in addition homeomorphic to $K$. If we denote by ${\mathcal G}_K$ the composition of the Gelfand transform of $A$ with the dual map to the homeomorphism $T_K\cong K$, we obtain a commutative diagram
$$
 \xymatrix @R=2pc @C=1.2pc
 {
 & A\ar[ld]_{\rho_{U_K}}\ar[dr]^{{\mathcal G}_{K}} &  \\
 A/U_K\ar[rr]_{} & & C(K)
 }.
$$

Together this proves \eqref{gamma_K-cong-C^*_A}, and \eqref{pi_K^L-cong-Bind(C^*_A)} and \eqref{gamma_A=lim-C^*_A} become its obvious corollaries. \epr

\blm\label{LM:spektr-rasshireniya} If the spectrum $M=\Spec(A)$ of a stereotype algebra $A$ is a $k$-space, then for each extension $\sigma:A\to C$ in the class $\Mor$ of all morphisms (in ${\tt InvSte}^\circledast$) with respect to the class of $C^*$-algebras the dual map of spectra
$$
\sigma^\star:\Spec(C)\to\Spec(A)=M\qquad\Big|\qquad
\sigma(s)=s\circ\sigma,\qquad s\in \Spec(C)
$$
is a homeomorphism of topological spaces.
 \elm
\bpr First, the map $\sigma^\star$ must be an injection, since if some characters $s\ne s'\in\Spec(C)$ have the same image under the action of $\sigma^\star$, i.e.
$$
s\circ\sigma=\sigma^\star(s)=\sigma^\star(s')=s'\circ\sigma,
$$
then this can be understood in such a way that the character $s\circ\sigma=s'\circ\sigma:A\to\C$ has two different continuations on $C$:
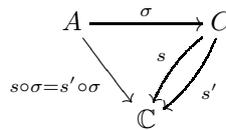
$$
 \xymatrix @R=2pc @C=1.2pc
 {
 A\ar[rr]^{\sigma}\ar[dr]_{s\circ\sigma=s'\circ\sigma} & & C\ar@/_1ex/[dl]_{s}\ar@/^1ex/[dl]^{s'}\\
  & \C &
 }
$$
This is impossible, since $\sigma$ is an extension, in particular, with respect to the $C^*$-algebra $\C$.

On the other hand, the map $\sigma^\star$ is a covering, i.e. for each compact set $K$ in $M$ there is a compact set $T$ in $\Spec(C)$ such that
$\sigma^\star(T)\supseteq K$. Indeed, if $K$ is a compact set in $M=\Spec(A)$, then, since $\sigma:A\to C$ is an extension with respect to the class of $C^*$-algebras, the natural homomorphism ${\mathcal G}_K:A\to C(K)$ into the $C^*$-algebra $C(K)$ have a continuation to $C$, i.e. a diagram arises:
$$
 
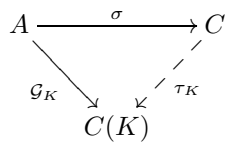
\xymatrix @R=2pc @C=1.2pc
 {
 A\ar[rr]^{\sigma}\ar[dr]_{{\mathcal G}_K} & & C\ar@{-->}[dl]^{\tau_K}\\
  & C(K) &
 }
$$
If we now put $T=\tau_K^\star(K)$, then
$$
\sigma^\star(T)=\sigma^\star\Big(\tau_K^\star(K)\Big)={\mathcal
G}_K^\star(K)=K.
$$
In addition, from the fact that $\sigma^\star$ is a covering, it follows that it is surjective. We obtain that $\sigma^\star:\Spec(C)\to\Spec(A)$ is a continuous bijective covering. Since $\Spec(A)$ is a $k$-space, the map $\sigma^\star$ is open, and thus, a homeomorphism. \epr

The following result supplements the results of Yu.~N.~Kuznetsova's paper \cite{Kuznetsova}:

\btm\label{EX:C^*-obolochka} If $A$ is a commutative involutive stereotype algebra with a paracompact locally compact involutive spectrum $M=\Spec(A)$, then its Gelfand transform ${\mathcal G}_A:A\to C(M)$ is its continuous envelope,  the Kuznetsova envelope, and the envelope in the classes of all morphisms and all epimorphisms in the category
${\tt InvSte}^\circledast$ with respect to the class of $C^*$-algebras:
$$
C(M)=A^\diamondsuit=\Env_{\tt C^*}^{\DEpi} A=\Env_{\tt C^*}^{\Epi} A=\Env_{\tt C^*}^{\Mor} A=\leftlim{\mathcal C}_A.
$$
 \etm
\bpr The equality $C(M)=\leftlim{\mathcal C}_A$ was already proved in Theorem \ref{LM:gamma_A=lim-C^*_A}.
On the other hand, by Proposition \ref{LM:rho-in-DEpi}, the morphism ${\mathcal G}_A:A\to C(M)$ is a dense epimorphism, and by \eqref{A^heartsuit=Env_BanAlg^DEpi-A=Env_BanAlg^Epi-A=A_B} we have the following chain
$$
C(M)=A^\diamondsuit=\Env_{\tt C^*}^{\DEpi} A=\Env_{\tt C^*}^{\Epi} A=\leftlim{\mathcal C}_A.
$$
It remains to prove the equality where the upper index is the class $\Mor$ of all morphisms in ${\tt InvSte}^\circledast$:
$$
\Env_{C^*}^{\Mor} A=C(M).
$$
Let us show first that ${\mathcal G}_A:A\to C(M)$ is an extension of $A$ with respect to the class of $C^*$-algebras. Let $\ph:A\to B$ be a morphism of $A$ into a $C^*$-algebra $B$. To construct a dotted arrow $\ph'$ for  \eqref{DEF:diagr-rasshirenie},
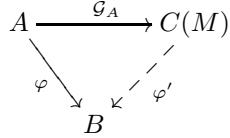
$$
 \xymatrix @R=2pc @C=1.2pc
 {
 A\ar[rr]^{{\mathcal G}_A}\ar[dr]_{\ph} & & C(M)\ar@{-->}[dl]^{\ph'}\\
  & B &
 }
$$
it is sufficient to think that $B$ is commutative and that $\ph(A)$ is dense in $B$ (since otherwise we can replace $B$ by the closure $\overline{\ph(A)}$ in $B$, and this is a commutative subalgebra in $B$). Then from the commutativity of $B$ it follows that $B$ has the form $C(K)$, and from the density of $\ph(A)$ in $B$ that the compact space $K$ is injectively embedded into $M=\Spec(A)$. Thus our diagram can be represented in the form
$$
 
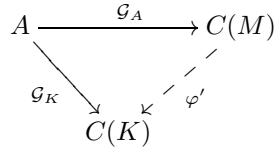
\xymatrix @R=2pc @C=1.2pc
 {
 A\ar[rr]^{{\mathcal G}_A}\ar[dr]_{{\mathcal G}_K} & & C(M)\ar@{-->}[dl]^{\ph'}\\
  & C(K) &
 }
$$
where $K$ is a compact set in $M$, and ${\mathcal G}_K$ is defined in \eqref{DEF:preobr-Gelfanda_K}. It is clear that $\ph'$ can be now defined as the restriction map $\pi_K$ from $M$ to $K$, which we considered above.
$$
 
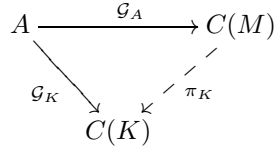
\xymatrix @R=2pc @C=1.2pc
 {
 A\ar[rr]^{{\mathcal G}_A}\ar[dr]_{{\mathcal G}_K} & & C(M)\ar@{-->}[dl]^{\pi_K}\\
  & C(K) &
 }
$$
And this dotted arrow is unique since by Proposition \ref{LM:rho-in-DEpi} ${\mathcal G}_A$ is an epimorphism.

Let us check now that ${\mathcal G}_A:A\to C(M)$ is a maximal extension, i.e. if we take another extension $\sigma:A\to C$, then there exists a morphism $\upsilon:C\to C(M)$ such that the following diagram is commutative:
 \beq\label{TOP:C->C(M)}
 \xymatrix @R=2pc @C=1.2pc
 {
 & A\ar[rd]^{{\mathcal G}_A}\ar[ld]_{\sigma} & \\
 C\ar@{-->}[rr]_{\upsilon} & & C(M)
 }
 \eeq
By Lemma \ref{LM:spektr-rasshireniya} the dual map of spectra
$\sigma^\star:\Spec(C)\to\Spec(A)=M$ is a homeomorphism. Therefore, the following map is defined:
$$
\upsilon:C\to C(M)\qquad \Big|\qquad
\upsilon(y)(t)=\underbrace{(\sigma^*)^{-1}(t)}_{\scriptsize\begin{matrix}\text{\rotatebox{90}{$\owns$}}
\\ \Spec(C)\end{matrix}}(y),\qquad y\in C,\ t\in M.
$$
It is trivially checked that this is a morphism of involutive stereotype algebras.
In addition \eqref{TOP:C->C(M)} will be commutative:
$$
\upsilon(\sigma(x))(t)=(\sigma^*)^{-1}(t)(\sigma(x))=\sigma^*\big((\sigma^*)^{-1}(t)\big)(x)=t(x)={\mathcal G}_A(x)(t),\qquad
x\in A,\ t\in M
$$
i.e. $\upsilon\circ\sigma={\mathcal G}_A$.

It remains to verify that the dotted arrow in \eqref{TOP:C->C(M)} is unique. Suppose that $\upsilon'$ is another dotted arrow with the same properties:
 \beq\label{upsilon-circ-sigma=rho=upsilon'-circ-sigma}
\upsilon\circ\sigma={\mathcal G}_A=\upsilon'\circ\sigma.
 \eeq
If $\upsilon$ and $\upsilon'$ are different, they do not coincide on some vector $y\in C$:
$$
\upsilon(y)\ne\upsilon'(y).
$$
Here in both sides there are functions on $M$, so the inequality means that they do not coincide in some point $t\in M$:
$$
\upsilon(y)(t)\ne\upsilon'(y)(t).
$$
Put
$$
s(z)=\upsilon(z)(t),\qquad s'(z)=\upsilon'(z)(t),\qquad z\in C,
$$
then we see that two different characters on $C$ give a same character in composition with $\sigma$:
$$
s(\sigma(x))=\upsilon(\sigma(x))(t)=\eqref{upsilon-circ-sigma=rho=upsilon'-circ-sigma}=
\upsilon'(\sigma(x))(t)=s'(\sigma(x)),\qquad x\in A.
$$
By Lemma \ref{LM:spektr-rasshireniya} this is impossible, so our initial supposition that $\upsilon\ne\upsilon'$ is also not true.
 \epr

\paragraph{Fourier transform on a commutative locally compact group.}

Let $G$ be a commutative locally compact group, $\mathcal{C}(G)$ the algebra of continuous functions on $G$, and $\mathcal{C}^\star(G)$ the algebra of measures with compact support on $G$ from Examples \ref{ex-10.3} and \ref{ex-10.7}. Let $G^\bullet$ be the dual group of characters on $G$, i.e. continuous homomorphisms $\chi:G\to\T$ into the circle $\T$
($G^\bullet$ is endowed with the pointwise algebraic operations and the topology of uniform convergence on compact sets in $G$), ${\mathcal F}_G:{\mathcal C}^\star(G)\to {\mathcal C}(G^\bullet)$ the Fourier transform on $G$, i.e. the homomorphism of algebras, acting by formula
$$
\overbrace{{\mathcal F}_G(\alpha)(\chi)}^{\scriptsize \begin{matrix}
\text{value of the function ${\mathcal F}_G(\alpha)\in {\mathcal C}(G^\bullet)$}\\
\text{in the point $\chi\in G^\bullet$} \\ \downarrow \end{matrix}}\kern-35pt=\kern-50pt\underbrace{\alpha(\chi)}_{\scriptsize \begin{matrix}\uparrow \\
\text{action of the functional $\alpha\in{\mathcal C}^\star(G)$}\\
\text{at the function $\chi\in G^\bullet\subseteq {\mathcal C}(G)$ }\end{matrix}}
\kern-50pt \qquad (\chi\in G^\bullet,\quad \alpha\in {\mathcal C}^\star(G))
$$

The following observation belongs to J.~N.~Kuznetsova \cite{Kuznetsova}:

\btm\label{TH:Fourier-LCA=obolochka} For each commutative locally compact group $G$ its Fourier transform ${\mathcal F}_G:{\mathcal
C}^\star(G)\to {\mathcal C}(G^\bullet)$ is a continuous envelope of the algebra ${\mathcal C}^\star(G)$, and it coincides with the Kuznetsova envelope and with the envelopes with respect to the class of $C^*$-algebras in the classes $\Mor$ of all morphisms and $\Epi$ of all epimorphisms (in the categories $\tt InvTopAlg$ and ${\tt InvSte}^\circledast$):
$$
{\mathcal F}_G=\diamondsuit_{{\mathcal C}^\star(G)}=\env_{\tt C^*}^{\DEpi}{\mathcal
C}^\star(G)=\env_{\tt C^*}^{\Epi}{\mathcal C}^\star(G)=\env_{\tt C^*}^{\Mor}{\mathcal
C}^\star(G)=\leftlim{\mathcal C}_{{\mathcal C}^\star(G)}.
$$
 \etm
\bpr The spectrum of the algebra ${\mathcal C}^\star(G)$ is homeomorphic to $G^\bullet$, so everything follows from Theorem \ref{EX:C^*-obolochka}. \epr

\section{Errata}

After the publication of this article I found several errors in it. I apologize for them and I give here the necessary corrections. The errors occur in Lemmas \ref{LM:obolochka-konusa} and \ref{LM:nachinka-konusa}, Theorems \ref{TH:funktorialnost-obolochki_F}, \ref{TH:funktorialnost-otpechatka_F}, \ref{TH:Arens-Michael=funktor} and \ref{TH:Kuznetsova=funktor}, and in Proposition \ref{lm-4.12}.

The error in Lemma \ref{LM:obolochka-konusa} is the following: if $\rho=\leftlim \rho^i:X\to\leftlim X^i$ is a projective limit of the projective cone $\{\rho^i:X\to X^i;\ i\in I\}$ from a given object $X$ into a covariant (or contravariant) system $\{X^i,\iota^j_i\}$, then in the diagram
\beq\label{proof-obolochka=lim-0-NEW}
\begin{diagram}
\node{X}\arrow[2]{e,t}{\rho}\arrow{se,b}{\rho^j}\node[2]{\leftlim X^i}\arrow{sw,b,--}{\pi^j}
\\
\node[2]{X^j}
\end{diagram}
\eeq
the morphism $\pi^j$ is not necessarily unique. As an example we can consider a projective limit in the category of vector spaces, for which some prolection $\pi^j$ is not zero, and if we take $X=0$, then  $\rho^j=0$, and in Diagram \eqref{proof-obolochka=lim-0} apart from the initial one $\pi^j\ne 0$ we can take $\pi^j=0$.

The right formulation of Lemma \ref{LM:obolochka-konusa} must be the following (the difference is that we add condition \eqref{rho-in-Epi}):

\blm\label{LM:obolochka-konusa-NEW} If the projective limit $\rho=\leftlim \rho^i:X\to\leftlim X^i$ of a projective cone $\{\rho^i:X\to X^i;\ i\in I\}$ from a given object $X$ into a covariant (or a contravariant) system $\{X^i,\iota^j_i\}$ is an epimorphism
\beq\label{rho-in-Epi}
\rho\in\Epi,
\eeq
then it is an envelope of the object $X$ in an arbitrary class $\varOmega$ that contains $\rho$ with respect to the system of morphisms $\{\rho^i; i\in I\}$:
 \beq\label{obolochka-konusa-Omega-NEW}
\rho=\leftlim \rho^i\in\varOmega
\quad\Longrightarrow\quad
\Env_{\{\rho^i;\ i\in I\}}^\varOmega X=\leftlim X^i
 \eeq
 In particular, this is always true for $\varOmega=\Mor({\tt K})$:
 \beq\label{obolochka-konusa-NEW}
\Env_{\{\rho^i;\ i\in I\}}^{\Mor({\tt K})}X=\leftlim X^i
 \eeq
\elm

A similar error occurs in the formulation of the dual Lemma \ref{LM:nachinka-konusa}. Its correct formulation must be the following (we add here condition \eqref{rho-in-Mono}):

\blm\label{LM:nachinka-konusa-NEW}
If an injective limit $\rho=\rightlim \rho^i:X\gets\rightlim X^i$ of an injective cone  $\{\rho^i:X\gets X^i;\ i\in I\}$ into a given object $X$ from a covariant (or a contravariant) system  $\{X^i,\iota^j_i\}$ is a monomorphism,
\beq\label{rho-in-Mono}
\rho\in\Mono,
\eeq
then it is a refinement if the object $X$ in an arbitrary class $\varGamma$ containing $\rho$ by means of the system of morphisms $\{\rho^i; i\in I\}$:
 \beq\label{otpechatok-konusa-Omega-NEW}
\rho=\rightlim \rho^i\in\varGamma
\quad\Longrightarrow\quad
\Rf_{\{\rho^i;\ i\in I\}}^\varGamma X=\rightlim X^i
 \eeq
 In particular, this is always true for $\varGamma=\Mor({\tt K})$:
 \beq\label{otpechatok-konusa-NEW}
\Rf_{\{\rho^i;\ i\in I\}}^{\Mor({\tt K})}X=\rightlim X^i
 \eeq
 \elm

These corrections imply four more corrections in the further text. First, in Theorem \ref{TH:funktorialnost-obolochki_F} we must add a new condition, \eqref{lim_N-subseteq-varOmega*}:

\btm\label{TH:funktorialnost-obolochki_F-NEW} Let ${\mathcal N}$ be a net of epimorphisms in the category ${\tt K}$. Then
 \bit{

\item[(a)] for each morphism $\alpha:X\to Y$ in ${\tt K}$ and any choice of local limits $\leftlim {\mathcal N}^X$ and $\leftlim {\mathcal N}^Y$ the formula
 \beq\label{DEF:alpha_F-NEW}
\alpha_{\mathcal N}=\leftlim_{\tau\in{\mathcal
N}_Y}\leftlim_{\sigma\in{\mathcal N}^X}\alpha_\sigma^\tau\circ\sigma_{\mathcal
N}
 \eeq
defines a morphism $\alpha_{\mathcal N}:X_{\mathcal N}\to Y_{\mathcal N}$ such that
 \beq\label{DIAGR:funktorialnost-lim-F-NEW}
\xymatrix @R=2.pc @C=10.0pc % @M=14pt
{
X\ar[d]^{\alpha}\ar[r]^{\leftlim {\mathcal N}^X} &  X_{\mathcal N}
\ar@{-->}[d]^{\alpha_{\mathcal N}} \\
Y\ar[r]^{\leftlim {\mathcal N}^Y} &  Y_{\mathcal N}
},
 \eeq

\item[(b)] the envelope $\Env_{\mathcal N}$ can be defined as a functor.
}\eit
If, in addition, all the local limits $\leftlim{\mathcal N}^X$ are epimorphisms
\beq\label{lim_N-subseteq-varOmega*}
\{\leftlim{\mathcal N}^X; \ X\in\Ob(\tt K)\}\subseteq\Epi(\tt
K),
\eeq
then
\bit{
\item[(c)]
for each object $X$ in ${\tt K}$ the local limit
$\leftlim{\mathcal N}^X:X\to X_{\mathcal N}$ is an envelope $\env_{\mathcal N} X$ of $X$ in the category  ${\tt K}$ with respect to the class of morphisms ${\mathcal N}$:
 \beq\label{lim F_X=env_F-X-NEW}
 \leftlim {\mathcal N}^X=\env_{\mathcal N} X,
 \eeq

 }\eit

\etm

And, second, in the dual Theorem \ref{TH:funktorialnost-otpechatka_F} we add a new condition, \eqref{lim_N-subseteq-varOmega**}:

\btm\label{TH:funktorialnost-otpechatka_F-NEW} Let ${\mathcal N}$ be a net of monomorphisms in a category ${\tt K}$. Then
 \bit{

\item[(a)] for each morphism $\alpha:X\to Y$ in
${\tt K}$ and for each choice of local limits $\rightlim {\mathcal N}_X$ and $\rightlim {\mathcal N}_Y$ the formula
 \beq\label{DEF:alpha_F-mono-NEW}
\alpha_{\mathcal N}=\rightlim_{\sigma\in{\mathcal
N}_X}\rightlim_{\tau\in{\mathcal N}_Y}\tau_{\mathcal N}\circ\alpha_\sigma^\tau
 \eeq
defines a morphism $\alpha_{\mathcal N}:X_{\mathcal N}\to Y_{\mathcal N}$ such that
 \beq\label{DIAGR:funktorialnost-lim-F-mono-NEW}
\xymatrix @R=2.pc @C=10.0pc % @M=14pt
{
X\ar[d]^{\alpha} &  X_{\mathcal N} \ar@{-->}[d]^{\alpha_{\mathcal N}}\ar[l]_{\rightlim {\mathcal N}_X} \\
Y &  Y_{\mathcal N} Y\ar[l]_{\rightlim {\mathcal N}_Y} \\
}
 \eeq

\item[(b)] the refinement $\Rf_{\mathcal N}$ can be defined as a functor.
 }\eit
If, in addition, all the local limits $\rightlim{\mathcal N}^X$ are monomorphisms
\beq\label{lim_N-subseteq-varOmega**}
\{\leftlim{\mathcal N}^X; \ X\in\Ob(\tt K)\}\subseteq\Mono(\tt
K),
\eeq
then
\bit{
\item[(c)] for each object $X$ in ${\tt K}$ the local limit $\rightlim{\mathcal N}_X:X_{\mathcal N}\to X$ is a refinement $\rf_{\mathcal N} X$ of $X$ in the category ${\tt K}$ by means of the class of morphisms ${\mathcal N}$:
 \beq\label{lim F_X=imp_F-X-NEW}
\rightlim {\mathcal N}_X=\rf_{\mathcal N} X,
 \eeq
}\eit

\etm

Further, in Theorem \ref{TH:Arens-Michael=funktor} the claim that the stereotype Arens---Michael envelope is an envelope with respect to the system of Banach quotient map remains unproved. As a corollary, the accurate formulation of this result must be as follows:

\btm\label{TH:Arens-Michael=funktor-NEW} The Arens---Michael envelope is a functor from ${\tt Ste}^\circledast$ into ${\tt Ste}^\circledast$. To any morphism $\ph:A\to B$ in ${\tt Ste}^\circledast$ it assigns by formula
 \beq\label{DEF:ph^heartsuit-NEW}
\ph_{\mathcal B}=\leftlim_{\tau\in{\mathcal B}_B}\leftlim_{\sigma\in{\mathcal
B}_A}\ph_\sigma^\tau\circ\sigma_{\mathcal B}
 \eeq
a morphism $\ph_{\mathcal B}:A_{\mathcal B}\to B_{\mathcal B}$ such that
 \beq\label{DIAGR:funktorialnost-heartsuit-NEW}
\xymatrix @R=2.pc @C=5.0pc % @M=14pt
{
A\ar[d]^{\ph}\ar[r]^{\leftlim{\mathcal B}_A} &  A_{\mathcal B}\ar@{-->}[d]^{\ph_{\mathcal B}} \\
B\ar[r]^{\leftlim{\mathcal B}_B} &  B_{\mathcal B} \\
}.
 \eeq
 \etm

In Theorem \ref{TH:Kuznetsova=funktor} the claim that the Kuznetsova envelope is an envelope with respect to the system of all $C^*$-quotient maps remains unproved. As a corollary, the accurate formulation of this result must be as follows:

\btm\label{TH:Kuznetsova=funktor} The Kuznetsova envelope is a functor from ${\tt Ste}^\circledast$ into ${\tt Ste}^\circledast$. To any morphism $\ph:A\to B$ in ${\tt Ste}^\circledast$ it assigns by formula
 \beq\label{DEF:ph^mathcal-C}
\ph_{\mathcal C}=\leftlim_{\tau\in{\mathcal C}_B}\leftlim_{\sigma\in{\mathcal
B}_A}\ph_\sigma^\tau\circ\sigma_{\mathcal C}
 \eeq
a morphism $\ph_{\mathcal C}:A_{\mathcal C}\to B_{\mathcal C}$ such that the following diagram is commutative,
 \beq\label{DIAGR:funktorialnost-mathcal-C}
\xymatrix @R=2.pc @C=5.0pc % @M=14pt
{
A\ar[d]^{\ph}\ar[r]^{\leftlim{\mathcal C}_A} &  A_{\mathcal C}\ar@{-->}[d]^{\ph_{\mathcal C}} \\
B\ar[r]^{\leftlim{\mathcal C}_B} &  B_{\mathcal C} \\
},
 \eeq
and the correspondence $(A,\ph)\mapsto(A_{\mathcal C},\ph_{\mathcal C})$ can be defined as a functor from
${\tt InvSte}^\circledast$ into ${\tt InvSte}^\circledast$.
 \etm

Finally, the error in Proposition \ref{lm-4.12} is that the statement (a) must from the very beginning assume the pseudocompleteness of the space $X$. The correct formulation of Proposition \ref{lm-4.12} must be the following:

\bprop \label{lm-4.12-NEW}
Let $E$ be a closed subspace in a locally convex space $X$, considered as a locally convex space with the topology induced from $X$, and let $E^\perp$ be the annihilator of the space $E$ in the dual space $X^\star$, again considered as a locally convex space with the topology induced from $X^\star$. Then
\begin{itemize}
\item[(a)] if $X$ is pseudocomplete, then there is a natural isomorphism of locally convex spaces
\beq\label{eq4.2-NEW}
  E^\star \cong X^\star/E^\perp,
\eeq
that generates isomorphisms of stereotype spaces
\beq\label{eq4.3-NEW}
  (E^\vartriangle)^\star \cong (X^\star/E^\perp)^\triangledown,  \quad\quad\quad
  E^\vartriangle \cong [(X^\star/E^\perp)^\triangledown]^\star;
\eeq
\item[(b)] if $X$ is stereotype then there is a natural isomorphism of locally convex spaces
\beq  \label{eq4.4-NEW}
  (E^\perp)^\star \cong X/E,
\eeq
that generates the isomorphisms of stereotype spaces
\beq
  ((E^\perp)^\vartriangle)^\star \cong (X/E)^\triangledown,  \quad\quad\quad
  (E^\perp)^\vartriangle \cong [(X/E)^\triangledown]^\star.
\label{eq4.5-NEW}
\eeq
\end{itemize}
\eprop

This error originates from an earlier author's work \cite{Akbarov}, where it is contained in Lemma 2.18. The exact formulation of this lemma must be the following:

\blm\label{lm-2.18-NEW} For each locally convex space $X$ and for each its closed subspace $E$ (with the topology induced from $X$) the mapping
$$
g\in E^\star \mapsto \widetilde{g}\in X^\star/E^\perp
$$
is bijective and open. If in addition $X$ is pseudocomplete, then this mapping is
(continuous and hence) an isomorphism of locally convex spaces:
\beq   \label{eq2.5-NEW}
  E^\star \cong X^\star/E^\perp
\eeq
\elm
\bpr  The bijectivity of this map is obvious.

1) Let us prove that it is open. For each basic neighbourhood of zero
$W=S^\circ_{E^\star}\in \mathcal{BU}(E^\star)$, $S\in \mathcal{BS}(E)$, we can consider the neighbourhood $U=S^\circ_{X^\star}\in \mathcal{BU}(X^\star)$ and put $V=\pi(U)$. Then  $\widetilde{W}=V$.

2) Suppose now that the space $X$ is pseudocomplete, and show that the map $g\mapsto \widetilde{g}$ is continuous. Take $V\in \mathcal{BU}(X^\star/E^\perp)$, and put $U=\pi^{-1}(V)$, where $\pi :X^\star\to X^\star/E^\perp$ is the quotient map. Then $V=\pi(U)$, $U\in \mathcal{BU}(X^\star)$ and $U+E^\perp=U$. Since $X$ is pseudocomplete, by \cite[Theorem 2.12]{Akbarov}, $X\cong X^{\star\star}$ (a bijection of sets), and this implies
$$
(^\circ_X{U})^\circ_{X^\star}=U.
$$
If we put $S={^\circ_X{U}}$, then we have $U=S^\circ_{X^\star}$, and by \cite[Theorem 2.1(d)]{Akbarov}, $S\in\mathcal{BS}(X)$. The condition $U+E^\perp=U$ is equivalent to the condition $S^\circ+E^\perp=S^\circ$ and means that $S\subseteq E$. Now we set  $W=S^\circ_{E^\star}\in\mathcal{BU}(E^\star)$, and we obtain that $W$ is the preimage of $V$ under the mapping $g\in E^\star \mapsto \widetilde{g}\in
X^\star/E^\perp$:
 \begin{multline*}
\widetilde{g}\in V=\pi(U)=\pi\left(S^\circ_{X^\star}\right)
\quad\Longleftrightarrow\quad \exists f\in S^\circ_{X^\star}:\ f|_E=g
\quad\Longleftrightarrow\quad \exists f\in X^\star:\ |f|_S\le 1,\ f|_E=g
\quad\Longleftrightarrow\\ \Longleftrightarrow\quad \exists f\in X^\star:\
|g|_S\le 1,\ f|_E=g \quad\Longleftrightarrow\quad |g|_S\le 1
\quad\Longleftrightarrow\quad g\in S^\circ_{E^\star}=W
 \end{multline*}
\epr

\bpr[Proof of Proposition \ref{lm-4.12-NEW}]
(a) If $X$ is pseudocomplete, then Equality \eqref{eq4.2-NEW} is just Equality \eqref{eq2.5-NEW} from Lemma  \ref{lm-2.18-NEW}. It implies
$$
(E^\vartriangle)^\star\cong\text{\cite[(3.18)]{Akbarov}}\cong (E^\star)^\triangledown \cong (X^\star/E^\perp)^\triangledown
$$
and (since $E^\vartriangle$ is a stereotype space)
$$
 E^\vartriangle \cong (E^\vartriangle)^{\star\star} \cong [(X^\star/E^\perp)^\triangledown]^\star.
$$
(b) If $X$ is stereotype, then $X=Y^\star$ for some stereotype space $Y$ (namely, for $Y=X^\star$). Hence from the already proven \eqref{eq4.2-NEW} we have \eqref{eq4.4-NEW}:
$$
(E^\perp)^\star\cong Y^\star/(E^\perp)^\perp\cong X/E
$$
Now we obtain
$$
((E^\perp)^\vartriangle)^\star\cong\text{\cite[(3.18)]{Akbarov}} \cong ((E^\perp)^\star)^\triangledown \cong (X/E)^\triangledown
$$
and (since $(E^\perp)^\vartriangle$ is a stereotype space)
$$
(E^\perp)^\vartriangle\cong ((E^\perp)^\vartriangle)^{\star\star} \cong [(X/E)^\triangledown]^\star
$$
\epr

The error in the formulation of Proposition \ref{lm-4.12} does not imply the further errors in the text of the paper.

\end{document}